\documentclass[bj,authoryear]{imsart}

\RequirePackage{amsthm,amsmath,amsfonts,amssymb}
\allowdisplaybreaks
\RequirePackage{graphicx}%% uncomment this for including figures

\usepackage{subcaption}
\usepackage{float} % [htbp] option for figure

\usepackage{caption}
\usepackage{booktabs}
\usepackage{float}
\usepackage{multirow}
\usepackage{diagbox}

\usepackage{xcite}
\usepackage{xr}
\makeatletter
\newcommand*{\addFileDependency}[1]{% argument=file name and extension
	\typeout{(#1)}% latexmk will find this if $recorder=0 (however, in that case, it will ignore #1 if it is a .aux or .pdf file etc and it exists! if it doesn't exist, it will appear in the list of dependents regardless)
	\@addtofilelist{#1}% if you want it to appear in \listfiles, not really necessary and latexmk doesn't use this
	\IfFileExists{#1}{}{\typeout{No file #1.}}% latexmk will find this message if #1 doesn't exist (yet)
}
\makeatother
\newcommand*{\myexternaldocument}[1]{%
	\externaldocument{#1}%
	\addFileDependency{#1.tex}%
	\addFileDependency{#1.aux}%
}
\myexternaldocument{CoDA-supp}
\startlocaldefs
\theoremstyle{plain}

\newtheorem{theorem}{Theorem}[section]
\newtheorem{lemma}[theorem]{Lemma}
\newtheorem{Corollary}[theorem]{Corollary}
\newtheorem{Proposition}[theorem]{Proposition}
\newtheorem{Condition}[theorem]{Condition}
\newtheorem{Assumption}[theorem]{Assumption}
\theoremstyle{remark}

\newtheorem{remark}{Remark}
\newcommand{\bA}{\boldsymbol{A}}
\newcommand{\bB}{\boldsymbol{B}}
\newcommand{\bC}{\boldsymbol{C}}
\newcommand{\bD}{\boldsymbol{D}}

\newcommand{\bI}{\boldsymbol{I}}

\newcommand{\bM}{\boldsymbol{M}}

\newcommand{\bQ}{\boldsymbol{Q}}
\newcommand{\bR}{\boldsymbol{R}}
\newcommand{\bS}{\boldsymbol{S}}

\newcommand{\bW}{\boldsymbol{W}}
\newcommand{\bX}{\boldsymbol{X}}
\newcommand{\bY}{\boldsymbol{Y}}
\newcommand{\bZ}{\boldsymbol{Z}}
\newcommand{\br}{\boldsymbol{r}}

\newcommand{\bv}{\boldsymbol{v}}
\newcommand{\bw}{\boldsymbol{w}}
\newcommand{\bx}{\boldsymbol{x}}

\newcommand{\bbR}{\mathbb{R}}

\newcommand{\bbC}{\mathbb{C}}
\newcommand{\bLambda}{\boldsymbol{\Lambda}}
\newcommand{\Prob}{\mathbb{P}}
\newcommand{\Expe}{\mathbb{E}}
\newcommand{\Var}{\mathrm{Var}}
\newcommand{\Cov}{\mathrm{Cov}}

\newcommand{\tr}{\mathrm{tr}}

\newcommand{\diag}{\mathrm{diag}}

\newcommand*{\dif}{\mathop{}\!\mathrm{d}}
\renewcommand{\bar}{\overline}
\usepackage{dsfont}
\newcommand{\indicator}{\mathds{1}}
\newcommand{\convas}{\overset{a.s.}{\longrightarrow}}

\newcommand{\convd}{\overset{d}{\longrightarrow}}

\endlocaldefs

\begin{document}

\begin{frontmatter}
%%%%%%%%%%%%%%%%%%%%%%%%%%%%%%%%%%%%%%%%%%%%%%
%%                                          %%
%% Enter the title of your article here     %%
%%                                          %%
%%%%%%%%%%%%%%%%%%%%%%%%%%%%%%%%%%%%%%%%%%%%%%
\title{On eigenvalues of sample covariance matrices based on high-dimensional compositional data}
%\title{A sample article title with some additional note\thanksref{T1}}
\runtitle{High-dimensional SCM based on compositional data}
%\thankstext{T1}{A sample of additional note to the title.}

\begin{aug}
%%%%%%%%%%%%%%%%%%%%%%%%%%%%%%%%%%%%%%%%%%%%%%%
%% ORCID can be inserted by command:         %%
%% \orcid{0000-0000-0000-0000}               %%
%%%%%%%%%%%%%%%%%%%%%%%%%%%%%%%%%%%%%%%%%%%%%%%
\author[A]{\inits{Q.}\fnms{Qianqian}~\snm{Jiang}\ead[label=e1]{jiangqq@sustech.edu.cn}}
\author[B]{\inits{J.}\fnms{Jiaxin}~\snm{Qiu}\ead[label=e2]{qiujx@connect.hku.hk}}
\author[C]{\inits{Z.}\fnms{Zeng}~\snm{Li}\ead[label=e3]{liz9@sustech.edu.cn}}
%%%%%%%%%%%%%%%%%%%%%%%%%%%%%%%%%%%%%%%%%%%%%%
%% Addresses                                %%
%%%%%%%%%%%%%%%%%%%%%%%%%%%%%%%%%%%%%%%%%%%%%%
\address[A]{Southern University of Science and Technology, Shenzhen, China\printead[presep={,\ }]{e1}}
\address[B]{The University of Hong Kong, Hong Kong, China\printead[presep={,\ }]{e2}}
\address[C]{Southern University of Science and Technology, Shenzhen, China\printead[presep={,\ }]{e3}}
\end{aug}

\begin{abstract}
This paper studies the asymptotic spectral properties of the sample covariance matrix for high-dimensional compositional data, including the limiting spectral distribution, the limit of extreme eigenvalues, and the central limit theorem for linear spectral statistics. All asymptotic results are derived under the high-dimensional regime where the data dimension increases to infinity proportionally with the sample size. The findings reveal that the limiting spectral distribution is the well-known Mar\v{c}enko-Pastur law. The largest (or smallest non-zero) eigenvalue converges almost surely to the left (or right) endpoint of the limiting spectral distribution, respectively. Moreover, the linear spectral statistics demonstrate a Gaussian limit. Simulation experiments demonstrate the accuracy of theoretical results.
\end{abstract}


\end{frontmatter}

%%%%%%%%%%%%%%%%%%%%%%%%%%%%%%%%%%%%%%%%%%%%%%
%%%% Main text entry area:

%\tableofcontents

\section{Introduction}
In recent years, there has been increasing interest in the analysis of high-dimensional compositional data (HCD), which arise in various fields including genomics, ecology, finance, and social sciences. Compositional data refers to observations whose sum is a constant, such as proportions or percentages. HCD often involve a large number of variables or features measured for each sample, posing unique challenges for analysis. In the field of genomics, HCD analysis plays a crucial role in studying the composition and abundance of microbial communities, such as the human gut microbiome. Understanding the microbial composition and its relationship with health and disease has significant implications for personalized medicine and therapeutic interventions.

Statistical inference in HCD involves microbial mean tests, covariance matrix structural tests, and linear regression hypothesis testing. These inferences are intricately linked to the statistical properties of the sample covariance matrix. Mean tests typically utilize sum-of-squares-type and maximum-type statistics for dense and sparse alternative hypotheses, respectively. \cite{Cao2018TwosampleTO} extended the maximum test framework by \cite{Cai2014TwosampleTO} for compositional data. However, there's a gap in having a suitable sum-of-squares-type statistic for dense alternatives in HCD mean tests. Many sum-of-squares-type statistics, like Hotelling's $T^2$-statistic, rely on the sample covariance matrix. For bacterial species correlation, \cite{Faust2012MicrobialCR} introduced the permutation-renormalization bootstrap (ReBoot), directly calculating correlations from compositional components. Shuffling is suggested due to compositional data's closure constraint, introducing negative correlations. Yet, compositional data's unique properties require an additional normalization step within the same sample post-shuffling, potentially impacting the theoretical validity of permutation and resampling methods. Additionally, resampling increases computational complexity for p-value calculation and confidence interval construction. To address these challenges, \cite{Wu2011LinkingLD} developed a covariance matrix element hypothesis testing method, allowing control over false discovery proportion (FDP) and false discovery rate (FDR). All these studies are closely related to the sample covariance matrix of HCD.

Current research predominantly focuses on sparse compositional data. In dense scenarios, researchers often turn to the spectral properties of sample covariance matrices. Despite this, there is a notable gap in the field of random matrices where specific attention to structures resembling compositional data, where row sum of the data matrix is constant, is lacking. Statistical inference for HCD encounters challenges arising not only from constraints but also from high dimensionality. Recognizing the crucial role of spectral theory in sample covariance matrices is also vital for addressing statistical challenges associated with high-dimensional data. Importantly, while previous research on statistical inference for HCD has overlooked studies under the spectral theory of sample covariance matrices, our work takes on these challenges from a Random Matrix Theory perspective. Existing literature extensively covers spectral properties of large-dimensional sample covariance matrices, but most results rely on independent component data structure, i.e. $\bZ=\boldsymbol{\Gamma}\bX$, where $\boldsymbol{\Gamma}$ is determined, and $\bX$ has independent and identically distributed (i.i.d.) components. Seminal works by \cite{marchenko1967distribution} and \cite{jonsson1982limit} established the limiting spectral distribution (LSD) of the sample covariance matrix $n^{-1}\bX\bX^{'}$, where $\bX$ is an i.i.d. data matrix with zero mean, leading to the well-known Mar\v{c}enko-Pastur law. Subsequent research by \cite{Y2003ALT} and \cite{silverstein1995empirical} extended these findings to the sample covariance matrix $n^{-1}\bX\boldsymbol{\Sigma}\bX^{'}$ for data with a linear dependence structure. \cite{Zhang2007spectral} extended to the general separable product form $n^{-1}\bA^{1/2}\bX\bB\bX^{'}\bA^{1/2}$, where $\bA$ is nonnegative definite, and $\bB$ is Hermitian. Another important area of interest is the investigation of extreme eigenvalues. \cite{johnstone2001distribution} explored the fluctuation of the extreme eigenvalues of the sample covariance matrix $n^{-1}\bX\bX^{'}$, proving that the standardized largest eigenvalue follows the Tracy-Widom law. Related extensions include sample covariance matrices with linear dependence structures \citep{karoui2007tracy}, Kendall rank correlation coefficient matrices \citep{Bao2017TracyWidomLF}, among others. Considerable attention has also been given to the study of linear functionals of eigenvalues. \cite{bai2004clt} established the Central Limit Theorem (CLT) for the Linear Spectral Statistics (LSS) of the sample covariance matrix $n^{-1}\bA^{1/2}\bX\bX^{'}\bA^{1/2}$, later extended to sample correlation coefficient matrices \citep{gao2017high}, and separable product matrices \citep{Bai2019central}. To summarize, existing results in spectral theory of large dimensional sample covariance matrix  predominantly rely on independent component data structure which, unfortunately, HCD does not fit in.

Specifically, current second-order limit theorems do not  apply to HCD, making the exploration of spectral theory for HCD with distinct constraints crucial. This paper delves into spectral theory for  sample covariance matrices of HCD, including LSD, extreme eigenvalues, and CLT for LSS. Analyzing HCD faces challenges due to compositional data's specific dependence structure, making existing techniques for i.i.d. observations less applicable. However, we can assume that HCD are generated from unobservable basis data, while the underlying basis data follow independent component model structure. In this way, spectral analysis of the sample covariance matrix of HCD can be approached through the basis data. In fact, the structure of the sample covariance matrix of HCD is similar to that of the Pearson sample correlation matrix in basis data. Therefore, we leverage the analysis methods of the spectral theory of the Pearson sample correlation matrix to study the spectral theory of the sample covariance matrix of HCD. In the field of random matrices, research on the spectral theory of the Pearson sample correlation matrix based on independent data is relatively mature. \cite{jiang2004limiting} demonstrated that the LSD of sample correlation matrix for i.i.d data is the well-known Mar\v{c}enko-Pastur law. \cite{gao2017high} derive the CLT for LSS of the Pearson sample correlation matrix. The derivation of spectral theory for the sample covariance matrix of HCD can benefit from methods in this context. The LSD of the sample covariance matrix for HCD in Theorem \ref{thm:CoDa-LSD} is established following the strategy in \cite{jiang2004limiting}, and we further investigate the extreme eigenvalues in Proposition \ref{extreigen}. The proof strategy of CLT for LSS in Theorem \ref{clt} follows the methodologies outlined in \cite{bai2004clt} for the sample covariance matrix and \cite{gao2017high} for the sample correlation matrix. However, due to the dependence inherent in HCD, certain tools from these works cannot be directly applied to the sample covariance matrix of HCD. In response, we introduce new techniques. Specifically, we establish concentration inequalities for compositional data. One of the central ideas of the paper, grounded in concentration phenomena, permeates the entire proof (details in Section \ref{importlem} and Section \ref{cltcenBp0}), where we develop three crucial technique lemmas (see Lemmas \ref{estsqu} - \ref{quadform}) essential for the proof. Finally, it is noteworthy that the mean and variance-covariance in Theorem \ref{clt} differ from those in \cite{bai2004clt}, and additional terms are present in both the mean and variance-covariance.

The paper is organized as follows. Section \ref{LSDExei} investigates the LSD and extreme eigenvalues of the sample covariance matrix for HCD. Section \ref{mrclt} establishes our main CLT for LSS of the sample covariance matrix for HCD. Section \ref{numexp} reports numerical studies. Technical proofs and
lemmas are relegated to Section \ref{prfcltall} and the supplementary document.

Before moving forward, let us introduce some notations that will be used throughout this paper. We adopt the convention of using regular letters for scalars and using bold-face letters for vectors or matrices.
For any matrix $\bA$, we denote its $(i,j)$-th entry by $A_{ij}$, its transpose by $\bA'$, its trace by $\tr(\bA)$, its $j$-th largest eigenvalue by $\lambda_j(\bA)$, its spectral norm by $\|\bA\|=\sqrt{\lambda_{1}(\bA\bA')}$. For a set of random variables $\{X_n\}_{n=1}^{\infty}$ and a corresponding set of nonnegative real numbers $\{a_n\}_{n=1}^{\infty}$, we write $X_n=O_P(a_n)$ if for any $\varepsilon>0$, there exists a constant $C>0$ and $N>0$ such that $\Prob(|X_n/a_n|\geq  C)\leq   \varepsilon$ holds for all $n\geq  N$; and we write $X_n=o_P(a_n)$ if $\lim_{n\to\infty} \Prob(|X_n/a_n|\geq  \varepsilon ) = 0$ holds for any $\varepsilon >0$; and we write $X_n\stackrel{a.s.}{\to}a$ ($X_n\stackrel{i.p.}{\to}a$, resp.) if $X_n$ converges almost surely (in probability, resp.) to $a$. We denote by $C$ and $K$ are constants, which may be different from line to line.

\section{Main Results}
\subsection{Preliminaries and Notations}
Let  $\bX_n =(\bx_1,\ldots,\bx_n)'$ denote the $n\times p$ observed data matrix, where each $\bx_i$
 represents compositions that lie in the $(p-1)$-dimensional simplex $\mathcal{S}^{p-1} =
\{(y_1, \ldots, y_p) :\, \sum_{j=1}^{p}y_j = 1,\, y_j \geq  0\}$. We assume that the compositional variables arise from a vector of latent variables, which we call the basis. Let
 $\bW_n = (w_{ij})_{n\times p}$ denote the $n\times p$ matrices of unobserved bases, where $w_{ij}$'s are positive and i.i.d. with mean $\mu>0$ and variance $\sigma^2$. The
observed compositional data is generated via the normalization
\begin{equation}
    x_{ij} = \frac{w_{ij}}{\sum_{\ell=1}^{p}w_{i\ell}},\qquad 1\leq   i \leq   n,\; 1\leq   j \leq   p.
\end{equation}
% The centered version of $\bX_n$ is given by $\bC_n\bX_n $, where $\bC_n = \bI_n-n^{-1}\boldsymbol{1}_n\boldsymbol{1}_n'$.
% The MLE (maximum likelihood estimate) sample covariance matrix of the compositional data $\bX_n$ is defined by
% \begin{equation*}
% 	\bS_n = \frac{1}{n}(\bC_n\bX_n)'(\bC_n\bX_n) = \frac{1}{n}\bX_n'\bC_n\bX_n .
% \end{equation*}
% The centralized sample covariance matrix of the compositional data $\bX_n$ is defined by
% \begin{equation*}
% 	\bS_n^0 = \frac{1}{n}(\bX_n-\Expe\bX_n)'(\bX_n-\Expe\bX_n).
% \end{equation*}
The unbiased sample covariance matrix of $\bX_n$ is defined by
$
	\bS_{n,N} = \frac{1}{N}(\bC_n\bX_n)'(\bC_n\bX_n),
$
where $\bC_n = \bI_n-(1/n)\boldsymbol{1}_n\boldsymbol{1}_n'$, $\boldsymbol{1}_n$ is a $n$-dimensional vector of all ones, and $N=n-1$ is the adjusted sample size.
We rescale $\bS_{n,N}$ as
\[
 \bB_{p,N}=p^2\bS_{n,N} =\frac{1}{N}(p\bX_n)'\bC_n(p\bX_n).
\]
For any $p \times p$ Hermitian matrix $\bB_p$ with eigenvalues $\lambda_1,\ldots,\lambda_p$, its \emph{empirical spectral distribution} (ESD) is defined by
\begin{equation}
        F^{\bB_p}(x)=\frac{1}{p}\sum^{p}_{i=1}I_{\{\lambda_i (\bB_p)\leq   x\}},
\end{equation}
where $I_{\{\cdot\}}$ denotes the indicator function. If $F^{\bB_p}(x)$ converges to a non-random limit $F(x)$ as $p\rightarrow\infty$, we call $F (x)$ the \emph{limiting spectral distribution} of $\bB_p$. The LSD of $\bB_p$ is described in terms of its Stieltjes transform. The Stieltjes transform of any cumulative distribution function $G$ is defined by
\begin{align}
    m_G(z)=\int\frac{1}{\lambda-z}\dif G(\lambda),\qquad z \in \mathbb{C}^+ := \{z: \Im(z)>0 \}.
\end{align}
Many classes of statistics related to the eigenvalues of the sample covariance matrix $\bB_{p,N}$ are important for multivariate inference, particularly functionals of the ESD. To explore this, for any function $f$ defined on $[0,\infty)$, we consider the \emph{linear spectral statistics} of $\bB_{p,N}$ given by
\begin{align}
    \int f(x)\dif F^{\bB_{p,N}}(x)=\frac{1}{p}\sum_{i=1}^p f\bigl(\lambda_i\bigr),
\end{align}
where $\lambda_i$, $i=1,\ldots,p$, are eigenvalues of $\bB_{p,N}$.

In this paper, we study the asymptotic spectral properties of $\bB_{p,N}$, including the LSD (see, Theorem \ref{thm:CoDa-LSD}), the behavior of extreme eigenvalues (see, Proposition \ref{extreigen}), and the CLT for LSS (see, Theorem \ref{clt}). 
%  Define $G_{p,N}(f)=p\int f(x)\dif \{F^{\bB_{p,N}}(x)- F^{c_N}(x)\}$,
% where $F^{c_N}(x)$ substitutes $c_N$ for $c$ in $F^{c}(x)$, the LSD of $\bB_{p,N}$. We prove that under the assumption  $\Expe\left|w_{11}-\mu\right|^{4}<\infty$ and analyticity of $f$, the rate $\int f(x)\dif \{F^{\bB_{p,N}}(x)- F^{c_N}(x)\}$, where $c_N=p/N$, approaches zero is essentially $1/n$ and $G_{p,N}(f)$ convergence weakly to a Gaussian vector, see Theorem \ref{clt}.  
% We studies the asymptotic spectral properties of $\bB_{p,N}$, including the LSD, the behavior of extreme eigenvalues, the CLT for LSS. We prove that under the assumption  $\Expe\left|w_{11}-\mu\right|^{4}<\infty$ and analyticity of $f$, the rate $\int f(x)\dif \{F^{\bB_{p,N}}(x)- F^{c_N}(x)\}$, where $c_N=p/N$, approaches zero is essentially $1/n$ and $G_{p,N}(f)$ convergence weakly to a Gaussian vector, see Theorem \ref{clt}.  

% Define $G_{p,N}(x)=p\left[F^{\bB_{p,N}}(x)-F^{c_N}(x)\right]$,
% where $F^{c_N}(x)$ denotes the function by substituting $c_N$ for $c$ in $F^{c}(x)$, the LSD of $\bB_{p,N}$. 

\subsection{Limiting spectral distribution and Extreme eigenvalues}\label{LSDExei}

Analyzing HCD poses challenges due to its unique dependence structure, making existing techniques for i.i.d. observations less applicable. To overcome this difficulty, we assume that the compositional data is generated from basis data and the basis data follows the commonly used independent component structure. Specifically, the unbiased sample covariance matrix of $\bX_n$ is defined by
$$
\bS_{n,N}=\frac{1}{N} \bX_n^{'} \bC_n \bX_n=\frac{1}{N}\bW_n^{'}\boldsymbol{\Lambda}_n\bC_n\boldsymbol{\Lambda}_n\bW_n, $$
% =\frac{1}{N}\left(\begin{array}{ccc}
% \frac{w_{11}}{\sum_{j=1}^p w_{1 j}} & \cdots & \frac{w_{1 p}}{\sum_{j=1}^p w_{1 j}} \\
% \vdots & \ddots & \vdots \\
% \frac{w_{n 1}}{\sum_{j=1}^p w_{n j}} & \cdots & \frac{w_{n p}}{\sum_{j=1}^p w_{n j}}
% \end{array}\right)^{'} \mathbf{C}_n\left(\begin{array}{ccc}
% \frac{w_{11}}{\sum_{j=1}^p w_{1 j}} & \cdots & \frac{w_{1 p}}{\sum_{j=1}^p w_{1 j}} \\
% \vdots & \ddots & \vdots \\
% \frac{w_{n 1}}{\sum_{j=1}^p w_{n j}} & \cdots & \frac{w_{n p}}{\sum_{j=1}^p w_{n j}}
% \end{array}\right),
where 
$$
\bX_n =\left(\begin{array}{ccc}
\frac{1}{\sum_{j=1}^p w_{1 j}} & \cdots & 0 \\
\vdots & \ddots & \vdots \\
0 & \cdots & \frac{1}{\sum_{j=1}^p w_{n j}}
\end{array}\right)_{n \times n}\left(\begin{array}{ccc}
w_{11} & \cdots & w_{1 p} \\
\vdots & \ddots & \vdots \\
w_{n 1} & \cdots & w_{n p}
\end{array}\right)_{n \times p}:=\boldsymbol{\Lambda}_n\bW_n.
$$
Here we assume $\bW_n$ has i.i.d. components $w_{ij}$ satisfying $\Expe(w_{ij})=\mu>0$, $\Var(w_{ij})=\sigma^2$. Recall that the Pearson sample correlation matrix for $\bW_n$ expressed as
$$
\bR_n = \frac{1}{n}\widetilde{\mathbf{X}}_n^{'} \mathbf{C}_n \widetilde{\mathbf{X}}_n=\frac{1}{n}\widetilde{\boldsymbol{\Lambda}}_p\bW_n^{'}\mathbf{C}_n \bW_n\widetilde{\boldsymbol{\Lambda}}_p,
% =\frac{1}{n}\left(\begin{array}{ccc}
% \frac{w_{11}}{\left\|\boldsymbol{w}_1\right\|_2} & \cdots & \frac{w_{1 p}}{\left\|\boldsymbol{w}_p\right\|_2} \\
% \vdots & \ddots & \vdots \\
% \frac{w_{n 1}}{\left\|\boldsymbol{w}_1\right\|_2} & \cdots & \frac{w_{n p}}{\left\|\boldsymbol{w}_p\right\|_2}
% \end{array}\right)^{'} \quad \mathbf{C}_n\left(\begin{array}{ccc}
% \frac{w_{11}}{\left\|\boldsymbol{w}_1\right\|_2} & \cdots & \frac{w_{1 p}}{\left\|\boldsymbol{w}_p\right\|_2} \\
% \vdots & \ddots & \vdots \\
% \frac{w_{n 1}}{\left\|\boldsymbol{w}_1\right\|_2} & \cdots & \frac{w_{n p}}{\left\|\boldsymbol{w}_p\right\|_2}
% \end{array}\right),
$$
where $\left\|\boldsymbol{w}_j\right\|_2=\left\{(1 / n) \sum_{i=1}^n\left(w_{i j}-\bar{\bar{w}}_j\right)^2\right\}^{1 / 2}, \bar{\bar{w}}_j=(1 / n) \sum_{i=1}^n w_{i j}, j=1, \cdots, p,$ and 
$$
\begin{gathered}
\widetilde{\mathbf{X}}_{n}=\left(\begin{array}{ccc}
w_{11} & \cdots & w_{1 p} \\
\vdots & \ddots & \vdots \\
w_{n 1} & \cdots & w_{n p}
\end{array}\right)_{n \times p}\left(\begin{array}{ccc}
\left\|\boldsymbol{w}_1\right\|_2^{-1} & \cdots & 0 \\
\vdots & \ddots & \vdots \\
0 & \cdots & \left\|\boldsymbol{w}_p\right\|_2^{-1}
\end{array}\right)_{p \times p}:= \bW_n\widetilde{\boldsymbol{\Lambda}}_p.
\end{gathered}
$$
It can be seen that the normalizing matrix $\boldsymbol{\Lambda}_n$ of $\bS_{n,N}$ is very similar to $\widetilde{\boldsymbol{\Lambda}}_p$ of $\bR_n$. The former uses $(\sum_{j=1}^pw_{ij})^{-1}$ for normalization, while the latter utilizes $\left\|\boldsymbol{w}_j\right\|_2^{-1}$. This allows us to leverage the techniques from the spectral theory of the Pearson sample correlation matrix in studying the asymptotic spectral properties of the sample covariance matrix for HCD.
% In the case of $\bX_n$, we employ $(\sum_{j=1}^pw_{ij})^{-1}$ for regularization, while for $\bW_n$, we use $\left\|\boldsymbol{w}_j\right\|_2^{-1}$. Hence, we can draw on the analysis methods of the spectral theory of the Pearson sample correlation matrix to study the asymptotic spectral theory of the sample covariance matrix for HCD.

Before diving into linear functionals of eigenvalues of $\bB_{p,N}$, we first explore its LSD and extreme eigenvalues. Specifically, suppose the following assumptions hold,
\begin{Assumption}\label{asforw}
       $\{w_{ij}>0,i=1,\ldots,n,j=1,\ldots,p\}$ are i.i.d. real random variables with  $\Expe w_{11}=\mu>0$, $\Expe (w_{11}-\mu)^2=\sigma^2$ and  $\Expe|w_{11}-\mu|^4<\infty$.
\end{Assumption}
\begin{Assumption}\label{asforcn}
    $c_N=p/N$ tends to a positive $c>0$ as $p,N\rightarrow\infty.$
\end{Assumption}
% In this section, let
% \begin{align*}  \bB_p&=p^2\bS_n, \ \ \bS_n=\frac{1}{n}\bX_n'\bC_n\bX_n,
% \end{align*}
% where $\bS_n$ is the ME (moment estimator) sample covariance matrix of the compositional data $\bX_n$.
% We prove that as $N\to \infty, p=p(N)\to\infty$ and $p/N=c_N\to c \in (0,\infty)$, $F^{\bB_{p,N}}$ almost surely converges to the standard Mar\v{c}enko-Pastur (MP) law, as stated in Theorem \ref{thm:CoDa-LSD}. The proof of Theorem \ref{thm:CoDa-LSD} and Proposition \ref{extreigen} are postponed to Sections \ref{prfLSD}-\ref{prfextreigen}.
\begin{theorem}\label{thm:CoDa-LSD}
	Under Assumptions \ref{asforw} and \ref{asforcn}, with probability one, the ESD of $\bB_{p,N}$ converges weakly to a deterministic probability distribution with a density function
	\begin{equation}\label{eq:MP_density}
	f(x) = \begin{cases}
		\frac{\mu^2}{2\pi c \sigma^2 x} \sqrt{(b-x)(x-a)},&\text{if }x\in[a,b],\\
		0,& \text{otherwise},
	\end{cases}
	\end{equation}
	and a point mass $1-1/c$ at $x=0$ if $c>1$, where $a=\tfrac{\sigma^2}{\mu^2}(1-\sqrt{c})^2$ and $b=\tfrac{\sigma^2}{\mu^2}(1+\sqrt{c})^2$.
\end{theorem}

The proof of Theorem \ref{thm:CoDa-LSD} is postponed to the supplementary file in Section \ref{prfLSD}. 

The LSD $F^{c}(x)$ has a Dirac mass $1-1/c$
at the origin when $c > 1$. We see that $m(\bar{z})=\bar{m(z)}$. For each $z\in \mathbb{C}^+=\{z: \Im(z)>0 \}$, by Theorem \ref{thm:CoDa-LSD} the Stieltjes transform $m(z)=m_{F^c}(z)$ is the unique solution of $m=\frac{1}{\sigma^2/\mu^2(1-c-czm)-z}$
% \begin{align*}
%     m=\frac{1}{\sigma^2/\mu^2(1-c-czm)-z},
% \end{align*}
in the set $\{m\in\mathbb{C}: \frac{1-c}{z}+\underline{m}(z)\in\mathbb{C}^+\}$.
Define $\underline{m}(z)$ to be the Stieltjes transform of the companion
LSD $\underline{F}^{c}(x)=(1-c)\delta_0+cF^{c}(x)$,
% \begin{align*}
%     \underline{F}^{c}(x)=(1-c)\delta_0+cF^{c}(x)
% \end{align*}
where $\delta_0$ is the point distribution at zero. Then $\underline{m}(z)$ is the unique solution in $\{\underline{m}\in\bbC:\frac{1-c}{z}+\underline{m}(z)\in\mathbb{C}^+\}$
of the equation 
\begin{align}
     z &= -\frac{1}{\underline{m}(z)}+\frac{c\sigma^2/\mu^2}{1+\sigma^2/\mu^2\underline{m}(z)},\ \ \ \ z\in \mathbb{C}^+ \label{zmbar}.
\end{align}
% Note that $\rank(\bS_n)\leq  \min\{\rank(\bX),\rank(\bC_n)\}\leq   \min\{p-1,n-1\}$. If $p >  n-1$, the $p-n+1$ smallest eigenvalues of $\bS_n$ are all zero. 
% We denote the smallest eigenvalue of the matrix $\bB_p$ as $\lambda_{\min}(\bB_p)$.
% \begin{equation}\label{eq:min-eig}
% 	\lambda_{\min}(\bB_p) =\begin{cases}
% 		\lambda_{p-1}(p^2\bS_n),&\text{if $p\leq   n-1$};\\
% 		\lambda_{n-1}(p^2\bS_n),&\text{if $p>n-1$}.
% 	\end{cases}
% \end{equation}
\begin{Proposition}\label{extreigen}
    Under Assumptions \ref{asforw} and \ref{asforcn}, we have
    \begin{equation}\label{eq:extreme-eig-CoDa}
		\lambda_{\max}(\bB_{p,N}) \convas \frac{\sigma^2}{\mu^2}(1+\sqrt{c})^2 \qquad \text{and}\qquad \lambda_{\min} (\bB_{p,N}) \convas \frac{\sigma^2}{\mu^2}(1-\sqrt{c})^2,
	\end{equation}
 where $\lambda_{\max}(\bB_{p,N})$ is the largest eigenvalue of $\bB_{p,N}$, and $\lambda_{\min}(\bB_{p,N})$ is the smallest non-zero eigenvalue of $\bB_{p,N}$. Furthermore, for any $\ell>0$, $\eta_1>\frac{\sigma^2}{\mu^2}(1+\sqrt{c})^2$ and $0<\eta_2<\frac{\sigma^2}{\mu^2}(1-\sqrt{c})^2\cdot I_{\{0<c<1\}}$, we have
	\[
		\Prob\bigl(\lambda_{\max}(\bB_{p,N})\geq  \eta_1\bigr) = o(n^{-\ell})
	\qquad\text{and}\qquad
		\Prob\bigl(\lambda_{\min}(\bB_{p,N})\leq   \eta_2\bigr) = o(n^{-\ell}).
	\]
\end{Proposition}

% \begin{lemma}\label{eigen}
% 	Under conditions and notations of Theorem \ref{thm:CoDa-LSD}, for any $\ell>0$, $\eta_1>\frac{\sigma^2}{\mu^2}(1+\sqrt{c})^2$ and $0<\eta_2<\indicator_{(0,1)}(c)\cdot\frac{\sigma^2}{\mu^2}(1-\sqrt{c})^2$, we have
% 	\[
% 		\Prob\bigl(\lambda_{\max}(\bB_{p,N})\geq  \eta_1\bigr) = o(n^{-\ell})
% 	\qquad\text{and}\qquad
% 		\Prob\bigl(\lambda_{\min}(\bB_{p,N})\leq   \eta_2\bigr) = o(n^{-\ell}).
% 	\]
% \end{lemma}
% \begin{align}\label{suppofLSD}
%     The \ LSD \ has \ support \ [a,b]=\left[\frac{\sigma^2}{\mu^2}(1-\sqrt{c})^2,\frac{\sigma^2}{\mu^2}(1+\sqrt{c})^2\right],
% \end{align}
The proof of Proposition \ref{extreigen} is postponed to the supplementary file in Section \ref{prfextreigen}.
\begin{remark}
The  LSD  has  support $\left[\frac{\sigma^2}{\mu^2}(1-\sqrt{c})^2,\frac{\sigma^2}{\mu^2}(1+\sqrt{c})^2\right]$, where it has a density function. The results of extreme eigenvalues find application in locating eigenvalues of the population covariance matrix and in proving the CLT for LSS. Proposition \ref{extreigen} shows that with probability 1, there are no eigenvalues of $\bB_{p,N}$ outside the support of LSD under Assumptions \ref{asforw}-\ref{asforcn}. These lemmas are crucial for applying the Cauchy integral formula (see, equation \eqref{chucform}) and proving tightness.
\end{remark}
\subsection{CLT for LSS}\label{mrclt}
We focus on linear functionals of eigenvalues of $\bB_{p,N}$, i.e. $\frac{1}{p}\sum_{i=1}^pf(\lambda_i)$. Naturally it converges to the functional integration of LSD of $\bB_{p,N}$, i.e. $\int f(x)\dif F^c(x)$. In this section, we explore second order fluctuation of $\frac{1}{p}\sum_{i=1}^pf(\lambda_i)$ describing how such LSS converges to its first order limit. Define 
$$G_{p,N}(f)=p\int f(x)\dif \{F^{\bB_{p,N}}(x)- F^{c_N}(x)\}$$
where $F^{c_N}(x)$ substitutes $c_N$ for $c$ in $F^{c}(x)$, the LSD of $\bB_{p,N}$. We show that under Assumptions \ref{asforw} -- \ref{asforcn} and the analyticity of $f$, the rate
% $$ G_{p,N}(f):= \int f(x)\dif G_{p,N}(x)=\sum_{i=1}^p f(\lambda_i)-p\int f(x)\dif F^{c_N}(x)$$ 
$\int f(x)\dif \{F^{\bB_{p,N}}(x)- F^{c_N}(x)\}$, 
approaching zero is essentially $1/n$ and $G_{p,N}(f)$ convergence weakly to a Gaussian variable. Before presenting the main result, we first recall some notation. Let
% $$G_{p,N}(f):=p\left[\int f(x)\dif F^{\bB_{p,N}}(x)-\int f(x)\dif F^{c_N}(x)\right]$$
$m(z)$ be the Stieltjes transform of the LSD $F^{c}(x)$ and $\underline{m}(z)$ be the Stieltjes transform of the companion LSD $\underline{F}^{c}(x)$. Furthermore, we define $m'(z)$ as the first derivative of $m(z)$ with respect to $z$ throughout the rest of this paper. The main result is stated in the following theorem.
% Note that when a finite-horizon proxy $F^{c_N}(x)$ is substituted for $F^c(x)$,
% the properties and relationships in Remark \ref{remarkLSD} hold with the parameters $c$ replaced by
% $c_N$. The main result is formulated in the theorem below.
\begin{theorem}\label{clt}
   Under Assumptions \ref{asforw} and \ref{asforcn}, let $f_1, f_2, \ldots, f_k$ be functions on $\mathbb{R}$ and analytic on an open interval containing
\begin{equation}\label{suppofLSD}
    \Bigl[\frac{\sigma^2}{\mu^2}(1-\sqrt{c})^2, \; \frac{\sigma^2}{\mu^2}(1+\sqrt{c})^2\Bigr].
\end{equation}
Then, the random vector $\left( G_{p,N}(f_1), \ldots, G_{p,N}(f_k)\right)$ forms a tight sequence in $p$ and converges weakly to a Gaussian vector $(X_{f_1}, \ldots, X_{f_k})$ with mean function
\begin{align}
    \Expe X_f %=& -\frac{1}{2\pi i}\oint_{\mathcal{C}} f(z)\Big(-\underline{m}(z)\left[1-c\frac{\sigma^{4}}{\mu^{4}}\underline{m}^{2}(z)(1+\frac{\sigma^{2}}{\mu^{2}}\underline{m}(z))^{-2}\right]^{-1}\\
    %&\times \biggl[\left(-z\underline{m}(z)+c\frac{\sigma^{2}}{\mu^{2}}z^{2}m(z)\underline{m}(z)^{2}\right)\left(1+\frac{\sigma^{2}}{\mu^{2}}\underline{m}(z)\right)^{-1} \times \left(h_{1}m(z)+\frac{\sigma^{2}}{\mu^{2}}m(z)+\frac{\sigma^{2}}{\mu^{2}}\frac{1}{z}\right)\\
    %&\quad\quad \times (h_{1}m(z)+\frac{\sigma^{2}}{\mu^{2}}m(z)+\frac{\sigma^{2}}{\mu^{2}}\frac{1}{z}) \\
    %&\qquad-cz^{2}\underline{m}^{2}(z)\left(1+\frac{\sigma^{2}}{\mu^{2}}\underline{m}(z)\right)^{-1}\left(h_{2}m^{2}(z)-2\frac{\sigma^{2}}{\mu^{2}}h_{1}m^{2}(z)+2\frac{\sigma^{4}}{\mu^{4}}m'(z)\right)\\
    %&\qquad+c\frac{\sigma^{4}}{\mu^{4}}\underline{m}^{2}(z)\left(1+\frac{\sigma^{2}}{\mu^{2}}\underline{m}(z)\right)^{-3}\left(1-c\frac{\sigma^{4}}{\mu^{4}}\underline{m}^{2}(z)(1+\frac{\sigma^{2}}{\mu^{2}}\underline{m}(z))^{-2}\right)^{-1}\biggr]\Big)\dif z\nonumber\\
    &= \frac{1}{2\pi i}\oint_{\mathcal{C}} c\frac{\sigma^{4}}{\mu^{4}}f(z)\underline{m}^{3}(z)\left\{1+\frac{\sigma^{2}}{\mu^{2}}\underline{m}(z)\right\}^{-3}\left[1-c\frac{\sigma^{4}}{\mu^{4}}\underline{m}^{2}(z)\left\{1+\frac{\sigma^{2}}{\mu^{2}}\underline{m}(z)\right\}^{-2}\right]^{-2}\label{cltmean}\\
    % &\qquad \times\left(1-c\frac{\sigma^{4}}{\mu^{4}}\underline{m}^{2}(z)(1+\frac{\sigma^{2}}{\mu^{2}}\underline{m}(z))^{-2}\right)^{-1}\dif z\nonumber\\
    &\quad - \frac{1}{2\pi i}\oint_{\mathcal{C}} f(z)\underline{m}(z)\left[1-c\frac{\sigma^{4}}{\mu^{4}}\underline{m}^{2}(z)\left\{1+\frac{\sigma^{2}}{\mu^{2}}\underline{m}(z)\right\}^{-2}\right]^{-1}\nonumber\\
    &\qquad\qquad\qquad\times z\underline{m}(z)\left\{1+\frac{\sigma^{2}}{\mu^{2}}\underline{m}(z)\right\}^{-1} \times \left\{h_1 m(z)+\frac{\sigma^{2}}{\mu^{2}}m(z)+\frac{\sigma^{2}}{\mu^{2}}\frac{1}{z}\right\}\dif z\nonumber\\
    %&\quad\quad \times (h_{1}m(z)+\frac{\sigma^{2}}{\mu^{2}}m(z)+\frac{\sigma^{2}}{\mu^{2}}\frac{1}{z}) \\
    &\quad-\frac{1}{2\pi i}\oint_{\mathcal{C}} cf(z)z^{2}\underline{m}^{3}(z)\left[1-c\frac{\sigma^{4}}{\mu^{4}}\underline{m}^{2}(z)\left\{1+\frac{\sigma^{2}}{\mu^{2}}\underline{m}(z)\right\}^{-2}\right]^{-1}\nonumber\\
    &\qquad\qquad\qquad\times\left\{1+\frac{\sigma^{2}}{\mu^{2}}\underline{m}(z)\right\}^{-1}\left\{(\alpha_1+\alpha_2)m^{2}(z)+2\frac{\sigma^{4}}{\mu^{4}}m'(z)\right\}\dif z,\nonumber
    %&\qquad+c\frac{\sigma^{4}}{\mu^{4}}\underline{m}^{2}(z)\left(1+\frac{\sigma^{2}}{\mu^{2}}\underline{m}(z)\right)^{-3}\left(1-c\frac{\sigma^{4}}{\mu^{4}}\underline{m}^{2}(z)(1+\frac{\sigma^{2}}{\mu^{2}}\underline{m}(z))^{-2}\right)^{-1}\biggr]\Big)\dif z
\end{align}
% and covariance function
% \begin{align}
%     \Cov(X_f, X_g) &= -\frac{1}{2\pi^2} \oint_{\mathcal{C}}\oint_{\mathcal{C}} \frac{f(z_1)g(z_2)}{[\underline{m}(z_1)-\underline{m}(z_2)]^2}\frac{\dif}{\dif z_1}\underline{m}(z_1)\frac{\dif}{\dif z_2} \underline{m}(z_2)\dif z_1\dif z_2\nonumber\\
%             &\quad -\frac{c}{4\pi^2}\left( \nu_4-3\nu_{12}+p(v_{12}-v_2^2)  \right)\nonumber\\
%             &\qquad\times\oint_{\mathcal{C}}\oint_{\mathcal{C}} f(z_1)g(z_2)\frac{1}{(1+\sigma^2/\mu^2\underline{m}(z_{1}))^{2}(1+\sigma^2/\mu^2\underline{m}(z_{2}))^{2}}\frac{\dif}{\dif z_1}\underline{m}(z_1)\frac{\dif}{\dif z_2} \underline{m}(z_2)\dif z_1\dif z_2.
% \end{align}
and covariance function
\begin{align}
    \Cov(X_f, X_g) &= -\frac{1}{2\pi^2} \oint_{\mathcal{C}_1}\oint_{\mathcal{C}_2} \frac{f(z_1)g(z_2)}{\left\{\underline{m}(z_1)-\underline{m}(z_2)\right\}^2}\dif \underline{m}(z_1)\dif \underline{m}(z_2)\label{cltcov}\\
            &\quad -\frac{c \left(\alpha_1+\alpha_2  \right)}{4\pi^2}\oint_{\mathcal{C}_1}\oint_{\mathcal{C}_2} \frac{f(z_1)g(z_2)}{\left\{1+\frac{\sigma^2}{\mu^2}\underline{m}(z_{1})\right\}^{2}\left\{1+\frac{\sigma^2}{\mu^2}\underline{m}(z_{2})\right\}^{2}}\dif \underline{m}(z_1)\dif \underline{m}(z_2),\nonumber
\end{align}
where 
\begin{align*}
   \alpha_1&= \lim_{p\rightarrow\infty}\left[\Expe\Bigl(\frac{w_{11}}{\bar{w}_1}-1\Bigr)^4 - 3\Expe\Bigl(\frac{w_{11}}{\bar{w}_1}-1\Bigr)^2\Bigl(\frac{w_{12}}{\bar{w}_1}-1\Bigr)^2\right],\\
   % = \frac{\Expe (w_{11}-\mu)^{4}}{\mu^{4}}-3\frac{(\Expe (w_{11}-\mu)^{2})^{2}}{\mu^{4}},\label{alpha1expr}\\
    % = \frac{\Expe (w_{11}-\mu)^{4}}{\mu^{4}}-3\frac{(\Expe (w_{11}-\mu)^{2})^{2}}{\mu^{4}},\ \
    %\label{alpha1expr}\\
    \alpha_2&= \lim_{p\rightarrow\infty} p\left[\Expe\Bigl(\frac{w_{11}}{\bar{w}_1}-1\Bigr)^2\Bigl(\frac{w_{12}}{\bar{w}_1}-1\Bigr)^2-\left\{\Expe\Bigl(\frac{w_{11}}{\overline{w}_1}-1\Bigr)^{2}\right\}^2\right],\\
    % =h_2-2\frac{\sigma^2}{\mu^2}h_1.\label{alpha2expr}\\
    % =h_2-2\frac{\sigma^2}{\mu^2}h_1,\label{alpha2expr}\\
    h_1&=\lim_{p\rightarrow\infty} p\left[\Expe\Bigl(\frac{w_{11}}{\overline{w}_1}-1\Bigr)^{2}-\frac{\sigma^2}{\mu^2}\right],
    % =-2\frac{\Expe w_{11}^3}{\mu^3}+3\Bigl(\frac{\sigma^2}{\mu^2}\Bigr)^2+5\frac{\sigma^2}{\mu^2}+2,\label{h1expr}\\
    % =-2\frac{\Expe w_{11}^3}{\mu^3}+3\Bigl(\frac{\sigma^2}{\mu^2}\Bigr)^2+5\frac{\sigma^2}{\mu^2}+2,\label{h1expr}
    % h_2&= -8\frac{\sigma^2}{\mu^2}\frac{\Expe w_{11}^3}{\mu^3}+10\Bigl(\frac{\sigma^2}{\mu^2}\Bigr)^3+22\Bigl(\frac{\sigma^2}{\mu^2}\Bigr)^2+8\frac{\sigma^2}{\mu^2}.\label{h2expr}
\end{align*}
and $\overline{w}_1=\sum_{j=1}^{p}w_{1j}/p$. The contours $\mathcal{C}, \mathcal{C}_1, \mathcal{C}_2$ in \eqref{cltmean} and \eqref{cltcov} are closed and taken in the positive direction in the complex plane, each enclosing the support of $F^{c}(x)$, i.e., $[\tfrac{\sigma^2}{\mu^2}(1-\sqrt{c})^2, \tfrac{\sigma^2}{\mu^2}(1+\sqrt{c})^2]$.
% where 
% \begin{align*}  
% %\quad  h_{1} = \frac{3\Expe w_{11}^{2}\Expe(w_{11}-\mu)^{2}}{\mu^{4}}-\frac{2\Expe w_{11}^{2}(w_{11}-\mu)}{\mu^{3}},\\
%    \beta_1 &= \frac{\Expe (w_{11}-\mu)^{4}}{\mu^{4}}-3\frac{(\Expe (w_{11}-\mu)^{2})^{2}}{\mu^{4}},
%    \nonumber\\
%    \beta_2&= h_2-2\frac{\sigma^2}{\mu^2}h_1=-4\frac{\Expe (w_{11}-\mu)^{3}\Expe (w_{11}-\mu)^{2}}{\mu^{5}}+4\frac{(\Expe (w_{11}-\mu)^{2})^{3}}{\mu^{6}}\nonumber\\
%    &\ \ \ \ \ \ \ \ \ \ \ \  \ \ \ \ \ \ \ \ \ \ +12\frac{\Expe (w_{11}-\mu)^{2}\Expe (w_{11}-\mu)}{\mu^{3}}-4,\nonumber\\
%    h_1&=-2\frac{\Expe(w_1-\mu)^3}{\mu^3}+3\frac{(\Expe(w_1-\mu)^2)^2}{\mu^4}-\frac{\Expe(w_1-\mu)^2}{\mu^2}-6\frac{\Expe(w_1-\mu)}{\mu}.\nonumber
%  %&h_{1}=\lim_{p\longrightarrow\infty}p(\sigma_{p}^{2}-\frac{\sigma^2}{\mu^2}),\ \ \ \     \beta_2= h_2-2\frac{\sigma^2}{\mu^2}h_1=\lim_{p\rightarrow\infty} p(v_{12}-v_2^2) .
% \end{align*}
\end{theorem}
\begin{remark}
% The emergence of parameters $h_1$ and $\alpha_2$ in our limiting mean and covariance functions may appear unconventional, but it stems from the unique aspects of our analysis. The reason is that, in the approximation of $\Expe\Bigl(\frac{w_{11}}{\bar{w}_1}-1\Bigr)^2$ and $\Expe\Bigl(\frac{w_{11}}{\bar{w}_1}-1\Bigr)^2\Bigl(\frac{w_{12}}{\bar{w}_1}-1\Bigr)^2$, the terms $h_1/p$ and $\alpha_2/p$ are non-negligible due to the multiplication by $p$ in the Central Limit Theorem (CLT) (see, Lemma \ref{estsqu} and Lemma \ref{quadform}). Additionally, our results introduce parameters $\alpha_1$ and $\alpha_2$ in place of traditional parameters like $\Expe |w_{11}|^4-3$ and $\Expe |w_{11}|^4-1$ in the limiting mean and covariance functions of the sample correlation matrix in \cite{gao2017high}. Notably, our results also introduce a novel parameter, $h_1$, in the mean function, setting our findings apart from conventional approaches.
The emergence of parameters $h_1$ and $\alpha_2$ in our limiting mean and covariance functions may appear unconventional, but it stems from the unique aspects of our analysis. This phenomenon arises from the non-negligible influence of terms $h_1/p$ and $\alpha_2/p$ in the approximation of $\Expe(\frac{w_{11}}{\bar{w}1}-1)^2$ and $\Expe(\frac{w_{11}}{\bar{w}_1}-1)^2(\frac{w_{12}}{\bar{w}_1}-1)^2$, driven by the multiplication by $p$ in the CLT (refer to Lemma \ref{estsqu} and Lemma \ref{quadform}). Furthermore, our results introduce parameters $\alpha_1$ and $\alpha_2$ in place of conventional parameters like $\Expe |w_{11}|^4-3$ and $\Expe |w_{11}|^4-1$ in the limiting mean and covariance functions of the sample correlation matrix in \cite{gao2017high}. 
Remarkably, our findings also bring forth a novel parameter, $h_1$, in the mean function, setting our results apart from conventional approaches.
\end{remark}
Applying Theorem \ref{clt} to three polynomial functions, we obtain the following corollary. The proof of Theorem \ref{clt} is postponed to Section \ref{prfcltall}, and detailed calculations in these applications are postponed to Section \ref{prfcorpoly} of the supplementary document.
\begin{Corollary}\label{corpoly}
Under conditions and notations in Theorem \ref{clt}, let $f_i=x^i$ for $i=1,2,3$, we have
% \begin{align*}
%     G_p(f_i) &= \tr (\bB_{p,N}^i) - s_i
%     \convd  \mathcal{N}(\mu_i,V_{i}),
% \end{align*}
\begin{align*}
    G_p(f_1) &= \tr (\bB_{p,N}) -p\frac{\sigma^2}{\mu^2}
    \convd  \mathcal{N}(\mu_1,V_{1}),\\
    G_p(f_2) &= \tr (\bB_{p,N}^2) -p(1+c_N)\Bigl(\frac{\sigma^2}{\mu^2} \Bigr)^2\convd  \mathcal{N}(\mu_2,V_{2}),\\
    G_p(f_3) &= \tr (\bB_{p,N}^3) -p(1+3c_N+c_N^2)\Bigl(\frac{\sigma^2}{\mu^2}\Bigr)^3\convd  \mathcal{N}(\mu_3,V_{3}),
\end{align*}
% $s_1=p\frac{\sigma^2}{\mu^2}$, $s_2=p(1+c_N)\Bigl(\frac{\sigma^2}{\mu^2} \Bigr)^2$, $s_3=p(1+3c_N+c_N^2)\Bigl(\frac{\sigma^2}{\mu^2}\Bigr)^3$,
where  $c_N=\frac{p}{N}$, and
% $\mu_1=h_1$, $\mu_2 
%     =  (1+c)\Bigl(\frac{\sigma^2}{\mu^2}\Bigr)^2+2(1+c)\frac{\sigma^2}{\mu^2}h_1+c (\alpha_1+\alpha_2)$, $\mu_3
%     = (2+6c+3c^2)\Bigl(\frac{\sigma^2}{\mu^2}\Bigr)^3+3(1+3c+c^2)\Bigl(\frac{\sigma^2}{\mu^2}\Bigr)^2h_1+3c(1+c)\frac{\sigma^2}{\mu^2}(\alpha_1+\alpha_2)$, and
%$\xi=\beta+h_2-2\frac{\sigma^2}{\mu^2}h_1$,
\begin{align}
    &\ \ \ \ \ \ \ \ \ \ \ \ \ \ \ \ \ \ \ \mu_1=h_1,\ \
    \mu_2 
    =  (1+c)\Bigl(\frac{\sigma^2}{\mu^2}\Bigr)^2+2(1+c)\frac{\sigma^2}{\mu^2}h_1+c (\alpha_1+\alpha_2 ),\nonumber\\
    %&=  (1+c)(\frac{\sigma^2}{\mu^2})^2+2(1+c)\frac{\sigma^2}{\mu^2}h_1+c\Big(\beta+h_2-2\frac{\sigma^2}{\mu^2}h_1\Big)\nonumber\\ %&=-c(\frac{\sigma^2}{\mu^2})^2+\frac{\sigma^2}{\mu^2}(\frac{\sigma^2}{\mu^2}+2ch_1+2h_1)+c\Big(\beta+h_2-2\frac{\sigma^2}{\mu^2}h_1\Big)+2c(\frac{\sigma^2}{\mu^2})^2\nonumber\\
    &\ \ \ \  \mu_3
    = (2+6c+3c^2)\Bigl(\frac{\sigma^2}{\mu^2}\Bigr)^3+3(1+3c+c^2)\Bigl(\frac{\sigma^2}{\mu^2}\Bigr)^2h_1+3c(1+c)\frac{\sigma^2}{\mu^2}(\alpha_1+\alpha_2),\nonumber\\
    %&= (2+6c+3c^2)(\frac{\sigma^2}{\mu^2})^3+3(1+3c+c^2)(\frac{\sigma^2}{\mu^2})^2h_1+3c(1+c)\frac{\sigma^2}{\mu^2}\Big(\beta+h_2-2\frac{\sigma^2}{\mu^2}h_1\Big)\nonumber
    %&\triangleq\Expe X_{f_3} %&= -3c(c+1)(\frac{\sigma^2}{\mu^2})^3+(\frac{\sigma^2}{\mu^2})^2\Big(3c\frac{\sigma^2}{\mu^2}+2\frac{\sigma^2}{\mu^2}+3c^2h_1+9ch_1+3h_1\Big)\nonumber\\
    %&\ \ \ \ \ \ +3c(c+1)\frac{\sigma^2}{\mu^2}\Big(\beta+h_2-2\frac{\sigma^2}{\mu^2}h_1\Big)+6c(c+1)(\frac{\sigma^2}{\mu^2})^3\nonumber\\
    V_{1} &= 2c\Bigl(\frac{\sigma^2}{\mu^2}\Bigr)^2+c(\alpha_1 + \alpha_2),\ \
    %&= 2c(\frac{\sigma^2}{\mu^2})^2+c\Big(\beta + h_2-2\frac{\sigma^2}{\mu^2}h_1\Big)\ \
    V_{2}=4c(2+c)(1+2c)\Bigl(\frac{\sigma^2}{\mu^2}\Bigr)^4+4c(1+c)^2\Bigl(\frac{\sigma^2}{\mu^2}\Bigr)^2(\alpha_1+\alpha_2),\nonumber\\
    %&=4c(2+c)(1+2c)(\frac{\sigma^2}{\mu^2})^4+4c(1+c)^2(\frac{\sigma^2}{\mu^2})^2\Big(\beta + h_2-2\frac{\sigma^2}{\mu^2}h_1\Big)
   &\ \ \ \ \ \ \ \ V_{3} =6c(1+6c+3c^2)(3+6c+c^2)\Bigl(\frac{\sigma^2}{\mu^2}\Bigr)^6 +9c(1+3c+c^2)^2\Bigl(\frac{\sigma^2}{\mu^2}\Bigr)^4(\alpha_1+\alpha_2).\nonumber%6c(1+6c+3c^2)(3+6c+c^2)(\frac{\sigma^2}{\mu^2})^6 +9c(1+3c+c^2)^2(\frac{\sigma^2}{\mu^2})^4(\beta + h_2-2\frac{\sigma^2}{\mu^2}h_1)
\end{align}
\end{Corollary}

\section{Numerical experiments}\label{numexp}
\subsection{Limiting spectral distribution}
% In this section, simulation experiments are conducted to examine the LSD of the sample covariance matrix $\bB_{p,N}$ generated by compositional data, as stated in Theorem \ref{thm:CoDa-LSD}. Recall that the compositional data $\{x_{ij}, 1\leq   i \leq   n, 1\leq   j \leq   p\}$ is generated by the normalization $x_{ij}=w_{ij}/\sum_{\ell=1}^p w_{i\ell}$. We generate the basis data $\{w_{ij}, 1\leq   i \leq   n, 1\leq   j \leq   p\}$ from three different populations, and draw the histograms of eigenvalues of $\bB_{p,N}$ and compare with their theoretical densities. Specifically, we consider three types of distributions for $w_{ij}$ as follows: 
In this section, simulation experiments are conducted to verify the LSD of the sample covariance matrix $\bB_{p,N}$ from compositional data, as stated in Theorem \ref{thm:CoDa-LSD}. Compositional data $\{x_{ij}\}_{1\leq i \leq n, 1\leq j \leq p}$ is generated by the normalization $x_{ij}=w_{ij}/\sum_{\ell=1}^p w_{i\ell}$. We generate basis data $w_{ij}$ from three populations, drawing histograms of eigenvalues of $\bB_{p,N}$ and comparing them with theoretical densities. Specifically, three types of distributions for $w_{ij}$ are considered:
\begin{enumerate}
    \item $w_{ij}$ follows the exponential distribution with rate parameter $5$;
    \item $w_{ij}$ follows the truncated standard normal distribution lying within the interval $(0,10)$, denoted by $\mathrm{TN}(0,1;0,10)$, where the first two parameters ($0$ and $1$) represent the mean and variance of the standard normal distribution;
    \item $w_{ij}$ follows the Poisson distribution with parameter $10$.
\end{enumerate}
% The pair of dimension and sample size, $(p, n)$, is taken to be $(500,500)$ or $(500,800)$.
% We list the histograms of eigenvalues of $\bB_{p,N}$ generated by three populations under different combinations of $(p, n)$, and compare them with their corresponding limiting densities in Figures \ref{lsdc1} -- \ref{lsdc2}. It can be seen from Figures \ref{lsdc1} -- \ref{lsdc2} that all the histograms conform to their theoretical limits, which proves the accuracy of our theoretical results.

\begin{figure}[htbp]
	\centering
	\begin{subfigure}{.32\textwidth}
		\centering		\includegraphics[width=.9\linewidth]{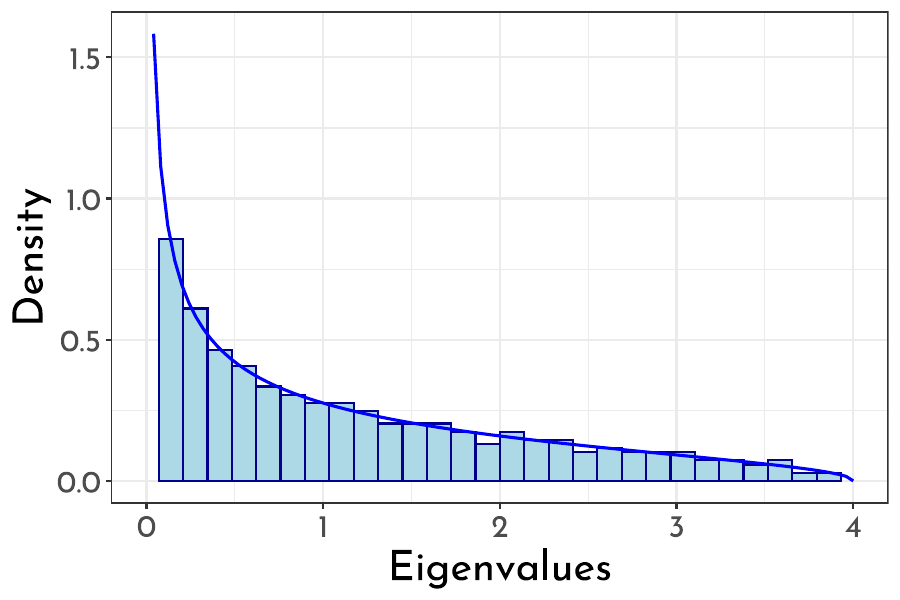}
		\caption{Exponential(5)}
	\end{subfigure}%
	\begin{subfigure}{.32\textwidth}
		\centering		\includegraphics[width=.9\linewidth]{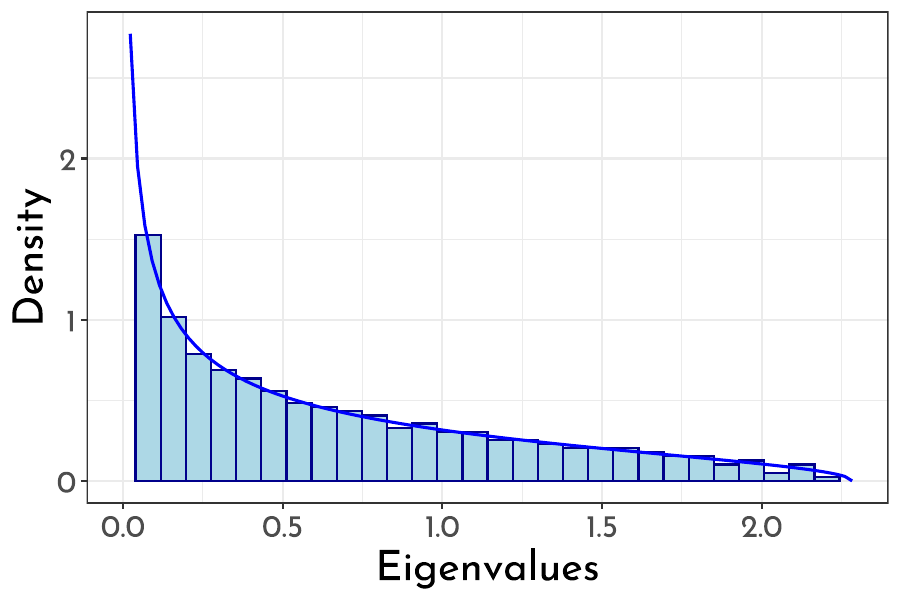}
		\caption{TN(0,1;0,10)}
	\end{subfigure}
	\begin{subfigure}{.32\textwidth}
		\centering		\includegraphics[width=.9\linewidth]{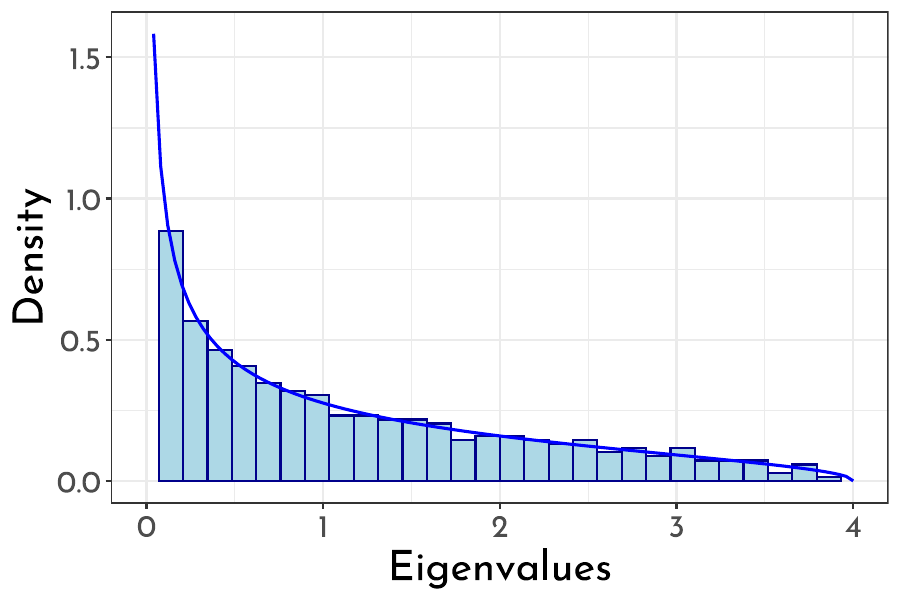}
		\caption{Poisson(10)}
	\end{subfigure}%
	\caption{Histograms of sample eigenvalues of $\bB_{p,N}$ with $(p,n)=(500,500)$.  The curves are density functions of their corresponding limiting spectral distribution.}\label{lsdc1}
\end{figure}

\begin{figure}[htbp]
	\centering
	\begin{subfigure}{.32\textwidth}
		\centering		\includegraphics[width=.9\linewidth]{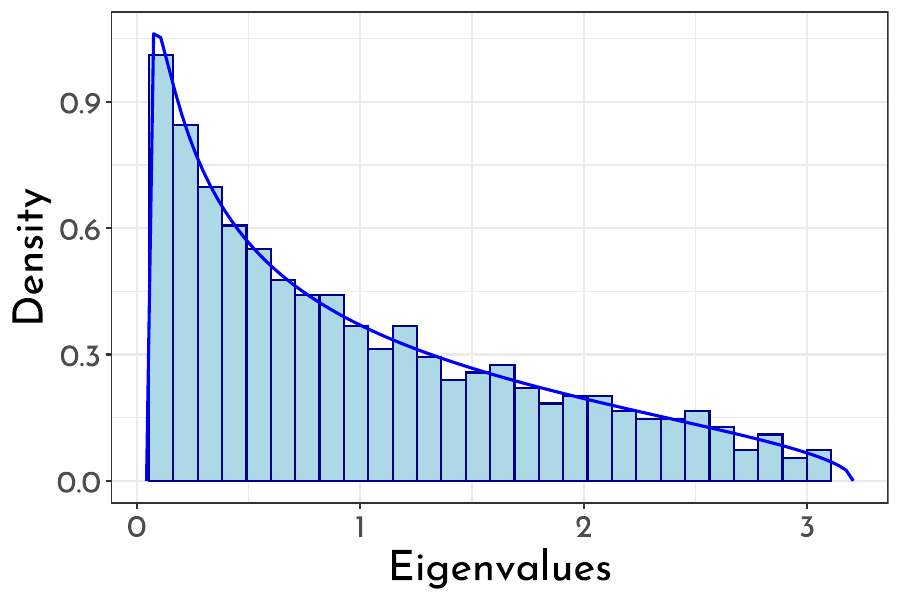}
		\caption{Exponential(5)}
	\end{subfigure}%
	\begin{subfigure}{.32\textwidth}
		\centering
        \includegraphics[width=.9\linewidth]{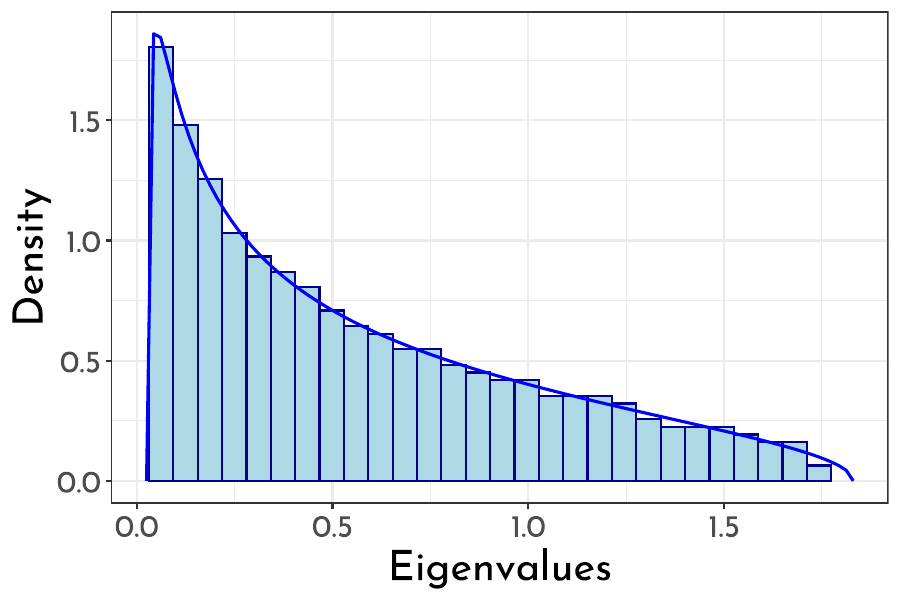}
		\caption{TN(0,1;0,10)}
	\end{subfigure}
	\begin{subfigure}{.32\textwidth}
		\centering
		\includegraphics[width=.9\linewidth]{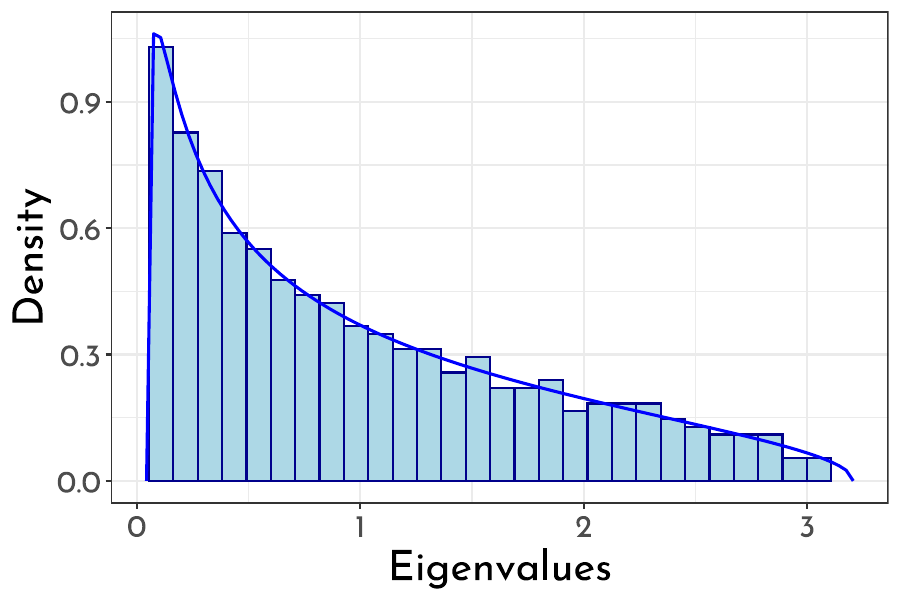}
		\caption{Poisson(10)}
	\end{subfigure}%
	\caption{Histograms of sample eigenvalues of $\bB_{p,N}$ with $(p,n)=(500,800)$. The curves are density functions of their corresponding limiting spectral distribution.}\label{lsdc2}
\end{figure}

The dimension and sample size pair, $(p, n)$, is set to $(500,500)$ or $(500,800)$. We display histograms of eigenvalues of $\bB_{p,N}$ generated by three populations under various $(p, n)$ combinations and compare them with their respective limiting densities in Figures \ref{lsdc1} -- \ref{lsdc2}. Figures \ref{lsdc1} -- \ref{lsdc2} reveal that all histograms align with their theoretical limits, affirming the accuracy of our theoretical results.

\subsection{CLT for LSS}
In this section, we implement some simulation studies to examine finite-sample properties of some LSS for $\bB_{p,N}$ by comparing their empirical means and variances with theoretical limiting values, as stated in Corollary \ref{corpoly}.
% \begin{enumerate}
%     \item finite-sample properties of some LSS for $\bB_{p,N}$ by comparing their empirical means and variances with theoretical limiting values, as stated in Corollary \ref{corpoly}.
%     % \item finite-sample properties of $M_p(z)$ by comparing their empirical means and variances with theoretical limiting values, as stated in Proposition \ref{mpzclt}.
% \end{enumerate} 

In the following, we present the numerical simulation of CLT for LSS.
% \subsubsection{CLT for LSS}
First, we compare the empirical mean and variance of $G_{p,N}(x^r) = \tr(\bB_{p,N}^r)-p\int x^r \dif F^{c_N}(x)$, $r = 1,2,3$, with their corresponding theoretical limits in Corollary \ref{corpoly}.
Two types of data distribution of $w_{ij}$ are consider:
\begin{enumerate}
    \item $w_{ij}$ follows the exponential distribution with rate parameter $5$;
    \item $w_{ij}$ follows the Chi-squared distribution with degree of freedom $1$. 
\end{enumerate}

\begin{figure}[htbp]
	\centering
	\begin{subfigure}{.32\textwidth}
		\centering		\includegraphics[width=.9\linewidth]{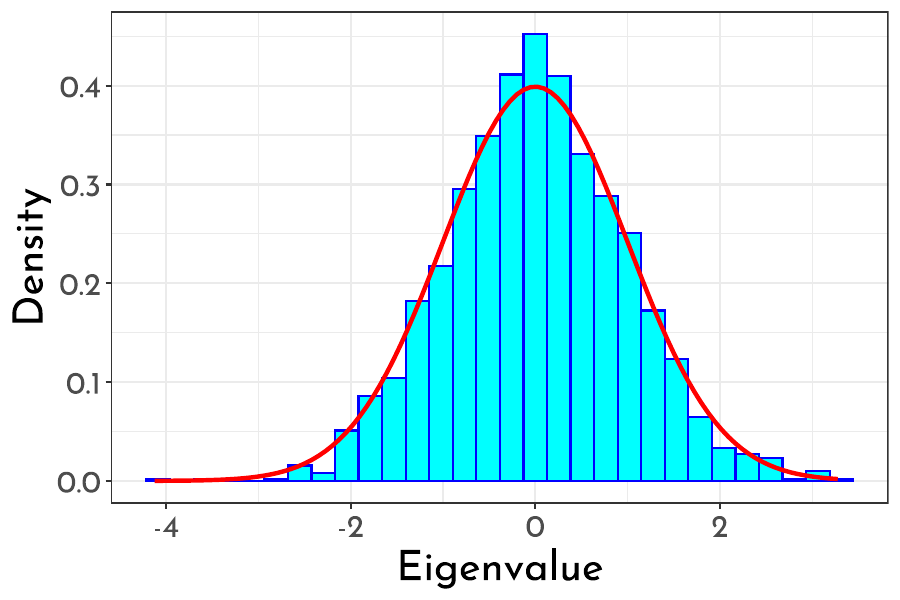}
		\caption{$\bar{G}_{p,N}(x)$}
	\end{subfigure}%
	\begin{subfigure}{.32\textwidth}
		\centering		\includegraphics[width=.9\linewidth]{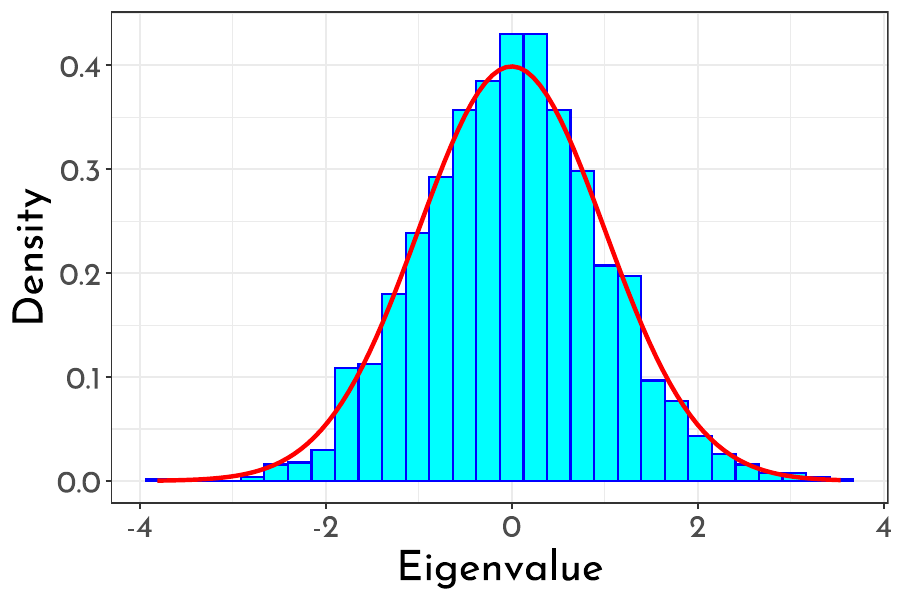}
		\caption{$\bar{G}_{p,N}(x^2)$}
	\end{subfigure}
	\begin{subfigure}{.32\textwidth}
		\centering		\includegraphics[width=.9\linewidth]{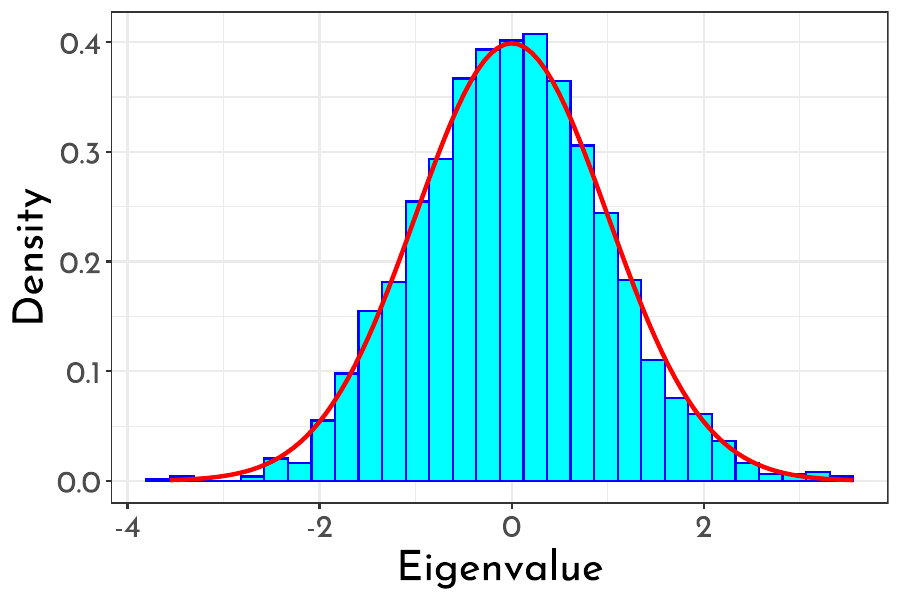}
		\caption{$\bar{G}_{p,N}(x^3)$}
	\end{subfigure}%
	\caption{Histograms of normalized LSS $\bar{G}_{p,N}(x^r)=(G_{p,N}(x^r)-\mu_r)/\sqrt{V_{r}}$, $r=1,2,3$, with $w_{ij}\sim\mathrm{Exp}(5)$ and $p=n=400$.  The curves are density functions of the standard normal distribution.}
        \label{fig:LSS-Exp}
\end{figure}

\begin{figure}[htbp]
	\centering
	\begin{subfigure}{.32\textwidth}
		\centering		\includegraphics[width=.9\linewidth]{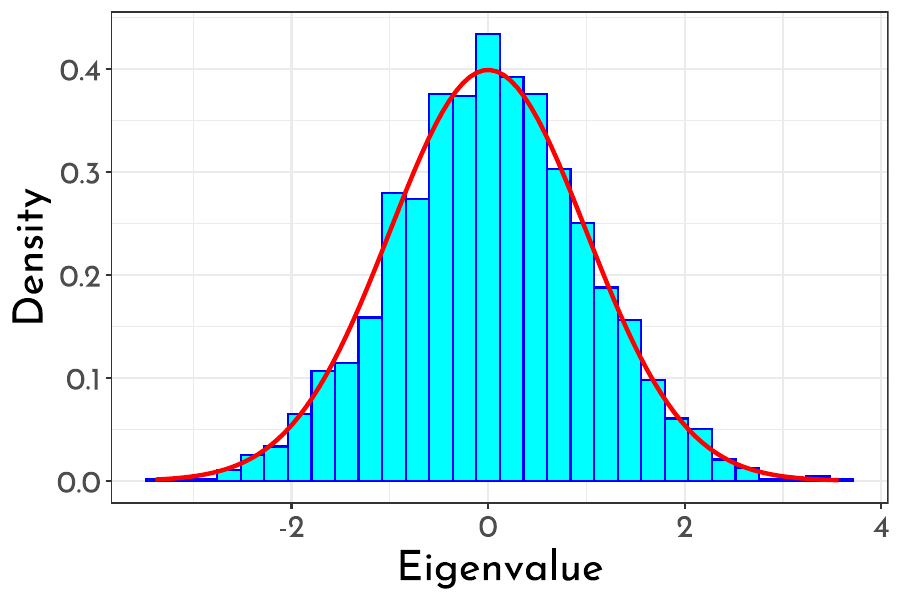}
		\caption{$\bar{G}_{p,N}(x)$}
	\end{subfigure}%
	\begin{subfigure}{.32\textwidth}
		\centering		\includegraphics[width=.9\linewidth]{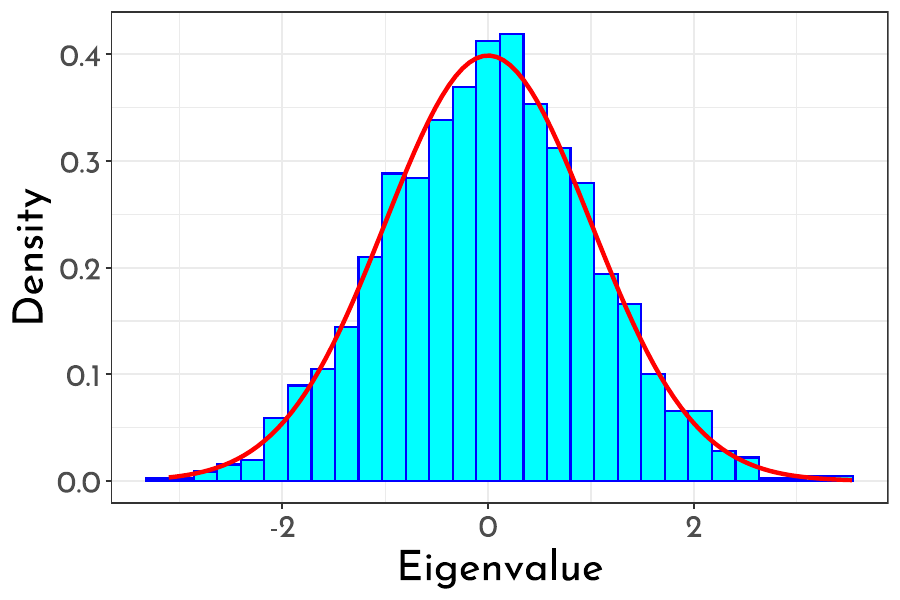}
		\caption{$\bar{G}_{p,N}(x^2)$}
	\end{subfigure}
	\begin{subfigure}{.32\textwidth}
		\centering		\includegraphics[width=.9\linewidth]{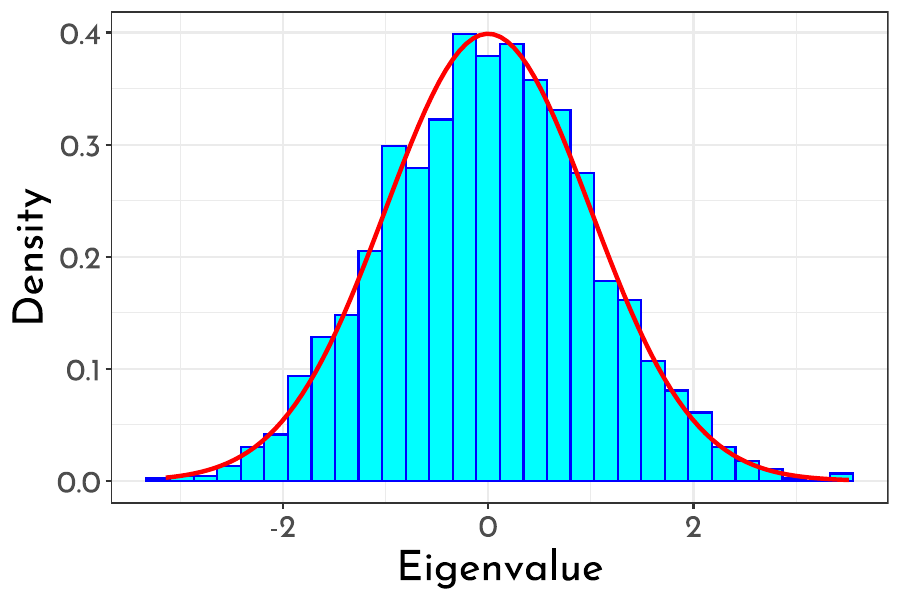}
		\caption{$\bar{G}_{p,N}(x^3)$}
	\end{subfigure}%
	\caption{Histograms of normalized LSS $\bar{G}_{p,N}(x^r)=(G_{p,N}(x^r)-\mu_r)/\sqrt{V_{r}}$, $r=1,2,3$, with $w_{ij}\sim\chi^2(1)$ and $p=300, n=400$.  The curves are density functions of the standard normal distribution.}
        \label{fig:LSS-Chisq}
\end{figure}

Empirical mean and variance of $\{G_{p,N}(x^r)\}$, $r=1,2,3$, are calculated for various combinations of $(p,n)$ with $p/n=3/4$ or $p/n=1$.
For each pair of $(p,n)$, $2000$ independent replications are used to obtain the empirical values. Tables \ref{tab:LSS-Exp} -- \ref{tab:LSS-chisq} report the empirical results for $\mathrm{Exp}(5)$ population and $\chi^2(1)$ population, respectively. As shown in Tables \ref{tab:LSS-Exp} -- \ref{tab:LSS-chisq}, the empirical mean and variance of $\{G_{p,N}(x^r)\}$ closely match their respective theoretical limits under all scenarios.
To verify the asymptotic normality of LSS, we draw the histogram of normalized LSS, $\bar{G}_{p,N}(x^r)=(G_{p,N}(x^r)-\mu_r)/\sqrt{V_{r}}$, $r=1,2,3$, where $\mu_r$ and $V_{r}$ are defined in Corollary \ref{corpoly}, and compare them with the standard normal density. 
Figures \ref{fig:LSS-Exp} and \ref{fig:LSS-Chisq}
depict the histograms of $\bar{G}_{p,N}(x^r)$ for $\mathrm{Exp}(5)$ population with $p/n=1$ and $\chi^2(1)$ population with $p/n=3/4$, respectively. The histograms for the cases of $\mathrm{Exp}(5)$ population with $p/n=3/4$ and $\chi^2(1)$ population with $p/n=1$ exhibit similar patterns and are omitted for brevity. It can be seen from Figures \ref{fig:LSS-Exp} -- \ref{fig:LSS-Chisq} that all
the histograms conform to the standard normal density, which fully supports our theoretical results.
\vspace{0.18cm}
% The histograms for the cases of $\mathrm{Exp}(5)$ population with $p/n=3/4$ and $\chi^2(1)$ population with $p/n=1$ exhibit similar patterns and are omitted for brevity.  
\begin{table}[htbp]
  \centering
  \caption{Empirical mean and variance of $G_{p,N}(x^r)$, $r=1,2,3$, with $w_{ij}\sim \mathrm{Exp}(5)$.}
    \begin{tabular}{ccccccccccc}
    \toprule
          &       &       & \multicolumn{2}{c}{$G_{p,N}(x)$} &       & \multicolumn{2}{c}{$G_{p,N}(x^2)$} &       & \multicolumn{2}{c}{$G_{p,N}(x^3)$} \\
\cmidrule{4-5}\cmidrule{7-8}\cmidrule{10-11}          & \multicolumn{1}{c}{$p/n$} & \multicolumn{1}{c}{$n$} & mean  & var   &       & mean  & var   &       & mean  & var \\
    \midrule
    \multicolumn{1}{c}{\multirow{4}[2]{*}{Emp}} & \multicolumn{1}{c}{\multirow{4}[2]{*}{$3/4$}} & \multicolumn{1}{c}{100} & -2.01 & 2.63  &       & -4    & 36.54 &       & -7.82 & 463.32 \\
          &       & \multicolumn{1}{c}{200} & -1.99 & 2.93  &       & -3.85 & 39.73 &       & -7.23 & 485.05 \\
          &       & \multicolumn{1}{c}{300} & -1.93 & 3.03  &       & -3.57 & 40.3  &       & -6.32 & 483.76 \\
          &       & \multicolumn{1}{c}{400} & -2.04 & 2.95  &       & -3.98 & 38.78 &       & -7.67 & 460.01 \\
    \cmidrule{2-11}
    Theo  &       &       & \textbf{-2} & \textbf{3} &       & \textbf{-3.75} & \textbf{39} &       & \textbf{-6.81} & \textbf{457} \\
    \midrule
    \multicolumn{1}{c}{\multirow{4}[2]{*}{Emp}} & \multicolumn{1}{c}{\multirow{4}[2]{*}{$1$}} & \multicolumn{1}{c}{100} & -1.91 & 3.61  &       & -3.83 & 64.09 &       & -6.56 & 1064.75 \\
          &       & \multicolumn{1}{c}{200} & -1.96 & 3.89  &       & -3.96 & 68.37 &       & -6.91 & 1090.14 \\
          &       & \multicolumn{1}{c}{300} & -2.01 & 3.97  &       & -4.06 & 68.7  &       & -7.16 & 1082.72 \\
          &       & \multicolumn{1}{c}{400} & -1.98 & 3.71  &       & -3.99 & 64.22 &       & -7.07 & 1010.09 \\
    \cmidrule{2-11}
    Theo  &       &       & \textbf{-2} & \textbf{4} &       & \textbf{-4} & \textbf{68} &       & \textbf{-7} & \textbf{1050} \\
    \bottomrule
    \end{tabular}%
  \label{tab:LSS-Exp}%
\end{table}%

\begin{table}[htbp]
  \centering
  \caption{Empirical mean and variance of $G_{p,N}(x^r)$, $r=1,2,3$, with $w_{ij}\sim \chi^2(1)$.}
    \begin{tabular}{ccccccccccc}
    \toprule
          &       &       & \multicolumn{2}{c}{$G_{p,N}(x)$} &       & \multicolumn{2}{c}{$G_{p,N}(x^2)$} &       & \multicolumn{2}{c}{$G_{p,N}(x^3)$} \\
\cmidrule{4-5}\cmidrule{7-8}\cmidrule{10-11}          & \multicolumn{1}{c}{$p/n$} & \multicolumn{1}{c}{$n$} & mean  & var   &       & mean  & var   &       & mean  & var \\
    \midrule
    \multicolumn{1}{c}{\multirow{4}[2]{*}{Emp}} & \multicolumn{1}{c}{\multirow{4}[2]{*}{$3/4$}} & \multicolumn{1}{c}{100} & -5.79 & 15.53 &       & -24.19 & 888.99 &       & -97.31 & 46790.03 \\
          &       & \multicolumn{1}{c}{200} & -5.96 & 16.74 &       & -24.39 & 920.63 &       & -96.17 & 45375.75 \\
          &       & \multicolumn{1}{c}{300} & -5.94 & 16.6  &       & -23.75 & 882.92 &       & -90.59 & 42487.68 \\
          &       & \multicolumn{1}{c}{400} & -5.88 & 17.51 &       & -22.68 & 912.28 &       & -81.2 & 42922.06 \\
    \cmidrule{2-11}
    Theo  &       &       & \textbf{-6} & \textbf{18} &       & \textbf{-23} & \textbf{918} &       & \textbf{-83} & \textbf{41806.12} \\
    \midrule
    \multicolumn{1}{c}{\multirow{4}[2]{*}{Emp}} & \multicolumn{1}{c}{\multirow{4}[2]{*}{$1$}} & \multicolumn{1}{c}{100} & -5.92 & 20.81 &       & -26.15 & 1563.02 &       & -102.73 & 107846.2 \\
          &       & \multicolumn{1}{c}{200} & -5.98 & 23.01 &       & -25.15 & 1639.95 &       & -90.25 & 105467.9 \\
          &       & \multicolumn{1}{c}{300} & -5.81 & 21.82 &       & -23.16 & 1526.34 &       & -74.54 & 96864.11 \\
          &       & \multicolumn{1}{c}{400} & -6.13 & 23.18 &       & -25.41 & 1599.96 &       & -90.31 & 99475.82 \\
    \cmidrule{2-11}
    Theo  &       &       & \textbf{-6} & \textbf{24} &       & \textbf{-24} & \textbf{1600} &       & \textbf{-80} & \textbf{96000} \\
    \bottomrule
    \end{tabular}%
  \label{tab:LSS-chisq}%
\end{table}%

\section{Proof of Theorem \ref{clt}}\label{prfcltall}
In this section, we first present the difference between the CLT for centralized sample covariance $\bB_p^0$ and unbiased sample covariance $\bB_{p,N}$ by substitution principle in Section \ref{subsprin}, where
\begin{align}
    \bB_{p}^0=p^2\bS_n^0=\frac{p^2}{n}(\bX_n-\Expe \bX_n)'(\bX_n-\Expe \bX_n)=\frac{1}{n}\bY_n'\bY_n,\ \
     \bB_{p,N}=p^2\bS_{n,N} =\frac{1}{N}(p\bX_n)'\bC_n(p\bX_n),
\end{align}
% \begin{align}
%     \bB_{p}^0&=p^2\bS_n^0=\frac{p^2}{n}(\bX_n-\Expe \bX_n)'(\bX_n-\Expe \bX_n)=\frac{1}{n}\bY_n'\bY_n,\\
%      \bB_{p,N}&=p^2\bS_{n,N} =\frac{1}{N}(p\bX_n)'\bC_n(p\bX_n),
% \end{align}
and $\bY_n=(y_{ij})_{n\times p}$, $y_{ij}=\frac{w_{ij}}{\bar{w}_i}-1$ and $\bar{w}_i=\frac{1}{p}\sum_{j=1}^p w_{ij}$. By substituting the adjusted sample size $N=n-1$ for the actual sample size $n$ in the centering term, the unbiased sample covariance matrix $\bB_{p,N}$ and the centralized sample covariance $\bB_{p}^0$ share the same CLT (see, Section \ref{subsprin}). The general strategy of the main proof of Theorem \ref{clt} is explained in the following and four major steps of the general strategy are presented in Section \ref{cltcenBp0}.

The general strategy of the proof follows the method established in \cite{bai2004clt} and \cite{gao2017high}, with necessary adjustments for handling the sample covariance matrix of HCD, where conventional tools are not directly applicable. Our novel techniques play a pivotal role in overcoming these challenges. To begin with, we follow the strategy in \cite{jiang2004limiting} to establish the LSD of $\bB_{p,N}$ in Theorem \ref{thm:CoDa-LSD}. Then, we develop Proposition \ref{extreigen} to find the extreme eigenvalues of $\bB_{p,N}$. Notably, these extreme eigenvalues are highly concentrated around two edges of the support, a crucial aspect for applying the Cauchy integral formula \eqref{chucform} and proving tightness. Given that compositional data $x_{ij}=w_{ij}/\sum_{j=1}^pw_{ij}$ are not i.i.d., dealing with the CLT for LSS of the unbiased sample covariance matrix $\bB_{p,N}$ presents challenges. To address this, we employ the substitution principle \citep{zheng2015substitution} to reduce the problem to the CLT for LSS of the centralized sample covariance $\bB_p^0$. By substituting the adjusted sample size $N=n-1$ for the actual sample size $n$ in the centering term, both the unbiased sample covariance matrix $\bB_{p,N}$ and the centralized sample covariance $\bB_{p}^0$ share the same CLT (see Section \ref{subsprin}). We then leverage the independence of samples to further study the CLT for LSS of $\bB_p^0$. Specifically, we exploit the independence of samples to establish independence for $\br_i=\frac{1}{\sqrt{n}}\bigl(\frac{w_{i1}}{\bar{w}_i}-1,\ldots,\frac{w_{ip}}{\bar{w}_i}-1\bigr)', i=1,2,\ldots,n$, and express $\bB_p^0$ as $\bB_p^0= \frac{1}{n}\bY_n'\bY_n = \sum_{i=1}^n \br_i\br_i'$. The ultimate goal is to establish the CLT for LSS of $\bB_p^0$.

 By the Cauchy integral formula, we have
\begin{equation}\label{chucform}
\int f(x)\dif G(x)=-\frac{1}{2\pi i}\oint_\mathcal{C}f(z)m_G(z)\dif z
\end{equation}
valid for any c.d.f $G$ and any analytic function $f$ on an open set containing the support of $G$, where $\oint_\mathcal{C}$ is the contour integration in the anti-clockwise direction. In our case, $G(x)=p(F^{\bB_p^0}(x)-F^{c_n}(x))$.  Therefore, the problem of finding the limiting distribution reduces to the study of $M_p(z)$ defined as follows:
\begin{align*}
    M_p(z)&=p\big[m_{p}(z)-m_{p}^{0}(z)\big]=n\big[\underline{m}_{p}(z)-\underline{m}_{p}^{0}(z)\big],\\
m_{p}(z)&=m_{F^{\bB_p^0}}(z)=\frac{1}{p}\tr\left[(\bB_p^0-z\bI_p)^{-1}\right],\qquad m_{p}^{0}(z)=m_{F^{c_n}}(z),\\ 
\underline{m}_{p}(z)&=m_{F^{\underline{\bB}_p^0}}(z)=\frac{1}{p}\tr\left[(\underline{\bB}_p^0-zI_n)^{-1}\right],\qquad\underline{m}_{p}^{0}(z)=\underline{m}_{F^{c_n}}(z),\\
\underline{\bB}_p^0&=p^2\underline{\bS}_n^0=\frac{p^2}{n}(\bX_n-\Expe\bX_n)(\bX_n-\Expe\bX_n)'.
\end{align*}

Note that the support of $F^{\bB_{p,N}}$ is random. Fortunately, we have shown that the extreme eigenvalues of $\bB_{p,N}$ are highly concentrated around two edges of the support of the limiting MP law $F^c(x)$ (see, Theorem \ref{thm:CoDa-LSD}, Proposition \ref{extreigen}). Then the contour $\mathcal{C}$ can be appropriately chosen. Moreover,  as in \cite{bai2004clt}, by Proposition \ref{extreigen}, we can replace the process $\{M_p(z),z\in\mathcal{C}\}$ by a slightly modified process $\{\widehat{M}_p(z), z\in\mathcal{C}\}$. Below we present the definitions of the contour $\mathcal{C}$ and the modified process $\widehat{M}_p(z)$. Let $x_r$ be any number greater than $\frac{\sigma^2}{\mu^2}(1+\sqrt{c})^2$. Let $x_l$ be any negative number if the left endpoint of \eqref{suppofLSD} is zero. Otherwise we choose $x_l\in(0,\frac{\sigma^2}{\mu^2}(1-\sqrt{c})^2)$. Now let $\mathcal{C}_u=\{x+iv_0:x\in[x_l,x_r]\}$. Then we define $\mathcal{C}^+ := \{x_l+iv:v\in[0,v_0]\}\cup\mathcal{C}_u\cup\{x_r+iv: v\in[0,v_0]\}$, and $\mathcal{C}=\mathcal{C}^+\cup \overline{\mathcal{C}^+}$. Now we define the subsets $\mathcal{C}_n$ of $\mathcal{C}$ on which $M_p(\cdot)$ equals to $\widehat{M}_p(\cdot)$. Choose sequence $\{\varepsilon_n\}$ decreasing to zero satisfying for some $\alpha\in(0,1)$, $\varepsilon_n\geq n^{-\alpha}$. Let
\begin{equation*}
\mathcal{C}_l=\left\{
\begin{array}{ccc}
\{x_l+iv:v\in[n^{-1}\varepsilon_n,v_0]\}~~~~ \mathrm{if}~x_l>0,\\
\{x_l+iv:v\in[0,v_0]\} ~~~~~~~~~~\mathrm{if}~x_l<0,
\end{array}
\right.
\end{equation*}
and $\mathcal{C}_r=\{x_r+iv: v\in[n^{-1}\varepsilon,v_0]\}$. Then $\mathcal{C}_n=\mathcal{C}_l\cup\mathcal{C}_u\cup\mathcal{C}_r$. For $z=x+iv$, we define
% \begin{equation*}
% \widehat{M}_p(z)=\left\{
% \begin{array}{ccc}
% &M_p(z)\mathrm{for}~ z\in \mathcal{C}_p,\\
% &M_p(x_r+in^{-1}\varepsilon_n)\ \ \ \ \mathrm{for}~ x=x_r,v\in[0,n^{-1}\varepsilon_n], and \ if x_l>0,\\
% &M_p(x_l+in^{-1}\varepsilon_n)\mathrm{for}~
% x=x_l,v\in[0,n^{-1}\varepsilon_n].
% \end{array}\right.
% \end{equation*}
\begin{equation*}
\widehat{M}_p(z)= \begin{cases}M_p(z), & \text { for } z \in \mathcal{C}_n \\ 
M_p(x_r+i n^{-1} \varepsilon_n), & \text { for } x=x_r, v \in[0, n^{-1} \varepsilon_n], \text { and if } x_l>0\\ 
M_p(x_l+i n^{-1} \varepsilon_n), & \text { for } x=x_l,v\in[0,n^{-1}\varepsilon_n],\end{cases}
\end{equation*}

Most of the paper will deal with proving the following proposition.
\begin{Proposition}\label{mpzclt}
  Under Assumption \ref{asforw} and \ref{asforcn}, then $\widehat{M}_{p}(\cdot)$ converges weakly to a two-dimensional Gaussian process
$M(\cdot)$ for $z \in \mathcal{C}$, with means
\begin{align}\label{mpzmean}
    \Expe M(z) =& -\underline{m}(z)\left[1-c\frac{\sigma^{4}}{\mu^{4}}\underline{m}^{2}(z)\Bigl(1+\frac{\sigma^{2}}{\mu^{2}}\underline{m}(z)\Bigr)^{-2}\right]^{-1}\nonumber\\
    &\times \biggl[-z\underline{m}(z)\left(1+\frac{\sigma^{2}}{\mu^{2}}\underline{m}(z)\right)^{-1} \times \left(h_{1}m(z)+\frac{\sigma^{2}}{\mu^{2}}m(z)+\frac{\sigma^{2}}{\mu^{2}}\frac{1}{z}\right)\nonumber\\
    &\qquad-cz^{2}\underline{m}^{2}(z)\left(1+\frac{\sigma^{2}}{\mu^{2}}\underline{m}(z)\right)^{-1}\left(\alpha_1 m^{2}(z)+\alpha_2 m^{2}(z)+2\frac{\sigma^{4}}{\mu^{4}}m'(z)\right)\nonumber\\
    &\qquad+c\frac{\sigma^{4}}{\mu^{4}}\underline{m}^{2}(z)\left(1+\frac{\sigma^{2}}{\mu^{2}}\underline{m}(z)\right)^{-3}\left(1-c\frac{\sigma^{4}}{\mu^{4}}\underline{m}^{2}(z)\left(1+\frac{\sigma^{2}}{\mu^{2}}\underline{m}(z)\right)^{-2}\right)^{-1}\biggr],
\end{align}
and covariance function
\begin{align}\label{mpzcov}
    \Cov(M(z_{1}),M(z_{2})) 
    &=2\left[\frac{\underline{m}'(z_{1})\underline{m}'(z_{2})}{(\underline{m}(z_{1})-\underline{m}(z_{2}))^{2}}-\frac{1}{(z_{1}-z_{2})^{2}}\right]\nonumber\\
    &\quad+ c\left(\alpha_1 +\alpha_2  \right)\times\frac{\underline{m}'(z_{1})\underline{m}'(z_{2})}{(1+\sigma^2/\mu^2\underline{m}(z_{1}))^{2}(1+\sigma^2/\mu^2\underline{m}(z_{2}))^{2}}.
\end{align}
\end{Proposition}
% \textbf{Step 4}: The limiting
% process of $M_p(z)$ on $\mathcal{C}$.
% To prove Theorem \ref{clt}, as in \cite{bai2004clt}, it suffices to prove the CLT for $\widehat{M}_p(z)$ with $z\in \mathcal{C}$. We then prove the limiting process of $M_p(z)$ on $\mathcal{C}$. We state the result in the following proposition and most of the remaining work will deal with proving this Proposition.

Now we explain how Theorem \ref{clt} follows from the above proposition. As in \cite{bai2004clt}, with probability 1, $\left|\int f(z)(M_p(z)-\widehat{M}_p(z))\dif z\right|\to 0$ as $n\to\infty$. Combining this observation with equation \eqref{chucform}, Theorem \ref{clt} follows from Proposition \ref{mpzclt}. To prove Proposition \ref{mpzclt}, we decompose $M_p(z)$ into a random part $M^{(1)}_p(z)$ and a deterministic part $M^{(2)}_p(z)$ for $z\in\mathcal{C}_n$, that is, $M_p(z)=M_p^{(1)}(z)+M_p^{(2)}(z)$, where  
\[
    M_p^{(1)}(z)=p\big[m_{p}(z)-\mathbb{E}m_{p}(z)\big] 
    \quad\text{and}\quad
    M_p^{(2)}(z)=p\big[\mathbb{E}m_{p}(z)-m_{p}^{0}(z)\big].
\]
The random part contributes to the covariance function and the deterministic part contributes
to the mean function. By Theorem 8.1 in Billingsley (1968), the proof of Proposition \ref{mpzclt} is
then complete if we can verify the following four steps:\\
  \textbf{Step 1} Truncation.\\
  \textbf{Step 2} Finite-dimensional convergence of $M_p^{(1)}(z)$ in distribution on $\mathcal{C}_n$  to a centered multivariate Gaussian random vector with covariance function given by \eqref{mpzcov}.\\
    \textbf{Step 3} Tightness of the $M_p^{(1)}(z)$ for $z\in\mathcal{C}_n$.\\
    \textbf{Step 4} Convergence of the non-random part $M_p^{(2)}(z)$ to \eqref{mpzmean} on $z\in\mathcal{C}_n$. 
    
The proof of these steps is presented in the coming sections. Before that, we introduce the substitution principle and crucial lemmas in Sections \ref{subsprin} and \ref{importlem} respectively. The former explains the reduction of problem of the CLT for LSS of $\bB_{p,N}$ to that of $\bB_{p}^0$, while the latter provides essential lemmas for these four steps in proving the CLT for LSS of $\bB_{p}^0$.
% Prior to that, we introduce the substitution principle and crucial lemmas in Sections \ref{subsprin} and \ref{importlem}. The former explains reducing the problem of the CLT for the LSS of $\bB_{p,N}$ to that of $\bB_{p}^0$, while the latter provides essential lemmas for these four steps in proving the CLT for LSS of $\bB_{p}^0$. 
% Before delving into the proof for these four steps, we introduce the substitution principle and some crucial lemmas in Sections \ref{subsprin} and \ref{importlem} respectively. The former explains how we can reduce the problem of the CLT for the LSS of the unbiased sample covariance matrix $\bB_{p,N}$ to that of the centralized sample covariance $\bB_{p}^0$, while the latter provides essential lemmas for these four steps in proving the CLT for LSS of the centralized sample covariance $\bB_{p}^0$. The detailed proofs for Sections \ref{subsprin} and \ref{importlem} are postponed to Section \ref{pfsubstiprintr}-\ref{prfquadform} of the supplementary document.
\subsection{Substitution principle}\label{subsprin}

% Let
% \begin{align}
%  % \bB_{p,N}&=p^2\bS_{n,N} =p^2\frac{1}{N}\bX_n'\bC_n\bX_n,\\
%     \bB_p^0=p^2\bS_n^0 &= p^2\frac{1}{n}(\bX_n-\Expe\bX_n)'(\bX_n-\Expe\bX_n)= \frac{1}{n}\bY_n'\bY_n
% \end{align}
% where $\bY_n=(y_{ij})_{n\times p}$, $y_{ij}=\frac{w_{ij}}{\bar{w}_i}-1$, $\bar{w}_i =\sum_{j=1}^pw_{ij}/p$.
By the Cauchy integral formula, we have
\begin{align}\label{chucformunbias}
G_{p,N}(f)
&=-\frac{1}{2\pi i}\oint_\mathcal{C}f(z)\left[\operatorname{tr}(\bB_{p,N}-z\bI_p)^{-1}-pm_N^0(z)\right]\dif z
%&=\sum_{i=1}^p f(\lambda_i)-p\int f(x)\dif F^{c_N}(x)\nonumber\\
\end{align}
% \begin{align}\label{chucformunbias}
% G_{p,N}(f):= \int f(x)\dif G_{p,N}(x)
% &=-\frac{1}{2\pi i}\oint_\mathcal{C}f(z)\left[\operatorname{tr}(\bB_{p,N}-z\bI_p)^{-1}-pm_N^0(z)\right]\dif z
% %&=\sum_{i=1}^p f(\lambda_i)-p\int f(x)\dif F^{c_N}(x)\nonumber\\
% \end{align}
valid for any function $f$ analytic on an open set containing the support of $G_{p,N}$, where
\begin{align}
    m_N^0(z) &\equiv m_{F^{c_N}}(z)= \frac{1}{\sigma^2/\mu^2(1-c_N-c_Nzm_N^0)-z},\label{munbiased}\\ 
    % \label{munbiased}\\
   \underline{m}_N^0(z)
    &\equiv\underline{m}_{F^{c_N}}(z)=  -\frac{1-c_N}{z} + cm_N^0(z),\label{mbarunbiased}\\ 
    z &=  -\frac{1}{\underline{m}_N^0(z)}+c_N\frac{\sigma^2/\mu^2}{1+\sigma^2/\mu^2\underline{m}_N^0(z)}\label{zmbarunbiased}%z&=&c\frac{\sigma^2/\mu^2{1+\sigma^2/\mu^2\underline_{m}}\nonumber%-\frac{1}{\underline_{m}}  \label{zmbar}
\end{align}
with $c_N=\frac{p}{N}$. To obtain the asymptotic distribution of $G_{p,N}(f)$, it is necessary to find the asymptotic distribution of $\operatorname{tr}(\bB_{p,N}-z\bI_p)^{-1}-pm_N^0(z)$.
% \begin{eqnarray}
%     \operatorname{tr}(\bB_{p,N}-z\bI_p)^{-1}-pm_N^0(z)
% \end{eqnarray}
To achieve this, we derive the following Lemma \ref{substiprintr} whose proof is postponed to Section \ref{pfsubstiprintr} of the supplementary document. 
\begin{lemma}\label{substiprintr}
    Under conditions and notations in Theorem \ref{clt}, as $n\rightarrow\infty$,
    \begin{eqnarray}
         \operatorname{tr}(\bB_{p,N}-z\bI_p)^{-1}-pm_N^0(z) =  \operatorname{tr}(\bB_{p}^0-z\bI_p)^{-1}-pm_n^0(z)+o_P(1).
    \end{eqnarray}
\end{lemma}
By Lemma \ref{substiprintr}, the asymptotic distribution of $G_{p,N}(f)$ is identical to that of $G_p^0(f)$, i.e.,
\begin{align}
    G_{p,N}(f) &=\sum_{i=1}^p f(\lambda_i(\bB_{p,N}))-p\int f(x)\dif F^{c_N}(x)\Longrightarrow \mathcal{N}(m(f),v(f)),\label{GpNlimdist}\\
    G_{p}^0(f) &=\sum_{i=1}^p f(\lambda_i(\bB_{p}^0))-p\int f(x)\dif F^{c_n}(x)\Longrightarrow \mathcal{N}(m(f),v(f)),\label{Gp0limdist}
\end{align}
a Guassian distribution whose parameters $m(f)$ and $v(f)$ depend only on the LSD $F^c(x)$ and $f$,
% \begin{align}
% G_{p,N}(f) &=\sum_{i=1}^p f(\lambda_i)-p\int f(x)dF^{c_N}(x)\Longrightarrow \mathbf{N}(\mu_1,\Sigma_1),\nonumber\\
%     G_{p}^0(f) &=\sum_{i=1}^p f(\lambda_i^0)-p\int f(x)dF^{c_n}(x)\Longrightarrow \mathbf{N}(\mu_1,\Sigma_1),\nonumber
% \end{align}
% to obtain the asmptotic distribution of $G_p(f)$ is equivalent to find  the asmptotic distribution of $G_p^0(f)$, 
where
\begin{align}\label{chucformcent}
G_p^0(f)
&=-\frac{1}{2\pi i}\oint_\mathcal{C}f(z)\left[\operatorname{tr}(\bB_{p}^0-z\bI_p)^{-1}-pm_n^0(z)\right]\dif z
\end{align}
% \begin{align}\label{chucformcent}
% G_p^0(f):= \int f(x)\dif G_p^0(x)
% &=-\frac{1}{2\pi i}\oint_\mathcal{C}f(z)\left[\operatorname{tr}(\bB_{p}^0-z\bI_p)^{-1}-pm_n^0(z)\right]\dif z
% \end{align}
and $c_n=\frac{p}{n}$, $m_n^{0}(z)=m_{F^{c_n}}(z)$ (note that we denote $m_n^{0}(z)$ as $m_p^{0}(z)$ in other sections except this subsection).

\subsection{Some important lemmas}\label{importlem}

Before delving into the proof of the CLT for LSS, it is crucial to introduce three pivotal lemmas, representing novel contributions to this paper, that unveil concentration phenomena. Lemma \ref{estsqu} is crafted to estimate essential parameters, facilitating the derivation of estimates of any order. Concerning $\nu_2$ and $\nu_{12}$, the terms $h_1/p$ and $h_2/p$ emerge as non-negligible due to the multiplication by $p$ in the CLT. To address these parameters, we establish that the probability of the event $B^c_{p}(\epsilon)$ decays polynomially to $0$ and leverage Taylor expansion on the event $B_{p}(\epsilon)=\{\omega: |\overline{w}_i-\mu|\leq \epsilon, \overline{w}_i=\sum_{j=1}^{p}w_{ij}/p\}$ to handle the issue of dependence. The proof of the CLT for LSS relies on two pivotal steps: the moment inequality for random quadratic forms and the precise estimation of the expectation of the product of two random quadratic forms. Lemma \ref{bound} establishes the former step, essential for converting them into the corresponding traces, while Lemma \ref{quadform} establishes the latter step, enabling the application of CLT for martingale differences. Both Lemma \ref{bound} and Lemma \ref{quadform} heavily hinge on the estimation of parameters $\nu_2$, $\nu_4$, and $\nu_{12}$ in Lemma \ref{estsqu}. The proof of Lemmas \ref{estsqu} -- \ref{quadform} are postponed to
Sections \ref{prfestsqu} -- \ref{prfquadform} of the supplementary document.

\begin{lemma}\label{estsqu}
Suppose that $\bw=\left(w_{1}, \ldots ,  w_{p}\right)'$ has i.i.d. entries with $\Expe w_1=\mu$, $\Expe (w_1-\mu)^2=\sigma^2$, and $\Expe\left|w_{1}-\mu\right|^{4}<\infty$, let $\overline{w}=\frac{1}{p}\sum_{j=1}^{p}w_{j}$, then there exists a constant $K>0$, such that for any $0<\epsilon<1/2$ and $p>0$,
\begin{align}
\nu_2:=\Expe\Bigl(\frac{w_{1}}{\overline{w}}-1\Bigr)^{2} &=\Expe\Bigl(\frac{w_{1}}{\mu}-1\Bigr)^{2}+\frac{1}{p}h_1+ o(p^{-1}),\label{eq:moment1}\\
%[\frac{3\Expe w_{1}^{2}\Expe(w_{1}-\mu)^{2}}{\mu^{4}}-\frac{2\Expe w_{1}^{2}(w_{1}-\mu)}{\mu^{3}}]\\
% + \frac{3}{p^{2}}\Bigl[\frac{\Expe w_{1}^{2}(w_{1}-\mu)^{2}}{\mu^{4}}-\frac{\Expe w_{1}^{2}\Expe(w_{1}-\mu)^{2}}{\mu^{4}}\Bigr]
\nu_{12}:=   \Expe\Bigl(\frac{w_1}{\bar{w}}-1\Bigr)^2\Bigl(\frac{w_2}{\bar{w}}-1\Bigr)^2&=\Bigl[\Expe\Bigl(\frac{w_1}{\mu}-1\Bigr)^2\Bigr]^2+\frac{1}{p}h_2+o(p^{-1}),\label{eq:moment2}\\
 \nu_4:=  \Expe\Bigl(\frac{w_1}{\bar{w}}-1\Bigr)^4 &= \Expe\Bigl(\frac{w_1}{\mu}-1\Bigr)^4 + o(1),\label{eq:moment3}
\end{align}
where  
\begin{align}
    h_1=-2\frac{\Expe w_{11}^3}{\mu^3}+3\Bigl(\frac{\sigma^2}{\mu^2}\Bigr)^2+5\frac{\sigma^2}{\mu^2}+2, \ \
    h_2= -8\frac{\sigma^2}{\mu^2}\frac{\Expe w_{11}^3}{\mu^3}+10\Bigl(\frac{\sigma^2}{\mu^2}\Bigr)^3+22\Bigl(\frac{\sigma^2}{\mu^2}\Bigr)^2+8\frac{\sigma^2}{\mu^2}.\label{hexpr}
\end{align}

% hence 
%     \begin{align}
%         \nu_4 -3\nu_{12}&= \alpha_1+ o(1),\label{beta1appr}\\
%         p(\nu_{12}-\nu_2^2)&=\alpha_2+o(1),\label{beta2appr}\\%= -4\frac{\sigma^2}{\mu^2}\frac{\Expe w_1^3}{\mu^3}+4(\frac{\sigma^2}{\mu^2})^3+12(\frac{\sigma^2}{\mu^2})^2+4\frac{\sigma^2}{\mu^2},\label{beta2appr}\\
%         %=-4\frac{\Expe (w_{11}-\mu)^{3}\Expe (w_{11}-\mu)^{2}}{\mu^{5}}\nonumber\\
%    %&\ \ \ \ \ \ \ \ \ \ \ \  \ \ \ \ \ \ \ \ \ \ +4\frac{(\Expe (w_{11}-\mu)^{2})^{3}}{\mu^{6}}+12\frac{\Expe (w_{11}-\mu)^{2}\Expe (w_{11}-\mu)}{\mu^{3}}-4,
%         p\Bigl(\nu_2-\frac{\sigma^2}{\mu^2}\Bigr)&=h_1+ o(1),\label{h1appr}  %\nu_2&=\frac{\sigma^2}{\mu^2}+o(1).\label{eigenappr}
%     \end{align} 
% where $h_1$ is in \eqref{h1expr} and $h_2$ is in \eqref{h2expr}.

% \begin{align*}
% |a|\leq K\epsilon^{3}, \quad |b|&\leq K [(\Prob(B^{c}_{p}(\epsilon)))^{1/2}+p^{-1}(\Prob(B^{c}_{p}(\epsilon)))^{1/4}+p^{-2}], \quad |c|\leq Kp^{2}\Prob(B^{c}_{p}(\epsilon)),\\
%    &\ \ \ \   B_{p}(\epsilon) =\{\omega: |\overline{w}-u|\leq \epsilon, \overline{w}=\sum_{j=1}^{p}w_{j}/p \},\\
%   \ \ \ \ \ & \Prob(B^{c}_{p}(\epsilon)) \leq K_{q_{1}}\epsilon^{-kq_{1}}(\sigma^{kq_{1}}p^{-kq_{1}/2}+p^{-kq_{1}+1}\Expe\left|w_{1}\right|^{kq_{1}}),
% \end{align*}
% in which $\epsilon>0$ is a constant and $k, q_{1} >0$ are constants.
\end{lemma}

\begin{lemma}\label{bound}
    Suppose that $\bw=\left(w_{1}, \ldots ,  w_{p}\right)'$ has i.i.d. entries with $\Expe w_1=\mu$, $\Expe (w_1-\mu)^2=\sigma^2$,  for any $p \times p$  matrix $\bA$ and $q \geq 2$, we have there is a positive constant $K_{q}$ depending on $q$ such that
    
\begin{align}
    \Expe\Bigl|\br'\bA \br-\frac{1}{n}\nu_2\tr \bA\Bigr|^{q} &\leq K_{q} \Bigl[n^{-q}\Bigl(\bigl(\Expe\left|w_{1}\right|^{4} \tr (\bA \bA') \bigr)^{q / 2}+\Expe|w_{1}|^{2 q} \tr(\bA \bA')^{q / 2}\Bigr)\nonumber\\
    &\qquad\qquad  + n^{q}  \Prob\bigl(B^{c}_{p}(\epsilon)\bigr) \|\bA\|^{q}+n^{-q}\|\bA^{q}\|h_1^q\Bigr],\label{eq:quad-form-moment}
\end{align}
where $\br = \frac{1}{\sqrt{n}}(w_{1}/\overline{w}-1,\dots,w_{p}/\overline{w}-1)'$, $h_1$ is in \eqref{hexpr},
$
B_{p}(\epsilon) =\{\omega: |\overline{w}-u|\leq \epsilon, \overline{w}=\sum_{j=1}^{p}w_{j}/p \},
$
and
\begin{equation}\label{eq:Bc_prob}
\Prob(B^{c}_{p}(\epsilon)) \leq K\epsilon^{-kq_{1}}(\sigma^{kq_{1}}p^{-kq_{1}/2}+p^{-kq_{1}+1}\Expe\left|w_{1}\right|^{kq_{1}}),
\end{equation}
in which $\epsilon, k, q_{1} >0$ are constants.
Furthermore, if $\|\bA\|\leq K$ and $|w_j-\mu|<\delta_n\sqrt{n}$ for all $j=1,\ldots,p$, then, for any $q\geq 2$,
\[
\Expe\Bigl|\br' \bA \br-\frac{1}{n}\nu_2\tr \bA\Bigr|^{q} \leq K_{q} n^{-1}\delta_{n}^{2q-4}.
\]
\end{lemma}

\begin{lemma}\label{quadform}
        Suppose that $\bw=\left(w_{1}, \ldots ,  w_{p}\right)'$  has i.i.d. entries with $\Expe w_1=\mu$, $\Expe (w_1-\mu)^2=\sigma^2$,  $\bA$ and $\bB$ are $p \times p$  matrices , if $\|\bA\|\leq K$ and $\|\bB\|\leq K$, 
    then
\begin{align}\label{quadformnonran}
    &\ \ \ \ \Expe \Bigl(\br'A\br-\frac{1}{n}\nu_2\tr \bA\Bigr)\Bigl(\br'\bB\br-\frac{1}{n}\nu_2\tr \bB\Bigr)\nonumber\\
    &=\frac{1}{n^{2}}(\nu_4-3\nu_{12} ) \sum_{i=1}^{p}A_{ii}B_{ii}+\frac{1}{n^{2}}\nu_{12}\bigl(\tr  (\bA\bB')+\tr (\bA\bB)\bigr) +\frac{1}{n^{2}} (\nu_{12}-\nu_{2}^{2}) \tr \bA \tr \bB 
    +o(n^{-1}).\nonumber
\end{align}
% where $\nu_2=\Expe(\frac{w_1}{\bar{w}}-1)^2$, $\nu_4=\Expe(\frac{w_1}{\bar{w}}-1)^4$, $\nu_{12}=\Expe(\frac{w_1}{\bar{w}}-1)^2(\frac{w_2}{\bar{w}}-1)^2$. 
%$B_{p}(\epsilon) =\{\omega: |\overline{w}-u|\leq \epsilon, \overline{w}=\frac{1}{p}\sum_{j=1}^{p}w_{j} \}$.
\end{lemma}

\subsection{CLT for
LSS of the centralized sample covariance $B_p^0$}\label{cltcenBp0}

% \subsubsection{Truncation}\label{pretrun}

\textbf{Step 1: Truncation.} We begin the proof of Proposition \ref{mpzclt} with the replacement of the entries of $\bW_n$ with truncated variables. We can choose a positive sequence of $\{\delta_n\}$ such that
\[
\delta_n\rightarrow0, 
\quad \delta_nn^{1/4}\rightarrow\infty,
\quad \delta_n^{-4} \Expe w_{11}^4I_{\{|w_{11}-\mu|\geq\delta_n\sqrt{n}\}}\rightarrow 0.
\]

% Let 
% \[
%     \bB_p^0=\frac{p^2}{n}(\bX_n-\Expe\bX_n)'(\bX_n-\Expe\bX_n).
% \]

Let $\hat\bB_p^0=\frac{p^2}{n}(\hat{\bX}_n-\Expe\hat{\bX}_n)'(\hat{\bX}_n-\Expe\hat{\bX}_n)$, where $\hat\bW_n$ is $n\times p$ matrix having $\hat{w}_{ij}=w_{ij}I_{\{|w_{ij}-\mu|<\delta_n\sqrt{n}\}}$.
% $\hat\bLambda_n$ is $n\times n$ diagonal matrix having $(\sum_{j=1}^{p}\hat{w}_{ij})^{-1}=(\sum_{j=1}^{p}w_{ij}I_{\{|w_{ij}-\mu|<\delta_n\sqrt{n}\}})^{-1}$.
We then have
\begin{align*}
\Prob(\bB_p^0\neq\hat\bB_p^0)&\leq \Prob\Big(\bigcup_{i\leq n,j\leq p}(|w_{ij}-\mu|\geq\delta_n\sqrt{n})\Big)\leq np\cdot \Prob(|w_{ij}-\mu|\geq\delta_n\sqrt{n})\\
&\leq K\delta_n^{-4}\int_{\{|w_{ij}-\mu|\geq\delta_n\sqrt{n}\}}|w_{11}|^4=o(1).
\end{align*}
Let $\hat G_p^0(x)$ be $G_p^0(x)$ with $\bB_p^0$ replaced by $\hat\bB_p^0$, then
$\Prob(\hat{G}_p^0(x)\neq G_p^0(x))\leq \Prob(\bB_p^0\neq\hat\bB_p^0)=o(1)$. In view of the above, we obtain 
$$\int f_j(x)\dif G_p^0(x)=\int f_j(x)\dif \hat{G}_p^0(x)+o_P(1).$$
To simplify notation, we below still use $w_{ij}$ instead of $\hat{w}_{ij}$, and assume that
\begin{align}\label{truncation}
    |w_{ij}-\mu|< \delta_n\sqrt{n},
    \quad \Expe w_{ij}=\mu>0, 
    \quad \Expe|w_{ij}-\mu|^2=\sigma^2,
    \quad \Expe|w_{ij}-\mu|^4<\infty.
\end{align}
\textbf{Step 2: Finite dimensional convergence of $M^{(1)}_p(z)$ in distribution}
\begin{lemma}
    Under conditions and notations in Theorem \ref{clt}, as $p\rightarrow\infty$, for any set of $r$ points $\{z_1, z_2, . . ., z_r\} \bigcup \mathcal{C}$, the random vector
$\bigl(M_p^{(1)}(z_1), \ldots ,M_p^{(1)}(z_r)\bigr)$ converges weakly to a $r$-dimensional centered Gaussian distribution with covariance function in \eqref{mpzcov}.
\end{lemma}
We now proceed to the proof of this lemma. By the fact that a random vector is multivariate normally distributed if and only if every linear combination of its components is normally distributed, we need only show that, for any positive integer $r$ and any complex sequence
${a_j}$, the sum 
\[
    \sum_{j=1}^ra_jM^{(1)}_p(z_j)
\]
converges weakly to a Gaussian random variable. To this end, we first decompose the random part $M^{(1)}_p(z)$ as a sum of martingale difference, which is given in \eqref{mdm1}. Then, we apply the martingale CLT (Lemma \ref{billin}) to obtain the asymptotic distribution of $M^{(1)}_p(z)$.  Details of these two steps are provided in the following two parts.
% Details of these two steps are provided in Sections \ref{mddmp1} and \ref{mcltmd}, respectively.

% \subsubsection{Martingale difference decomposition of $M^1_p(z)$}\label{mddmp1}
\textbf{Part 1: Martingale difference decomposition of $M^{(1)}_p(z)$.} First, we introduce some notations. In the following proof, we assume that $v=\Im z\ge v_0>0$. Moreover, for $j=1,2,\ldots,n$, let
\begin{align}
&\br_j=\frac{1}{\sqrt{n}}\left(\frac{w_{j1}}{\bar{w}_j}-1,\ldots,\frac{w_{jp}}{\bar{w}_j}-1\right)',\;\bD(z)=\bB_p^0-z\bI_p, \;\bD_j(z)=\bD(z)-\br_j\br_j', \nonumber\\
&\beta_j(z)=\frac{1}{1+\br_j'\bD_j^{-1}(z)\br_j}, \ \ \bar{\beta}_j(z)=\frac{1}{1+\frac{1}{n}\nu_2\tr\bD^{-1}_j(z)},\;b_p(z)=\frac{1}{1+\frac{1}{n}\nu_2\Expe \tr\bD_1^{-1}(z)},\nonumber
\end{align}
$\varepsilon_j(z)=\br_j'\bD_j^{-1}(z)\br_j-\frac{1}{n}\nu_{2}\tr\bD_j^{-1}(z)$ and  $\delta_j(z)=\br_j'\bD_j^{-2}(z)\br_j-\frac{1}{n}\nu_2\tr\bD_j^{-2}(z)=\frac{d}{\dif z}\varepsilon_j(z)$. By Lemma \ref{bound}, we have, for any $r\geq 2$,
\begin{align}
\Expe|\varepsilon_j(z)|^r\leq \frac{K}{v^{2r}}n^{-1}\delta_n^{2r-4}\ \ and \ \
\Expe|\delta_j(z)|^r\leq \frac{K}{v^{2r}}n^{-1}\delta_n^{2r-4}.\label{epsirmboun}
% \label{deltrmboun}
\end{align}
It is easy to see that
\begin{align}\label{dminusdj}
\bD^{-1}(z)-\bD_{j}^{-1}(z)=-\bD_{j}^{-1}(z)\br_j\br_j'\bD_{j}^{-1}(z)\beta_{j}(z),
\end{align}
where we use the formula that $\bA_1^{-1}-\bA_2^{-1}=\bA_2^{-1}(\bA_2-\bA_1)\bA_1^{-1}$ holds for any two invertible matrices $\bA_1$ and $\bA_2$. Note that $|\beta_j(z)|$, $|\bar{\beta}_j(z)|$ and $|b_n(z)|$ are bounded by $\frac{|z|}{v}$. We also get that  for any $j$,
\begin{align}\label{covmatrixmain}
    \Expe (\br_j\br_j') = \frac{1}{n}\nu_{2} \left(-\frac{1}{p-1}\boldsymbol{1}_{p}\boldsymbol{1}'_{p} + \frac{1}{p-1}\bI_p + \bI_p \right).
\end{align}

Let $\Expe_0(\cdot)$ denote expectation and $\Expe_j(\cdot)$ denote conditional expectation with respect to the $\sigma$-field generated by $\br_1, \br_2, \ldots, \br_j$, where $j=1,2,\ldots,n$. Next, we write $M_p^{(1)}(z)$ as a sum of martingale difference sequences (MDS), and then utilize the CLT of MDS (Lemma \ref{billin}) to derive the asymptotic distribution of $M_p^{(1)}(z)$, which can be written as
\begin{align}\label{m1}
% M_p^{(1)}(z)
% &= 
p[m_{p}(z)-\Expe M_{p}(z)]
% & =\tr[\bD^{-1}(z)-\Expe\bD^{-1}(z)]\nonumber\\
= \sum^{n}_{j=1}[\tr(\Expe_j-\Expe_{j-1})\bD^{-1}(z)]= -\sum^{n}_{j=1}(\Expe_j-\Expe_{j-1})\beta_j(z)\br_j'\bD_j^{-2}(z)\br_j.
% = \sum^{n}_{j=1}[\tr\Expe_j\bD^{-1}(z)-\tr\Expe_{j-1}\bD^{-1}(z)]= -\sum^{n}_{j=1}(\Expe_j-\Expe_{j-1})\beta_j(z)\br_j'\bD_j^{-2}(z)\br_j.
% &= \ \sum^{n}_{j=1}\Big(\tr\Expe_j[\bD^{-1}(z)-\bD^{-1}_j(z)]-\tr\Expe_{j-1}[\bD^{-1}(z)-\bD_j^{-1}(z)]\Big)\nonumber\\
% &= \ -\sum^{n}_{j=1}(\Expe_j-\Expe_{j-1})\beta_j(z)\br_j'\bD_j^{-2}(z)\br_j.
\end{align}
Write
$\beta_j(z)=\bar{\beta}_j(z)-\beta_j(z)\bar{\beta}_j(z)\varepsilon_j(z)
=\bar{\beta}_j(z)-\bar{\beta}_j^2(z)\varepsilon_j(z)+\bar{\beta}_j^2(z)\beta_j(z)\varepsilon_j^2(z)$.
From this and the definition of $\delta_j(z)$, (\ref{m1}) has the following expression
\begin{align}\label{m1ejej-1brdr}
&(\Expe_j-\Expe_{j-1})\beta_j(z)\br_j'\bD_j^{-2}(z)\br_j = (\Expe_j-\Expe_{j-1})\big[\big(\bar{\beta}_j(z)-\bar{\beta}_j^2(z)\varepsilon_j(z)
\nonumber\\
&\ \ \ \ \ \ \ \ \ \ \ \ \ \ \ +\bar{\beta}_j^2(z)\beta_j(z)\varepsilon_j^2(z)\big)
\big(\delta_j(z)+\frac{1}{n}\nu_2\tr\bD_j^{-2}(z)\big)\big] \nonumber\\
% =&  (\Expe_j-\Expe_{j-1})\big[\bar{\beta}_j(z)\delta_j(z)
% \nonumber\\
% &\ \ \ \ \ \ \ \ \ - \ \bar{\beta}_j^2(z)\varepsilon_j(z)\delta_j(z)-\bar{\beta}_j^2(z)\varepsilon_j(z)\frac{1}{n}\nu_2tr\bD_j^{-2}(z)+\bar{\beta}_j^2(z)\beta_j(z)\varepsilon_j^2(z)\br_j'\bD_j^{-2}(z)\br_j\big]\nonumber\\
% =& \ \Expe_j\big[\bar{\beta}_j(z)\delta_j(z)-\bar{\beta}_j^2(z)\varepsilon_j(z)\frac{1}{n}\nu_2tr\bD_j^{-2}(z)\big] -  \Expe_{j-1}\big[\bar{\beta}_j(z)\delta_j(z)-\bar{\beta}_j^2(z)\varepsilon_j(z)\frac{1}{n}\nu_2tr\bD_j^{-2}(z)\big]
% \nonumber\\
% & \ \ \ \ \ \ \ \ \ - \ (\Expe_j-\Expe_{j-1})\big[\bar{\beta}_j^2(z)\big(\varepsilon_j(z)\delta_j(z)-\beta_j(z)\br_j'\bD_j^{-2}(z)\br_j\varepsilon_j^2(z)\big)\big],\nonumber\\
=& -Y_j(z) + \Expe_{j-1} Y_j(z)  - (\Expe_j-\Expe_{j-1})\big[\bar{\beta}_j^2(z)\big(\varepsilon_j(z)\delta_j(z)-\beta_j(z)\br_j'\bD_j^{-2}(z)\br_j\varepsilon_j^2(z)\big)\big],
\end{align}
where the second equality uses the fact that $(\Expe_j-\Expe_{j-1})\bar{\beta}_j(z)\tr\bD_j^{-2}(z)=0$, and
\begin{eqnarray*}
Y_j(z)=-\Expe_j\Big(\bar{\beta}_j(z)\delta_j(z)-\bar{\beta}_j^2(z)\varepsilon_j(z)\frac{1}{n}\nu_2\tr\bD_j^{-2}(z)\Big)=-\Expe_j\frac{\dif}{\dif z}\big(\bar{\beta}_j(z)\varepsilon_j(z)\big).
\end{eqnarray*}
By (\ref{epsirmboun}), we have
\begin{align}\label{m1smallterm1}
\Expe\Big|\sum^{n}_{j=1}(\Expe_j-\Expe_{j-1})\bar{\beta}_j^2(z)\varepsilon_j(z)\delta_j(z)\Big|^2
% &=\sum^{n}_{j=1}\Expe\big|(\Expe_j-\Expe_{j-1})\bar{\beta}_j^2(z)\varepsilon_j(z)\delta_j(z)\big|^2\nonumber\\
&\leq 4\sum^{n}_{j=1}\Expe\big|\bar{\beta}_j^2(z)\varepsilon_j(z)\delta_j(z)\big|^2=o(1),
\end{align}
here we leverage the the martingale difference property of $(\Expe_j-\Expe_{j-1})\bar{\beta}_j^2(z)\varepsilon_j(z)\delta_j(z)$.
Thus, $\sum^{n}_{j=1}(\Expe_j-\Expe_{j-1})\bar{\beta}_j^2(z)\varepsilon_j(z)\delta_j(z)$ converges to $0$ in probability. By the same argument, we have
\begin{eqnarray}\label{m1smallterm2}
\sum^{n}_{j=1}(\Expe_j-\Expe_{j-1})\bar{\beta}_j^2(z)\beta_j(z)\br_j'\bD_j^{-2}(z)\br_j\varepsilon_j^2(z)\stackrel{i.p.}{\rightarrow}0.
\end{eqnarray}
Then, equations \eqref{m1} -- \eqref{m1smallterm2} imply that
\begin{align}\label{mdm1}
    M^{(1)}_p(z) = \sum_{j=1}^n \{Y_j(z) - \Expe_{j-1} Y_j(z)\}+o_P(1),
\end{align}
where $\left\{Y_j(z) - \Expe_{j-1} Y_j(z),j=1,\ldots,n\right\}$ is a sequence of martingale difference.

% \subsubsection{Application of martingales CLT to \eqref{mdm1}}\label{mcltmd}

\textbf{Part 2: Application of martingales CLT to \eqref{mdm1}.} 
To prove finite-dimensional convergence
of $M_p^{(1)}(z)$, $z \in \mathcal{C}$, we need only to consider the limit of the following martingale difference decomposition:
\begin{eqnarray*}
\sum_{i=1}^r\alpha_i M^{(1)}_p(z_i) =\sum^{r}_{i=1}\alpha_i\sum^{n}_{j=1}(Y_j(z_i)-\Expe_{j-1} Y_j(z_i))+o(1)=\sum^{n}_{j=1}\sum^{r}_{i=1}\alpha_i(Y_j(z_i)-\Expe_{j-1} Y_j(z_i))+o(1),
\end{eqnarray*}
where $\Im(z_i)\neq 0$, $\{\alpha_i: i=1,2,\ldots,r\}$ are constants. We apply the martingale
CLT (Lemma \ref{billin}) to this martingale difference decomposition of $\sum_{i=1}^r\alpha_i M^{(1)}_p(z_i)$. To this
end, we need to check two conditions:
\begin{Condition}\label{condi1}
\begin{eqnarray}\label{condi1eq}
&&\qquad\sum^{n}_{j=1}\Expe\left(\left|\sum^{r}_{i=1}\alpha_i\left
(Y_j(z_i)-\Expe_{j-1}Y_j(z_i)\right)\right|^2I_{\{|\sum^{r}_{i=1}\alpha_i\left(Y_j(z_i)+\Expe_{j-1}Y_j(z_i)\right)|\geq\varepsilon\}}\right)\rightarrow 0.
\end{eqnarray}
\end{Condition}
\begin{Condition}\label{condi2}
\begin{eqnarray}\label{condi2eq}
&\sum^{n}_{j=1}\Expe_{j-1}\bigl[\bigl(Y_j(z_1)-\Expe_{j-1} Y_j(z_1)\bigr)\bigl(Y_j(z_2)-\Expe_{j-1} Y_j(z_2)\bigr)\bigr]
\end{eqnarray}
converges in probability to a constant.
\end{Condition}
First, we verify Condition \ref{condi1}. By Lemma \ref{bound}, we obtain
\begin{eqnarray}\label{y4bound}
\Expe|Y_j(z)|^4\leq K
\Expe|\varepsilon_j(z)|^4=o\left(\frac{1}{p}\right).
\end{eqnarray}
Furthermore, by Jensen's inequality and \eqref{y4bound},
\begin{align}\label{condy4bound}
    \Expe|\Expe_{j-1}Y_j(z)|^4 \leq \Expe(\Expe_{j-1}|Y_j(z)|^4)= \Expe|Y_j(z)|^4=o\left(\frac{1}{p}\right).
\end{align}
It follows from (\ref{y4bound}) and \eqref{condy4bound} that
\begin{align*}
\text{the left hand side of}\ \ \eqref{condi1eq}
% \leq &\; \frac{1}{\varepsilon^2}\sum^{n}_{j=1}\Expe\left|\sum^{r}_{i=1}\alpha_i\left(Y_j(z_i)-\Expe_{j-1}Y_j(z_i)\right)\right|^4\nonumber\\
\leq &\; K\left(\frac{1}{\varepsilon^2}\sum^{n}_{j=1}\Expe\left|\sum^{r}_{i=1}\alpha_iY_j(z_i)\right|^4+\frac{1}{\varepsilon^2}\sum^{n}_{j=1}\Expe\left|\sum^{r}_{i=1}\alpha_i\Expe_{j-1}Y_j(z_i)\right|^4\right)\rightarrow 0.
\end{align*}
Then, we verify Condition \ref{condi2}.
Since 
\begin{align*}
    \eqref{condi2eq}&=\sum^{n}_{j=1}\Expe_{j-1}[Y_j(z_1)Y_j(z_2)]-\sum^{n}_{j=1}[\Expe_{j-1}Y_j(z_1)][\Expe_{j-1}Y_j(z_2)]\\
    &=\frac{\partial^2}{\partial z_1\partial z_2}\Big(\sum^{n}_{j=1}\Expe_{j-1}\big[\Expe_j\big(\bar{\beta}_j(z_1)\varepsilon_j(z_1)\big)\Expe_j\big(\bar{\beta}_j(z_2)\varepsilon_j(z_2)\big)\big]\Big)\\
    &\ \ \ \ -\frac{\partial^2}{\partial z_1\partial z_2}\Big(\sum^{n}_{j=1}\big[\Expe_{j-1}\bar{\beta}_j(z_1)\varepsilon_j(z_1)][\Expe_{j-1}\bar{\beta}_j(z_2)\varepsilon_j(z_2)\big]\Big),
\end{align*}
% \begin{eqnarray*}
% \frac{\partial^2}{\partial z_1\partial z_2}\Big(\sum^{n}_{j=1}\Expe_{j-1}\big[\Expe_j\big(\bar{\beta}_j(z_1)\varepsilon_j(z_1)\big)\Expe_j\big(\bar{\beta}_j(z_2)\varepsilon_j(z_2)\big)\big]\Big)=\sum^{n}_{j=1}\Expe_{j-1}[Y_j(z_1)Y_j(z_2)],
% \end{eqnarray*}
% and
% \begin{eqnarray*}
% \frac{\partial^2}{\partial z_1\partial z_2}\Big(\sum^{n}_{j=1}\big[\Expe_{j-1}\bar{\beta}_j(z_1)\varepsilon_j(z_1)][\Expe_{j-1}\bar{\beta}_j(z_2)\varepsilon_j(z_2)\big]\Big)=\sum^{n}_{j=1}[\Expe_{j-1}Y_j(z_1)][\Expe_{j-1}Y_j(z_2)],
% \end{eqnarray*}
% and by the same arguments as those on page 571 of \cite{bai2004clt},
it is enough to consider the limits of
\begin{eqnarray}\label{precltlim}
\sum^{n}_{j=1}\Expe_{j-1}\big[\Expe_j\big(\bar{\beta}_j(z_1)\varepsilon_j(z_1)\big)\Expe_j\big(\bar{\beta}_j(z_2)\varepsilon_j(z_2)\big)\big]
\end{eqnarray}
and
\begin{align}\label{smallprecltlim}
    \sum^{n}_{j=1}\big[\Expe_{j-1}\bar{\beta}_j(z_1)\varepsilon_j(z_1)][\Expe_{j-1}\bar{\beta}_j(z_2)\varepsilon_j(z_2)\big].
\end{align}
The limit of \eqref{smallprecltlim} is provided in the following lemma.
\begin{lemma}\label{ej-1term0}
Under conditions and notations in Theorem \ref{clt}, then
\begin{align*}
     \sum^{n}_{j=1}\big[\Expe_{j-1}\bar{\beta}_j(z_1)\varepsilon_j(z_1)][\Expe_{j-1}\bar{\beta}_j(z_2)\varepsilon_j(z_2)\big]\stackrel{i.p.}{\rightarrow}0.
\end{align*}
%converges in probability to 0.
\end{lemma}
The proof of Lemma \ref{ej-1term0} is postponed to Section \ref{prfej-1term0} of the supplementary document. By Lemma \ref{ej-1term0}, the remaining work is to consider the limit of \eqref{precltlim}. Since the following inequalities hold:
\begin{align}
    \big|\tr\big(\bD^{-1}(z)-\bD_j^{-1}(z)\big)\bA\big|&\leq\frac{\|\bA\|}{\Im(z)},\label{trda}\\
    \Expe\big|\frac{1}{n}\nu_2\tr\bD^{-1}(z)-\Expe\frac{1}{n}\nu_2\tr\bD^{-1}(z)\big|^{q}&\leq C_qn^{-q/2}v_0^{-q},\label{ded}\\
    \mathbb{E}\big|\bar{\beta}_j(z_i)-b_n(z_i)\big|^2&\leq \frac{K|z_i|^4}{nv_0^6},\label{betabpbound}
\end{align}
% For any $p\times p$ matrix $\bA$, 
% \begin{equation}\label{trda}
% \big|\tr\big(\bD^{-1}(z)-\bD_j^{-1}(z)\big)\bA\big|\leq\frac{\|\bA\|}{\Im(z)}.
% \end{equation}
% By Lemma \ref{burk} and \eqref{trda}, similarly to (4.3) in \cite{Bai1998No}, 
% \begin{align}\label{ded}
%     \Expe\big|\frac{1}{n}\nu_2\tr\bD^{-1}(z)-\Expe\frac{1}{n}\nu_2\tr\bD^{-1}(z)\big|^{q}\leq C_qn^{-q/2}v_0^{-q},
% \end{align}
% thus, by \eqref{ded} we obtain
% \begin{eqnarray}\label{betabpbound}
% \mathbb{E}\big|\bar{\beta}_j(z_i)-b_n(z_i)\big|^2\leq \frac{K|z_i|^4}{nv_0^6}.
% \end{eqnarray}
% By (\ref{betabpbound}), we have
% \begin{eqnarray*}
% \Expe\big|
% \Expe_{j-1}[\mathbb{E}_j\big(\bar{\beta}_j(z_1)\varepsilon_j(z_1)\big)\Expe_j\big(\bar{\beta}_j(z_2)\varepsilon_j(z_2)\big)]
% -\Expe_{j-1}[\Expe_j\big(b_p(z_1)\varepsilon_j(z_1)\big)\Expe_j\big(b_p(z_2)\varepsilon_j(z_2)\big)]\big|
% =O(n^{-3/2}).
% \end{eqnarray*}
% From this, it follows that
% \begin{eqnarray*}
% \sum^{n}_{j=1}\Expe_{j-1}\big[\Expe_j\big(\bar{\beta}_j(z_1)\varepsilon_j(z_1)\big)\Expe_j\big(\bar{\beta}_j(z_2)\varepsilon_j(z_2)\big)\big]
% -b_p(z_1)b_p(z_2)\sum^{n}_{j=1}\Expe_{j-1}\big[\Expe_j\big(\varepsilon_j(z_1)\big)\Expe_j\big(\varepsilon_j(z_2)\big)\big]
% \stackrel{i.p.}{\rightarrow}0.
% \end{eqnarray*}
it is enough to prove that
\begin{eqnarray}\label{bprecltlim}
b_p(z_1)b_p(z_2)\sum^{n}_{j=1}\Expe_{j-1}\big[\Expe_j\big(\varepsilon_j(z_1)\big)\Expe_j\big(\varepsilon_j(z_2)\big)\big]
\end{eqnarray}
converges to a constant in probability, which further gives the limit of (\ref{precltlim}). By Lemma \ref{quadform}, we have
\begin{align}
     % &\qquad b_p(z_1)b_p(z_2)\sum_{j=1}^{n}\Expe_{j-1}\Bigl[ \Expe_j(\varepsilon_j(z_1))\Expe_j(\varepsilon_j(z_2))\Bigr]\nonumber\\
    \eqref{bprecltlim}&= b_p(z_1)b_p(z_2)\sum_{j=1}^{n}\biggl[\sum_{i=1}^{p}\frac{1}{n^{2}}(\nu_4 - 3\nu_{12} ) \Expe_j(\bD_j^{-1}(z_1))_{i i}\Expe_j(\bD_j^{-1}(z_2))_{i i}\nonumber\\
    &+\frac{1}{n^{2}}\nu_{12}\Bigl(\tr \left[\Expe_j\bD_j^{-1}(z_1)\Expe_j(\bD_j^{-1}(z_2))'\right]+\tr \left[\Expe_j\bD_j^{-1}(z_1)\Expe_j\bD_j^{-1}(z_2)\right]\Bigr)\nonumber\\
    &+\frac{1}{n^{2}}(\nu_{12}-\nu_2^2) \tr  \left[\Expe_j\bD_j^{-1}(z_1)\right] \tr \left[\Expe_j\bD_j^{-1}(z_2)\right] \biggr]+o(1)
    =:I_1 + I_2 + I_3+ I_4 + o(1),
\end{align}
where
\begin{align}
    I_1 &= \frac{1}{n^{2}}(\nu_4-3\nu_{12})b_p(z_1)b_p(z_2)\sum_{j=1}^{n}\sum_{i=1}^{p} \Expe_j(\bD_j^{-1}(z_1))_{i i}\Expe_j(\bD_j^{-1}(z_2))_{i i},\nonumber\\
    I_2 &= \frac{1}{n^{2}}\nu_{12}b_p(z_1)b_p(z_2)\sum_{j=1}^{n}\tr \left[\Expe_j\bD_j^{-1}(z_1)\Expe_j(\bD_j^{-1}(z_2))'\right],\nonumber\\
    I_3 &= \frac{1}{n^{2}}\nu_{12}b_p(z_1)b_p(z_2)\sum_{j=1}^{n}\tr \left[\Expe_j\bD_j^{-1}(z_1)\Expe_j\bD_j^{-1}(z_2)\right],\nonumber\\
    I_4 &= \frac{1}{n^{2}}(\nu_{12}-\nu_2^2) b_p(z_1)b_p(z_2)\sum_{j=1}^{n}\tr \left[\Expe_j\bD_j^{-1}(z_1)\right] \tr \left[\Expe_j\bD_j^{-1}(z_2)\right].\nonumber
\end{align}
In the following Steps (i)-(iii), we derive as $p\to\infty$,
\begin{align}
&\frac{\partial^2}{\partial z_1\partial z_2}I_1\stackrel{i.p.}{\rightarrow}c\alpha_1\frac{\underline{m}'(z)\underline{m}'(z_2)}{\big(1+\frac{\sigma^2}{\mu^2}\underline{m}(z_1)\big)^2\big(1+\frac{\sigma^2}{\mu^2}\underline{m}(z_2)\big)^2},\label{parI1lim}\\
&
\frac{\partial^2}{\partial z_1\partial z_2}I_2,\ \frac{\partial^2}{\partial z_1\partial z_2}I_3\stackrel{i.p.}{\rightarrow}\frac{\partial}{\partial z_2}\big(\frac{\partial a(z_1,z_2)/\partial z_1}{1-a(z_1,z_2)}\big)= \frac{\underline{m}'(z_1)\underline{m}'(z_2)}{(\underline{m}(z_1)-\underline{m}(z_2))^2} - \frac{1}{(z_1-z_2)^2},\label{parI3lim}\\
&\frac{\partial^2}{\partial z_1\partial z_2}I_4\stackrel{i.p.}{\rightarrow}c\Bigl(h_2-2\frac{\sigma^2}{\mu^2}h_1\Bigr)\frac{\underline{m}'(z_1)\underline{m}'(z_2)}{\big(1+\frac{\sigma^2}{\mu^2}\underline{m}(z_1)\big)^2\big(1+\frac{\sigma^2}{\mu^2}\underline{m}(z_2)\big)^2}.\label{parI4lim}
\end{align}

\underline{Step (i)}: \textbf{Consider $I_2$ and $I_3$.} Let $\bD_{ij}(z)=\bD(z)-\br_i\br_i'-\br_j\br_j'$, $
    b_1(z)=\frac{1}{1+\frac{1}{n}\nu_2\Expe \tr\bD_{12}^{-1}(z)}$ and $\beta_{ij}(z)=\frac{1}{1+\br_i'\bD_{ij}^{-1}(z)\br_i}$.
We have the equality 
    $\bD_j(z_1)+z_1\bI_p-\frac{n-1}{n}\nu_2b_1(z_1)\bI_p=\sum^{n}_{i\neq j}\br_i\br_i'-\frac{n-1}{n}\nu_2b_1(z_1)\bI_p$.
Multiplying by $(z_1\bI_p-\frac{n-1}{n}\nu_2b_1(z)\bI_p)^{-1}$ on the left-hand side and $\bD_j^{-1}(z_1)$ on the right-hand side, and using
$\br_i'\bD_j^{-1}(z_1)=\beta_{ij}(z_1)\br_i'\bD_{ij}^{-1}(z_1)$, we get
\begin{align}
\bD_j^{-1}(z_1)&=-\bQ_p(z_1)+\sum_{i\neq j}^{n}\beta_{ij}(z_1)\bQ_p(z_1)\br_i\br_i'\bD_{ij}^{-1}(z_1)
-\frac{n-1}{n}\nu_2b_1(z_1)\bQ_p(z_1)\bD_j^{-1}(z_1)\nonumber\\
% &=-\bQ_p(z_1)+\sum_{i\neq j}^{n}b_1(z_1)\bQ_p(z_1)\br_i\br_i'\bD_{ij}^{-1}(z_1)+\sum_{i\neq j}^{n}(\beta_{ij}(z_1)-b_1(z_1))\bQ_p(z_1)\br_i\br_i'\bD_{ij}^{-1}(z_1)\nonumber\\
% &&-\frac{n-1}{n}\nu_2b_1(z_1)\bQ_p(z_1)\bD_{ij}^{-1}(z_1)+\frac{n-1}{n}\nu_2b_1(z_1)\bQ_p(z_1)(\bD_{ij}^{-1}(z_1)-\bD_{j}^{-1}(z_1))\nonumber\\
&=-\bQ_p(z_1)+b_1(z_1)\bA(z_1)+\bB(z_1)+\bC(z_1),\label{djdecomp}
\end{align}
where $\bQ_p(z_1)=\big(z_1\bI_p-\frac{n-1}{n}\nu_2b_1(z_1)\bI_p\big)^{-1}$,
$
\bA(z_1)=\sum_{i\neq j}^{n}\bQ_p(z_1)\big(\br_i\br_i'-\frac{1}{n}\nu_2\bI_p\big)\bD_{ij}^{-1}(z_1)$,
$\bB(z_1)=\sum_{i\neq j}^{n}\big(\beta_{ij}(z_1)-b_1(z_1)\big)\bQ_p(z_1)\br_i\br_i'\bD_{ij}^{-1}(z_1)$, $\bC(z_1)=\frac{1}{n}\nu_2b_1(z_1)\bQ_p(z_1)\sum_{i\neq j}^{n}\big(\bD_{ij}^{-1}(z_1)-\bD_j^{-1}(z_1)\big)$.
For any real $t$, 
$
\biggl|1-\frac{t}{z(1+n^{-1}\nu_2\Expe \tr\bD_{12}^{-1}(z))}\biggr|^{-1} 
\leq\frac{\big|z\big(1+n^{-1}\nu_2\Expe \tr\bD_{12}^{-1}(z)\big)\big|}{\Im\big(z(1+n^{-1}\nu_2\Expe \tr\bD_{12}^{-1}(z))\big)}\leq \frac{|z|(1+p/(nv_0))}{v_0}$.
Thus,
\begin{eqnarray}\label{Qbound}
\|\bQ_p(z_1)\|\leq\frac{1+p/(nv_0)}{v_0}.
\end{eqnarray}
For any random matrix $\bM$, denote its nonrandom bound on the spectrum norm of $\bM$ by $|||\bM|||$.
From \eqref{betabpbound}, Lemma \ref{bound}, \eqref{Qbound} and (\ref{trda}), we get, for any $\bM$,
\begin{align}
\Expe\big|\tr\bB(z_1)\bM\big|
% &\leq n\Expe^{1/2}\big(\big|\beta_{12}(z_1)-b_1(z_1)\big|^2\big)\cdot\Expe^{1/2}\big(\big|\br_i'\bD_{ij}^{-1}(z_1)\bM \bQ_p(z_1)\br_i\big|^2\big)\nonumber\\
&\leq K|||\bM|||\frac{|z_1|^2(1+p/(nv_0))}{v_0^5}n^{1/2},
% \label{bmbound},
\ \
\big|\tr\bC(z_1)\bM\big| 
% &\leq n^{-1}\frac{\sigma^2}{\mu^2}\sum_{i\neq j}^{n}\frac{\|\bM b_1(z_1)\bQ_p(z_1)\|}{v_0}\nonumber\\
\leq |||\bM|||\frac{|z_1|(1+p/(nv_0))}{v_0^3},\label{cmbound}\\
\Expe\big|\tr\bA(z_1)\bM\big|
% &\leq(\Expe\big|\tr\bA(z_1)\bM\big|^2)^{1/2}\nonumber\\
% &\leq (n\Expe\sum_{i\neq j}^{n}\big|\operatorname{tr}\bQ_p(z_1)\big(\br_i\br_i'-\frac{1}{n}\frac{\sigma^2}{\mu^2}\bI_p\big)\bD_{ij}^{-1}(z_1)\bM\big|^2)^{1/2}\nonumber\\
% &=(n\sum_{i\neq j}^{n}\Expe\big|\br_i'\bD_{ij}^{-1}(z_1)\bM\bQ_p(z_1)\br_i-\frac{1}{n}\frac{\sigma^2}{\mu^2}\operatorname{tr}\bD_{ij}^{-1}(z_1)\bM\bQ_p(z_1)\big|^2)^{1/2}\nonumber\\
% &\leq Kn^{1/2}\Expe^{1/2}\big(\tr\bD_{ij}^{-1}(z_1)\bM\bQ_p(z_1)\bQ_p(\bar{z}_1)\bM'\bD_{ij}^{-1}(\bar z_1)\big)\nonumber\\
&\leq K\|\bM\|\frac{1+p/(nv_0)}{v_0^2}n^{1/2}.\label{ambound}
\end{align}
Note that
\begin{align}   \operatorname{tr}\Expe_j\big(\bA(z_1)\big)\bD_j^{-1}(z_2) 
% &= \operatorname{tr}\Bigl[\Expe_j\Big(\sum_{i\neq j}^{n}\bQ_p(z_1)\big(\br_i\br_i'-\frac{1}{n}\nu_2\bI_p\big)\bD_{ij}^{-1}(z_1)\Big)\bD_j^{-1}(z_2)\Bigr]\nonumber\\
% &= \operatorname{tr}\Bigl[\Expe_j\Big(\sum_{i< j}^{n}\bQ_p(z_1)\big(\br_i\br_i'-\frac{1}{n}\nu_2\bI_p\big)\bD_{ij}^{-1}(z_1)\Big)\bD_j^{-1}(z_2)\Bigr]\nonumber\\
% &\ \ \ \ +\operatorname{tr}\Bigl[\Expe_j\Big(\sum_{i> j}^{n}\bQ_p(z_1)\big(\br_i\br_i'-\frac{1}{n}\nu_2\bI_p\big)\bD_{ij}^{-1}(z_1)\Big)\bD_j^{-1}(z_2)\Bigr]\nonumber\\
=&\operatorname{tr}\sum_{i<j}^{n}\bQ_p(z_1)\big(\br_i\br_i'-n^{-1}\nu_2\bI_p\big)\Expe_j\big(\bD_{ij}^{-1}(z_1)\big)\bD_{ij}^{-1}(z_2)\nonumber\\
&+ \operatorname{tr}\sum_{i<j}^{n}\bQ_p(z_1)\big(\br_i\br_i'-n^{-1}\nu_2\bI_p\big)\Expe_j\big(\bD_{ij}^{-1}(z_1)\big)(\bD_{j}^{-1}(z_2)-\bD_{ij}^{-1}(z_2))\nonumber\\
&+\operatorname{tr}\,\Expe_j\Big(\sum_{i> j}^{n}\bQ_p(z_1)\big(\br_i\br_i'-\frac{1}{n}\nu_2\bI_p\big)\bD_{ij}^{-1}(z_1)\Big)\bD_j^{-1}(z_2),\nonumber
\end{align}
therefore, by using \eqref{dminusdj},  we can write 
\begin{align}  
\tr\Expe_j\big(\bA(z_1)\big)\bD_j^{-1}(z_2)=A_1(z_1,z_2)+A_2(z_1,z_2)+A_3(z_1,z_2)+R(z_1,z_2)\label{trejadjdecomp},
\end{align}
where
\begin{align}
A_1(z_1,z_2)
% &=-\operatorname{tr}\sum_{i<j}^{n}\bQ_p(z_1)\br_i\br_i'\Expe_j\big(\bD_{ij}^{-1}(z_1)\big)\beta_{ij}(z_2)\bD_{ij}^{-1}(z_2)\br_i\br_i'\bD_{ij}^{-1}(z_2)\nonumber\\
&=-\sum_{i<j}^{n}\beta_{ij}(z_2)\br_i'\Expe_j\big(\bD_{ij}^{-1}(z_1)\big)\bD_{ij}^{-1}(z_2)\br_i\br_i'\bD_{ij}^{-1}(z_2)\bQ_p(z_1)\br_i,\label{eq:A1z1z2}\\
A_2(z_1,z_2)&=-\operatorname{tr}\sum_{i<j}^{n}\bQ_p(z_1)n^{-1}\nu_2\Expe_j\big(\bD_{ij}^{-1}(z_1)\big)\big(\bD_j^{-1}(z_2)-\bD_{ij}^{-1}(z_2)\big), \nonumber \\
A_3(z_1,z_2)&=\operatorname{tr}\sum_{i<j}^{n}\bQ_p(z_1)\big(\br_i\br_i'-n^{-1}\nu_2\bI_p\big)\Expe_j\big(\bD_{ij}^{-1}(z_1)\big)\bD_{ij}^{-1}(z_2),\nonumber\\ 
R(z_1 , z_2)&=\operatorname{tr}\sum_{i>j}^{n}\bQ_p(z_1)\Bigl(-\frac{1}{n(p-1)}\nu_2\boldsymbol{1}_{p}\boldsymbol{1}_{p}'+\frac{1}{n(p-1)}\nu_2\bI_p\Bigr)\Expe_{j}\bigl(\bD_{i j}^{-1}(z_1)\bigr)\bD_{j}^{-1}(z_2),\label{Rterm}
\end{align}
and $\boldsymbol{1}_{p}$ is a $p$-dimensional vector with all elements being $1$.
It is easy to see that $R(z_1 , z_2)=O(1)$.
We get from (\ref{trda}) and (\ref{Qbound}) that $|A_2(z_1,z_2)|\leq\frac{1+p/(nv_0)}{v_0^2}$. Similar to (\ref{ambound}), we have $\Expe|A_3(z_1,z_2)|\leq\frac{1+p/(nv_0)}{v_0^3}n^{1/2}$. By similar calculation of \cite{bai2004clt}, we get the following lemma and its proof is postponed to Section \ref{prftrejdjdjeqmain} of the supplementary document. 
\begin{lemma}\label{trejdjdjeqmain}
Under conditions and notations in Theorem \ref{clt}, for any $1\leq j \leq n$,
\begin{align}
&\;\tr\big(\Expe_j\big(\bD_j^{-1}(z_1)\big)\bD_j^{-1}(z_2)\big)\nonumber\\
&\qquad\times\biggl[1-\frac{j-1}{n}\nu_2^2\underline{m}_{p}^{0}(z_1)\underline{m}_{p}^{0}(z_2)\frac{c_n}{(1+\frac{n-1}{n}\nu_2\underline{m}_p^0(z_1))(1+\frac{n-1}{n}\nu_2\underline{m}_p^0(z_2))}\biggr]\nonumber\\
=&\;\frac{nc_n}{z_1z_2}\frac{1}{(1+\frac{n-1}{n}\nu_2\underline{m}_p^0(z_1))(1+\frac{n-1}{n}\nu_2\underline{m}_p^0(z_2))}+S(z_1,z_2),\label{eq:trD1D2_limit}
\end{align}
where $\Expe|S(z_1,z_2)|\leq Kn^{1/2}$.
\end{lemma}
By Lemma \ref{trejdjdjeqmain}, $I_3$ can be written as 
\begin{equation}\label{trdjdjeq2}
    I_3=a_p(z_1,z_2)\frac{1}{n}\sum^{n}_{j=1}\frac{1}{1-\frac{j-1}{n}a_p(z_1,z_2)}+A_4(z_1,z_2),
\end{equation}
where
$a_p(z_1,z_2)=\nu_2^2\frac{c_n\underline{m}_{p}^{0}(z_1)\underline{m}_{p}^{0}(z_2)}{(1+\frac{n-1}{n}\nu_2\underline{m}_p^0(z_1))(1+\frac{n-1}{n}\nu_2\underline{m}_p^0(z_2))}$ and  $\Expe|A_4(z_1,z_2)|\leq Kn^{-1/2}$.
By Lemma \ref{estsqu}, the limit of $a_p(z_1, z_2)$ is 
$a(z_1,z_2)=\frac{\sigma^4}{\mu^4}\frac{c \underline{m}(z_1)\underline{m}(z_2)}{\big(1+\frac{\sigma^2}{\mu^2}\underline{m}(z_1)\big)\big(1+\frac{\sigma^2}{\mu^2}\underline{m}(z_2)\big)}$.
Thus, by (\ref{trdjdjeq2}), the in probability (i.p.) limit of $\frac{\partial^2}{\partial z_2\partial z_1}I_3$ is in \eqref{parI3lim}.
% \begin{align}\label{parI3lim}
% \frac{\partial^2}{\partial z_2\partial z_1}\int^{a(z_1,z_2)}_{0}\frac{1}{1-z}\dif z=\frac{\partial}{\partial z_2}\big(\frac{\partial a(z_1,z_2)/\partial z_1}{1-a(z_1,z_2)}\big)= \frac{\underline{m}'(z_1)\underline{m}'(z_2)}{(\underline{m}(z_1)-\underline{m}(z_2))^2} - \frac{1}{(z_1-z_2)^2}.
% \end{align}
Similarly, we get the i.p. limit of $\frac{\partial^2}{\partial z_2\partial z_1}I_2$, which is also given by \eqref{parI3lim}.

\underline{Step (ii)}: \textbf{Consider $I_1$.} 
It is enough to find the limit of $\frac{1}{n^{2}}\sum_{j=1}^{n} \sum_{i=1}^{p} \Expe_j(D_j^{-1}(z_1))_{i i}\Expe_j(D_j^{-1}(z_2))_{i i}$. By similar calculation of \cite{gao2017high}, we get the following lemma and its proof is postponed to Section \ref{prfm1I1ejdj} of the supplementary document.
\begin{lemma}\label{m1I1ejdj}
Under conditions and notations in Theorem \ref{clt}, for any $1 \leq j\leq n$,
\begin{align*}
    \frac{1}{p}\sum_{i=1}^{p} \Expe_j(D_j^{-1}(z_1))_{i i}\Expe_j(D_j^{-1}(z_2))_{i i}\stackrel{i.p.}{\rightarrow}m(z_1)m(z_2).
\end{align*}
\end{lemma}
By \eqref{ded}, the formula (2.2) of \cite{Silverstein1995StrongCO}, $
\underline{m}_p(z)=-\frac{1}{zn}\sum^{n}_{j=1}\beta_j(z)
$, and Lemma \ref{bound}, we have 
\begin{align}
&|b_p(z)-\Expe\beta_1(z)|\leq Kn^{-1/2},\ \ 
% \label{bpbetaboun}\\
\Expe\beta_1(z)
=-z\mathbb{E}\underline{m}_p(z),\ \
% \label{betazm}\\
|b_p(z)+z\underline{m}_{p}^{0}(z)|\leq Kn^{-1/2}.\label{b1zm}
\end{align}
Thus, by \eqref{b1zm}, Lemma \ref{estsqu} and Lemma \ref{m1I1ejdj}, we have
\begin{align*}
    I_1\stackrel{i.p.}{\rightarrow} &\;c\alpha_1z_1z_2m(z_1)m(z_2)\underline{m}(z_1)\underline{m}(z_2)
    =\;c\alpha_1\frac{\underline{m}(z_1)\underline{m}(z_2)}{\big(1+\frac{\sigma^2}{\mu^2}\underline{m}(z_1)\big)\big(1+\frac{\sigma^2}{\mu^2}\underline{m}(z_2)\big)},
\end{align*}
% \begin{align}
%     I_1\stackrel{i.p.}{\rightarrow} &\;c\left(\frac{\Expe |w_{11}-\mu|^{4}}{\mu^{4}}-3\frac{|\Expe |w_{11}-\mu|^{2}|^{2}}{\mu^{4}} \right)z_1z_2m(z_1)m(z_2)\underline{m}(z_1)\underline{m}(z_2)\nonumber\\
%     =&\;c\left(\frac{\Expe |w_{11}-\mu|^{4}}{\mu^{4}}-3\frac{|\Expe |w_{11}-\mu|^{2}|^{2}}{\mu^{4}} \right)\frac{\underline{m}(z_1)\underline{m}(z_2)}{\big(1+\frac{\sigma^2}{\mu^2}\underline{m}(z_1)\big)\big(1+\frac{\sigma^2}{\mu^2}\underline{m}(z_2)\big)},
% \end{align}
where the equality above follows from $m(z)=-z^{-1}\big(1+\frac{\sigma^2}{\mu^2}\underline{m}(z)\big)^{-1}$. Thus, the in probability (i.p.) limit of $\frac{\partial^2}{\partial z_2\partial z_1}I_1$ is in \eqref{parI1lim}.
% \begin{eqnarray*}
% \frac{\partial^2}{\partial z_1\partial z_2}I_1\stackrel{i.p.}{\rightarrow}c\alpha_1\frac{\underline{m}'(z)\underline{m}'(z_2)}{\big(1+\frac{\sigma^2}{\mu^2}\underline{m}(z_1)\big)^2\big(1+\frac{\sigma^2}{\mu^2}\underline{m}(z_2)\big)^2}.
% \end{eqnarray*}

\underline{Step (iii)}: \textbf{Consider $I_4$}.
We have
$\Expe\left|\frac{1}{p}\tr\Expe_j\bD_j^{-1}(z_1)\frac{1}{p}\tr\Expe_j\bD_j^{-1}(z_2)-m_{p}^{0}(z_1) m_{p}^{0}(z_2)\right|=o(1)$. By Lemma \ref{estsqu}, we get 
% \eqref{beta2appr}.
\begin{align}\label{nu12minusn2sq}   \lim_{p\rightarrow\infty}p(\nu_{12}-\nu_2^2) = h_2-2\frac{\sigma^2}{\mu^2}h_1.
\end{align}
By \eqref{b1zm} and \eqref{nu12minusn2sq}, we have
\begin{align*}
I_4\stackrel{i.p.}{\rightarrow}c\Bigl(h_2-2\frac{\sigma^2}{\mu^2}h_1\Bigr)z_1m(z_1)z_2m(z_2)\underline{m}(z_1)\underline{m}(z_2)=c\Bigl(h_2-2\frac{\sigma^2}{\mu^2}h_1\Bigr)\frac{\underline{m}(z_1)\underline{m}(z_2)}{\big(1+\frac{\sigma^2}{\mu^2}\underline{m}(z_1)\big)\big(1+\frac{\sigma^2}{\mu^2}\underline{m}(z_2)\big)}.
\end{align*}
% where the equality above uses the relation between $m(z)$ and $\underline{m}(z)$: $m(z)=-\frac{1}{z\big(1+\frac{\sigma^2}{\mu^2}\underline{m}(z)\big)}$.
Then the in probability (i.p.) limit of the second derivative $\frac{\partial^2}{\partial z_2\partial z_1}I_4$ is in \eqref{parI4lim}.
% \begin{eqnarray*}
% \frac{\partial^2}{\partial z_1\partial z_2}I_4\stackrel{i.p.}{\rightarrow}c\Bigl(h_2-2\frac{\sigma^2}{\mu^2}h_1\Bigr)\frac{\underline{m}'(z_1)\underline{m}'(z_2)}{\big(1+\frac{\sigma^2}{\mu^2}\underline{m}(z_1)\big)^2\big(1+\frac{\sigma^2}{\mu^2}\underline{m}(z_2)\big)^2}.
% \end{eqnarray*}

% \subsection{Tightness of $M_p^{(1)}(z)$}\label{prftigmp1}

\textbf{Step 3: Tightness of $M_p^{(1)}(z)$.} To prove tightness of $M_p^{(1)}(z)$, it is sufficient to prove the moment condition of \cite{billingsley1968convergence}, i.e., $\sup_{n; z_1, z_2\in\mathcal{C}_n}\frac{\Expe|M_p^{(1)}(z_1)-M_p^{(1)}(z_2)|^2}{|z_1-z_2|^2}$
% \begin{eqnarray}\label{momentbound}
% \sup_{n; z_1, z_2\in\mathcal{C}_n}\frac{\Expe|M_p^{(1)}(z_1)-M_p^{(1)}(z_2)|^2}{|z_1-z_2|^2}
% \end{eqnarray}
is finite. Its proof exactly follows \cite{bai2004clt}, and is postponed to Section \ref{prftight} of the supplementary document.% The tightness of $M_p^{(1)}(z)$ is similar to that provided in \cite{bai2004clt}.

% \subsection{Convergence of $M_p^{(2)}(z)$}\label{prfconmp2}

\textbf{Step 4: Convergence of $M_p^{(2)}(z)$.} Similar to \cite{bai2004clt}, one can prove the inequality:
% \begin{eqnarray*}
% \sup_{z\in\mathcal{C}_n}|\underline{m}_p^{0}(z)-\underline{m}(z)|\rightarrow 0, \ \ as \ \ n\rightarrow\infty,
% \end{eqnarray*}
% \[
% \sup_{n;z\in\mathcal{C}_n}\Bigl\|\big(\frac{\sigma^2}{\mu^2}\Expe\underline{m}_p(z)\bI_p+\bI_p\big)^{-1}\Bigr\|<\infty,
% \]
% \begin{eqnarray*}
% \sup_{z\in\mathcal{C}_n}\big|\frac{c_n\frac{\sigma^2}{\mu^2}\Expe\big(\underline{m}_p^2(z)\big)}{\big(1+\frac{\sigma^2}{\mu^2}\Expe\underline{m}_p(z)\big)^2}\big|<\xi, \ \ \xi\in(0,1),
% \end{eqnarray*}
\begin{eqnarray}\label{dmedm}
\Expe\big|\tr\bD^{-1}_1(z)\bM-\Expe\tr\bD^{-1}_1(z)\bM\big|^2\leq K\|\bM\|^2.
\end{eqnarray}

We first present the following equations for later use,
$ M_{p}^{(2)}(z) = p (\Expe m_{F^{p^2 S_{n}^0}}(z)-m_{F^{c_n}}(z))= n(\Expe\underline{m}_p(z)-\underline{m}_p^0(z)),\ \ \underline{m}(z)
    =  -\frac{1-c}{z} + cm(z)$. The next step is to find $\Expe\underline{m}_p(z)$. From the identity \eqref{zmbar}, which is the inverse of $\underline{m}(z)$, we define 
\begin{align*}
	R_p(z) &:= \frac{1}{\Expe \underline{m}_p(z)}+z-c_n\frac{\sigma^2/\mu^2}{1+\sigma^2/\mu^2\Expe\underline{m}_p(z)}=\frac{1}{\Expe\underline{m}_p(z)}\Bigl(1-c_n+z\Expe \underline{m}_p(z)+\frac{c_n}{1+\sigma^2/\mu^2\Expe\underline{m}_p(z)}\Bigr),\label{eq:An_def}
\end{align*}
thus,
\begin{equation}
	\Expe \underline{m}_p(z)
 % &= \Bigl(-z+c_n\frac{\sigma^2/\mu^2}{1+\sigma^2/\mu^2\Expe\underline{m}_p(z)}+R_p(z)\Bigr)^{-1}\nonumber\\
	= \Bigl(-z+c_n\frac{\sigma^2/\mu^2}{1+\sigma^2/\mu^2\Expe\underline{m}_p(z)}+A_p(z)\big/\Expe\underline{m}_p(z)\Bigr)^{-1},\label{eq:Expe_mnz}
\end{equation}
where $
    A_p(z) = \frac{c_n} {1+\sigma^2/\mu^2\Expe\underline{m}_p(z)} + zc_n\Expe M_p(z)$.
Note that
\begin{align}\label{eq:mn0z}
    \underline{m}_{p}^{0} = \Bigl(-z+c_n\frac{\sigma^2/\mu^2}{1+\sigma^2/\mu^2 \underline{m}_{p}^{0}}\Bigr)^{-1}.
\end{align}
From \eqref{eq:Expe_mnz} and \eqref{eq:mn0z}, we get
% \begin{equation*}
% 	\Expe \underline{m}_p(z) - \underline{m}_p^0(z) = \underline{m}_p^0(z)\Expe\underline{m}_p(z)\Biggl[ c_n\frac{(\sigma^2/\mu^2)^2(\Expe \underline{m}_p(z) - \underline{m}_p^0(z))}{(1+\sigma^2/\mu^2\Expe\underline{m}_p(z))(1+\sigma^2/\mu^2\underline{m}_p^0(z))} -\frac{A_p(z)}{\Expe \underline{m}_p(z)} \biggr].
% \end{equation*}
% Thus,
\begin{align}%\label{empminusmp}
\Expe \underline{m}_p(z) - \underline{m}_p^0(z)
 =-\underline{m}_p^0(z) A_p(z)\Biggl[ 1 - c_n\underline{m}_p^0(z)\Expe \underline{m}_p(z)\frac{(\sigma^2/\mu^2)^2}{\bigl(1+\frac{\sigma^2}{\mu^2}\Expe\underline{m}_p(z)\bigr)\bigl(1+\frac{\sigma^2}{\mu^2}\underline{m}_p^0(z)\bigr)} \biggr]^{-1}.\label{eq:Mn2_An}
\end{align}
Our next task is to investigate the limiting behavior of $nA_p$. Let $\widetilde{\bQ}_p(z)=\bI_{p}+\frac{\sigma^2}{\mu^2}\Expe\underline{m}_p(z)\bI_{p}$, then
\begin{align}\label{prenAp}
    nA_p &= p\frac{1} {1+\sigma^2/\mu^2\Expe\underline{m}_p(z)} + pz\Expe M_p(z)= \Expe(\beta_1P_1(z)) + \Expe(\beta_1P_2(z)),
\end{align}
where
\begin{align*}
P_1(z)&=\Bigl[n\br_1' \bD_1^{-1}(z)\widetilde{\bQ}_p^{-1}(z)\br_1 -\frac{\sigma^2}{\mu^2}\operatorname{tr}\left(\widetilde{\bQ}_p^{-1}(z)\Expe\bD_1^{-1}(z)\right)\Bigr],\\
P_2(z)&=\Bigl[\frac{\sigma^2}{\mu^2}\operatorname{tr}(\bQ_{p}^{-1}(z) \Expe\bD_1^{-1}(z))-\frac{\sigma^2}{\mu^2}\operatorname{tr}\left(\bQ_{p}^{-1}(z)\Expe\bD^{-1}(z)\right)\Bigr].
\end{align*}
Since $\beta_1 =b_{p}-b_{p}^2\gamma_1+\beta_1b_{p}^2\gamma_1^2$, where $\gamma_1(z)=\br_1'\bD_1^{-1}(z)\br_1-\frac{1}{n}\nu_2\Expe\operatorname{tr}\bD_1^{-1}(z)$, we have
\begin{align}\label{betap1}
 \Expe(\beta_1P_1(z)) 
% &=\Expe\left(\beta_1\Bigl[n\br_1' \bD_1^{-1}(z)\widetilde{\bQ}_p^{-1}(z)\br_1 -\frac{\sigma^2}{\mu^2}\tr \left(\widetilde{\bQ}_p^{-1}(z)\Expe\bD_1^{-1}(z)\right)\Bigr]\right)\nonumber\\ 
= b_p(z)\Expe P_1(z) - b_p^{2}(z)\Expe(\gamma_1P_1(z)) + b_p^{2}(z)\Expe(\beta_1\gamma_1^{2}P_1(z)).
\end{align}
For $\Expe P_1(z)$, by \eqref{covmatrixmain}, we get
\begin{align}\label{p1}
    \Expe P_1(z) 
    % &=  \Expe\Bigl[n\br_1' \bD_1^{-1}(z)\widetilde{\bQ}_p^{-1}(z)\br_1 -\frac{\sigma^2}{\mu^2}\operatorname{tr}\left(\widetilde{\bQ}_p^{-1}(z)\Expe\bD_1^{-1}(z)\right)\Bigr]\nonumber\\
    % &= n\Expe\operatorname{tr}\Bigl[\frac{1}{1+\frac{\sigma^2}{\mu^2}\Expe\underline{m}_p(z)}\br_1\br_1'\bD_1^{-1}(z)\Bigr]  - \frac{\sigma^2}{\mu^2}\Expe\operatorname{tr}\Bigl[\frac{1}{1+\frac{\sigma^2}{\mu^2}\Expe\underline{m}_p(z)}\ExpeD_1^{-1}(z)\Bigr]\nonumber\\
    % &= \frac{n}{1+\frac{\sigma^2}{\mu^2}\Expe\underline{m}_p(z)}\Expe\operatorname{tr}\Bigl[\br_1'\bD_1^{-1}(z)\br_1-\frac{1}{n}\frac{\sigma^2}{\mu^2}\Expe\bD_1^{-1}(z)\Bigr] \nonumber\\
    &= \frac{n}{1+\frac{\sigma^2}{\mu^2}\Expe\underline{m}_p(z)}\Big(\Expe\gamma_1(z)+\frac{1}{n}\Bigl(\nu_2-\frac{\sigma^2}{\mu^2}\Bigr)\Expe\operatorname{tr}\bD_1^{-1}(z)\Big).
\end{align}
The estimates for $\Expe(\gamma_1P_1(z))$, $\Expe(\beta_1\gamma_1^{2}P_1(z))$, and $\Expe(\beta_1P_2(z))$ are provided in the following lemma.
\begin{lemma}\label{estbeta1gamma1p1p2}
Under conditions and notations in Theorem \ref{clt}, we have
\begin{align}\label{gammp1}
     & \Expe(\gamma_1P_1(z)) 
    % &=  \Expe\Bigl[\gamma_1(z)\Bigl(n\br_1' \bD_1^{-1}(z)\widetilde{\bQ}_p^{-1}(z)\br_1 -\frac{\sigma^2}{\mu^2}\operatorname{tr}\left(\widetilde{\bQ}_p^{-1}(z)\Expe\bD_1^{-1}(z)\right)\Bigr)\Bigr]\nonumber\\
    % = n\Expe\Bigl[\gamma_1(z)\br_1'\bD_1^{-1}(z)\widetilde{\bQ}_p^{-1}(z)\br_1\Bigr]  - \frac{\sigma^2/\mu^2}{1+\frac{\sigma^2}{\mu^2}\Expe\underline{m}_p(z)}\operatorname{tr}[\Expe\bD_1^{-1}(z)]\Expe\gamma_1(z)\nonumber\\
    =n\Expe\Bigl[\Bigl(  \br_1'\bD_1^{-1}(z)\br_1-\frac{1}{n}\nu_2\operatorname{tr}\bD_1^{-1}(z)   \Bigr)\times\Bigl(\br_1'\bD_1^{-1}(z)\widetilde{\bQ}_p^{-1}(z)\br_1-\frac{1}{n}\nu_2\operatorname{tr}\left[\bD_1^{-1}(z)\widetilde{\bQ}_p^{-1}(z)\right]\Bigr)\Bigr] \nonumber\\
    &+\frac{1}{n(p-1)}\nu_2^2\Expe\left(\operatorname{tr}\bD_1^{-1}(z)\operatorname{tr}[\bD_1^{-1}(z)\widetilde{\bQ}_p^{-1}(z)]\right)-\frac{1}{n(p-1)}\nu_2^2\Expe\left(\boldsymbol{1}'_{p}\bD_1^{-1}(z)\boldsymbol{1}_{p}\operatorname{tr}[\bD_1^{-1}(z)\widetilde{\bQ}_p^{-1}(z)]\right)\nonumber\\
    &\quad- \frac{\sigma^2/\mu^2}{1+\frac{\sigma^2}{\mu^2}\Expe\underline{m}_p(z)}\operatorname{tr}[\Expe\bD_1^{-1}(z)]\Expe\gamma_1(z)+o(1),
\end{align}
and
\begin{align}
   &\Expe(\beta_1(z)\gamma_1^{2}(z)P_1(z))=\;\Expe\Bigl(n\beta_1(z)\gamma_1^{2}(z)\br_1' \bD_1^{-1}(z)\widetilde{\bQ}_p^{-1}(z)\br_1 \Bigr)-\Expe\left(\beta_1(z)\gamma_1^{2}(z)\operatorname{tr}\left[\frac{\sigma^2}{\mu^2}\widetilde{\bQ}_p^{-1}(z)\bD_1^{-1}(z)\right]\right)\nonumber\\
&\qquad\qquad\qquad\qquad\ \ \ \ \ \ +\Cov\left(\beta_1(z)\gamma_1^{2}(z),\operatorname{tr}\left[\frac{\sigma^2}{\mu^2}\widetilde{\bQ}_p^{-1}(z)\bD_1^{-1}(z)\right]\right)=\;O(\delta_n^2),\label{begap1}
\end{align}
and
\begin{align}\label{betap2}
    \Expe(\beta_1(z)P_2(z))
     = \frac{p}{n(p-1)}\nu_2\frac{\sigma^2}{\mu^2}b_p^2(z)\Expe\operatorname{tr}\Bigl[\bD_1^{-1}(z)\widetilde{\bQ}_p^{-1}(z)\bD_1^{-1}(z)\Bigr] + O(n^{-1/2}).
\end{align}
\end{lemma}
The proof of Lemma \ref{estbeta1gamma1p1p2} is postponed to Section \ref{prfeqs} of the supplementary document. Therefore, from \eqref{prenAp} -- \eqref{betap2}, we get
\begin{align}\label{nAp}
      nA_p = J_1 + J_2 + J_3 + o(1),
\end{align}
where
\begin{align}
    J_1&= \frac{nb_p(z)}{1+\frac{\sigma^2}{\mu^2}\Expe\underline{m}_p(z)}\Big(\Expe\gamma_1(z)+\frac{1}{n}(\nu_2-\frac{\sigma^2}{\mu^2})\Expe\operatorname{tr}\bD_1^{-1}(z)\Big)\nonumber\\
&\ \ \ \ \ \ \ \ \ \ \ \ \ \ +\left(b_p^2(z)\frac{\sigma^2}{\mu^2}\frac{1}{1+\frac{\sigma^2}{\mu^2}\Expe\underline{m}_p(z)}\right)\operatorname{tr}[\Expe\bD_1^{-1}(z)]\Expe\gamma_1(z)\nonumber\\
 &\quad-b_p^2(z)\biggl[\frac{1}{n(p-1)}\nu_2^2\Expe\left(\operatorname{tr}\bD_1^{-1}(z)\operatorname{tr}[\bD_1^{-1}(z)\widetilde{\bQ}_p^{-1}(z)]\right)\label{J13}\\
   &\quad\quad\quad\quad-\frac{1}{n(p-1)}\nu_2^2\Expe\left(\boldsymbol{1}'_{p}\bD_1^{-1}(z)\boldsymbol{1}_{p}\operatorname{tr}[\bD_1^{-1}(z)\widetilde{\bQ}_p^{-1}(z)]\right)\biggr],\label{J14}\\
   J_2&=-nb_p^2(z)\Expe\left(\br_1'\bD_1^{-1}(z)\br_1-\frac{1}{n}\nu_2\operatorname{tr} \bD_1^{-1}(z)\right)\nonumber\\
&\quad\qquad\times\left(\br_1'\bD_1^{-1}(z)\widetilde{\bQ}_p^{-1}(z)\br_1-\frac{1}{n}\nu_2\operatorname{tr} [\bD_1^{-1}(z)\widetilde{\bQ}_p^{-1}(z)]\right),\nonumber\\
J_3&= \frac{p}{n(p-1)}b_p^2(z)\frac{\sigma^2}{\mu^2}\nu_2 \Expe \tr[\bD_1^{-1}(z)\widetilde{\bQ}_p^{-1}(z)\bD_1^{-1}(z)].\nonumber
\end{align}
The limits of $J_1$, $J_2$ and $J_3$ are provided in the following lemma. The proof of Lemma \ref{J123limit} is postponed to Section \ref{prfJ123limit} of the supplementary document.
\begin{lemma}\label{J123limit}
    Under conditions and notation in Theorem \ref{clt}, as $n\rightarrow\infty$,
    \begin{align*}
         &J_1 \stackrel{i.p.}{\rightarrow} \biggl(\frac{-z\underline{m}(z)}{1+\frac{\sigma^2}{\mu^2}\underline{m}(z)}\biggr)\Bigl(\frac{\sigma^2}{\mu^2}m(z)+\frac{\sigma^2}{\mu^2}\frac{1}{z}+h_1m(z)\Bigr),\\
    &J_{2} \stackrel{i.p.}{\rightarrow} -\frac{c\alpha_1z^2m^2(z)\underline{m}^2(z)}{1+\frac{\sigma^2}{\mu^2}\underline{m}(z)}-\frac{2cz^2m'(z)\underline{m}^2(z)}{1+\frac{\sigma^2}{\mu^2}\underline{m}(z)}\frac{|\Expe |w_{1}-\mu|^{2}|^{2}}{\mu^{4}}-c(h_2-2\frac{\sigma^2}{\mu^2}h_1)\frac{z^2m^2(z)\underline{m}^2(z)}{1+\frac{\sigma^2}{\mu^2}\underline{m}(z)},\\
    &J_3 \stackrel{i.p.}{\rightarrow} c\frac{\sigma^4}{\mu^4}\underline{m}^2(z)(1+\frac{\sigma^2}{\mu^2}\underline{m}(z))^{-3}\Bigl[1-c\frac{\sigma^4}{\mu^4}\underline{m}^2(z)(1+\frac{\sigma^2}{\mu^2}\underline{m}(z))^{-2}\Bigr]^{-1}.
    \end{align*}
\end{lemma}
From \eqref{eq:Mn2_An}, \eqref{nAp}, and Lemma \ref{J123limit}, we get \eqref{mpzmean}. The proof is completed.

\begin{supplement}
\stitle{Supplement to ``On eigenvalues of sample covariance matrices based on high-dimensional compositional data''.}
\sdescription{
This supplementary document contains some technical lemmas and their proofs, including proofs of Theorem \ref{thm:CoDa-LSD},  Proposition \ref{extreigen}, Lemmas \ref{estsqu} -- \ref{quadform}, Lemmas \ref{ej-1term0} --  \ref{J123limit}, Lemma \ref{xdxz}, Corollary \ref{corpoly}, the tightness of $M_p^{(1)}(z)$. We also report the numerical simulation of CLT for $M_p(z)$ in Section \ref{simuforcltmpz}.}
% We also report some additional simulation results in this document.
\end{supplement}

%%%%%%%%%%%%%%%%%%%%%%%%%%%%%%%%%%%%%%%%%%%%%%%%%%%%%%%%%%%%%
%%                  The Bibliography                       %%
%%                                                         %%
%%  imsart-???.bst  will be used to                        %%
%%  create a .BBL file for submission.                     %%
%%                                                         %%
%%  Note that the displayed Bibliography will not          %%
%%  necessarily be rendered by Latex exactly as specified  %%
%%  in the online Instructions for Authors.                %%
%%                                                         %%
%%  MR numbers will be added by VTeX.                      %%
%%                                                         %%
%%  Use \cite{...} to cite references in text.             %%
%%                                                         %%
%%%%%%%%%%%%%%%%%%%%%%%%%%%%%%%%%%%%%%%%%%%%%%%%%%%%%%%%%%%%%

%% if your bibliography is in bibtex format, uncomment commands:
\bibliographystyle{imsart-nameyear} % Style BST file (imsart-number.bst or imsart-nameyear.bst)
\bibliography{reference}       % Bibliography file (usually '*.bib')

%% or include bibliography directly:
% \begin{thebibliography}{}
% \bibitem[\protect\citeauthoryear{???}{???}]{b1}
% \end{thebibliography}

\end{document}

\begin{frontmatter}
\title{Supplement to ``On eigenvalues of sample covariance matrices based on high-dimensional compositional data''}
\runtitle{High-dimensional SCM based on compositional data}
%\thankstext{T1}{A sample additional note to the title.}
\end{frontmatter}

%%%%%%%%%%%%%%%%%%%%%%%%%%%%%%%%%%%%%%%%%%%%%%
%%%% Main text entry area:

\tableofcontents

This supplementary document contains some technical lemmas and their proofs, including proofs of Theorem \ref{thm:CoDa-LSD},  Proposition \ref{extreigen}, Lemmas \ref{estsqu} -- \ref{quadform}, Lemmas \ref{ej-1term0} --  \ref{J123limit}, Lemma \ref{xdxz}, Corollary \ref{corpoly}, the tightness of $M_p^{(1)}(z)$. We also report some additional simulation results in this document.

\section{Some technical lemmas}
\begin{lemma}[Weyl's inequality, Corollary 7.3.5 in \cite{horn2012matrix}]\label{lem:Weyl_ineq}
	Let $\bA$ and $\bB$ be two $p\times n$ matrices and let $r=\min\{p,n\}$. Let $\sigma_1(\bA)\geq   \cdots\geq   \sigma_r(\bA)$ and $\sigma_1(\bB)\geq   \cdots\geq   \sigma_r(\bB)$ be the nonincreasingly ordered singular values of $\bA$ and $\bB$, respectively. Then
	\[
	\max_{1\leq   i\leq   r}|\sigma_{i}(\bA)-\sigma_i(\bB)| \leq   \|\bA-\bB\|,
	\]
	where $\|\bA-\bB\|$ denotes the spectral norm of $\bA-\bB$.
\end{lemma}

\begin{lemma}[\cite{Burkholder1973Distribution}]\label{burk}
Let $\left\{X_k\right\}$ be a complex martingale difference sequence with respect to the increasing $\sigma$-field $\{\mathcal F_k\}$. Then for $q > 1$,
\[\Expe\left|\sum X_k\right|^q \le K_q \left\{\Expe\left(\sum \Expe_{k-1}|X_k|^2\right)^{q/2}+\Expe\left(\sum|X_k|^q\right)\right\}.
\]
\end{lemma}

\begin{lemma}[Theorem 35.12 of \cite{Billingsley1995Probability}]\label{billin}
Suppose for each $n$, $\{Y_{n1}, Y_{n2}, \ldots, Y_{nr_n}\}$ is a real martingale difference sequence with respect to the increasing $\sigma$-field $\{\mathcal{F}_{nj}\}$ having second moments. If as $n\rightarrow\infty$,
\begin{eqnarray}
\sum^{r_n}_{j=1}\Expe(Y_{nj}^2|\mathcal{F}_{n,j-1})\stackrel{i.p.}{\rightarrow}\sigma^2,
\end{eqnarray}
where $\sigma^2$ is positive constant, and for each $\varepsilon>0$,
\begin{eqnarray}
\sum^{r_n}_{j=1}\Expe(Y_{nj}^2I_{\{|Y_{nj}|\geq\varepsilon\}})\rightarrow 0,
\end{eqnarray}
then $\sum^{r_n}_{j=1}Y_{nj}\stackrel{D}{\rightarrow}\mathcal{N}(0,\sigma^2)$.
\end{lemma}

Lemma \ref{xdxz} is used in the proof of Lemma \ref{ej-1term0}, and its proof is provided in Section \ref{sec:lem-xdxz-proof}.
\begin{lemma}\label{xdxz}
    Suppose that $\bx_{p}=\frac{1}{\sqrt{p}}(1,1,\ldots,1)'$ is a $p$-dimensional normalized all-one vector, then for the truncated random variable satisfying \eqref{truncation}, we have
$\Expe|\bx_{p}'\bD^{-1}(z)\bx_{p} +\frac{1}{z}|^2\rightarrow 0$.
\end{lemma}

\section{Proofs}
This section contains proofs of Theorem \ref{thm:CoDa-LSD},  Proposition \ref{extreigen}, Lemmas \ref{estsqu} -- \ref{quadform}, Lemmas \ref{ej-1term0} --  \ref{J123limit}, Lemma \ref{xdxz}, Corollary \ref{corpoly}, tightness of $M_p^{(1)}(z)$.

% \subsection{Proof of Theorem \ref{thm:CoDa-LSD}, Lemmas \ref{extreigen} and \ref{eigen}}\label{prflsdandext}
\subsection{Proof of Theorem \ref{thm:CoDa-LSD}}\label{prfLSD}
We write
\[
	\bX_{n}  = \begin{pmatrix}
		\frac{1}{\sum_{j=1}^p w_{1j}} & \cdots & 0 \\
		\vdots & \ddots & \vdots \\
		0 & \cdots & \frac{1}{\sum_{j=1}^p w_{nj}}
	\end{pmatrix}
	\begin{pmatrix}
		w_{11} & \cdots & w_{1p} \\
		\vdots & \ddots & \vdots \\
		w_{n1} & \cdots & w_{np} 
	\end{pmatrix} =: \bLambda_n \bW_n,
\]
and let $\widetilde{\bW}_n= \tfrac{\bC_n\bW_n}{\sqrt{N}\mu}$. By Theorem 3.6 in \cite{BSbook}, the LSD of $\widetilde{\bW}_n'\widetilde{\bW}_n$ has a density function given by \eqref{eq:MP_density}. 
Let $L(\cdot,\cdot)$ be the Levy distance. To prove Theorem~\ref{thm:CoDa-LSD}, we need to show that $L(F^{\bB_{p,N}}, F^{\widetilde{\bW}_n'\widetilde{\bW}_n})\convas 0$.
Let $\bY_n = \tfrac{p\bC_n\bX_n}{\sqrt{N}}=\tfrac{p\bC_n\bLambda_n\bW_n}{\sqrt{N}}$, then we can write $\bB_{p,N} = \bY_n'\bY_n$. By Corollary A.42 in \cite{BSbook}, we have
\begin{align*}
    L^4(F^{\bB_{p,N}}, F^{\widetilde{\bW}_n'\widetilde{\bW}_n}) & = L^4(F^{\bY_n'\bY_n}, F^{\widetilde{\bW}_n'\widetilde{\bW}_n})\\
    & \leq   \frac{2}{p^2} \tr \bigl[ (\bY_n- \widetilde{\bW}_n )' (\bY_n - \widetilde{\bW}_n) \bigr]\cdot \tr\bigl(\bY_n'\bY_n + \widetilde{\bW}_n' \widetilde{\bW}_n\bigr).
\end{align*}
By using Von Neumann’s trace inequality and the law of large numbers, we have 
\begin{align*}
    \frac{1}{p}\tr(\bY_n'\bY_n) & = \frac{p}{N}\tr(\bW_n'\bLambda_n'\bC_n\bLambda_n\bW_n)\\
    & \leq   \frac{p\|\bC_n\|}{N}\tr(\bW_n'\bLambda_n'\bLambda_n\bW_n) = \frac{1}{N} \sum_{i=1}^n \frac{\sum_{j=1}^pw_{ij}^2/p}{\bigl(\sum_{j=1}^p w_{ij}/p\bigr)^2} =O_P(1),\\
    \frac{1}{p} \tr( \widetilde{\bW}_n' \widetilde{\bW}_n) & = \frac{1}{Np\mu^2} \tr(\bW_n'\bC_n\bW_n) \leq   \frac{\|\bC_n\|}{Np\mu^2} \tr(\bW_n'\bW_n) = \frac{1}{Np\mu^2} \sum_{i=1}^n\sum_{j=1}^p w_{ij}^2 = O_P(1).
\end{align*}
Thus, to complete the proof, it is sufficient to prove that 
\[
\frac{1}{p}\tr \bigl[ (\bY_n- \widetilde{\bW}_n )' (\bY_n - \widetilde{\bW}_n) \bigr] = o_P(1).
\]
Let $\bv_i$ be the $i$-th column of $\bY_n-\widetilde{\bW}_n$, that is, $\bv_i =\bC_n\Bigl( \tfrac{\bw_i}{\sqrt{N}(\sum_{j=1}^p w_{ij}/p)} - \tfrac{\bw_i}{\sqrt{N}\mu}\Bigr)$. Then,
\begin{align*}
    &\; \frac{1}{p} \tr \bigl[ (\bY_n-\widetilde{\bW}_n)'(\bY_n-\widetilde{\bW}_n) \bigr]  \\
    =& \; \frac{1}{p}\sum_{i=1}^{n} \|\bv_i\|^2 \leq      \frac{\|\bC_n\|}{Np}\sum_{i=1}^n\sum_{j=1}^p \Bigl(\frac{1}{\mu}-\frac{1}{\sum_{j=1}^p w_{ij}/p}\Bigr)^2 w_{ij}^2 \\
    = & \; \frac{1}{\mu^2}\Bigl( \frac{1}{Np} \sum_{i=1}^n\sum_{j=1}^p w_{ij}^2 \Bigr) - \frac{2}{N\mu}\sum_{i=1}^n\frac{\sum_{j=1}^p w_{ij}^2/p}{\sum_{j=1}^p w_{ij}/p} + \frac{1}{N}\sum_{i=1}^n\frac{\sum_{j=1}^p w_{ij}^2/p}{(\sum_{j=1}^p w_{ij}/p)^2}\\
    =&\; o_P(1).		
\end{align*}
This completes the proof.

\subsection{Proof of Proposition \ref{extreigen}}\label{prfextreigen}

First, we prove the convergence of extreme eigenvalues.
Recall that $\widetilde{\bW}_n= \tfrac{\bC_n\bW_n}{\sqrt{N}\mu}$. Let $\lambda_{\max}(\widetilde{\bW}_n'\widetilde{\bW}_n)$ and $\lambda_{\min}(\widetilde{\bW}_n'\widetilde{\bW}_n)$ be the largest eigenvalue and the smallest non-zero of $\widetilde{\bW}_n'\widetilde{\bW}_n$. By the equation (2.7) in \cite{jiang2004limiting} and Theorem 1.4 in \cite{xiao2010almost}, we have
\begin{equation}\label{eq:extreme-eig-SCM}
    \lim_{n\to\infty} \lambda_{\max}(\widetilde{\bW}_n'\widetilde{\bW}_n) = \frac{\sigma^2}{\mu^2}(1+\sqrt{c})^2\qquad \text{and}\qquad 
    \lim_{n\to\infty} \lambda_{\min}(\widetilde{\bW}_n'\widetilde{\bW}_n) = \frac{\sigma^2}{\mu^2}(1-\sqrt{c})^2.
\end{equation}
By Lemma \ref{lem:Weyl_ineq}, we have
\begin{align}
    \max_{1\leq   i\leq   p}\Bigl| \sqrt{\lambda_{i}(\bB_{p,N})}-\sqrt{\lambda_{i}(\widetilde{\bW}_n'\widetilde{\bW}_n)} \Bigr| & \leq    \Bigl\| \bC_n(p\bLambda_n-\bI)\tfrac{\bW_n}{\mu\sqrt{N}} \Bigr\|  \leq    \| p\mu\bLambda_n-\bI \|\cdot \Bigl\|\tfrac{\bW_n}{\mu\sqrt{N}}\Bigr\|.\label{eq:largest-eigval-sqrt}
\end{align}
By Lemma 2 in \citet{Bai1993limit}, we have 
$
\max_{1\leq   i \leq   n} |\sum_{j=1}^p w_{ij}/(p\mu) - 1 | \convas 0,
$
which implies that 
$
\|p\mu\bLambda-\bI\|=\max_{1\leq   i \leq   n} |\frac{\mu}{\sum_{j=1}^p w_{ij}/p} - 1|\convas 0
$.
From the Theorem 2.9 in \cite{benaych2012singular}, we conclude that $\|\bW_n/(\mu\sqrt{N})\|<\infty$. Combining these facts with 
\eqref{eq:extreme-eig-SCM} and \eqref{eq:largest-eigval-sqrt}, we have
\begin{equation*}
    \sqrt{\lambda_{\max}(\bB_{p,N})}-\sqrt{\lambda_{\max}(\widetilde{\bW}_n'\widetilde{\bW}_n)}   \convas 0\quad \text{and}\quad \sqrt{\lambda_{\min}(\bB_{p,N})}-\sqrt{\lambda_{\min}(\widetilde{\bW}_n'\widetilde{\bW}_n)}   \convas 0 ,
\end{equation*}
which together with \eqref{eq:extreme-eig-SCM} implies that
\begin{align*}
     \sqrt{\lambda_{\max}(\bB_{p,N})}+\sqrt{\lambda_{\max}(\widetilde{\bW}_n'\widetilde{\bW}_n)}  \convas \frac{2\sigma}{\mu}|1+\sqrt{c}|,\\[0.5em]
      \sqrt{\lambda_{\min}(\bB_{p,N})}+\sqrt{\lambda_{\min}(\widetilde{\bW}_n'\widetilde{\bW}_n)}  \convas \frac{2\sigma}{\mu}|1-\sqrt{c}|.
\end{align*}
Therefore,
\[
    \lambda_{\max}(\bB_{p,N})-\lambda_{\max}(\widetilde{\bW}_n'\widetilde{\bW}_n)  \convas 0 \quad\text{and}\quad \lambda_{\min}(\bB_{p,N})-\lambda_{\min}(\widetilde{\bW}_n'\widetilde{\bW}_n)  \convas 0 .
\]
This, together with \eqref{eq:extreme-eig-SCM}, completes the proof of \eqref{eq:extreme-eig-CoDa}. 

Now, we show that, with probability one, no eigenvalues of $\bB_{p,N}$ outside the support of LSD. 
For any positive constant $\varepsilon$ small enough such that $\eta_1-\varepsilon>\frac{\sigma^2}{\mu^2}(1+\sqrt{c})^2$ and $\eta_2+\varepsilon<\indicator_{(0,1)}(c)\cdot\frac{\sigma^2}{\mu^2}(1-\sqrt{c})^2$, we have 
\begin{align*}
    &\;\Prob\bigl(\lambda_{\max}(\bB_{p,N})\geq   \eta_1\bigr) \\ 
    =&\;  \Prob\bigl(\lambda_{\max}(\bB_{p,N})\geq   \eta_1, \lambda_{\max}(\widetilde{\bW}_n'\widetilde{\bW}_n)\geq   \eta_1-\varepsilon \bigr)\\
    &\qquad + \Prob\bigl(\lambda_{\max}(\bB_{p,N})\geq   \eta_1, \lambda_{\max}(\widetilde{\bW}_n'\widetilde{\bW}_n)< \eta_1-\varepsilon \bigr)\\
    \leq  &\; \Prob\bigl( \lambda_{\max}(\widetilde{\bW}_n'\widetilde{\bW}_n)\geq   \eta_1-\varepsilon \bigr) +  \Prob\bigl(|\lambda_{\max}(\bB_{p,N}) - \lambda_{\max}(\widetilde{\bW}_n'\widetilde{\bW}_n)|\geq   \varepsilon \bigr)
\end{align*}
and
\begin{align*}
    &\;\Prob\bigl(\lambda_{\min}(\bB_{p,N})\leq   \eta_2\bigr) \\ 
    =&\;  \Prob\bigl(\lambda_{\min}(\bB_{p,N})\leq   \eta_1, \lambda_{\min}(\widetilde{\bW}_n'\widetilde{\bW}_n)\leq   \eta_2-\varepsilon \bigr)\\
    &\qquad + \Prob\bigl(\lambda_{\min}(\bB_{p,N})\leq   \eta_1, \lambda_{\min}(\widetilde{\bW}_n'\widetilde{\bW}_n) > \eta_1-\varepsilon \bigr)\\
    \leq  &\; \Prob\bigl( \lambda_{\min}(\widetilde{\bW}_n'\widetilde{\bW}_n)\leq   \eta_2-\varepsilon \bigr) +  \Prob\bigl(|\lambda_{\min}(\bB_{p,N}) - \lambda_{\min}(\widetilde{\bW}_n'\widetilde{\bW}_n)|\geq   \varepsilon \bigr).
\end{align*} 
To prove this theorem, it suffices to give the following three estimations:
\begin{align}
    \Prob\Bigl(\max_{1\leq   i \leq   p} |\lambda_{i}(\bB_{p,N}) - \lambda_{i}(\widetilde{\bW}_n'\widetilde{\bW}_n)|\geq   \varepsilon \Bigr) & = o(n^{-\ell}),\label{eq:eig_diff_small_prob} \\
    \Prob\bigl( \lambda_{\max}(\widetilde{\bW}_n'\widetilde{\bW}_n)\geq   \eta_1-\varepsilon \bigr) &= o(n^{-\ell}),\label{eq:largest_eig_small_prob} \\
    \Prob\bigl( \lambda_{\min}(\widetilde{\bW}_n'\widetilde{\bW}_n)\leq   \eta_2-\varepsilon \bigr) &= o(n^{-\ell}).\label{eq:smallest_eig_small_prob} 
\end{align}

First, we prove \eqref{eq:eig_diff_small_prob}. In view of the inequalities \eqref{eq:largest-eigval-sqrt} and \eqref{eq:largest_eig_small_prob} (will be proved below), it suffices to show that, for any $\ell>0$ and $\varepsilon>0$, 
\begin{equation}\label{eq:chenpan2012_max_small_prob_1}
\Prob\biggl(\max_{1\leq   i \leq   n} \biggl|\frac{\sum_{j=1}^p w_{ij}/p}{\mu} - 1\biggr|\geq   \varepsilon \biggr) = o(n^{-\ell}),
\end{equation}
which follows from the equation (9) in \citet{chen2012convergence}. (Although \citet{chen2012convergence}'s ultra-high dimensional setting ($p/n\to\infty$ as $p\wedge n\to\infty$) is different from our setting, their equation (9) still holds when $p$ and $n$ are of the same order.) 

Second, we prove \eqref{eq:largest_eig_small_prob} and \eqref{eq:smallest_eig_small_prob}. Let $\bW_n^{\circ}=\frac{1}{\mu\sqrt{N}}(\bW_n-\mu\boldsymbol{1}_n\boldsymbol{1}_n')$. From \citet{bai2004clt}, we have
\begin{equation}\label{eq:SCM_largest_eig_small_prob}
    \Prob\Bigl(\lambda_{\max}\bigl((\bW_n^{\circ})'\bW_n^{\circ}\bigr)\geq   \eta_1-\varepsilon/2\Bigr) = o(n^{-\ell})
\end{equation}
and
\begin{equation}\label{eq:SCM_smallest_eig_small_prob}
    \Prob\Bigl(\lambda_{\min}\bigl((\bW_n^{\circ})'\bW_n^{\circ}\bigr)\leq   \eta_2-\varepsilon/2\Bigr) = o(n^{-\ell}).
\end{equation}
By using Lemma \ref{lem:Weyl_ineq}, we have 
\begin{align}
    &\;\max_{1\leq   i\leq   p}\Bigl| \lambda_i(\widetilde{\bW}_n'\widetilde{\bW}_n)-\lambda_i\bigl((\bW_n^{\circ})'\bW_n^{\circ}\bigr) \Bigr|\nonumber\\
    \leq  &\; \| \widetilde{\bW}_n-\bW_n^{\circ} \|
    = \biggl\{\frac{1}{N}\sum_{j=1}^p \biggl(\frac{\sum_{i=1}^n w_{ij}/n}{\mu}-1\biggr)^2\biggr\}^{1/2}\nonumber\\
    \leq  &\;  \biggl\{\frac{p}{N}\cdot \max_{1\leq   j\leq   p} \biggl(\frac{\sum_{i=1}^n w_{ij}/n}{\mu}-1\biggr)^2\biggr\}^{1/2}.\label{eq:eig_diff_Wtilde_Wo}
\end{align}
Similarly to \eqref{eq:chenpan2012_max_small_prob_1}, for any $\ell>0$ and $\varepsilon>0$,
\begin{equation}\label{eq:chenpan2012_max_small_prob_2}
 \Prob\Biggl(\max_{1\leq   j\leq   p} \biggl|\frac{\sum_{i=1}^n w_{ij}/n}{\mu}-1\biggr|\geq   \varepsilon\Biggr) = o(n^{-\ell}).
\end{equation}
This, together with \eqref{eq:SCM_largest_eig_small_prob} $\sim$ \eqref{eq:eig_diff_Wtilde_Wo}, completes the proof of \eqref{eq:largest_eig_small_prob} and \eqref{eq:smallest_eig_small_prob}.

\subsection{Proof of Lemma \ref{substiprintr}}\label{pfsubstiprintr}
   The proof of Lemma \ref{substiprintr} is quite similar to Sections 5.3.1, 5.3.2, and 5.5 of \cite{zheng2015substitution}, it is then omitted. For readers' convenience, we present the outline of the proof for Lemma \ref{substiprintr}. In this situation, $\bB_p^0= \frac{1}{n}\bY_n'\bY_n = \sum_{i=1}^n \br_i\br_i'$, $\br_i=\frac{1}{\sqrt{n}}(\frac{w_{i1}}{\bar{w}_i}-1,\ldots,\frac{w_{ip}}{\bar{w}_i}-1)'=\frac{1}{\sqrt{n}}(y_{i1},\ldots,y_{ip})', i=1,2,\ldots,n$ and $\Expe \br_i\br_i^{'}=\frac{1}{n}\Sigma$, where $\Sigma$ is the population covariance matrix of $p$ mutiple of the compositional data (i.e., $p\bX_n$) and $\Sigma=\nu_{2} \left(-\frac{1}{p-1}\boldsymbol{1}_{p}\boldsymbol{1}'_{p} + \frac{1}{p-1}\bI_p + \bI_p \right)$. As for moments of $y_{ij}$, by \eqref{eq:Bc_prob}, for any $q>0$, we have
\begin{align*}
       \Expe \left(y_{ij}\right)^q&=\Expe \left[\left(\frac{w_{ij}}{\bar{w_i}}-1\right)^qI(B_{p}(\epsilon))\right]+\Expe \left[\left(\frac{w_{ij}}{\bar{w_i}}-1\right)^qI(B^c_{p}(\epsilon))\right]\leq K\Expe \left(w_{ij}-\mu\right)^q.
   \end{align*}
% where $B_{p}(\epsilon) =\{\omega: |\overline{w_i}-u|\leq \epsilon, \overline{w_i}=\sum_{j=1}^{p}w_{ij}/p \}$. 
Therefore in the following proof, the requirement of truncation of $y_{ij}$ reduces to truncation of $w_{ij}$. First, we get that  \begin{align}\label{substiptineq1}
      &\operatorname{tr}(\bB_{p,N}-z\bI_p)^{-1}-pm_N^0(z) = \operatorname{tr}(\bB_{p}^0-z\bI_p)^{-1}-pm_n^0(z)+ p(m_n^0(z)-m_N^0(z))\nonumber\\
     & \ \ \ \ \ \ \ \ +\operatorname{tr}\bA^{-2}(z)\boldsymbol{\Delta}+\operatorname{tr}\bA^{-1}(z)(\boldsymbol{\Delta}\bA^{-1}(z))^2+\operatorname{tr}(\bA(z)-\boldsymbol{\Delta})^{-1}(\boldsymbol{\Delta}\bA^{-1}(z))^3,
    \end{align}
where $\bA(z)=\bB_{p}^0-z\bI_p$ and $\boldsymbol{\Delta}=\bB_{p}^0-\bB_{p,N}$.  Moreover, after truncation and normalization, for every $z\in\mathbb{C}^+=\{z:\Im z>0\}$,
    \begin{align}
         &p(m_n^0(z)-m_N^0(z))\stackrel{a.s.}{\rightarrow}(1+z\underline{m}(z))\frac{\underline{m}(z)+z\underline{m}^{'}(z)}{z\underline{m}(z)},\ \ \Expe |\operatorname{tr}(\bA^{-2}(z)\boldsymbol{\Delta})|^2=o(1),\label{substiptineq2}\\
        % &\Expe |\operatorname{tr}(\bA^{-2}(z)\boldsymbol{\Delta})|^2=o(1),\label{subsprin_small}\\
         & \operatorname{tr}(\bA^{-2}(z)\boldsymbol{\Delta}\bA^{-1}(z)\boldsymbol{\Delta})\stackrel{L_2}{\rightarrow}(\underline{m}(z)+z\underline{m}^{'}(z))(1+z\underline{m}(z)),\label{substiprineq4}\\
        & \operatorname{tr}(\bA^{-1}(z)\boldsymbol{\Delta})^3(\bA(z)-\boldsymbol{\Delta})^{-1} =(1+z\underline{m}(z)) \operatorname{tr}((\bA^{-1}(z)\boldsymbol{\Delta})^2(\bA(z)-\boldsymbol{\Delta})^{-1})+o_P(1).\label{substiprineq5}
      \end{align}
% Hence by  \eqref{subsprin_small}-\eqref{substiprineq5}
%     \begin{align}\label{substiptineq6}
%         &\ \ \ \ \ \ \operatorname{tr}\bA^{-2}(z)\boldsymbol{\Delta}+\operatorname{tr}\bA^{-1}(z)(\boldsymbol{\Delta}\bA^{-1}(z))^2+\operatorname{tr}(\bA(z)-\Delta)^{-1}(\Delta\bA^{-1}(z))^3 \nonumber\\
%         &=\frac{(\underline{m}(z)+z\underline{m}^{'}(z))(1+z\underline{m}(z))}{-z\underline{m}(z)}+o_P(1).
%     \end{align}
Thus, by \eqref{substiptineq1} -- \eqref{substiprineq5}, Lemma \ref{substiprintr} is obtained.

Note that, we also need to check the tightness of $\operatorname{tr}(\bB_{p,N}-z\bI_p)^{-1}-pm_N^0(z)$. Since
\begin{align*}
     \operatorname{tr}(\bB_{p,N}-z\bI_p)^{-1}-pm_N^0(z)& = \operatorname{tr}(\bB_{p,N} -z\bI_p)^{-1}-\operatorname{tr}(\bB_{p}^0-z\bI_p)^{-1}\\
     &\ \ +\operatorname{tr}(\bB_{p}^0-z\bI_p)^{-1}-pm_n^0(z)+pm_n^0(z)-pm_N^0(z),
\end{align*}
and the tightness of $\operatorname{tr}(\bB_{p}^0-z\bI_p)^{-1}-pm_n^0(z)$ is proved in Step 2 of Section \ref{cltcenBp0}, it suffices to prove tightness of $\operatorname{tr}(\bB_{p,N} -z\bI_p)^{-1}-\operatorname{tr}(\bB_{p}^0-z\bI_p)^{-1}$. It can be obtained from similar arguments in Section 5.3.2 of \cite{zheng2015substitution} and we omit the details. Finally, the proof is completed.

\subsection{Proof of Lemma \ref{estsqu}}\label{prfestsqu}

Note that, by Taylor expansion, there exist $C_{1}>0$ such that, for any $-1/2\leq x\leq1/2$,
\begin{equation*}
    \frac{1}{(1+x)^{2}} = 1-2x+3x^{2}+a(x), |a(x)|\leq C_{1}x^{3}.
\end{equation*}
Hence, there exist $C_{1}>0$ such that, for any $0<\epsilon<1/2$, on the event $B_{p}(\epsilon)=\{\omega: |\overline{w}-\mu|\leq \epsilon, \overline{w}=\sum_{j=1}^{p}w_{j}/p \}$,
\[
    \frac{1}{\overline{w}^{2}} = \frac{1}{\mu^{2}(\frac{\overline{w}-\mu}{\mu}+1)^{2}}
    = \frac{1}{\mu^{2}}\Big[1-\frac{2(\overline{w}-\mu)}{\mu}+\frac{3(\overline{w}-\mu)^{2}}{\mu^{2}}+a\Bigl(\frac{\overline{w}-\mu}{\mu}\Bigr)\Big],
\]
where $|a((\overline{w}-\mu)/\mu)|<C_{1}\epsilon^{3}$.
Hence, we have
\begin{align}
    \frac{w_{1}^{2}}{\overline{w}^{2}}I(B_{p}(\epsilon)) 
    % &=  \frac{1}{\mu^{2}}\Big[w_{1}^{2}-\frac{2w_{1}^{2}(\overline{w}-\mu)}{\mu}+\frac{3w_{1}^{2}(\overline{w}-\mu)^{2}}{\mu^{2}}+w_{1}^{2}a(\frac{\overline{w}-\mu}{\mu})\Big]I(B_{p}(\epsilon))\nonumber\\
    &=  \Big[\frac{w_{1}^{2}}{\mu^{2}}-\frac{2w_{1}^{2}(\overline{w}-\mu)}{\mu^{3}}+\frac{3w_{1}^{2}(\overline{w}-\mu)^{2}}{\mu^{4}}\Big]I(B_{p}(\epsilon))+\frac{w_{1}^{2}}{\mu^{2}}a\Bigl(\frac{\overline{w}-\mu}{\mu}\Bigr)I(B_{p}(\epsilon)),\label{Bintaylor}
\end{align}
where $\bigl|\frac{w_{1}^{2}}{\mu^{2}}a(\frac{\overline{w}-\mu}{\mu})I(B_{p}(\epsilon))\bigr|<C_{1}\frac{w_{1}^{2}}{\mu^{2}}\epsilon^{3}$.
Since $I(B_{p}(\epsilon))=1-I(B_{p}^{c}(\epsilon))$, 
\begin{align}
    \text{the right hand side of \eqref{Bintaylor}} =&  \frac{w_{1}^{2}}{\mu^{2}}-\frac{2w_{1}^{2}(\overline{w}-\mu)}{\mu^{3}}+\frac{3w_{1}^{2}(\overline{w}-\mu)^{2}}{\mu^{4}}\nonumber\\
    &- \Big[\frac{w_{1}^{2}}{\mu^{2}}-\frac{2w_{1}^{2}(\overline{w}-\mu)}{\mu^{3}}+\frac{3w_{1}^{2}(\overline{w}-\mu)^{2}}{\mu^{4}}\Big]I(B_{p}^{c}(\epsilon))+ a, \nonumber
\end{align}
where $|a|\leq C_{1}\frac{w_{1}^{2}}{\mu^{2}}\epsilon^{3}$. Therefore,
\begin{equation}\label{Binexpan}
\begin{aligned}
     \frac{w_{1}^{2}}{\overline{w}^{2}}I(B_{p}(\epsilon)) - \frac{w_{1}^{2}}{\mu^{2}} =&  -\frac{2w_{1}^{2}(\overline{w}-\mu)}{\mu^{3}}+\frac{3w_{1}^{2}(\overline{w}-\mu)^{2}}{\mu^{4}}\\
    &- \Big[\frac{w_{1}^{2}}{\mu^{2}}-\frac{2w_{1}^{2}(\overline{w}-\mu)}{\mu^{3}}+\frac{3w_{1}^{2}(\overline{w}-\mu)^{2}}{\mu^{4}}\Big]I(B_{p}^{c}(\epsilon))+ a,
\end{aligned}
\end{equation}
where $|a|\leq C_{1}\frac{w_{1}^{2}}{\mu^{2}}\epsilon^{3}$.
Taking expectation for \eqref{Binexpan} yields that
\begin{align}
     \Expe\Big(\frac{w_{1}^{2}}{\overline{w}^{2}}I(B_{p}(\epsilon))\Big) - \Expe\frac{w_{1}^{2}}{\mu^{2}} =&  -\frac{2\Expe w_{1}^{2}(\overline{w}-\mu)}{\mu^{3}}+\frac{3\Expe w_{1}^{2}(\overline{w}-\mu)^{2}}{\mu^{4}}\nonumber\\
    &- \Expe\Big[\frac{w_{1}^{2}}{\mu^{2}}-\frac{2w_{1}^{2}(\overline{w}-\mu)}{\mu^{3}}+\frac{3w_{1}^{2}(\overline{w}-\mu)^{2}}{\mu^{4}}\Big]I(B_{p}^{c}(\epsilon))+ \Expe a\nonumber\\
     =&  -\frac{2\Expe w_{1}^{2}(\overline{w}-\mu)}{\mu^{3}}+\frac{3\Expe w_{1}^{2}(\overline{w}-\mu)^{2}}{\mu^{4}}
     -\Expe b+\Expe a,\nonumber
\end{align}
where $|\Expe a|\leq C_{1}\frac{w_{1}^{2}}{\mu^{2}}\epsilon^{3}$.
Thus, 
\begin{align}\label{estrate}
     \Expe\frac{w_{1}^{2}}{\overline{w}^{2}} - \Expe\frac{w_{1}^{2}}{\mu^{2}}
     =  -\frac{2\Expe w_{1}^{2}(\overline{w}-\mu)}{\mu^{3}}+\frac{3\Expe w_{1}^{2}(\overline{w}-\mu)^{2}}{\mu^{4}} 
     -\Expe b+\Expe a+\Expe\Big(\frac{w_{1}^{2}}{\overline{w}^{2}}I(B^{c}_{p}(\epsilon))\Big).
\end{align}
Note that 
\begin{align}\label{cbound}
   |\Expe c|= \Bigl|\Expe\Big(\frac{w_{1}^{2}}{\overline{w}^{2}}I(B^{c}_{p}(\epsilon))\Big)\Bigr| \leq p^{2}\Prob(B^{c}_{p}(\epsilon)).
\end{align}
Next, we bound $\Expe b$. In save of notation, we denote by
\[
    y_{1}= \frac{w_{1}^{2}}{\mu^{2}},
    \qquad
    y_{2}= \frac{2w_{1}^{2}(\overline{w}-\mu)}{\mu^{3}},
    \qquad
    y_{3}= \frac{3w_{1}^{2}(\overline{w}-\mu)^{2}}{\mu^{4}}.
\]
It is obvious that
\begin{align}
    \Bigl|\Expe(y_{1}I(B_{p}^{c}(\epsilon)))\Bigr| &\leq C_2\Bigl|\Expe w_{1}^{2}I(B_{p}^{c}(\epsilon))\Bigr|\leq C_2\Prob(B^{c}_{p}(\epsilon))^{1/2}.\label{b1est}
\end{align}
Note that 
\begin{align}
    \Expe\frac{w_{1}^{2}(\overline{w}-\mu)}{\mu^{3}}
    &= \frac{1}{p}\sum_{i=1}^{p}\frac{\Expe w_{1}^{2}(w_{i}-\mu)}{\mu^{3}}= \frac{1}{p}\frac{\Expe w_{1}^{2}(w_{1}-\mu)}{\mu^{3}}\label{y2},
\end{align}
and
\begin{align}
   \Expe\frac{w_{1}^{2}(\overline{w}-\mu)^{2}}{\mu^{4}}&= \frac{1}{p^{2}}\sum_{i,j=1}^{p}\frac{\Expe w_{1}^{2}(w_{i}-\mu)(w_{j}-\mu)}{\mu^{4}}\nonumber\\
   &= \frac{1}{p^{2}}\sum_{i=1}^{p}\frac{\Expe w_{1}^{2}(w_{i}-\mu)^{2}}{\mu^{4}}\nonumber\\
   &= \frac{1}{p^{2}}\frac{\Expe w_{1}^{2}(w_{1}-\mu)^{2}}{\mu^{4}}+\frac{p-1}{p^{2}}\frac{\Expe w_{1}^{2}\Expe(w_{1}-\mu)^{2}}{\mu^{4}}.\label{y3}
\end{align}
By \eqref{y2}, we get
\begin{align}
    \Bigl|\Expe(y_{2}I(B_{p}^{c}(\epsilon)))\Bigr| &=\frac{2}{p\mu^{3}}\Bigl|\Expe w_{1}^{2}(w_{1}-\mu)I(B_{p}^{c}(\epsilon))\Bigr|\nonumber\\
    &\leq \frac{C_3}{p}\Big(\Bigl|\Expe w_{1}^{2}I(B_{p}^{c}(\epsilon))\Bigr|+\Bigl|\Expe w_{1}^{3}I(B_{p}^{c}(\epsilon))\Bigr|\Big)\nonumber\\
    &\leq \frac{C_3}{p}\Big(\Prob(B^{c}_{p}(\epsilon))^{1/2}+\Bigl|(\Expe w_{1}^{4})^{1/2}(\Expe w_{1}^{2}I(B_{p}^{c}(\epsilon)))^{1/2}\Bigr|\Big)\nonumber\\
    &\leq \frac{C_3}{p}\Big(\Prob(B^{c}_{p}(\epsilon))^{1/2}+\Prob(B^{c}_{p}(\epsilon))^{1/4}\Big).\label{b2est}
\end{align}
By \eqref{y3}, we get
\begin{align}
    \Bigl|\Expe(y_{3}I(B_{p}^{c}(\epsilon)))\Bigr| &\leq  \frac{1}{p^{2}\mu^{4}}\Bigl|\Expe w_{1}^{2}(w_{1}-\mu)^{2}I(B_{p}^{c}(\epsilon))\Bigr|+\frac{p-1}{p^{2}\mu^{4}}\Bigl|\bigl(\Expe w_{1}^{2}I(B_{p}^{c}(\epsilon)\bigr)\bigl(\Expe(w_{1}-\mu)^{2}I(B_{p}^{c}(\epsilon)\bigr)\Bigr|\nonumber\\
    &\leq \frac{C_4}{p^{2}}+\frac{C_4}{p}\bigl(\Prob(B^{c}_{p}(\epsilon))^{1/2}\cdot\Prob(B^{c}_{p}(\epsilon))^{1/2}\bigr)\nonumber\\
    &\leq C_4\Big(\frac{1}{p^{2}}+\frac{1}{p}\Prob(B^{c}_{p}(\epsilon))\Big).\label{b3est}
\end{align}
By \eqref{b2est}--\eqref{b1est}, we have
\begin{align}\label{bbound}
    |\Expe b|&= \Bigl|\Expe(y_1+y_2+y_3)I(B_{p}^{c}(\epsilon))\Bigr|\nonumber\\
    &\leq  \Bigl|\Expe(y_{1}I(B_{p}^{c}(\epsilon)))\Bigr|+\Bigl|\Expe(y_{2}I(B_{p}^{c}(\epsilon)))\Bigr| + \Bigl|\Expe(y_{3}I(B_{p}^{c}(\epsilon)))\Bigr|\nonumber\\
    &\leq C_5\Bigl(\Prob(B^{c}_{p}(\epsilon))^{1/2}+\frac{1}{p}\Prob(B^{c}_{p}(\epsilon))^{1/4}+p^{-2}\Bigr).
\end{align}
By using \eqref{eq:Bc_prob} and $\Expe\left|w_{1}-\mu\right|^{4}<\infty$ and $|\frac{w_{j}-\mu}{\sigma}|\leq \sqrt{n}\delta_{n}$, we have $\Expe a - \Expe b + \Expe c=o(p^{-1})$.

Thus, by \eqref{estrate}, \eqref{cbound} and \eqref{bbound}, we get
\begin{equation}\label{lem2expan1}
    \Expe\frac{w_{1}^{2}}{\overline{w}^{2}} - \Expe\frac{w_{1}^{2}}{\mu^{2}}
     =  -\frac{2\Expe w_{1}^{2}(\overline{w}-\mu)}{\mu^{3}}+\frac{3\Expe w_{1}^{2}(\overline{w}-\mu)^{2}}{\mu^{4}}+o(p^{-1}).
\end{equation}
Plugging \eqref{y2}--\eqref{y3} into \eqref{lem2expan1}, we get
\begin{align}
\Expe(\frac{w_{1}}{\overline{w}})^{2}-\Expe(\frac{w_{1}}{\mu})^{2} &=\frac{1}{p}\Bigr[\frac{3\Expe w_{1}^{2}\Expe(w_{1}-\mu)^{2}}{\mu^{4}}-\frac{2\Expe w_{1}^{2}(w_{1}-\mu)}{\mu^{3}}\Bigr]+o(p^{-1})\nonumber\\
& = \frac{1}{p}h_1 + o(p^{-1}),
\end{align}
which is \eqref{eq:moment1}.

Similar to the previous calculation, we obtain
% \begin{align}
%  &\Bigl(\frac{w_1}{\bar{w}}-1\Bigr)^2\Bigl(\frac{w_2}{\bar{w}}-1\Bigr)^2I(B_{p}(\epsilon))=(\frac{w_1}{\mu}-1)^2(\frac{w_2}{\mu}-1)^2I(B_{p}(\epsilon))\nonumber\\
% &\ \ \ \ \ \ +(\frac{\bar{w}}{\mu}-1)\Big[-2\frac{w_1}{\mu}(\frac{w_1}{\mu}-1)(\frac{w_2}{\mu}-1)^2-2\frac{w_2}{\mu}(\frac{w_2}{\mu}-1)(\frac{w_1}{\mu}-1)^2\Big]I(B_{p}(\epsilon))\nonumber\\
% &\ \ \ \ \ \ +(\frac{\bar{w}}{\mu}-1)^2\Big[4\frac{w_1}{\mu}\frac{w_2}{\mu}(\frac{w_1}{\mu}-1)(\frac{w_2}{\mu}-1)\nonumber\\
%     &\ \ \ \ \ \ \ \ \ \ +(\frac{w_1}{\mu}-1)^2(\frac{w_2}{\mu}-1)^2(\frac{w_1^2/\mu^2}{(w_1/\mu-1)^2}+2\frac{w_1/\mu}{w_1/\mu-1})\nonumber\\
%     &\ \ \ \ \ \ \ \ \ \ +(\frac{w_1}{\mu}-1)^2(\frac{w_2}{\mu}-1)^2(\frac{w_2^2/\mu^2}{(w_2/\mu-1)^2}+2\frac{w_2/\mu}{w_2/\mu-1})\Big]I(B_{p}(\epsilon))+o(p^{-1}),
% \end{align}
\begin{align*}
 &\;\Bigl(\frac{w_1}{\bar{w}}-1\Bigr)^2\Bigl(\frac{w_2}{\bar{w}}-1\Bigr)^2I(B_{p}(\epsilon))\\
 =&\;\Bigl(\frac{w_1}{\mu}-1\Bigr)^2\Bigl(\frac{w_2}{\mu}-1\Bigr)^2I(B_{p}(\epsilon))+\Bigl(\frac{\bar{w}}{\mu}-1\Bigr)f_1(w_1,w_2,\mu)I(B_{p}(\epsilon))\\
 &\qquad +\Bigl(\frac{\bar{w}}{\mu}-1\Bigr)^2f_2(w_1,w_2,\mu)I(B_{p}(\epsilon))+o(p^{-1}),
\end{align*}
where
\begin{align}
    f_1(w_1,w_2,\mu)&=-2\frac{w_1}{\mu}\Bigl(\frac{w_1}{\mu}-1\Bigr)\Bigl(\frac{w_2}{\mu}-1\Bigr)^2-2\frac{w_2}{\mu}\Bigl(\frac{w_2}{\mu}-1\Bigr)\Bigl(\frac{w_1}{\mu}-1\Bigr)^2,\nonumber\\
    f_2(w_1,w_2,\mu)&=4\frac{w_1}{\mu}\frac{w_2}{\mu}\Bigl(\frac{w_1}{\mu}-1\Bigr)\Bigl(\frac{w_2}{\mu}-1\Bigr)\nonumber\\
    &\quad +\Bigl(\frac{w_1}{\mu}-1\Bigr)^2\Bigl(\frac{w_2}{\mu}-1\Bigr)^2\Bigl(\frac{w_1^2/\mu^2}{(w_1/\mu-1)^2}+2\frac{w_1/\mu}{w_1/\mu-1}\Bigr)\nonumber\\
    &\quad +\Bigl(\frac{w_1}{\mu}-1\Bigr)^2\Bigl(\frac{w_2}{\mu}-1\Bigr)^2\Bigl(\frac{w_2^2/\mu^2}{(w_2/\mu-1)^2}+2\frac{w_2/\mu}{w_2/\mu-1}\Bigr).
\end{align}
Similar to \eqref{Bintaylor}--\eqref{lem2expan1}, we get
% \begin{align}\label{v12exp1}
%  &\Expe(\frac{w_1}{\bar{w}}-1)^2(\frac{w_2}{\bar{w}}-1)^2=\Expe(\frac{w_1}{\mu}-1)^2(\frac{w_2}{\mu}-1)^2\nonumber\\
% &\ \ \ \ \ \ +\Expe(\frac{\bar{w}}{\mu}-1)\Big[-2\frac{w_1}{\mu}(\frac{w_1}{\mu}-1)(\frac{w_2}{\mu}-1)^2-2\frac{w_2}{\mu}(\frac{w_2}{\mu}-1)(\frac{w_1}{\mu}-1)^2\Big]\nonumber\\
% &\ \ \ \ \ \ +\Expe(\frac{\bar{w}}{\mu}-1)^2\Big[4\frac{w_1}{\mu}\frac{w_2}{\mu}(\frac{w_1}{\mu}-1)(\frac{w_2}{\mu}-1)\nonumber\\
%     &\ \ \ \ \ \ \ \ \ \ +(\frac{w_1}{\mu}-1)^2(\frac{w_2}{\mu}-1)^2(\frac{w_1^2/\mu^2}{(w_1/\mu-1)^2}+2\frac{w_1/\mu}{w_1/\mu-1})\nonumber\\
%     &\ \ \ \ \ \ \ \ \ \ +(\frac{w_1}{\mu}-1)^2(\frac{w_2}{\mu}-1)^2(\frac{w_2^2/\mu^2}{(w_2/\mu-1)^2}+2\frac{w_2/\mu}{w_2/\mu-1})\Big]+o(p^{-1})\nonumber\\
%     &\ \ \ \ \ \ =T_1+T_2+T_3.
% \end{align}
\begin{align}\label{v12exp1}
 &\;\Expe(\frac{w_1}{\bar{w}}-1)^2(\frac{w_2}{\bar{w}}-1)^2\nonumber\\
 =&\;\Expe\Bigl(\frac{w_1}{\mu}-1\Bigr)^2\Bigl(\frac{w_2}{\mu}-1\Bigr)^2 +\Expe\Bigl(\frac{\bar{w}}{\mu}-1\Bigr)f_1(w_1,w_2,\mu) +\Expe\Bigl(\frac{\bar{w}}{\mu}-1\Bigr)^2f_2(w_1,w_2,\mu)+o(p^{-1})\nonumber\\
    =:&\;T_1+T_2+T_3+o(p^{-1}).
\end{align}
Similar to \eqref{y2} and \eqref{y3}, we obtain
\begin{align}\label{v12t2}
    T_2&=\frac{1}{p}\Expe\sum_{i=1}^{p}\Bigl(\frac{w_i}{\mu}-1\Bigr)f_1(w_1,w_2,\mu)=\frac{1}{p}\Expe\Bigl(\frac{w_1}{\mu}+\frac{w_2}{\mu}-2\Bigr)f_1(w_1,w_2,\mu)
\end{align}
and
\begin{align}
    T_3&=\frac{1}{p^2}\Expe\sum_{i,j=1}^{p}\Bigl(\frac{w_i}{\mu}-1\Bigr)\Bigl(\frac{w_j}{\mu}-1\Bigr)f_2(w_1,w_2,\mu)\nonumber\\
    &=\frac{1}{p^2}\Expe\sum_{i=1}^{p}\Bigl(\frac{w_i}{\mu}-1\Bigr)^2f_2(w_1,w_2,\mu)\nonumber\\
     &=\frac{1}{p^2}\Expe\sum_{i=1}^{2}\Bigl(\frac{w_i}{\mu}-1\Bigr)^2f_2(w_1,w_2,\mu)+\frac{p-2}{p^2}\Expe\Bigl(\frac{w_1}{\mu}-1\Bigr)^2\Expe f_2(w_1,w_2,\mu).\label{v12t3}
\end{align}
Thus, by \eqref{v12exp1} -- \eqref{v12t3}, we get
\begin{align}
     \Expe\Bigl(\frac{w_1}{\bar{w}}-1\Bigr)^2\Bigl(\frac{w_2}{\bar{w}}-1\Bigr)^2
     &=\Expe\Bigl(\frac{w_1}{\mu}-1\Bigr)^2\Bigl(\frac{w_2}{\mu}-1\Bigr)^2+\frac{1}{p}\Expe\Bigl(\frac{w_1}{\mu}+\frac{w_2}{\mu}-2\Bigr)f_1(w_1,w_2,\mu)\nonumber\\
     &\qquad +\frac{1}{p}\frac{\sigma^2}{\mu^2}\Expe f_2(w_1,w_2,\mu)+o(p^{-1})\nonumber\\
     & = \Expe\Bigl(\frac{w_1}{\mu}-1\Bigr)^2\Bigl(\frac{w_2}{\mu}-1\Bigr)^2 + \frac{1}{p}h_2 + o(p^{-1}).
\end{align}
which is \eqref{eq:moment2}.

Similarly, we get
\[
    \Expe\Bigl(\frac{w_1}{\bar{w}}-1\Bigr)^4 = \Expe\Bigl(\frac{w_1}{\mu}-1\Bigr)^4 + o(1),
\]
which is \eqref{eq:moment3}.

\subsection{Proof of Lemma \ref{bound}}\label{prfbound}

First, we prove \eqref{eq:Bc_prob}. By Markov's inequality and Burkh\"{o}lder inequality, we get
 \begin{align}
    \Prob(B^{c}_{p}(\epsilon)) &= \Prob(| \overline{w}-\mu | \geq \epsilon)\nonumber\\
    % &\leq \frac{\Expe|\overline{w}-\mu|^{kq_{1}}}{\epsilon^{kq_{1}}}\nonumber\\
    % &= \frac{\Expe|\sum_{j=1}^{p}w_{j}/p-\mu|^{kq_{1}}}{\epsilon^{kq_{1}}}\nonumber\\
    &\leq \epsilon^{-kq_{1}}\Expe|\sum_{j=1}^{p}(w_{j}-\mu)/p|^{kq_{1}}\nonumber\\
    % &= \epsilon^{-kq_{1}}p^{-kq_{1}}\Expe\Bigl|\sum_{j=1}^{p}(w_{j}-\mu)\Bigr|^{kq_{1}}\nonumber\\
     &\leq K_{kq_{1}}\epsilon^{-kq_{1}}p^{-kq_{1}}\Big[\Expe(\sum_{j=1}^{p}\Expe_{j-1}|w_{j}-\mu|^{2})^{kq_{1}/2}+\Expe\sum_{j=1}^{p}|w_{j}-\mu|^{kq_{1}}\Big]\nonumber\\
     &= K_{kq_{1}}\epsilon^{-kq_{1}}p^{-kq_{1}}\Big[(\sum_{j=1}^{p}\Expe|w_{j}-\mu|^{2})^{kq_{1}/2}+\Expe\sum_{j=1}^{p}|w_{j}-\mu|^{kq_{1}}\Big]\nonumber\\
     &\leq K_{kq_{1}}\epsilon^{-kq_{1}}p^{-kq_{1}}\Big[(p\sigma^{2})^{kq_{1}/2}+p\Expe|w_{j}-\mu|^{kq_{1}}\Big]\nonumber\\
     &= K_{kq_{1}}\sigma^{kq_{1}}\epsilon^{-kq_{1}}\Big[p^{-kq_{1}/2}+p^{-kq_{1}+1}\Expe|\frac{w_{j}-\mu}{\sigma}|^{kq_{1}}\Big]\label{proest}.
\end{align}

Next, we prove \eqref{eq:quad-form-moment}. For any $q\geq 2$,
\begin{align}
    \Expe\Bigl|\br' \bA \br-\frac{1}{n}\nu_2\operatorname{tr} \bA \Bigr|^{q} \leq K_q\left(\Expe\Bigl|\br' \bA \br-\frac{1}{n}\frac{\sigma^{2}}{\mu^{2}}\operatorname{tr} \bA\Bigr|^{q} + \Expe\Bigl|\frac{1}{n}\frac{\sigma^{2}}{\mu^{2}}\operatorname{tr} \bA - \frac{1}{n}\nu_2\operatorname{tr} \bA\Bigr|^{q}\right).
\end{align}
By Lemma \ref{estsqu}, we have
\begin{align}\label{firest}
    \Expe\biggl|\frac{1}{n}\frac{\sigma^{2}}{\mu^{2}}\operatorname{tr} \bA - \frac{1}{n}\nu_2\operatorname{tr} \bA\biggr|^{q} \leq K_q p^{-q}\|\bA\|^{q}h_1^{q}.
\end{align}
Next, consider $\Expe\bigl|\br' \bA \br-\frac{1}{n}\frac{\sigma^{2}}{\mu^{2}}\operatorname{tr} \bA\bigr|^{q} $.
% \[
% \br' \bA \br-\frac{1}{n}\frac{\sigma^{2}}{\mu^{2}}\operatorname{tr} \bA = \Bigl(\br' \bA \br-\frac{1}{n}\frac{\sigma^{2}}{\mu^{2}}\operatorname{tr} \bA\Bigr)I(B_{p}(\epsilon)) + \Bigl(\br' \bA \br-\frac{1}{n}\frac{\sigma^{2}}{\mu^{2}}\operatorname{tr} \bA\Bigr)I(B^{c}_{p}(\epsilon)),
% \]
% where $I(\cdot)$ denotes the indicator function.
There exists a positive constant $K_{q}$ such that 
\begin{align}
    &\;\Expe \biggl|\br' \bA \br-\frac{1}{n}\frac{\sigma^{2}}{\mu^{2}}\operatorname{tr} \bA \biggr| \nonumber\\ 
    \leq &\; K_{q}\Big(\Expe \Bigl| \Bigl(\br' \bA \br-\frac{1}{n}\frac{\sigma^{2}}{\mu^{2}}\operatorname{tr} \bA\Bigr)I(B_{p}(\epsilon)) \Bigr|^{q} + \Expe \Bigl| \Bigl(\br' \bA \br-\frac{1}{n}\frac{\sigma^{2}}{\mu^{2}}\operatorname{tr} \bA\Bigr)I(B_{p}(\epsilon))\Bigr|^{q}\Big).\label{Bcest}
\end{align}

\noindent\textbf{Step 1}: Since $\bigl|\br' \bA \br-\frac{1}{n}\frac{\sigma^{2}}{\mu^{2}}\operatorname{tr} \bA \bigr| \leq \|\br\|\|\bA\|+C\|\bA\|\leq C\cdot p\|\bA\|$, we have 
\begin{equation}\label{eq:quad-form-moment-1}
\Expe \Bigl| \Bigl(\br' \bA \br-\frac{1}{n}\frac{\sigma^{2}}{\mu^{2}}\operatorname{tr} \bA\Bigr)I(B^{c}_{p}(\epsilon)) \Bigr|^{q}\leq C\cdot p^{q}\|\bA\|^{q} \Prob(B^{c}_{p}(\epsilon)).
\end{equation}

\noindent\textbf{Step 2}: Estimating $\Expe \biggl| [\br' \bA \br-\frac{1}{n}\frac{\sigma^{2}}{\mu^{2}}\operatorname{tr} \bA]I(B_{p}(\epsilon)) \biggr|^{q}$.
Write
\begin{align}
    \br' A \br-\frac{1}{n}\frac{\sigma^{2}}{\mu^{2}}\operatorname{tr} \bA &= \frac{1}{n\overline{w}^{2}} (\bw-\overline{w}\boldsymbol{1}_p)'\bA(\bw-\overline{w}\boldsymbol{1}_p)-\frac{1}{n}\frac{\sigma^{2}}{\mu^{2}}\operatorname{tr} \bA\nonumber\\
    &= \frac{\sigma^{2}}{n\overline{w}^{2}}\Big[ \frac{(\bw-\overline{w}\boldsymbol{1}_p)'}{\sigma}\bA \frac{(\bw-\overline{w}\boldsymbol{1}_p)}{\sigma}-\operatorname{tr} \bA\Big]+\Big[\frac{1}{n}\frac{\sigma^{2}}{\overline{w}^{2}}\operatorname{tr} \bA-\frac{1}{n}\frac{\sigma^{2}}{\mu^{2}}\operatorname{tr} \bA\Big]\nonumber\\
    &=: v_{1}+v_{2},\label{eq:quad-form-moment-2}
\end{align}
where $\boldsymbol{1}_p=(1,1,\ldots,1)'\in\bbR^p$.
For $0<\epsilon<1/2$, there exists a positive constant $K_{q}$ such that 
\[
  \Expe \biggl| [\br' \bA \br-\frac{1}{n}\frac{\sigma^{2}}{\mu^{2}}\operatorname{tr} \bA]I(B_{p}(\epsilon)) \biggr|^{q} \leq  K_{q} \Big( \Expe|v_{1}I(B_{p}(\epsilon))|^{q} + \Expe|v_{2}I(B_{p}(\epsilon))|^{q}\Big).
\]
By \eqref{proest}, we also have 
\begin{equation}\label{qjuest}
    \Expe |\bar{w} -\mu|^q \leq K_q\sigma^q(p^{-q/2}+p^{-q+1}\Expe \left|\frac{w_j-\mu}{\sigma}\right|^q).
\end{equation}
On the event $B_{p}(\epsilon)$, we have $-\epsilon\leq \overline{w}-\mu\leq \epsilon$, by this fact and \eqref{qjuest}, we get
\begin{align}\label{Bv2est}
    \Expe|v_{2}I(B_{p}(\epsilon))|^{q} &= \frac{1}{n}\frac{\sigma^{2}}{\overline{w}^{2}}\operatorname{tr} \bA-\frac{1}{n}\frac{\sigma^{2}}{\mu^{2}}\operatorname{tr} \bA\nonumber\\
    &= \sigma^{2q}\biggl|\frac{\operatorname{tr} \bA}{n}\biggr|^{q}  \Expe \biggl|(\frac{1}{\overline{w}^{2}}-\frac{1}{\mu})I(B_{p}(\epsilon))\biggr|^{q}\nonumber\\
    &= \sigma^{2q}\biggl|\frac{\operatorname{tr} \bA}{n}\biggr|^{q}  \Expe \biggl|(\frac{(\overline{w}-\mu)(\overline{w}+\mu)}{\overline{w}^{2}\mu^{2}})I(B_{p}(\epsilon))\biggr|^{q}\nonumber\\
    &= K_{\mu,\sigma,q} \biggl|\frac{\operatorname{tr} \bA}{n}\biggr|^{q} \Expe|\overline{w}-\mu|^{q}\nonumber\\
    & \leq K_{\mu,\sigma,q} \biggl|\frac{\operatorname{tr} \bA}{n}\biggr|^{q} \Big[p^{-q/2}+p^{-q+1}\Expe|\frac{w_{1}-\mu}{\sigma}|^{q}\Big]\nonumber\\
    & \leq K_{\mu,\sigma,q} \| \bA\|^{q} \Big[p^{-q/2}+p^{-q+1}\Expe|\frac{w_{1}-\mu}{\sigma}|^{q}\Big].
\end{align}
Next, we consider $v_{1}$. Note that $\frac{\bw-\overline{w}\boldsymbol{1}_p}{\sigma}=\frac{\bw-\mu\boldsymbol{1}_p}{\sigma}-\frac{\overline{w}\boldsymbol{1}_p-\mu\boldsymbol{1}_p}{\sigma}$, we get
\begin{align}\label{v1Best}
    \Expe|v_{1}I(B_{p}(\epsilon))|^{q} &= \Expe\biggl|\frac{\sigma^{2}}{n\overline{w}^{2}}\Big[ \frac{(\bw-\overline{w}\boldsymbol{1}_p)'}{\sigma}\bA \frac{(\bw-\overline{w}\boldsymbol{1}_p)}{\sigma}-\operatorname{tr} \bA\Big]I(B_{p}(\epsilon))\biggr|^{q}\nonumber\\
    &\leq K_{\mu,\sigma,q}n^{-q}\Expe\biggl|\Big[ \frac{(\bw-\overline{w}\boldsymbol{1}_p)'}{\sigma}\bA \frac{(\bw-\overline{w}\boldsymbol{1}_p)}{\sigma}-\operatorname{tr} \bA\Big]I(B_{p}(\epsilon))\biggr|^{q}\nonumber\\
    &\leq K_{\mu,\sigma,q}n^{-q}\Expe\biggl|\Big[\frac{(\bw-\mu\boldsymbol{1}_p)'}{\sigma}\bA \frac{(\bw-\mu\boldsymbol{1}_p)}{\sigma}-\operatorname{tr} \bA\Big]I(B_{p}(\epsilon))\biggr|^{q}\nonumber\\
    &\qquad + K_{\mu,\sigma,q}n^{-q}\Expe\biggl|\Big[\frac{(\overline{w}\boldsymbol{1}_p-\mu\boldsymbol{1}_p)'}{\sigma}\bA \frac{(\bw-\mu\boldsymbol{1}_p)}{\sigma}\Big]I(B_{p}(\epsilon))\biggr|^{q}\nonumber\\
    &\qquad + K_{\mu,\sigma,q}n^{-q}\Expe\biggl|\Big[\frac{(\bw-\mu\boldsymbol{1}_p)'}{\sigma}\bA \frac{(\overline{w}\boldsymbol{1}_p-\mu\boldsymbol{1}_p)}{\sigma}\Big]I(B_{p}(\epsilon))\biggr|^{q}\nonumber\\
    &\qquad + K_{\mu,\sigma,q}n^{-q}\Expe\biggl|\Big[\frac{(\overline{w}\boldsymbol{1}_p-\mu\boldsymbol{1}_p)'}{\sigma}\bA \frac{(\overline{w}\boldsymbol{1}_p-\mu\boldsymbol{1}_p)}{\sigma}\Big]I(B_{p}(\epsilon))\biggr|^{q}\nonumber\\
    &=: K_{\mu,\sigma,q}n^{-q}(V_{11}+V_{12}+V_{13}+V_{14}).
\end{align}
By Lemma 2.2 in \cite{bai2004clt}, we have
\begin{align}\label{v11Best}
    V_{11}\leq K_{q}\Big[\Big(\Expe\Bigl|\frac{w_{1}-\mu}{\sigma}\Bigr|^{4}\operatorname{tr} (\bA\bA')\Big)^{q/2}+\Expe\Bigl|\frac{w_{1}-\mu}{\sigma}\Bigr|^{2q}\operatorname{tr} (\bA\bA')^{q/2}\Big].
\end{align}
In the following, we consider $V_{12}$, $V_{13}$ and $V_{14}$. Note that $\frac{(\overline{w}-\mu)\boldsymbol{1}'_{p}}{\sigma}=\frac{1}{p}\boldsymbol{1}_{p}\boldsymbol{1}'_{p}\frac{(\bw-\mu\boldsymbol{1}_p)}{\sigma}$ 
and
\[\Expe\frac{1}{p}\operatorname{tr}(\boldsymbol{1}_{p}\boldsymbol{1}'_{p}\bA)=\Expe\frac{1}{p}\operatorname{tr}(\bA\boldsymbol{1}_{p}\boldsymbol{1}'_{p})=\Expe\frac{1}{p^{2}}\operatorname{tr}(\boldsymbol{1}_{p}\boldsymbol{1}'_{p}\bA\boldsymbol{1}_{p}\boldsymbol{1}'_{p})\leq \sqrt{p}\|\bA\|,
\]we get
\begin{align}\label{v12Best}
    V_{12}&= K_{\mu,\sigma,q}n^{-q}\Expe\biggl|\Big[\frac{1}{p}\frac{(\bw-\mu\boldsymbol{1}_p)'}{\sigma}\boldsymbol{1}_{p}\boldsymbol{1}'_{p}\bA \frac{(\bw-\mu\boldsymbol{1}_p)}{\sigma}\Big]I(B_{p}(\epsilon))\biggr|^{q}\nonumber\\
    &\leq K_{q}\Big(\Expe\biggl|\frac{(\bw-\mu\boldsymbol{1}_p)'}{\sigma}(\frac{1}{p}\boldsymbol{1}_{p}\boldsymbol{1}'_{p}\bA) \frac{(\bw-\mu\boldsymbol{1}_p)}{\sigma}-\operatorname{tr}(\frac{1}{p}\boldsymbol{1}_{p}\boldsymbol{1}'_{p}\bA)\biggr|^{q}+\Expe\biggl|\operatorname{tr}(\frac{1}{p}\boldsymbol{1}_{p}\boldsymbol{1}'_{p}\bA)\biggr|^{q}\Big)\nonumber\\
    &\leq K_{q}\Big(\Expe\biggl|\frac{(\bw-\mu\boldsymbol{1}_p)'}{\sigma}(\frac{1}{p}\boldsymbol{1}_{p}\boldsymbol{1}'_{p}\bA) \frac{(\bw-\mu\boldsymbol{1}_p)}{\sigma}-\operatorname{tr}(\frac{1}{p}\boldsymbol{1}_{p}\boldsymbol{1}'_{p}\bA)\biggr|^{q}+p^{q/2}\|\bA\|^{q}\Big)\nonumber\\
    & \leq V_{11}.
\end{align}
Similarly, we get
\begin{equation}\label{v134Best}
    V_{13}\leq V_{11}, \qquad V_{14}\leq V_{11}.
\end{equation}
By \eqref{v1Best}--\eqref{v134Best}, we get
\begin{align}\label{Bv1est}
    \Expe|v_{1}I(B_{p}(\epsilon))|^{q}\leq K_{q}n^{-q}\Big[\Big(\Expe\biggl|\frac{w_{1}-\mu}{\sigma}\biggr|^{4}\operatorname{tr} (\bA\bA')\Big)^{q/2}+\Expe\biggl|\frac{w_{1}-\mu}{\sigma}\biggr|^{2q}\operatorname{tr} (\bA\bA')^{q/2}\Big].
\end{align}
From \eqref{Bcest} --\eqref{Bv2est} and \eqref{Bv1est}, we have
\begin{align}\label{secest}
    &\;\Expe \Bigl|\br' \bA \br-\frac{1}{n}\frac{\sigma^{2}}{\mu^{2}}\operatorname{tr} \bA \Bigr|^{q}\nonumber\\
    \leq &\;K_{q}\Bigl(n^{-q}\Big[\Big(\Expe\Bigl|\frac{w_{1}-\mu}{\sigma}\Bigr|^{4}\operatorname{tr} (\bA\bA')\Big)^{q/2}+\Expe\Bigl|\frac{w_{1}-\mu}{\sigma}\Bigr|^{2q}\operatorname{tr} (\bA\bA')^{q/2}\Big]\nonumber\\
    &\qquad+\| \bA\|^{q} \Big[p^{-q/2}+p^{-q+1}\Expe|\frac{w_{j}-\mu}{\sigma}|^{q}\Big]+p^{q}\|\bA\|^{q} \Prob(B^{c}_{p}(\epsilon))\biggr)\nonumber\\
    \leq&\; K_{q}\Bigl(n^{-q}\Big[\Big(\Expe\Bigl|\frac{w_{1}-\mu}{\sigma}\Bigr|^{4}\operatorname{tr} (\bA\bA')\Big)^{q/2}+\Expe\Bigl|\frac{w_{1}-\mu}{\sigma}\Bigr|^{2q}\operatorname{tr} (\bA\bA')^{q/2}\Big]+n^{q}\|\bA\|^{q} \Prob(B^{c}_{p}(\epsilon))\Bigr).
\end{align}
By \eqref{firest} and \eqref{secest}, we get
\begin{align}\label{finest}
    &\;\Expe\biggl|\br' \bA r - \frac{1}{n}\nu_2\operatorname{tr} \bA\biggr|^{q} \nonumber\\
    \leq&\; K_{q}\Bigl(n^{-q}\Big[\Big(\Expe\Bigl|\frac{w_{1}-\mu}{\sigma}\Bigr|^{4}\operatorname{tr} (\bA\bA')\Big)^{q/2}+\Expe\Bigl|\frac{w_{1}-\mu}{\sigma}\Bigr|^{2q}\operatorname{tr} (\bA\bA')^{q/2}\Big]\nonumber\\
    &\qquad\qquad +n^{q}\|\bA\|^{q} \Prob(B^{c}_{p}(\epsilon))+n^{-q}\|\bA\|^{q}h_1^{q}\Bigr).
\end{align}

If $\Expe\left|w_{1}-\mu\right|^{4}<\infty$, $\|\bA\|\leq K$ and $|\frac{w_{j}-\mu}{\sigma}|\leq \sqrt{n}\delta_{n}$, then, for any $q\geq 2$,
\begin{align}
n^{-q}\Big[\Big(\Expe\biggl|\frac{w_{1}-\mu}{\sigma}\biggr|^{4}\operatorname{tr} (\bA\bA')\Big)^{q/2}+\Expe\biggl|\frac{w_{1}-\mu}{\sigma}\biggr|^{2q}\operatorname{tr} (\bA\bA')^{q/2}\Big] \leq K_q n^{-1}\delta_{n}^{2q-4}
\end{align}
and
\begin{align}
    n^{-q}\|\bA\|^{q}h_1^{q} \leq K_q n^{-1}\delta_{n}^{2q-4}.
\end{align}
Taking $\epsilon = n^{-\alpha}$, $0<\alpha<1/2$ and $kq_{1}>\frac{3q}{1-2\alpha}$ yields that
\[
    n^{q}\|\bA\|^{q} \Prob(B^{c}_{p}(\epsilon)) \leq K_q n^{-1}\delta_{n}^{2q-4}.
\]
Hence, we have
\[
    \Expe \biggl|r' \bA r-\frac{1}{n}\nu_2\operatorname{tr} \bA \biggr|^{q} \leq K_{q}n^{-1}\delta_{n}^{2q-4},
\]
which is derived by the following calculation,
\begin{align*}
    \Prob(B^{c}_{p}(\epsilon)) &\leq K_{kq_{1}}\sigma^{kq_{1}}\epsilon^{-kq_{1}}\Big[p^{-kq_{1}/2}+p^{-kq_{1}+1}\Expe\biggl|\frac{w_{j}-\mu}{\sigma}\biggr|^{kq_{1}}\Big]\nonumber\\
    &=K_{kq_{1}}\sigma^{kq_{1}}(P_1+P_2),
\end{align*}
and
\begin{align*}
    P_2 &= \epsilon^{-kq_{1}}p^{-kq_{1}+1}\Expe\biggl|\frac{w_{j}-\mu}{\sigma}\biggr|^{kq_{1}}\nonumber\\
    &\leq \epsilon^{-kq_{1}}p^{-kq_{1}+1}(n^{1/2}\delta_{n})^{kq_{1}-4}\nonumber\\
    &\leq \epsilon^{-kq_{1}}n^{-kq_{1}/2-1}\nonumber\\
    &\leq P_1.
\end{align*}
Therefore, 
\begin{align}\label{prordset}
    \Prob(B^{c}_{p}(\epsilon)) &\leq 2K_{kq_{1}}\sigma^{kq_{1}}P_1.
\end{align}
Take $\epsilon = n^{-\alpha}$ $(0<\alpha<1/2)$ and $kq_{1}>\frac{3q}{1-2\alpha}$ into $\eqref{prordset}$, we have
\[
    n^q\Prob(B^{c}_{p}(\epsilon)) 
    \leq K_{q}n^qn^{-(\frac{1}{2}-\alpha)kq_{1}} \leq K_{q}n^{-q/2}\leq K_{q}n^{-1}\delta_{n}^{2q-4}.
\]

\subsection{Proof of Lemma \ref{quadform}}\label{prfquadform}

It is obvious that
\begin{align}
    &\Expe \Bigl(\br'\bA\br-\frac{1}{n}\nu_2\operatorname{tr} \bA\Bigr) \Bigl(\br'\bB\br-\frac{1}{n}\nu_2\operatorname{tr} \bB\Bigr)\nonumber\\
    =& \Expe\Bigl(\br' \bA \br \br' \bB \br\Bigr) -\Bigl(\frac{1}{n}\nu_2\operatorname{tr} \bB \Bigr)\Expe\Bigl(\br'\bA\br-\frac{1}{n}\nu_2\operatorname{tr} \bA\Bigr)\nonumber\\
    &\quad-\Bigl(\frac{1}{n}\nu_2\operatorname{tr} \bA\Bigr) \Expe\Bigl(\br'\bB\br-\frac{1}{n}\nu_2\operatorname{tr} \bB\Bigr)-\frac{1}{n^{2}}\nu_2^2\operatorname{tr} \bA\operatorname{tr} \bB.
\end{align}
\noindent\textbf{Step 1}: Consider $\Expe\left(\br'\bA\br-\frac{1}{n}\nu_2\operatorname{tr} \bA\right)$.

Recall that $\br = \frac{1}{\sqrt{n}}(\frac{w_{1}}{\overline{w}}-1,\ldots,\frac{w_{p}}{\overline{w}}-1)'$. 
Note that $\{w_j,1\leq j\leq p\}$ are i.i.d., and
\begin{align}\label{summean}
    \sum_{j=1}^{p}\frac{w_j}{\overline{w}} = p,
\end{align}
and
\[
    \sum_{j=1}^{p} \left(\frac{w_j}{\overline{w}}-1\right)^{2} + \sum_{j_1\neq j_2}\left(\frac{w_{j_1}}{\overline{w}}-1\right)\left(\frac{w_{j_2}}{\overline{w}}-1\right) = 0.
\]
Taking expectation on the above two identities yields that $\Expe \frac{w_j}{\overline{w}} = 1$,
and, for $1\leq j_1\neq j_2\leq p$, 
\[
    \Expe\Bigl(\frac{w_{j_1}}{\overline{w}}-1\Bigr)\Bigl(\frac{w_{j_2}}{\overline{w}}-1\Bigr) = -\frac{\nu_2}{p-1},
\]
where $\nu_{2} = \Var(\frac{w_j}{\overline{w}})$.
This implies that
\begin{equation}\label{covmatrix}
    \Expe (\br\br') = \frac{\nu_{2}}{n} \left(-\frac{1}{p-1}\boldsymbol{1}_{p}\boldsymbol{1}'_{p} + \frac{1}{p-1}\bI_p + \bI_p \right).
\end{equation}
By \eqref{covmatrix}, we have
\begin{align}
    \Expe\left(\br'\bA\br-\frac{\nu_2}{n}\operatorname{tr} \bA\right) &= \Expe\operatorname{tr}\left(\bA\br\br'\right) - \frac{\nu_2}{n}\operatorname{tr} \bA \nonumber\\
    &= \operatorname{tr}(\bA \Expe rr') - \frac{\nu_2}{n}\operatorname{tr} \bA\nonumber\\
    &= -\frac{\nu_2}{n(p-1)}\sum_{k\neq l}\bA_{kl},\nonumber
\end{align}
thus, 
\begin{equation}\label{secoterm}
    \left(\frac{1}{n}\nu_2\operatorname{tr} \bB \right)\Expe\left(\br'\bA\br-\frac{1}{n}\nu_2\operatorname{tr} \bA\right) 
    = -\frac{1}{n^{2}(p-1)}\nu_{2}^{2}\operatorname{tr} \bB\sum_{k\neq l}\bA_{kl}.
\end{equation}
Similarly, 
\begin{equation}\label{thirterm}
    \left(\frac{1}{n}\nu_2\operatorname{tr} \bA \right)\Expe\left(\br'\bB\br-\frac{1}{n}\nu_2\operatorname{tr} \bB\right) 
    = -\frac{1}{n^{2}(p-1)}\nu_2^{2}\operatorname{tr} \bA\sum_{k\neq l}\bB_{kl}.
\end{equation}

\noindent\textbf{Step 2}: Consider $\Expe\left(\br' \bA \br \br' \bB \br\right)$.

Let $R_j = \frac{1}{\sqrt{n}}\bigl(\frac{w_{j}}{\overline{w}}-1\bigr)$, $1\leq j \leq p$, then we have $\sum_{j=1}^p R_j = 0$, and  
\begin{align*}
\nu_4 &:= \Expe\left(\frac{w_{1}}{\overline{w}}-1 \right)^{4} = n^2\Expe R_1^4,\\
\nu_{12} &:= \Expe\bigl[\bigl(\frac{w_{1}}{\overline{w}}-1 \bigr)^{2}\bigl(\frac{w_{2}}{\overline{w}}-1 \bigr)^{2}\bigr] = n^2\Expe \left(R_1^2R_2^2\right).
\end{align*}
It follows from \eqref{summean} that
\begin{equation}\label{r1112}
    \Expe R_1^3 R_2 
    =\frac{1}{p-1} \Expe \Bigl[R_{1}^{3}\Bigl(\sum_{j=1}^{p}R_j-R_1\Bigr)\Bigr] = -\frac{1}{p-1}\nu_4,
\end{equation}
\begin{align}\label{r1123}
    \Expe R_1^2 R_2 R_3 &= \frac{1}{p-2}\Expe \Bigl[R_1^2 R_2 \Bigl(\sum_{j=1}^{p}R_j-R_1-R_2\Bigr)\Bigr]\nonumber\\
     &= -\frac{1}{p-2}\Expe (R_1^3 R_2)-\frac{1}{p-2}\Expe (R_1^2 R_2^2)\nonumber\\
    &= \frac{1}{(p-1)(p-2)}\nu_4 - \frac{1}{p-2}\nu_{12},
\end{align}
\begin{align}\label{r1234}
    \Expe R_1 R_2 R_3 R_4 &= \frac{1}{p-3}\Expe \Bigl[R_1 R_2 R_3\Bigl(\sum_{j=1}^{p}R_j-R_1-R_2-R_3\Bigr)\Bigr]\nonumber\\
    &= -\frac{3}{p-3}\Expe (R_1^2 R_2 R_3)\nonumber\\
    &= -\frac{3}{(p-1)(p-2)(p-3)}\nu_4 + \frac{3}{(p-2)(p-3)}\nu_{12}.
\end{align}

To calculate  $\Expe\left(\br' \bA \br \br' \bB \br\right)$, we expand it as
\begin{align}
\Expe\left(\br' \bA \br \br' \bB \br\right) 
=  \Expe\Bigl(\sum_{i, j} R_{i} A_{i j} R_{j} \sum_{k, l} R_{k} B_{k l} R_{l}\Bigr)=\sum_{i, j, k, l} \Expe \left(R_{i}R_{j}R_{k}R_{l}  A_{i j} B_{k l} \right).\label{quadexpa}
\end{align}
To calculate \eqref{quadexpa}, we split it into the following cases:
\begin{itemize}
	\item[1] \[i=j=k=l, \; \sum_{i}\left(R_{i}^{4}\right) A_{i i} B_{i i};\]
	\item[2] \[i=j, k=l, i \neq k,  \sum_{\substack{i, k \\ i \neq k}} \left(R_{i}^{2}R_{k}^{2}\right)A_{i i} B_{k k}; \]
	\item[3]  \[i=j, k \neq l, \sum_{\substack{i, k, l \\ k \neq l}} \left(R_{i}^{2}R_{k}R_{l}\right) A_{i i} B_{k l}; \]
	\item[4]  \[i \neq j, k=l, \sum_{\substack{i, j, k \\ i \neq j}} \left(R_{i}R_{j}R_{k}^{2}\right)A_{i j} B_{k k}; \]
	\item[5]   \[i \neq j, k \neq l, i=k, j=l, \sum_{\substack{i, j \\ i \neq j}}\left(R^{2}_{i} R^{2}_{j}\right) A_{i j} B_{i j} ;\]
	\item[6]   \[i \neq j, k \neq l, i=l, j=k, \sum_{\substack{i, j\\ i \neq j}}\left(R^{2}_{i} R^{2}_{j}\right) A_{i j} B_{j i} ;\]
	\item[7]   \[i \neq j, k \neq l, i=k, l \neq j, \sum_{\substack{i, j, l \\ i \neq j \neq l}}\left(R_{i}R_{j} R_{i}R_{l}\right) A_{i j} B_{i l} ;\]
	\item[8]   \[i \neq j, k \neq l, l=j, i \neq k, \sum_{\substack{i, j, k \\ i \neq j \neq k}}\left(R_{i}R_{j} R_{k}R_{j}\right) A_{i j} B_{k j};\]
	\item[9]   \[i \neq j, k \neq l, k=j, i \neq l, \sum_{\substack{i, j, l \\ i \neq j \neq l}}\left(R_{i}R_{j} R_{j}R_{l}\right) A_{i j} B_{j l};\]
	\item[10]   \[i \neq j, k \neq l, i=l, k \neq j, \sum_{\substack{i, j, k \\ i \neq j \neq k}}\left(R_{i}R_{j} R_{k}R_{i}\right) A_{i j} B_{k i} ;\]
	\item[11]   \[i \neq j, k \neq l, l \neq j, i \neq k, \sum_{\substack{i, j, k, l \\ i \neq j \neq k \neq l}}\left(R_{i}R_{j} R_{k}R_{l}\right)A_{i j} B_{k l}.\]
\end{itemize}
For ease of presentation, we still keep  $\nu_{4}$ in the expectations although we have obtained its value.
The expectations of all cases are listed as follows.

\noindent \underline{Case 1}:  
\[
\Expe\sum_{i}\left(R_{i}^{4}\right) A_{i i} B_{i i} = \Expe \left(R_{i}^{4}\right) \sum_{i} A_{i i} B_{i i} = \frac{1}{n^2}\nu_{4} \sum_{i} A_{i i} B_{i i}.
\]
    
\noindent \underline{Case 2}: 
    \begin{align}
        \Expe\sum_{\substack{i, k \\ i \neq k}} \left(R_{i}^{2}R_{k}^{2}\right)A_{i i} B_{k k} &= \Expe \left(R_{i}^{2}R_{k}^{2}\right)\sum_{\substack{i, k \\ i \neq k}}A_{i i}B_{k k} = \frac{\nu_{12}}{n^2}\sum_{\substack{i, k \\ i \neq k}}A_{i i}B_{k k} \nonumber\\
        &= \frac{\nu_{12}}{n^2}\Bigl(\operatorname{tr} \bA\operatorname{tr} \bB-\sum_{i}A_{i i} B_{i i}\Bigr).\nonumber
    \end{align}
    
\noindent \underline{Case 3}: From \eqref{r1112} and \eqref{r1123}, we obtain
\begin{align}\label{precase3}
\Expe\sum_{\substack{i, k, l \\ k \neq l}} \left(R_{i}^{2}R_{k}R_{l}\right) A_{i i} B_{k l}
=&  \Expe R_{1}^{2}R_{2}R_{3} \sum_{\substack{i, k, l \\
i \neq k \neq l}} A_{i i} B_{k l}+\Expe R_{1}^{3}R_{2} \sum_{\substack{i, l \\
l \neq i}} A_{i i} B_{i l}+\Expe R_{1}^{3}R_{2}\sum_{\substack{i, k \\
k \neq i}} A_{i i} B_{k i} \nonumber\\
=&  \left(\frac{\nu_{4}}{n^2(p-1)(p-2)} -\frac{\nu_{12}}{n^2(p-2)}\right) \sum_{\substack{i, k, l \\
i \neq k \neq l}} A_{i i} B_{k l} \nonumber\\
&\quad-\frac{\nu_{4}}{n^2(p-1)}  \Biggl(\sum_{\substack{i, l \\
l \neq i}} A_{i i}B_{i l}+\sum_{\substack{i, k \\
k \neq i}}A_{i i}B_{k i}\Biggr).
\end{align}

Note that
\begin{align}
\Expe\operatorname{tr}\bA\leq \Expe(p\sum_iA_{ii}^2)^{1/2}\leq \Expe(p\operatorname{tr}(\bA\bA'))^{1/2}\leq\Expe (p^2\|\bA\|^2)^{1/2}=p\Expe\|\bA\|,
\end{align}
and
\begin{align}\label{triiklall}
     \Expe \Bigl|\sum_{\substack{i, k, l }} A_{i i} B_{k l}\Bigr| &= \Expe \Bigl|\operatorname{tr}\bA \boldsymbol{1}_{p}'\bB\boldsymbol{1}_{p}\Bigr| \nonumber\\
     &\leq p\Expe(\|\bA\|\cdot |\boldsymbol{1}_{p}'\bB\boldsymbol{1}_{p}|)\nonumber\\
     % &\leq p\Expe\|\bA\|\cdot\|\boldsymbol{1}_{p}'\|\cdot\|\bB\boldsymbol{1}_{p}\|\nonumber\\
     % &\leq p\Expe\|\bA\|\cdot\|\bB\|\cdot\|\boldsymbol{1}_{p}\|^2\nonumber\\
     &\leq p^2\Expe(\|\bA\|\cdot\|\bB\|)\nonumber\\
      &\leq p^2(\Expe\|\bA\|^{q_1})^{1/q_1}(\Expe\|\bB\|)^{1/q_2}\nonumber\\
     &= O(p^{2}),
\end{align}
and, for $q_1, q_2 >0$ and $1/q_1+1/q_2=1$,
\begin{align}\label{triiil}
    \Expe\Bigl|\sum_{\substack{i, l}} A_{i i}B_{i l}\Bigr|% &=\Expe\left|\sum_{i}A_{i i}\sum_{l}B_{i,l}\right|\nonumber\\
   &= \Expe\Bigl|\sum_{\substack{i}} A_{i i}(\boldsymbol{1}_{p}^{i})'\bB\boldsymbol{1}_{p}\Bigr|\nonumber\\
   &\leq \Expe\Bigl(\Bigl|\sum_{\substack{i}} A_{i i}^{2}\Bigr|\Bigr)^{1/2}\Bigl(\Bigl[\sum_{i}\left((\boldsymbol{1}_{p}'\bB\boldsymbol{1}_{p}^{i})\right)\left((\boldsymbol{1}_{p}^{i})'\bB\boldsymbol{1}_{p}\right)\Bigr]\Bigr)^{1/2}\nonumber\\
   &\leq p\Expe\|\bA\|\cdot\|(\boldsymbol{1}_{p}'\bB\boldsymbol{1}^{1}_{p})((\boldsymbol{1}^{1}_{p})'\bB\boldsymbol{1}_{p})\|^{1/2}\nonumber\\
   &\leq p^{3/2}\Expe\|\bA\|\cdot\|\bB\|\nonumber\\
   &\leq p^{3/2}(\Expe\|\bA\|^{q_1})^{1/q_1}(\Expe\|\bB\|)^{1/q_2}\nonumber\\
   &= O(p^{3/2}),
\end{align}
where $\boldsymbol{1}_p^{i}$ is the $p$-dimensional vector with all components being $0$ except for the $i$-th component being $1$,
and
\begin{align}\label{triill}
     \Expe\sum_{\substack{i, l}} A_{i i}B_{l l} =  \Expe\operatorname{tr} \bA\operatorname{tr}\bB \leq p^2\Expe(\|\bA\|\cdot\|\bB\|)=O(p^2),
\end{align}
and
\begin{align}\label{triiii}
    \Expe\sum_{i} A_{i i}B_{i i} \leq p  \Expe\|\diag(\bA)\diag(\bB)\|\leq p\Expe(\|\bA\|\cdot\|\bB\|).
\end{align}
Thus, by \eqref{triiklall}--\eqref{triiii}, we have
\begin{align}\label{triikl}
     \Expe\sum_{\substack{i, k, l \\
i \neq k \neq l}} A_{i i} B_{k l} = \Expe\left(\sum_{\substack{i, k, l }} A_{i i} B_{k l} - \sum_{\substack{i, l}} A_{i i}B_{i l} - \sum_{\substack{i, l}} A_{i i}B_{l l} + \sum_{i} A_{i i} B_{i i}\right) = O(p^{2}).
\end{align}

It follows from \eqref{precase3}--\eqref{triikl} that
\[
\Expe\sum_{\substack{i, k, l \\ k \neq l}} \left(R_{i}^{2}R_{k}R_{l}\right) A_{i i} B_{k l}
=  -\frac{\nu_{12}}{n^2(p-2)} \sum_{\substack{i, k, l \\
i \neq k \neq l}} A_{i i} B_{k l} + o(n^{-1}).
\]

\noindent \underline{Case 4}: Similarly to Case 3, one can conclude that
% \begin{align}\label{precase4}
%     \Expe \sum_{\substack{i, j, k \\ i \neq j}} \left(R_{i}R_{j}R_{k}^{2}\right)A_{i j} B_{k k} 
% =&  \Expe R_{1}^{2}R_{2}R_{3} \sum_{\substack{i, j, k \\
% i \neq j \neq k}} A_{i j} B_{k k}+\Expe R_{1}^{3}R_{2} \sum_{\substack{j, k \\
% j \neq k}} A_{k j} B_{k k}+\Expe R_{1}^{3}R_{2}\sum_{\substack{i, k \\
% k \neq i}} A_{i k} B_{k k} \nonumber\\
% =&  \left(\frac{\nu_{4}}{n^2(p-1)(p-2)} -\frac{\nu_{12}}{n^2(p-2)}\right) \sum_{\substack{i, j, k \\
% i \neq j \neq k}} A_{i j} B_{k k} \nonumber\\
% &\quad-\frac{\nu_{4}}{n^2(p-1)}  \Biggl(\sum_{\substack{k, j \\
% k \neq j}} A_{k j}B_{k k}+\sum_{\substack{i, k \\
% k \neq i}}A_{i k}B_{k k}\Biggr).
% \end{align}

% It follows from \eqref{precase4} and \eqref{triiil}--\eqref{triikl} that
\[
    \Expe \sum_{\substack{i, j, k \\ i \neq j}} \left(R_{i}R_{j}R_{k}^{2}\right)A_{i j} B_{k k} 
=  -\frac{\nu_{12}}{n^2(p-2)} \sum_{\substack{i, j, k \\
i \neq j \neq k}} A_{i j} B_{k k} + o(n^{-1}).
\]

\noindent \underline{Case 5}:
\[
     \Expe \sum_{\substack{i, j \\ i \neq j}}\left(R^{2}_{i} R^{2}_{j}\right) A_{i j} B_{i j} 
    = \Expe \left(R^{2}_{1} R^{2}_{2}\right)  \sum_{\substack{i, j \\ i \neq j}} A_{i j} B_{i j} = \frac{\nu_{12}}{n^2} \bigl(\operatorname{tr} (\bA\bB')-\sum_{i} A_{i i} B_{i i}\bigr).
\]

\noindent \underline{Case 6}: 
\[
   \Expe \sum_{\substack{i, j \\ i \neq j}}\left(R^{2}_{i} R^{2}_{j}\right) A_{i j} B_{j i} =  \Expe \left(R^{2}_{1} R^{2}_{2}\right)  \sum_{\substack{i, j \\ i \neq j}} A_{i j} B_{j i} = \frac{\nu_{12}}{n^2} \bigl(\operatorname{tr} (\bA\bB)-\sum_{i} A_{i i} B_{i i}\bigr).
\]

\noindent \underline{Case 7}: By \eqref{r1123}, we have
\begin{align}\label{precase7}
    \Expe \sum_{\substack{i, j, l \\ i \neq j \neq l}} \left(R_{i}R_{j} R_{i}R_{l}\right) A_{i j} B_{i l} &= \Expe \left(R_{1}^2R_{2} R_{3}\right) \sum_{\substack{i, j, l \\ i \neq j \neq l}}A_{i j} B_{i l}\nonumber\\
    &= \biggl(\frac{1}{n^2(p-1)(p-2)} \nu_{4}-\frac{1}{n^2(p-2)}\nu_{12}\biggr) \sum_{\substack{i, j, l \\
i \neq j \neq l}} A_{i j} B_{i l}.
\end{align}

Note that,  $\Expe\boldsymbol{1}_{p}'\bA'\bB\boldsymbol{1}_{p}\leq p \Expe(\|\bA\|\cdot\|\bB\|)$, $\Expe\operatorname{tr} (\bA'\bB)\leq p \Expe(\|\bA\|\cdot\|\bB\|)$,
and by \eqref{triiil} and \eqref{triiii}, we have
\begin{align}
    \Expe\sum_{i\neq l} A_{i i} B_{i l}=\Expe\sum_{i, l} A_{i i} B_{i l}-\Expe\sum_{i} A_{i i} B_{i i}=O(p^{3/2}),
\end{align}
thus,
\begin{align}\label{trikil}
     \Expe \sum_{\substack{i, j, l \\
i \neq j \neq l}} A_{i j} B_{i l} =  \Expe \Bigl(\boldsymbol{1}_{p}'\bA'\bB\boldsymbol{1}_{p} - \operatorname{tr} (\bA'\bB) - 2\sum_{i\neq l} A_{i i} B_{i l}\Bigr)= O(p^{3/2}).
\end{align}

It follows from \eqref{precase7}--\eqref{trikil} that
\[
    \Expe \sum_{\substack{i, j, l \\ i \neq j \neq l}} \left(R_{i}R_{j} R_{i}R_{l}\right) A_{i j} B_{i l} = o(n^{-1}).
\]

\noindent \underline{Case 8}: Similarly to Case 7, we have
% \begin{align}
%     \Expe \sum_{\substack{i, j, k \\ i \neq j \neq k}} \left(R_{i}R_{j} R_{k}R_{j}\right) A_{i j} B_{k j} &= \Expe \left(R_{i}R_{j} R_{k}R_{j}\right)  \sum_{\substack{i, j, k \\ i \neq j \neq k}}^{p} A_{i j} B_{k j} \nonumber\\
%     &= \left(\frac{1}{n^2(p-1)(p-2)} \nu_{4}-\frac{1}{n^2(p-2)}\nu_{12}\right)  \sum_{\substack{i, j, k \\ i \neq j \neq k}}^{p} A_{i j} B_{k j}.
% \end{align}
% Similarly to \eqref{trikil}
% \begin{equation}
%     \sum_{\substack{i, j, k \\ i \neq j \neq k}}^{p} A_{i j} B_{k j} = O(p).
% \end{equation}
% Thus,
\[
    \Expe \sum_{\substack{i, j, k \\ i \neq j \neq k}} \left(R_{i}R_{j} R_{k}R_{j}\right) A_{i j} B_{k j} = o(n^{-1}).
\]

\noindent \underline{Case 9}: Similarly to Case 7, we have 

\[
    \Expe \sum_{\substack{i, j, l \\ i \neq j \neq l}}\left(R_{i}R_{j} R_{j}R_{l}\right) A_{i j} B_{j l} = \left(\frac{\nu_{4}}{n^2(p-1)(p-2)} -\frac{\nu_{12}}{n^2(p-2)}\right) \sum_{\substack{i, j, l \\ i \neq j \neq l}} A_{i j} B_{j l} = o(n^{-1}).
\]

\noindent \underline{Case 10}: Similarly to Case 7, we have 

\[
    \Expe \sum_{\substack{i \neq j, k \\ i \neq j \neq k}} \left(R_{i}R_{j} R_{k}R_{i}\right) A_{i j} B_{k i} = \left(\frac{\nu_{4}}{n^2(p-1)(p-2)} -\frac{\nu_{12}}{n^2(p-2)}\right) \sum_{\substack{i, j, k \\ i \neq j \neq k}} A_{i j} B_{k i} = o(n^{-1}).
\]

\noindent \underline{Case 11}: By \eqref{r1234}, we have
\begin{align}
    &\quad\Expe \sum_{\substack{i, j, k, l \\ i \neq j \neq k \neq l}}\left(R_{i}R_{j} R_{k}R_{l}\right) A_{i j} B_{k l} \nonumber\\
    &= \Expe \left(R_{1}R_{2} R_{3}R_{4}\right) \sum_{\substack{i, j, k, l \\ i \neq j \neq k \neq l}} A_{i j} B_{k l} \nonumber\\
    &= \left(-\frac{3\nu_{4}}{n^2(p-1)(p-2)(p-3)} +\frac{3\nu_{12}}{n^2(p-2)(p-3)}\right) \sum_{\substack{i, j, k, l \\ i \neq j \neq k \neq l}} A_{i j} B_{k l}.\nonumber
\end{align}
Note that, by \eqref{trikil} we have
\begin{align}
    \Expe\sum_{\substack{i, j, k, l \\ i \neq j \neq k \neq l}} A_{i j} B_{k l}& = \Expe\left(\boldsymbol{1}_{p}'\bA\boldsymbol{1}_{p} - \operatorname{tr}\bA\right)\left(\boldsymbol{1}_{p}'\bB\boldsymbol{1}_{p} - \operatorname{tr}\bB\right) - 2\Expe\sum_{\substack{i, k, l \\ i\neq k \neq l}} A_{i k} B_{i l} \nonumber\\
    &\leq Kp^2 \Expe(\|\bA\|\cdot\|\bB\|)-2\Expe\sum_{\substack{i, k, l \\ i\neq k \neq l}} A_{i k} B_{i l}= O(p^{2}).\nonumber
\end{align}
Thus,
\[
   \Expe \sum_{\substack{i, j, k, l \\ i \neq j \neq k \neq l}}\left(R_{i}R_{j} R_{k}R_{l}\right) A_{i j} B_{k l} = o(n^{-1}).
\]

%Therefore, we only need to consider terms including $\nu_{4}$, $\nu_{12}$, $\frac{1}{p-1}\nu_{4}$, $\frac{1}{p-1}\nu_{12}$ and constant terms.
Combining Cases 1 -- 11 gives us 
\begin{align}
    \Expe\left(\br' \bA \br \br' \bB \br\right)  
%     =& \sum_{i=1}^{p}\frac{1}{n^{2}}\Big(\nu_4-3\nu_{12} \Big) A_{ii}B_{ii}\nonumber\\
%     &+\frac{1}{n^{2}}\nu_{12}\Big(\operatorname{tr}\bA\operatorname{tr}\bB + \operatorname{tr} (\bA\bB')+\operatorname{tr} (\bA\bB)\Big)\nonumber\\
%     &-\frac{1}{n^2(p-2)}\nu_{12} \sum_{\substack{i, k, l \\
% i \neq k, k\neq l, i\neq l}} A_{i i} B_{k l}-\frac{1}{n^2(p-2)}\nu_{12} \sum_{\substack{i, j, k \\
% i \neq j, j \neq k, i\neq k}} A_{i j} B_{k k}+ o(n^{-1})\nonumber\\
=& \frac{\nu_4-3\nu_{12}}{n^{2}} \sum_{i=1}^{p}A_{ii}B_{ii}+\frac{\nu_{12}}{n^{2}}\Big(\operatorname{tr}\bA\operatorname{tr}\bB + \operatorname{tr} (\bA\bB')+\operatorname{tr} (\bA\bB)\Big)\nonumber\\
    &-\frac{1}{n^2(p-2)}\nu_{12} \Biggl(\sum_{\substack{i, k, l \\
i \neq k, k \neq l, i\neq l}} A_{i i} B_{k l} +\sum_{\substack{i, j, k \\
i \neq j, j \neq k, i\neq k}} A_{i j} B_{k k}\Biggr)+ o(n^{-1}).\nonumber
\end{align}

\noindent\textbf{Step 3}: By \eqref{secoterm}--\eqref{thirterm}, we have
\begin{align}
    &-\left(\frac{\nu_2}{n}\operatorname{tr} \bB \right)\Expe\left(\br'\bA\br-\frac{\nu_2}{n}\operatorname{tr} \bA\right) - \left(\frac{\nu_2}{n}\operatorname{tr} \bA \right)\Expe\left(\br'\bB\br-\frac{\nu_2}{n}\operatorname{tr} \bB\right) \nonumber\\
    =& \frac{\nu_{2}^{2}}{n^{2}(p-1)}\operatorname{tr} \bB\sum_{k\neq l}\bA_{kl} +\frac{\nu_{2}^{2}}{n^{2}(p-1)}\operatorname{tr} \bA\sum_{k\neq l}\bB_{kl}\nonumber\\
    =& \frac{\nu_{2}^{2}}{n^{2}(p-1)}\nonumber\\
    &\quad\times \left(\sum_{\substack{k, l \\ k \neq l}} B_{k k}A_{k l} + \sum_{\substack{k, l \\ k \neq l}} A_{k k} B_{k l}+\sum_{\substack{k, l \\ k \neq l}} B_{l l}A_{k l} + \sum_{\substack{k, l \\ k \neq l}} A_{l l} B_{k l}+\sum_{\substack{i, k, l \\
i \neq k, k \neq l, i\neq l}} B_{i i}A_{k l} + \sum_{\substack{i, k, l \\
i \neq k, k \neq l, i\neq l}} A_{i i} B_{k l}\right)\nonumber\\
    =& \frac{\nu_{2}^{2}}{n^{2}(p-1)}\left(\sum_{\substack{i, k, l \\
i \neq k, k \neq l,i\neq l}} B_{i i}A_{k l} + \sum_{\substack{i, k, l \\
i \neq k, k \neq l, i\neq l}} A_{i i} B_{k l}\right) + o(n^{-1}).\nonumber
\end{align}

Following the above calculations, we have 
\begin{align*}
    &\Expe \Big(\br'A\br-\frac{1}{n}\nu_2\operatorname{tr} \bA\Big)\Big(\br'\bB\br-\frac{1}{n}\nu_2\operatorname{tr} \bB\Big)\nonumber\\
    =& \frac{\nu_4-3\nu_{12}}{n^{2}} \sum_{i=1}^{p}A_{ii}B_{ii}+\frac{\nu_{12}}{n^{2}}\Big(\operatorname{tr}\bA\operatorname{tr}\bB + \operatorname{tr} (\bA\bB')+\operatorname{tr} (\bA\bB)\Big)\nonumber\\
    &-\frac{\nu_{12}}{n^{2}(p-2)} \sum_{\substack{i, k, l \\
i \neq k, k \neq l, i\neq l}} A_{i i} B_{k l}-\frac{\nu_{12}}{n^2(p-2)} \sum_{\substack{i, j, k \\
i \neq j, j \neq k, i\neq k}} A_{i j} B_{k k}\nonumber\\
& -\frac{\nu_2^2}{n^2}\operatorname{tr} \bA\operatorname{tr} \bB+ \frac{\nu_2^2}{n^{2}(p-1)}\left(\sum_{\substack{i, k, l \\
i \neq k, k \neq l, i\neq l}} B_{i i}A_{k l} + \sum_{\substack{i, k, l \\
i \neq k, k \neq l, i\neq l}} A_{i i} B_{k l}\right)+o(n^{-1}) \nonumber\\ 
    =&\frac{\nu_4-3\nu_{12}}{n^{2}} \sum_{i=1}^{p}A_{ii}B_{ii}+\frac{\nu_{12}}{n^{2}}\Big(\operatorname{tr} (\bA\bB')+\operatorname{tr} (\bA\bB)\Big)+\frac{\nu_{12}-\nu_{2}^{2}}{n^{2}} \operatorname{tr} \bA \operatorname{tr} \bB\\
    &+\frac{\nu_2^2-\nu_{12}}{n^{2}(p-1)}\left(\sum_{\substack{i, k, l \\
i \neq k, k \neq l, i\neq l}} B_{i i}A_{k l} + \sum_{\substack{i, k, l \\
i \neq k \neq l}} A_{i i} B_{k l}\right)\nonumber\\
    &- \frac{\nu_{12}}{n^{2}(p-1)(p-2)}\left(\sum_{\substack{i, k, l \\
i \neq k, k \neq l, i\neq l}} B_{i i}A_{k l} + \sum_{\substack{i, k, l \\
i \neq k, k \neq l, i\neq l}} A_{i i} B_{k l}\right)+o(n^{-1})\nonumber\\ %+ K\epsilon n^{-1} + n^{3} \Prob(B^{c}_{p}(\epsilon))\nonumber\\
     =&\frac{\nu_4-3\nu_{12}}{n^{2}} \sum_{i=1}^{p}A_{ii}B_{ii}
    +\frac{\nu_{12}}{n^{2}}\Bigl(\operatorname{tr} (\bA\bB')+\operatorname{tr} (\bA\bB)\Bigr)+\frac{\nu_{12}-\nu_{2}^{2}}{n^{2}}  \operatorname{tr} \bA \operatorname{tr} \bB  +o(n^{-1}).
\end{align*} 
%\begin{align}
 %   \Expe \Big(\br'A\br-\frac{1}{n}\nu_2\operatorname{tr} \bA\Big)\Big(\br'\bB\br-\frac{1}{n}\nu_2\operatorname{tr} \bB\Big)
%    &=\sum_{i=1}^{p}\frac{1}{n^{2}}\Big(\nu_4-3\nu_{12} \Big) A_{ii}B_{ii}\nonumber\\
 %   &+\frac{1}{n^{2}}\nu_{12}\Big(\operatorname{tr} (\bA\bB')+\operatorname{tr} (\bA\bB)\Big)\nonumber\\
  %  &+\frac{1}{n^{2}} \Big(\nu_{12}-\nu_{2}^{2}\Big) \operatorname{tr} \bA \operatorname{tr} \bB \nonumber\\
%  &+\frac{1}{n^{2}(p-1)}(\nu_2^{2}-\nu_{12})\left(\sum_{\substack{i, k, l \\
%i \neq k, k \neq l, i\neq l}} B_{i i}A_{k l} + \sum_{\substack{i, k, l \\
%i \neq k \neq l}} A_{i i} B_{k l}\right)\nonumber\\
%    &+o(n^{-1}) + K\epsilon n^{-1} + n^{3} \Prob(B^{c}_{p}(\epsilon)),
%\end{align}

\subsection{Proof of Lemma \ref{ej-1term0}}\label{prfej-1term0}

For any $p\times p$ matrix $\bA$, we have 
\begin{eqnarray}\label{trdasupl}
\big|\tr\big(\bD^{-1}(z)-\bD_j^{-1}(z)\big)\bA\big|\leq\frac{\|\bA\|}{\Im(z)}.
\end{eqnarray}
By Lemma \ref{burk} and equation  \eqref{trdasupl}, similarly to (4.3) in \cite{Bai1998No}, 
\begin{align}\label{dedsupl}
    \Expe\big|\frac{\nu_2}{n}\tr\bD^{-1}(z)-\Expe\frac{\nu_2}{n}\tr\bD^{-1}(z)\big|^{q}\leq C_qn^{-q/2}v_0^{-q},
\end{align}
which implies that
\begin{eqnarray}\label{betabpboundsupl}
\Expe\big|\bar{\beta}_j(z)-b_p(z)\big|^2\leq \frac{K|z|^4}{nv_0^6}.
\end{eqnarray}

Now, we prove that the difference between 
\[
    \sum^{n}_{j=1}\big[\Expe_{j-1}\bar{\beta}_j(z_1)\varepsilon_j(z_1)][\Expe_{j-1}\bar{\beta}_j(z_2)\varepsilon_j(z_2)\big] 
\]
and
\[
b_p(z_1)b_p(z_2)\sum^{n}_{j=1}\big[\Expe_{j-1}\varepsilon_j(z_1)][\Expe_{j-1}\big(\varepsilon_j(z_2)]
\]
converges to zero. 
We write
\begin{align}
  &\;\Expe\Big|\big[\Expe_{j-1}\bar{\beta}_j(z_1)\varepsilon_j(z_1)][\Expe_{j-1}\bar{\beta}_j(z_2)\varepsilon_j(z_2)\big]\nonumber\\
  &\qquad-b_p(z_1)b_p(z_2)\big[\Expe_{j-1}\varepsilon_j(z_1)][\Expe_{j-1}\big(\varepsilon_j(z_2)]\Big|\nonumber\\
  \leq&\; \Expe\Big|\big[\Expe_{j-1}(\bar{\beta}_j(z_1)\varepsilon_j(z_1)-b_p(z_1)\varepsilon_j(z_1))][\Expe_{j-1}\bar{\beta}_j(z_2)\varepsilon_j(z_2)\big]
  \Big|\nonumber\\  &\qquad+\Expe\Big|\big[\Expe_{j-1}b_p(z_1)\varepsilon_j(z_1)][\Expe_{j-1}\bar{\beta}_j(z_2)\varepsilon_j(z_2)-\Expe_{j-1}(b_p(z_2)\varepsilon_j(z_2))\big]\Big|\nonumber\\
  =:&\; I_1 + I_2.\nonumber
\end{align}
Note that, for any $q_1,q_2\geq 1$ with $1/q_1+1/q_2=1$, we have
\begin{align}
    I_1 &=\Expe\Big|\bigl[\Expe_{j-1}(\bar{\beta}_j(z_1)-b_p(z_1))\varepsilon_j(z_1)\bigr]\bigl[\Expe_{j-1}\bar{\beta}_j(z_2)\varepsilon_j(z_2)\bigr]\Big|\nonumber\\
    &\leq \left(\Expe\Big|\Expe_{j-1}(\bar{\beta}_j(z_1)-b_p(z_1))\varepsilon_j(z_1)\Big|^{q_1}\right)^{1/q_1}\left(\Expe\Big|\Expe_{j-1}\bar{\beta}_j(z_2)\varepsilon_j(z_2)\Big|^{q_2}\right)^{1/q_2}\nonumber\\
    &\leq \left(\Expe\Big|(\bar{\beta}_j(z_1)-b_p(z_1))\varepsilon_j(z_1)\Big|^{q_1}\right)^{1/q_1}\left(\Expe\Big|\bar{\beta}_j(z_2)\varepsilon_j(z_2)\Big|^{q_2}\right)^{1/q_2}\nonumber\\
    &=: I_{11} \times I_{12}.
\end{align}
From Lemma \ref{bound} and \eqref{betabpboundsupl}, for any $1<q_1<2$, $q_2>2$, we have
\begin{align}
    (I_{11})^{q_1} &=\Expe\Big|\big[(\bar{\beta}_j(z_1)-b_p(z_1))\varepsilon_j(z_1)]\Big|^{q_1}\nonumber\\
    &\leq \Bigl(\Expe\bigl|\bar{\beta}_j(z_1)-b_p(z_1)\bigr|^{q_1\frac{2}{q_1}}\Bigr)^{q_1/2}\Bigl(\Expe\bigl|\varepsilon_j(z_1)\bigr|^{q_1\frac{2}{2-q_1}}\Bigr)^{\frac{2-q_1}{2}}\nonumber\\
    &\leq Kn^{-\frac{q_1}{2}}n^{-\frac{2-q_1}{2}}\delta_n^{4q_1-4} = n^{-1}\delta_n^{4q_1-4}, 
\end{align}
and
\[
    I_{12} = \Expe\bigl|\bar{\beta}_j(z_2)\varepsilon_j(z_2)\bigr|^{q_2}\leq K\frac{|z|^{q_2}}{v^{q_2}}n^{-1}\delta_n^{2q_2-4}.
\]
Thus,
\begin{align}
    I_1 \leq Kn^{-1}\delta_n^{\frac{4q_1-4}{q_1}+\frac{2q_2-4}{q_2}} = Kn^{-1}\delta_n^{2}.
\end{align}
Similarly, we can obtain
$I_2 \leq Kn^{-1}\delta_n^{2}.$
Thus, 
\begin{align}\label{yj-1btran}
  &\;\Expe\Big|\big[\Expe_{j-1}\bar{\beta}_j(z_1)\varepsilon_j(z_1)][\Expe_{j-1}\bar{\beta}_j(z_2)\varepsilon_j(z_2)\big]\nonumber\\
  &\qquad-b_p(z_1)b_p(z_2)\big[\Expe_{j-1}\varepsilon_j(z_1)][\Expe_{j-1}\big(\varepsilon_j(z_2)]\Big|\nonumber\\
  \leq&\; Kn^{-1}\delta_n^{2}.
\end{align}
From \eqref{yj-1btran}, we get
\begin{align}
   \sum^{n}_{j=1}\big[\Expe_{j-1}&\bar{\beta}_j(z_1)\varepsilon_j(z_1)][\Expe_{j-1}\bar{\beta}_j(z_2)\varepsilon_j(z_2)\big]\nonumber\\
   &-b_p(z_1)b_p(z_2)\sum^{n}_{j=1}\big[\Expe_{j-1}\varepsilon_j(z_1)][\Expe_{j-1}\big(\varepsilon_j(z_2)] \stackrel{i.p.}{\rightarrow}0.\label{eq:diff-betaeps-beps}
\end{align}

Next, we prove that 
\begin{eqnarray}
b_p(z_1)b_p(z_2)\sum^{n}_{j=1}\big[\Expe_{j-1}\varepsilon_j(z_1)][\Expe_{j-1}\big(\varepsilon_j(z_2)] \stackrel{i.p.}{\rightarrow}0.
\end{eqnarray}
By \eqref{covmatrix}, we get
\begin{align}
\Expe_{j-1}\varepsilon_j(z) &= \Expe_{j-1}\left(\br_j'\bD_j^{-1}(z)\br_j-\frac{1}{n}\nu_2tr\bD_j^{-1}(z)\right)\nonumber\\
% &= \Expe_{j-1}\left(\operatorname{tr}[\bD_j^{-1}(z)\br_j\br_j']-\frac{1}{n}\nu_2tr\bD_j^{-1}(z)\right)\nonumber\\
% &= \operatorname{tr}\left[\Expe_{j-1}\bD_j^{-1}(z)\Expe_{j-1}\br_j\br_j'\right]-\Expe_{j-1}\left(\frac{1}{n}\nu_2tr\bD_j^{-1}(z)\right)\nonumber\\
&= \operatorname{tr}\bigl[\Expe_{j-1}\bD_j^{-1}(z)\Expe\br_j\br_j'\bigr]-\frac{1}{n}\nu_2\Expe_{j-1}\tr\bD_j^{-1}(z)\nonumber\\
&= \frac{1}{n}\nu_2 \operatorname{tr}\Bigl[\Expe_{j-1}\bD_j^{-1}(z) \Bigl(-\frac{1}{p-1}\boldsymbol{1}_{p}\boldsymbol{1}'_{p} + \frac{1}{p-1}\bI_p + \bI_p \Bigr)\Bigr]-\frac{1}{n}\nu_2\Expe_{j-1}\bigl(\tr\bD_j^{-1}(z)\bigr)\nonumber\\
% &=  -\frac{1}{n(p-1)}\nu_2\operatorname{tr}\left[\Expe_{j-1}\bD_j^{-1}(z)\boldsymbol{1}_{p}\boldsymbol{1}'_{p}\right]
% +\frac{1}{n(p-1)}\nu_2\operatorname{tr}\left[\Expe_{j-1}\bD_j^{-1}(z)\right]\nonumber\\
% &\qquad+\frac{1}{n}\nu_2 \operatorname{tr}\left[\Expe_{j-1}\bD_j^{-1}(z)\right] -\Expe_{j-1}\left(\frac{1}{n}\nu_2tr\bD_j^{-1}(z)\right)\nonumber\\
&=  -\frac{1}{n(p-1)}\nu_2\Expe_{j-1}\left[\boldsymbol{1}'_{p}\bD_j^{-1}(z)\boldsymbol{1}_{p}\right]
+\frac{1}{n(p-1)}\nu_2\Expe_{j-1}\bigl[\operatorname{tr}\bD_j^{-1}(z)\bigr].\nonumber
%&\qquad+\frac{1}{n}\textcolor{red}{(\sigma_{p}^{2} - \frac{\sigma^2}{\mu^2})}\Expe_{j-1}\operatorname{tr}[\bD_j^{-1}(z)], \nonumber
\end{align}
Thus,
\begin{align}\label{epjepj}
&\;\sum^{n}_{j=1}\bigl[\Expe_{j-1}\varepsilon_j(z_1)\bigr]\bigl[\Expe_{j-1}\varepsilon_j(z_2)\bigr]\nonumber\\
=&\;\frac{\nu_2^2}{n^2(p-1)^2}\sum^{n}_{j=1} \Expe_{j-1}\Bigl(\operatorname{tr}\bD_j^{-1}(z_1)-\boldsymbol{1}'_{p}\bD_j^{-1}(z_1)\boldsymbol{1}_{p}\Bigr)\nonumber\\
&\qquad\qquad\qquad\qquad \times \Expe_{j-1}\Bigl(\operatorname{tr}\bD_j^{-1}(z_2)-\boldsymbol{1}'_{p}\bD_j^{-1}(z_2)\boldsymbol{1}_{p}\Bigr).
\end{align}
By Lemma 2.3 in \cite{bai2004clt}, we have
\begin{align}\label{trdidilimsupl}
   \frac{1}{p}\operatorname{tr}\left[\bD_j^{-1}(z)\right]\stackrel{L^2}{\rightarrow} m(z),\ as\ p \rightarrow \infty,
\end{align}
and by Lemma \ref{xdxz}, 
\begin{align}\label{epdepsupl}
    \frac{1}{p}\boldsymbol{1}'_{p}\bD_j^{-1}(z)\boldsymbol{1}_{p} \stackrel{L^2}{\rightarrow} 
    -\frac{1}{z}, \ as\ p \rightarrow \infty.
\end{align}
By \eqref{dedsupl}, we have 
\begin{eqnarray}\label{bpbetabounsupl}
|b_p(z)-\Expe\beta_1(z)|\leq Kn^{-1/2}.
\end{eqnarray}
From the formula (2.2) of \cite{Silverstein1995StrongCO}, 
$
\underline{m}_p(z)=-\frac{1}{zn}\sum^{n}_{j=1}\beta_j(z),
$
we have
\begin{eqnarray}\label{betazmsupl}
\Expe\beta_1(z)=-z\Expe\underline{m}_p(z).
\end{eqnarray}
From \eqref{bpbetabounsupl}, \eqref{betazmsupl}, and Lemma \ref{bound}, we have
\begin{eqnarray}\label{b1zmsupl}
|b_p(z)+z\underline{m}_{p}^{0}(z)|\leq Kn^{-1/2}.
\end{eqnarray}
From \eqref{epjepj} -- \eqref{b1zmsupl},
%\begin{align}
 %   |b_1(z)+z\underline{m}_{p}^{0}(z)|\leq Kn^{-1/2},
%\end{align}
we get
\[
b_p(z_1)b_p(z_2)\sum^{n}_{j=1}\big[\Expe_{j-1}\varepsilon_j(z_1)][\Expe_{j-1}\big(\varepsilon_j(z_2)] \stackrel{i.p.}{\rightarrow}0.
\]
This, together with \eqref{eq:diff-betaeps-beps}, completes the proof of Lemma \ref{ej-1term0}.

\subsection{Proof of Lemma \ref{trejdjdjeqmain}}\label{prftrejdjdjeqmain}

Using Lemma \ref{bound} and equation \eqref{betabpboundsupl}, we have, for $i<j$,
\begin{align}\label{betabdd}
&\;\Expe\Bigl|\beta_{ij}(z_2)\br_i'\Expe_j\big(\bD_{ij}^{-1}(z_1)\big)\bD_{ij}^{-1}(z_2)\br_i\br_i'\bD_{ij}^{-1}(z_2)\bQ_p(z_1)\br_i\nonumber\\
&\qquad -n^{-2}\nu_2^2b_1(z_2)\operatorname{tr}\Bigl(\Expe_j\big(\bD_{ij}^{-1}(z_1)\big)\bD_{ij}^{-1}(z_2)\Bigr)\operatorname{tr}\Bigl(\bD_{ij}^{-1}(z_2)\bQ_p(z_1)\Bigr)\Bigr|\nonumber\\
\leq&\;\Expe\Bigl|\beta_{ij}(z_2)\br_i'\Expe_j\big(\bD_{ij}^{-1}(z_1)\big)\bD_{ij}^{-1}(z_2)\br_i\br_i'\bD_{ij}^{-1}(z_2)\bQ_p(z_1)\br_i\nonumber\\
&\qquad -n^{-2}\nu_2^2\beta_{ij}(z_2)\operatorname{tr}\Bigl(\Expe_j\bigl(\bD_{ij}^{-1}(z_1)\bigr)\bD_{ij}^{-1}(z_2)\Bigr)\operatorname{tr}\Bigl(\bD_{ij}^{-1}(z_2)\bQ_p(z_1)\Bigr)\Bigr|\nonumber\\
&\quad +\Expe\Bigl|n^{-2}\nu_2^2\bigl(\beta_{ij}(z_2)-b_1(z_2)\bigr)\operatorname{tr}\Bigl(\Expe_j\bigl(\bD_{ij}^{-1}(z_1)\bigr)\bD_{ij}^{-1}(z_2)\Bigr)\operatorname{tr}\Bigl(\bD_{ij}^{-1}(z_2)\bQ_p(z_1)\Bigr)\Bigr|\nonumber\\
\leq&\;\Expe\Bigl|\beta_{ij}(z_2)\Bigl[\br_i'\Expe_j\big(\bD_{ij}^{-1}(z_1)\big)\bD_{ij}^{-1}(z_2)\br_i-\frac{1}{n}\nu_2\operatorname{tr}\Bigl(\Expe_j\big(\bD_{ij}^{-1}(z_1)\big)\bD_{ij}^{-1}(z_2)\Bigr)\Bigr]\br_i'\bD_{ij}^{-1}(z_2)\bQ_p(z_1)\br_i\nonumber\\
&\qquad +\frac{1}{n}\nu_2\beta_{ij}(z_2)\operatorname{tr}\Bigl(\Expe_j\bigl(\bD_{ij}^{-1}(z_1)\bigr)\bD_{ij}^{-1}(z_2)\Bigr)\Bigl[\br_i'\bD_{ij}^{-1}(z_2)\bQ_p(z_1)\br_i-\frac{1}{n}\nu_2\operatorname{tr}\Bigl(\bD_{ij}^{-1}(z_2)\bQ_p(z_1)\Bigr)\Bigr]\Bigr|\nonumber\\
&\quad +\Expe\Bigl|n^{-2}\nu_2^2\bigl(\beta_{ij}(z_2)-b_1(z_2)\bigr)\operatorname{tr}\Bigl(\Expe_j\big(\bD_{ij}^{-1}(z_1)\big)\bD_{ij}^{-1}(z_2)\Bigr)\operatorname{tr}\Bigl(\bD_{ij}^{-1}(z_2)\bQ_p(z_1)\Bigr)\Bigr|\nonumber\\
%&\ \ \ \ \ -n^{-2}\nu_2^2b_1(z_2)tr\big(\Expe_j\big(\bD_{ij}^{-1}(z_1)\big)\bD_{ij}^{-1}(z_2))\operatorname{tr}(\bD_{ij}^{-1}(z_2)\bQ_p(z_1))\big|\nonumber\\
\leq&\; Kn^{-1/2}.
\end{align}
By \eqref{trdasupl}, we have
\begin{eqnarray}\label{dijdj}
&&\Bigl|\operatorname{tr}\Bigl(\Expe_j\bigl(\bD_{ij}^{-1}(z_1)\bigr)\bD_{ij}^{-1}(z_2)\Bigr)\tr\Bigl(\bD_{ij}^{-1}(z_2)\bQ_p(z_1)\Bigr)\nonumber\\
&&\qquad -\operatorname{tr}\Bigl(\Expe_j\bigl(\bD_{j}^{-1}(z_1)\bigr)\bD_{j}^{-1}(z_2)\Bigr)\tr\Bigl(\bD_{j}^{-1}(z_2)\bQ_p(z_1)\Bigr)\Bigr|\leq Kn.
\end{eqnarray}
It follows from \eqref{betabdd} and \eqref{dijdj} that
\begin{align}\label{A1sim}
\Expe\Bigl|A_1(z_1,z_2)+\frac{j-1}{n^2}\nu_2^2b_1(z_2)\operatorname{tr}\Bigl(\Expe_j\big(\bD_j^{-1}(z_1)\big)\bD_j^{-1}(z_2)\Bigr)\operatorname{tr}\Bigl(\bD_j^{-1}(z_2)\bQ_p(z_1)\Bigr)\Bigr|\leq Kn^{1/2},
\end{align}
where $A_1(z_1,z_2)$ is defined in \eqref{eq:A1z1z2}.
Note that 
\begin{align}
\operatorname{tr}\bigl(\Expe_j\big(\bD_j^{-1}(z_1)\big)\bD_j^{-1}(z_2)\bigr)&=\operatorname{tr}\Big(\Expe_j\big[-\bQ_p(z_1)+b_1(z)\bA(z_1)\big]\bD_j^{-1}(z_2)\Big) + O(n^{1/2})\nonumber\\
    &= -\operatorname{tr}\big(\bQ_p(z_1)\bD_j^{-1}(z_2)\big)+b_1(z_1)\tr\Expe_j\big(\bA(z_1)\big)\bD_j^{-1}(z_2)+O(n^{1/2})\nonumber\\
    &=-\operatorname{tr}\big(\bQ_p(z_1)\bD_j^{-1}(z_2)\big)+b_1(z_1)A_1(z_1,z_2)+O(n^{1/2}).\nonumber
\end{align}
Therefore, from \eqref{djdecomp}--\eqref{trejadjdecomp} and \eqref{betabdd}--\eqref{A1sim}, we can write
\begin{align*}
&\;\operatorname{tr}\big(\Expe_j\big(\bD_j^{-1}(z_1)\big)\bD_j^{-1}(z_2)\big)
\Bigl[1+\frac{j-1}{n^2}\nu_2^2b_1(z_1)b_1(z_2)\operatorname{tr}\big(\bD_j^{-1}(z_2)\bQ_p(z_1)\big)\Bigr]\nonumber\\
=&\;-\operatorname{tr}\big(\bQ_p(z_1)\bD_j^{-1}(z_2)\big)+A_4(z_1,z_2),
\end{align*}
where 
$
    \Expe|A_4(z_1,z_2)|\leq Kn^{1/2}.
$
Using the expression for $\bD_j^{-1}(z_2)$ in \eqref{djdecomp}, \eqref{trda}, and \eqref{cmbound}-\eqref{Rterm}, we have
\begin{align}\label{trdjdj}
&\;\tr\big(\Expe_j\big(\bD_j^{-1}(z_1)\big)\bD_j^{-1}(z_2)\big)
\Bigl[1-\frac{j-1}{n^2}\nu_2^2b_1(z_1)b_1(z_2)\tr\big(\bQ_p(z_2)\bQ_p(z_1)\big)\Bigr]\nonumber\\
=&\;\tr\Bigl(\bQ_p(z_2)\bQ_p(z_1)\Bigr)+A_5(z_1,z_2),
\end{align}
where 
$
    \Expe|A_5(z_1,z_2)|\leq Kn^{1/2}.
$
By \eqref{b1zmsupl} and \eqref{trdjdj}, we get
\begin{align}\label{trdjdjeq}
&\;\operatorname{tr}\Bigl(\Expe_j\big(\bD_j^{-1}(z_1)\big)\bD_j^{-1}(z_2)\Bigr)
\biggl\{1-\frac{j-1}{n^2}\nu_2^2\underline{m}_{p}^{0}(z_1)\underline{m}_{p}^{0}(z_2)\nonumber\\
&\qquad\qquad \times\operatorname{tr}\biggl[\Bigl(\bI_p+\frac{n-1}{n}\nu_2\underline{m}_p^0(z_1)\Bigr)^{-1}\Bigl(\bI_p+\frac{n-1}{n}\nu_2\underline{m}_p^0(z_2)\Bigr)^{-1}\biggr]\biggr\}\nonumber\\
=&\;\frac{1}{z_1z_2}\operatorname{tr}\left[\Bigl(\bI_p+\frac{n-1}{n}\nu_2\underline{m}_p^0(z_1)\Bigr)^{-1}\Bigl(\bI_p+\frac{n-1}{n}\nu_2\underline{m}_p^0(z_2)\Bigr)^{-1}\right]+A_6(z_1,z_2),
\end{align}
where $
    \Expe|A_6(z_1,z_2)|\leq Kn^{1/2}.
$
Thus, we rewrite \eqref{trdjdjeq} as \eqref{eq:trD1D2_limit} in Lemma \ref{trejdjdjeqmain}.
% \begin{eqnarray}\label{trdjdjeq2supl}
% &&\operatorname{tr}\big(\Expe_j\big(\bD_j^{-1}(z_1)\big)\bD_j^{-1}(z_2)\big)\nonumber\\
% &&\qquad\times\big[1-\frac{j-1}{n}\nu_2^2\underline{m}_{p}^{0}(z_1)\underline{m}_{p}^{0}(z_2)\frac{c_n}{(1+\frac{n-1}{n}\nu_2\underline{m}_p^0(z_1))(1+\frac{n-1}{n}\nu_2\underline{m}_p^0(z_2))}\big]\nonumber\\
% &&=\frac{nc_n}{z_1z_2}\frac{1}{(1+\frac{n-1}{n}\nu_2\underline{m}_p^0(z_1))(1+\frac{n-1}{n}\nu_2\underline{m}_p^0(z_2))}+A_6(z_1,z_2).\nonumber
% \end{eqnarray}

\subsection{Proof of Lemma \ref{m1I1ejdj}}\label{prfm1I1ejdj}

By Burkholder's inequality, the inequality $|\beta_{ij}(z)|\leq \frac{|z|}{v_0}$, and Lemma \ref{bound}, we get
\begin{align}
    &\;\Expe\bigl|(\boldsymbol{1}_{p}^{i})'\bigl(\bD_{1}^{-1}(z_1)-\Expe\bD_{1}^{-1}(z_1)\bigr)\boldsymbol{1}_{p}^{i}\bigr|^2\nonumber\\
 &=\Expe\biggl|\sum_{j=1}^{n}\Bigl(\Expe_j(\boldsymbol{1}_{p}^{i})'\bD_{1}^{-1}(z_1)\boldsymbol{1}_{p}^{i}-\Expe_{j-1}(\boldsymbol{1}_{p}^{i})'\bD_{1}^{-1}(z_1)\boldsymbol{1}_{p}^{i}\Bigr)\biggr|^2\nonumber\\
 &=\Expe\biggl|\sum_{j=1}^{n}\Bigl[\Expe_j(\boldsymbol{1}_{p}^{i})'\Bigl(\bD_{1}^{-1}(z_1)-\bD_{1j}^{-1}(z_1)\Bigr)\boldsymbol{1}_{p}^{i}-\Expe_{j-1}(\boldsymbol{1}_{p}^{i})'\Bigl(\bD_{1}^{-1}(z_1)-\bD_{1j}^{-1}(z_1)\Bigr)\boldsymbol{1}_{p}^{i}\Bigr]\biggr|^2\nonumber\\
 &=\Expe\biggl|\sum_{j=1}^{n}(\Expe_j-\Expe_{j-1})(\boldsymbol{1}_{p}^{i})'\Bigl(\bD_{1}^{-1}(z_1)-\bD_{1j}^{-1}(z_1)\Bigr)\boldsymbol{1}_{p}^{i}\biggr|^2\nonumber\\
    &\leq K\sum_{j=1}^{n} \Expe\Bigl|(\Expe_j-\Expe_{j-1})\beta_{1j}(z_1)\br_j'\bD_{1j}^{-1}\boldsymbol{1}_{p}^{i}(\boldsymbol{1}_{p}^{i})'\bD_{1j}^{-1}\br_j\Bigr|^2\nonumber\\
    &\leq K\sum_{j=1}^{n} \Expe\bigl|\beta_{1j}(z_1)\br_j'\bD_{1j}^{-1}\boldsymbol{1}_{p}^{i}(\boldsymbol{1}_{p}^{i})'\bD_{1j}^{-1}\br_j\bigr|^2\nonumber\\
    &\leq Kn^{-1},
\end{align}
where $\boldsymbol{1}_{p}^{i}$ is the $p$-dimensional vector with all components being $0$ except for the $i$-th component being $1$.
Hence, we have
\begin{align}
&\;\Expe\biggl|\frac{1}{n^{2}}\sum_{j=1}^{n} \sum_{i=1}^{p} \Expe_j\bigl(\bD_j^{-1}(z_1)-\Expe \bD_j^{-1}(z_1)\bigr)_{i i}\Expe_j\bigl(\bD_j^{-1}(z_2)\bigr)_{i i}\biggr|
\nonumber\\
\leq&\;\frac{nK}{n^2 v_0} \sum^{p}_{i=1}\Expe\Bigl|(\boldsymbol{1}_{p}^{i})'(\bD_{1}^{-1}(z_1)-\Expe\bD_{1}^{-1}(z_1))\boldsymbol{1}_{p}^{i}\Bigr| \leq K n^{-1/2},\nonumber
\end{align}
and thus
\begin{eqnarray}\label{ejdjedj}
\frac{1}{n^{2}}\sum_{j=1}^{n} \sum_{i=1}^{p} \Expe_j\bigl(\bD_j^{-1}(z_1)-\Expe \bD_j^{-1}(z_1)\bigr)_{i i}\Expe_j\bigl(\bD_j^{-1}(z_2)\bigr)_{i i}=O_P(n^{-1/2}).
\end{eqnarray}
With \eqref{ejdjedj}, it remains to find the limit of
\begin{eqnarray}\label{I1tr}
\frac{1}{p}\sum^{p}_{i=1}\Expe\big(\bD_j^{-1}(z_1)\big)_{ii}\Expe\big(\bD_j^{-1}(z_2)\big)_{ii}.
\end{eqnarray}
It is easy to see that the sum of expectations in \eqref{I1tr} is exactly the same for any $j$. Moreover, we have
\[
\frac{1}{p}\sum^{p}_{i=1}\Expe\big(\bD_j^{-1}(z_1)\big)_{ii}\Expe\big(\bD_j^{-1}(z_2)\big)_{ii}\stackrel{i.p.}{\rightarrow}m(z_1)m(z_2).
\]
This completes the proof of Lemma \ref{m1I1ejdj}.

\subsection{Proof of Lemma \ref{estbeta1gamma1p1p2}}\label{prfeqs}

First, we prove \eqref{gammp1}. For $\Expe(\gamma_1(z)P_1(z))$, we have
\begin{align}
    &\;\Expe(\gamma_1(z)P_1(z)) \nonumber\\  =&\;\Expe\Bigl[\gamma_1(z)\Bigl(n\br_1' \bD_1^{-1}(z)\widetilde{\bQ}_p^{-1}(z)\br_1 -\frac{\sigma^2}{\mu^2}\operatorname{tr}\left(\widetilde{\bQ}_p^{-1}(z)\Expe\bD_1^{-1}(z)\right)\Bigr)\Bigr]\nonumber\\
    =&\; n\Expe\Bigl[\gamma_1(z)\br_1'\bD_1^{-1}(z)\widetilde{\bQ}_p^{-1}(z)\br_1\Bigr]  - \frac{\sigma^2/\mu^2}{1+\frac{\sigma^2}{\mu^2}\Expe\underline{m}_p(z)}\operatorname{tr}[\Expe\bD_1^{-1}(z)]\Expe\gamma_1(z)\nonumber\nonumber\\
    =&\; n\Expe\Bigl[\Bigl(  \br_1'\bD_1^{-1}(z)\br_1-\frac{1}{n}\nu_2\operatorname{tr}\bD_1^{-1}(z)+\frac{1}{n}\nu_2\operatorname{tr}\bD_1^{-1}(z)-\frac{1}{n}\nu_2\Expe\operatorname{tr}\bD_1^{-1}(z)   \Bigr)\nonumber\\
 &\qquad\times\Bigl(\br_1'\bD_1^{-1}(z)\widetilde{\bQ}_p^{-1}(z)\br_1-\frac{1}{n}\nu_2\operatorname{tr}\left[\bD_1^{-1}(z)\widetilde{\bQ}_p^{-1}(z)\right]+\frac{1}{n}\nu_2\operatorname{tr}\left[\bD_1^{-1}(z)\widetilde{\bQ}_p^{-1}(z)\right]\Bigr)\Bigr] \nonumber\\
    &\quad- \frac{\sigma^2/\mu^2}{1+\frac{\sigma^2}{\mu^2}\Expe\underline{m}_p(z)}\operatorname{tr}[\Expe\bD_1^{-1}]\Expe\gamma_1(z)\nonumber\nonumber\\
     =&\;n\Expe\Bigl[\Bigl(  \br_1'\bD_1^{-1}(z)\br_1-\frac{1}{n}\nu_2\operatorname{tr}\bD_1^{-1}(z)   \Bigr)\Bigl(\br_1'\bD_1^{-1}(z)\widetilde{\bQ}_p^{-1}(z)\br_1-\frac{1}{n}\nu_2\operatorname{tr}\left[\bD_1^{-1}(z)\widetilde{\bQ}_p^{-1}(z)\right]\Bigr)\Bigr] \nonumber\\
    &\quad+\frac{1}{n(p-1)}\nu_2^2\Expe\left(\operatorname{tr}\bD_1^{-1}(z)\operatorname{tr}[\bD_1^{-1}(z)\widetilde{\bQ}_p^{-1}(z)]\right)-\frac{1}{n(p-1)}\nu_2^2\Expe\left(\boldsymbol{1}'_{p}\bD_1^{-1}(z)\boldsymbol{1}_{p}\operatorname{tr}[\bD_1^{-1}(z)\widetilde{\bQ}_p^{-1}(z)]\right)\nonumber\\
    &\quad+\frac{1}{n}(1+\frac{1}{p-1})\nu_2^2\Cov\left(\operatorname{tr}\bD_1^{-1}(z),\operatorname{tr}[\bD_1^{-1}(z)\widetilde{\bQ}_p^{-1}(z)]\right)\nonumber\\
    &\quad-\frac{1}{n(p-1)}\nu_2^2\Big[\Expe\left(\operatorname{tr}\bD_1^{-1}(z)[\boldsymbol{1}'_{p}\bD_1^{-1}(z)\widetilde{\bQ}_p^{-1}(z)\boldsymbol{1}_{p}]\right)  -\Expe\left(\operatorname{tr}\bD_1^{-1}(z)\right)\Expe\left(\boldsymbol{1}'_{p}\bD_1^{-1}(z)\widetilde{\bQ}_p^{-1}(z)\boldsymbol{1}_{p}\right)\Big]\nonumber\\
    &\quad- \frac{\sigma^2/\mu^2}{1+\frac{\sigma^2}{\mu^2}\Expe\underline{m}_p(z)}\operatorname{tr}[\Expe\bD_1^{-1}(z)]\Expe\gamma_1(z)\nonumber\\
        =&\;n\Expe\Bigl[\Bigl(  \br_1'\bD_1^{-1}(z)\br_1-\frac{1}{n}\nu_2\operatorname{tr}\bD_1^{-1}(z)   \Bigr)\Bigl(\br_1'\bD_1^{-1}(z)\widetilde{\bQ}_p^{-1}(z)\br_1-\frac{1}{n}\nu_2\operatorname{tr}\left[\bD_1^{-1}(z)\widetilde{\bQ}_p^{-1}(z)\right]\Bigr)\Bigr] \nonumber\\
    &\quad+\frac{1}{n(p-1)}\nu_2^2\Expe\left(\operatorname{tr}\bD_1^{-1}(z)\operatorname{tr}[\bD_1^{-1}(z)\widetilde{\bQ}_p^{-1}(z)]\right)-\frac{1}{n(p-1)}\nu_2^2\Expe\left(\boldsymbol{1}'_{p}\bD_1^{-1}(z)\boldsymbol{1}_{p}\operatorname{tr}[\bD_1^{-1}(z)\widetilde{\bQ}_p^{-1}(z)]\right)\nonumber\\
    &\quad- \frac{\sigma^2/\mu^2}{1+\frac{\sigma^2}{\mu^2}\Expe\underline{m}_p(z)}\operatorname{tr}[\Expe\bD_1^{-1}(z)]\Expe\gamma_1(z)+o(1),\nonumber
\end{align}
which is the equation \eqref{gammp1}.
Below are some interpretations of the above equalities:
\begin{enumerate}
    \item The fourth equality uses the following derivation: By \eqref{covmatrix}, i.e., $\Expe (\br\br') = \frac{1}{n}\nu_{2} \bigl(-\frac{1}{p-1}\boldsymbol{1}_{p}\boldsymbol{1}'_{p} + \frac{1}{p-1}\bI_p + \bI_p \bigr)$, we get
\begin{align}
     &\ \ \ \  n\Expe\Bigl[\Bigl(  \br_1'\bD_1^{-1}(z)\br_1-\frac{1}{n}\nu_2\operatorname{tr}\bD_1^{-1}(z)\Bigr)\Bigl(\frac{1}{n}\nu_2\operatorname{tr}\left[\bD_1^{-1}(z)\widetilde{\bQ}_p^{-1}(z)\right]\Bigr)\Bigr]\nonumber\\
 &\ \ \ \ \ \ + n\Expe\Bigl[\Bigl(\frac{1}{n}\nu_2\operatorname{tr}\bD_1^{-1}(z)-\frac{1}{n}\nu_2\Expe\operatorname{tr}\bD_1^{-1}(z)   \Bigr)\Bigl(\br_1'\bD_1^{-1}(z)\widetilde{\bQ}_p^{-1}(z)\br_1\Bigr)\Bigr] \nonumber\\
  &=\nu_2\Expe\Bigl( \br_1'\bD_1^{-1}(z)\br_1\operatorname{tr}\left[\bD_1^{-1}(z)\widetilde{\bQ}_p^{-1}(z)\right]\Bigr)-\frac{1}{n}\nu_2^2\Expe\Bigl( \operatorname{tr}\bD_1^{-1}(z)\operatorname{tr}\left[\bD_1^{-1}(z)\widetilde{\bQ}_p^{-1}(z)\right]\Bigr)\nonumber\\
 &\ \ \ \ \ \ + \nu_2\Expe\Bigl(\operatorname{tr}\bD_1^{-1}(z) \br_1'\bD_1^{-1}(z)\widetilde{\bQ}_p^{-1}(z)\br_1\Bigr) -\nu_2\Expe\operatorname{tr}\bD_1^{-1}(z)   \Expe\Bigl(\br_1'\bD_1^{-1}(z)\widetilde{\bQ}_p^{-1}(z)\br_1\Bigr) \nonumber\\
  &=\nu_2\operatorname{tr}\Bigl\{\Expe\br_1\br_1'\Expe\Bigl[\bD_1^{-1}(z)\operatorname{tr}\left(\bD_1^{-1}(z)\widetilde{\bQ}_p^{-1}(z)\right)\Bigr]\Bigr\}-\frac{1}{n}\nu_2^2\Expe\Bigl( \operatorname{tr}\bD_1^{-1}(z)\operatorname{tr}\left[\bD_1^{-1}(z)\widetilde{\bQ}_p^{-1}(z)\right]\Bigr)\nonumber\\
 &\ \ \ \ \ \ + \nu_2\operatorname{tr}\Bigl[\Expe(\bD_1^{-1}(z)\widetilde{\bQ}_p^{-1}(z)\operatorname{tr}\bD_1^{-1}(z)) \Expe\br_1\br_1'\Bigr] -\nu_2\Expe\operatorname{tr}\bD_1^{-1}(z)   \operatorname{tr}\Bigl[\Expe\Bigl(\bD_1^{-1}(z)\widetilde{\bQ}_p^{-1}(z)\Bigr)\Expe\br_1\br_1'\Bigr] \nonumber\\
  &=\frac{1}{n}\nu_2^2\Expe\Bigl( \operatorname{tr}\bD_1^{-1}(z)\operatorname{tr}\left[\bD_1^{-1}(z)\widetilde{\bQ}_p^{-1}(z)\right]\Bigr) +\frac{1}{n(p-1)}\nu_2^2\Expe\left(\operatorname{tr}\bD_1^{-1}(z)\operatorname{tr}[\bD_1^{-1}(z)\widetilde{\bQ}_p^{-1}(z)]\right)\nonumber\\
   &\quad-\frac{1}{n(p-1)}\nu_2^2\Expe\left(\boldsymbol{1}'_{p}\bD_1^{-1}(z)\boldsymbol{1}_{p}\operatorname{tr}[\bD_1^{-1}(z)\widetilde{\bQ}_p^{-1}(z)]\right) -\frac{1}{n}\nu_2^2\Expe\Bigl( \operatorname{tr}\bD_1^{-1}(z)\operatorname{tr}\left[\bD_1^{-1}(z)\widetilde{\bQ}_p^{-1}(z)\right]\Bigr)\nonumber\\
    &\quad+\frac{1}{n}(1+\frac{1}{p-1})\nu_2^2\Cov\left(\operatorname{tr}\bD_1^{-1}(z),\operatorname{tr}[\bD_1^{-1}(z)\widetilde{\bQ}_p^{-1}(z)]\right)\nonumber\\
    &\quad-\frac{1}{n(p-1)}\nu_2^2\Big[\Expe\left(\operatorname{tr}\bD_1^{-1}(z)[\boldsymbol{1}'_{p}\bD_1^{-1}(z)\widetilde{\bQ}_p^{-1}(z)\boldsymbol{1}_{p}]\right)  -\Expe\left(\operatorname{tr}\bD_1^{-1}(z)\right)\Expe\left(\boldsymbol{1}'_{p}\bD_1^{-1}(z)\widetilde{\bQ}_p^{-1}(z)\boldsymbol{1}_{p}\right)\Big]\nonumber\\
      &=\frac{\nu_2^2}{n(p-1)}\Expe\left(\operatorname{tr}\bD_1^{-1}(z)\operatorname{tr}[\bD_1^{-1}(z)\widetilde{\bQ}_p^{-1}(z)]\right)-\frac{\nu_2^2}{n(p-1)}\Expe\left(\boldsymbol{1}'_{p}\bD_1^{-1}(z)\boldsymbol{1}_{p}\operatorname{tr}[\bD_1^{-1}(z)\widetilde{\bQ}_p^{-1}(z)]\right)\nonumber\\
    &\quad+\frac{1}{n}(1+\frac{1}{p-1})\nu_2^2\Cov\left(\operatorname{tr}\bD_1^{-1}(z),\operatorname{tr}[\bD_1^{-1}(z)\widetilde{\bQ}_p^{-1}(z)]\right)\nonumber\\
    &\quad-\frac{\nu_2^2}{n(p-1)}\Big[\Expe\left(\operatorname{tr}\bD_1^{-1}(z)[\boldsymbol{1}'_{p}\bD_1^{-1}(z)\widetilde{\bQ}_p^{-1}(z)\boldsymbol{1}_{p}]\right)-\Expe\left(\operatorname{tr}\bD_1^{-1}(z)\right)\Expe\left(\boldsymbol{1}'_{p}\bD_1^{-1}(z)\widetilde{\bQ}_p^{-1}(z)\boldsymbol{1}_{p}\right)\Big].\nonumber
\end{align}

\item The last equality is due to
\[
    \frac{1}{n}\Cov\bigl(\operatorname{tr}\bD_1^{-1}(z),\, \operatorname{tr}\bigl[\bD_1^{-1}(z)\widetilde{\bQ}_p^{-1}(z)\bigr]\bigr) = O(n^{-1}),
\]
which follows from \eqref{dmedm}.
\end{enumerate}

Next, we prove \eqref{begap1}.
For $\Expe(\beta_1(z)\gamma_1^{2}(z)P_1(z))$, we have
\begin{align}
   &\;\Expe\bigl(\beta_1(z)\gamma_1^{2}(z)P_1(z)\bigr)\nonumber\\
   % =&\; \Expe\Bigl[\beta_1(z)\gamma_1^{2}(z)\Bigl(n\br_1' \bD_1^{-1}(z)\widetilde{\bQ}_p^{-1}(z)\br_1 -\frac{\sigma^2}{\mu^2}\operatorname{tr}\left[\widetilde{\bQ}_p^{-1}(z)\Expe\bD_1^{-1}(z)\right]\Bigr)\Bigr]\nonumber\\
   =&\;  \Expe\Bigl(n\beta_1(z)\gamma_1^{2}(z)\br_1' \bD_1^{-1}(z)\widetilde{\bQ}_p^{-1}(z)\br_1 \Bigr)-\Expe\left(\beta_1(z)\gamma_1^{2}(z)\right)\Expe\left(\frac{\sigma^2}{\mu^2}\operatorname{tr}\left[\widetilde{\bQ}_p^{-1}(z)\bD_1^{-1}(z)\right]\right)\nonumber\\
    =&\;  \Expe\Bigl(n\beta_1(z)\gamma_1^{2}(z)\br_1' \bD_1^{-1}(z)\widetilde{\bQ}_p^{-1}(z)\br_1 \Bigr)-\Expe\left(\beta_1(z)\gamma_1^{2}(z)\operatorname{tr}\left[\frac{\sigma^2}{\mu^2}\widetilde{\bQ}_p^{-1}(z)\bD_1^{-1}(z)\right]\right)\nonumber\\
   &\qquad+\Cov\left(\beta_1(z)\gamma_1^{2}(z),\operatorname{tr}\left[\frac{\sigma^2}{\mu^2}\widetilde{\bQ}_p^{-1}(z)\bD_1^{-1}(z)\right]\right).\label{begap1-eq1}
\end{align}
From Lemma \ref{bound} and equation \eqref{dmedm}, we have
\begin{align}
    &\;\Expe\Bigl(n\beta_1(z)\gamma_1^{2}(z)\br_1' \bD_1^{-1}(z)\widetilde{\bQ}_p^{-1}(z)\br_1 \Bigr)-\Expe\left(\beta_1(z)\gamma_1^{2}(z)\operatorname{tr}\left[\frac{\sigma^2}{\mu^2}\widetilde{\bQ}_p^{-1}(z)\bD_1^{-1}(z)\right]\right)\nonumber\\
    \leq&\; n\Bigl[\Expe|\gamma_1^{2}(z)\beta_1(z)|^2\Bigr]^{1/2}\biggl[\Expe\Bigl|\br_1' \bD_1^{-1}(z)\widetilde{\bQ}_p^{-1}(z)\br_1 -\frac{1}{n}\frac{\sigma^2}{\mu^2}\operatorname{tr}\Bigl(\widetilde{\bQ}_p^{-1}(z)\bD_1^{-1}(z)\Bigr)\Bigr|^2\biggr]^{1/2}\nonumber\\
    \leq&\; Kn(n^{-1}\delta_n^4)^{1/2}n^{-1/2}=K\delta_n^2,\label{begap1-eq2}
\end{align}
and
\begin{align}
    &\;\Cov\left(\beta_1(z)\gamma_1^{2}(z),\operatorname{tr}\left[\frac{\sigma^2}{\mu^2}\widetilde{\bQ}_p^{-1}(z)\bD_1^{-1}(z)\right]\right)\nonumber\\
    \leq&\; \bigl(\Expe|\beta_1(z)|^4\bigr)^{1/4}\bigl(\Expe|\gamma_1^{2}(z)|^4\bigr)^{1/4}\nonumber\\
&\qquad \times\biggl(\Expe\biggl|\operatorname{tr}\Bigl(\frac{\sigma^2}{\mu^2}\widetilde{\bQ}_p^{-1}(z)\bD_1^{-1}(z)\Bigr)-\Expe\operatorname{tr}\Bigl(\frac{\sigma^2}{\mu^2}\widetilde{\bQ}_p^{-1}(z)\bD_1^{-1}(z)\Bigr)\biggr|^2\biggr)^{1/2}\nonumber\\
    \leq&\; Kn^{-1/4}\delta_n^3.\label{begap1-eq3}
\end{align}
From \eqref{begap1-eq1}, \eqref{begap1-eq2}, and \eqref{begap1-eq3}, we complete the proof of equation \eqref{begap1}.

Finally, we consider \eqref{betap2}. From \eqref{covmatrix}, $\bD_1^{-1}(z)-\bD^{-1}(z)=\beta_1(z)\bD_1^{-1}(z)\br_1\br_1'\bD_1^{-1}(z)$, $\beta_1(z)=b_p(z) - b_p(z)\gamma_1(z)\beta_1(z)$ and $ \Expe\beta_1(z)=b_p(z)+o(n^{-1/2})$, we have
\begin{align}\label{betap2supl}
    &\;\Expe(\beta_1(z)P_2(z))\nonumber\\
    % =&\; \Expe\left(\beta_1(z)\Bigl[\frac{\sigma^2}{\mu^2}\operatorname{tr}(\widetilde{\bQ}_p^{-1}(z)\Expe\bD_1^{-1}(z))-\frac{\sigma^2}{\mu^2}\operatorname{tr}\left(\widetilde{\bQ}_p^{-1}(z)\Expe\bD^{-1}(z)\right)\Bigr]\right)\nonumber\\
    =&\; \frac{\sigma^2}{\mu^2}\Expe(\beta_1(z))\operatorname{tr}\Bigl[\widetilde{\bQ}_p^{-1}(z)\left(\Expe\bD_1^{-1}(z)-\Expe\bD^{-1}(z)\right)\Bigr]\nonumber\\
      =&\; \frac{\sigma^2}{\mu^2}\Expe(\beta_1(z))\operatorname{tr}\Bigl[\widetilde{\bQ}_p^{-1}(z)\Expe\left(\beta_1(z)\bD_1^{-1}(z)\br_1\br_1'\bD_1^{-1}(z)\right)\Bigr]\nonumber\\
       =&\; \frac{\sigma^2}{\mu^2}\Expe(\beta_1(z))\Expe\Bigl[(b_p(z) - b_p(z)\beta_1(z)\gamma_1(z))\br_1'\bD_1^{-1}(z)\widetilde{\bQ}_p^{-1}(z)\bD_1^{-1}(z)\br_1\Bigr]\nonumber\\
        =&\; b_p(z)\frac{\sigma^2}{\mu^2}\Expe(\beta_1(z))\Bigl[\Expe\br_1'\bD_1^{-1}(z)\widetilde{\bQ}_p^{-1}(z)\bD_1^{-1}(z)\br_1\nonumber\\
        &\qquad\qquad\qquad\qquad- \Expe\beta_1(z)\gamma_1(z)\br_1'\bD_1^{-1}(z)\widetilde{\bQ}_p^{-1}(z) \bD_1^{-1}(z)\br_1\Bigr]\nonumber\\
        =&\; b_p^2(z)\frac{\sigma^2}{\mu^2}\Bigl[\Expe\br_1'\bD_1^{-1}(z)\widetilde{\bQ}_p^{-1}(z)\bD_1^{-1}(z)\br_1\Bigr] + O(n^{-1/2})\nonumber\\
        =&\; b_p^2(z)\frac{\sigma^2}{\mu^2}\operatorname{tr}\Bigl[\Expe\bD_1^{-1}(z)\widetilde{\bQ}_p^{-1}(z)\bD_1^{-1}(z)\Expe\br_1\br_1'\Bigr] + O(n^{-1/2})\nonumber\\
        =&\; b_p^2(z)\frac{\sigma^2}{\mu^2}\operatorname{tr}\Bigl[\Expe\bD_1^{-1}(z)\widetilde{\bQ}_p^{-1}(z)\bD_1^{-1}(z)\nonumber\\
        &\qquad\qquad\qquad\qquad\times  \frac{1}{n}\nu_2 \Bigl(-\frac{1}{p-1}\boldsymbol{1}_{p}\boldsymbol{1}'_{p} + \frac{1}{p-1}\bI_p + \bI_p \Bigr)  \Bigr] + O(n^{-1/2})\nonumber\\
        =&\; \frac{p}{n(p-1)}\nu_2\frac{\sigma^2}{\mu^2}b_p^2(z)\Expe\operatorname{tr}\Bigl[\bD_1^{-1}(z)\widetilde{\bQ}_p^{-1}(z)\bD_1^{-1}(z)\Bigr] \nonumber\\
        &\qquad-\frac{1}{n(p-1)}\nu_2\frac{\sigma^2}{\mu^2}b_p^2(z)\Expe\operatorname{tr}\Bigl[  \bD_1^{-1}(z)\widetilde{\bQ}_p^{-1}(z)\bD_1^{-1}(z)\boldsymbol{1}_{p}\boldsymbol{1}'_{p} \Bigr] + O(n^{-1/2})\nonumber\\
     =&\; \frac{p}{n(p-1)}\nu_2\frac{\sigma^2}{\mu^2}b_p^2(z)\Expe\operatorname{tr}\Bigl[\bD_1^{-1}(z)\widetilde{\bQ}_p^{-1}(z)\bD_1^{-1}(z)\Bigr] + O(n^{-1/2}),\nonumber
\end{align}
which is the equation \eqref{betap2}.
Below are some interpretations of the above equalities:
\begin{enumerate}
    \item The fifth equality follows from
    \[
    \Expe(\beta_1(z)\gamma_1(z)\br_1'\bD_1^{-1}(z)\widetilde{\bQ}_p^{-1}(z) \bD_1^{-1}(z)\br_1)\leq Kn^{-1/2},
    \]
    which can proved by using Lemma \ref{bound}.
    \item The last equality follows from
    \[
        \frac{1}{n(p-1)}\nu_2\frac{\sigma^2}{\mu^2}b_p^2(z)\Expe\operatorname{tr}\Bigl[  \bD_1^{-1}(z)\widetilde{\bQ}_p^{-1}(z)\bD_1^{-1}(z)\boldsymbol{1}_{p}\boldsymbol{1}'_{p} \Bigr] = O(n^{-1}).
    \]
    This can be proved by using Lemma \ref{estsqu} and Lemma \ref{xdxz}.
\end{enumerate}

\subsection{Proof of Lemma \ref{J123limit}}\label{prfJ123limit}
\textbf{Step 1: Consider $J_1$.} By \eqref{covmatrixmain}, we get
\begin{align}\label{ngamma}
    n\Expe\gamma_1(z) 
    % &= n\Expe\left( \br_1'\bD_1^{-1}(z)\br_1-\frac{1}{n}\nu_2\Expe\operatorname{tr}\bD_1^{-1}(z)  \right)\nonumber\\
    % &= n\left( \operatorname{tr}\Bigl[\Expe\bD_1^{-1}(z)\Expe\br_1\br_1'\Bigr]-\frac{1}{n}\nu_2\Expe\operatorname{tr}\bD_1^{-1}(z)  \right)\nonumber\\
    &= n\left( \frac{1}{n}\nu_2\Expe \tr \Bigl[\Bigl(-\frac{1}{p-1}\boldsymbol{1}_{p}\boldsymbol{1}'_{p} + \frac{1}{p-1}\bI_p + \bI_p \Bigr)\bD_1^{-1}(z)  \Bigr]-\frac{1}{n}\nu_2\Expe\operatorname{tr}\bD_1^{-1}(z)  \right)\nonumber\\
     &= -\frac{1}{p-1}\nu_2\Expe\operatorname{tr}\bigl(\boldsymbol{1}'_{p}\bD_1^{-1}(z)\boldsymbol{1}_{p} \bigr) + \frac{1}{p-1}\nu_2\Expe\operatorname{tr}\bD_1^{-1}(z).
\end{align}
By Lemma \ref{xdxz}, we have
\begin{align}\label{j12limit}
    &\;-b_p^2(z)\biggl[\frac{1}{n(p-1)}\nu_2^2\Expe\Bigl(\operatorname{tr}\bD_1^{-1}(z)\operatorname{tr}\bigl(\bD_1^{-1}(z)\widetilde{\bQ}_p^{-1}(z)\bigr)\Bigr)\nonumber\\
   &\quad\quad\quad\quad-\frac{1}{n(p-1)}\nu_2^2\Expe\Bigl(\boldsymbol{1}'_{p}\bD_1^{-1}(z)\boldsymbol{1}_{p}\operatorname{tr}\bigl(\bD_1^{-1}(z)\widetilde{\bQ}_p^{-1}(z)\bigr)\Bigr)\biggr]\nonumber\\
   \stackrel{i.p.}{\rightarrow}&\;\frac{-cz^2\underline{m}^2(z)\sigma^4/\mu^4}{1+\frac{\sigma^2}{\mu^2}\underline{m}(z)}\Bigl\{m^2(z)+\frac{m(z)}{z}\Bigr\}.
\end{align}
By Lemma 2.3 in \cite{bai2004clt}, we have
\begin{align}\label{trdidilim}
   \frac{1}{p}\operatorname{tr}\left[\bD_j^{-1}(z)\right]\stackrel{L^2}{\rightarrow} m(z),\ as\ p \rightarrow \infty.
\end{align}
By Lemma \ref{estsqu}, equations \eqref{ngamma}, \eqref{j12limit},
 \eqref{trdidilim}, and \eqref{b1zm}, we have
\begin{align}\label{J1meanlim}
    % J_1 \stackrel{i.p.}{\rightarrow} &h_1m(z)\frac{-z\underline{m}(z)}{1+\frac{\sigma^2}{\mu^2}\underline{m}(z)}+\Bigl(\frac{-z\underline{m}(z)}{1+\frac{\sigma^2}{\mu^2}\underline{m}(z)}+\frac{c\sigma^2/\mu^2z^2\underline{m}^{2}(z)m(z)}{1+\frac{\sigma^2}{\mu^2}\underline{m}(z)} \Bigr)\Bigl(\frac{\sigma^2}{\mu^2}m(z)+\frac{\sigma^2}{\mu^2}\frac{1}{z}\Bigr)\nonumber\\
    % &\ \ \ \ -\frac{c\sigma^4/\mu^4z^2\underline{m}^2(z)m(z)}{1+\frac{\sigma^2}{\mu^2}\underline{m}(z)}(m(z)+\frac{1}{z})\nonumber\\
     J_1 \stackrel{i.p.}{\rightarrow} \biggl(\frac{-z\underline{m}(z)}{1+\frac{\sigma^2}{\mu^2}\underline{m}(z)}\biggr)\Bigl(\frac{\sigma^2}{\mu^2}m(z)+\frac{\sigma^2}{\mu^2}\frac{1}{z}+h_1m(z)\Bigr).
\end{align}
\textbf{Step 2: Consider $J_2$.} By Lemma \ref{quadform}, we have
\begin{align}\label{J2mean}
  J_2 = J_{21} + J_{22} + J_{23} + J_{24} + o(1),
  % &= -nb_p^2(z)\Expe\left(\br_1'\bD_1^{-1}(z)\br_1-\frac{1}{n}\nu_2\operatorname{tr} \bD_1^{-1}(z)\right)\nonumber\\  &\qquad\qquad\times\left(\br_1'\bD_1^{-1}(z)\widetilde{\bQ}_p^{-1}(z)\br_1-\frac{1}{n}\nu_2\operatorname{tr} [\bD_1^{-1}(z)\widetilde{\bQ}_p^{-1}(z)]\right)\nonumber\\
    % &= -nb_p^2(z)\Bigl[\frac{1}{n^2}\Big(\nu_4-3\nu_{12} \Big)  \sum_{i=1}^{p}\Expe(\bD_1^{-1}(z))_{ii}(\bD_1^{-1}(z)\widetilde{\bQ}_p^{-1}(z))_{ii}\nonumber\\
    %  &\quad+\frac{1}{n^2}\nu_{12}\left(\Expe\operatorname{tr}\bD_1^{-1}(z)\widetilde{\bQ}_p^{-1}(z)\bD_1^{-1}(z)+\Expe\operatorname{tr}\bD_1^{-1}(z)\widetilde{\bQ}_p^{-1}(z)(\bD_1^{-1}(z))'\right)\nonumber\\
    %  &\quad+\frac{1}{n^2}\textcolor{red}{(\nu_{12}-\nu_2^2)}E\left(\operatorname{tr}\bD_1^{-1}(z)\operatorname{tr}\bD_1^{-1}(z)\widetilde{\bQ}_p^{-1}(z)\right)\Bigr]+o(1)\nonumber\\
\end{align}
where
\begin{align}
    J_{21}  &= -nb_p^2(z)\Bigl[\frac{1}{n^2}(\nu_4-3\nu_{12} )  \sum_{i=1}^{p}\Expe\bigl[\bigl(\bD_1^{-1}(z)\bigr)_{ii}\bigl(\bD_1^{-1}(z)\widetilde{\bQ}_p^{-1}(z)\bigr)_{ii}\bigr]\Bigr],\nonumber\\
    J_{22}&=-nb_p^2(z)\Bigl[\frac{1}{n^2}\nu_{12}\Expe\operatorname{tr}\bigl(\bD_1^{-1}(z)\widetilde{\bQ}_p^{-1}(z)\bD_1^{-1}(z)\bigr)\Bigr],\nonumber\\
    J_{23}&=-nb_p^2(z)\Bigl[\frac{1}{n^2}\nu_{12}\Expe\operatorname{tr}\bigl(\bD_1^{-1}(z)\widetilde{\bQ}_p^{-1}(z)(\bD_1^{-1}(z))'\bigr)\Bigr],\nonumber\\
    J_{24}&=-nb_p^2(z)\Bigl[\frac{1}{n^2}(\nu_{12}-\nu_2^2)\Expe\bigl[\operatorname{tr}\bD_1^{-1}(z)\operatorname{tr}\bigl(\bD_1^{-1}(z)\widetilde{\bQ}_p^{-1}(z)\bigr)\bigr]\Bigr].\nonumber
\end{align}
Since
$\frac{1}{p}\sum_{i=1}^{p}\Expe(\bD_1^{-1}(z))_{ii}(\bD_1^{-1}(z)\widetilde{\bQ}_p^{-1}(z))_{ii}\stackrel{i.p.}{\rightarrow} \frac{m^2(z)}{1+\frac{\sigma^2}{\mu^2}\underline{m}(z)}$,
\begin{align}\label{J21mean}
    J_{21} \stackrel{i.p.}{\rightarrow} \frac{-cz^2m^2(z)\underline{m}^2(z)}{1+\frac{\sigma^2}{\mu^2}\underline{m}(z)}\Big(\frac{\Expe |w_{1}-\mu|^{4}}{\mu^{4}}-3\frac{|\Expe |w_{1}-\mu|^{2}|^{2}}{\mu^{4}} \Big).
\end{align}
Note that
\begin{align}
    \frac{1}{p}\Expe\operatorname{tr}\bD_1^{-2}(z) \stackrel{i.p.}{\rightarrow} m'(z),
\end{align}
% and
% \begin{align}
%    J_{22} &= -nb_p^2(z)\times\frac{1}{n^2}\nu_{12}\Expe\operatorname{tr}\bD_1^{-1}(z)\widetilde{\bQ}_p^{-1}(z)\bD_1^{-1}(z)\nonumber\\
%    &= -\frac{b_p^2(z)}{n}\nu_{12}\frac{1}{1+\frac{\sigma^2}{\mu^2}\underline{m}(z)}\Expe\operatorname{tr}\bD_1^{-2}(z)
% \end{align}
thus we get
\begin{align}\label{J22mean}
    J_{22}\stackrel{i.p.}{\rightarrow} \frac{-cz^2m'(z)\underline{m}^2(z)}{1+\frac{\sigma^2}{\mu^2}\underline{m}(z)}\frac{|\Expe |w_{1}-\mu|^{2}|^{2}}{\mu^{4}}.
\end{align}
Similarly, we have
\begin{align}\label{J23mean}
    J_{23}\stackrel{i.p.}{\rightarrow} \frac{-cz^2m'(z)\underline{m}^2(z)}{1+\frac{\sigma^2}{\mu^2}\underline{m}(z)}\frac{|\Expe |w_{1}-\mu|^{2}|^{2}}{\mu^{4}}.
\end{align}
% Since
% \begin{align}
%     J_{24} &=-nb_p^2(z)\times \frac{1}{n^2}\textcolor{red}{(\nu_{12}-\nu_2^2)}\Expe\left(\operatorname{tr}\bD_1^{-1}(z)\operatorname{tr}\bD_1^{-1}(z)\widetilde{\bQ}_p^{-1}(z)\right)\nonumber\\
%     &=  -\frac{1}{n}\textcolor{red}{(\nu_{12}-\nu_2^2)}b_p^2(z)\frac{1}{1+\frac{\sigma^2}{\mu^2}\underline{m}(z)}\Expe\left(\operatorname{tr}\bD_1^{-1}(z)\operatorname{tr}\bD_1^{-1}(z)\right)\nonumber
% \end{align}
By Lemma \ref{estsqu}, we get
\begin{align}\label{J24mean}
     J_{24} \stackrel{i.p.}{\rightarrow} -c\Bigl(h_2-2\frac{\sigma^2}{\mu^2}h_1\Bigr)\frac{z^2m^2(z)\underline{m}^2(z)}{1+\frac{\sigma^2}{\mu^2}\underline{m}(z)}.
\end{align}
From \eqref{J2mean}, \eqref{J21mean}, \eqref{J22mean}, \eqref{J23mean}, and \eqref{J24mean}, we have
\begin{align}\label{J2meanlim}
    J_{2} \stackrel{i.p.}{\rightarrow} &-\frac{cz^2m^2(z)\underline{m}^2(z)}{1+\frac{\sigma^2}{\mu^2}\underline{m}(z)}\Big(\frac{\Expe |w_{1}-\mu|^{4}}{\mu^{4}}-3\frac{|\Expe |w_{1}-\mu|^{2}|^{2}}{\mu^{4}} \Big)\nonumber\\
    &-\frac{2cz^2m'(z)\underline{m}^2(z)}{1+\frac{\sigma^2}{\mu^2}\underline{m}(z)}\frac{|\Expe |w_{1}-\mu|^{2}|^{2}}{\mu^{4}}-c\Bigl(h_2-2\frac{\sigma^2}{\mu^2}h_1\Bigr)\frac{z^2m^2(z)\underline{m}^2(z)}{1+\frac{\sigma^2}{\mu^2}\underline{m}(z)}.
\end{align}
\textbf{Step 3: Consider $J_3$.} To calculate the limit of $J_3$, we can expand $\bD_1^{-1}(z)$ like \eqref{djdecomp} and find the limit of $J_3$ using the method similarly to \cite{bai2004clt}. The limit of $J_3$ is
\begin{align}\label{J3meanlim}
    J_3 \stackrel{i.p.}{\rightarrow} c\frac{\sigma^4}{\mu^4}\underline{m}^2(z)(1+\frac{\sigma^2}{\mu^2}\underline{m}(z))^{-3}\Bigl[1-c\frac{\sigma^4}{\mu^4}\underline{m}^2(z)(1+\frac{\sigma^2}{\mu^2}\underline{m}(z))^{-2}\Bigr]^{-1}.
\end{align}
% One can also find the limit of $J_3$ by the following calculation, note that
% \begin{align}
%      J_3 &= \frac{p}{n(p-1)}b_p^2(z)\frac{\sigma^2}{\mu^2}\nu_2 \frac{1}{1+\frac{\sigma^2}{\mu^2}\underline{m}(z)}\Expe \operatorname{tr}[\bD_1^{-2}(z)]\nonumber,
% \end{align}
% and
% \begin{align}
%     \frac{1}{p}\operatorname{tr}[\bD_1^{-2}(z)]\stackrel{L^2}{\rightarrow} m'(z),
% \end{align}
% therefore,
% \begin{align}
%     J_3 \stackrel{i.p.}{\rightarrow} c\frac{\sigma^4}{\mu^4}\frac{z^2\underline{m}^2(z)m'(z)}{1+\frac{\sigma^2}{\mu^2}\underline{m}(z)}.
% \end{align}
From \eqref{J1meanlim}, \eqref{J2meanlim}, and \eqref{J3meanlim}, the proof is completed.

\subsection{Proof of Lemma \ref{xdxz}}\label{sec:lem-xdxz-proof}

By Lemma \ref{bound}, we obtain, for any $2\leq r\in \mathbb{N}^+$,
\begin{align}\label{eq:rdxxdr}
	&\;\Expe|\br_{j}'\bD_{j}^{-1}(z)\bx_{p}\bx_{p}'\bD_{j}^{-1}(z)\br_{j}|^{r}\nonumber\\
    \leq&\; K_r\Bigl(\Expe\Bigl|\br_{j}'\bD_{j}^{-1}(z)\bx_{p}\bx_{p}'\bD_{j}^{-1}(z)\br_{j}-\frac{1}{n}\nu_2\operatorname{tr}\bigl[\bD_{j}^{-1}(z)\bx_{p}\bx_{p}'\bD_{j}^{-1}(z)\bigr]\Bigr|^{r}\nonumber\\
	&\qquad \qquad + \frac{1}{n}\nu_2\Expe\bigl|\operatorname{tr}[\bD_{j}^{-1}(z)\bx_{p}\bx_{p}'\bD_{j}^{-1}(z)]\bigr|^{r}\Bigr)\nonumber\\  
	\leq&\; K_r(n^{-2}\delta_n^{2r-4}+n^{-r}),
\end{align}
where $\frac{1}{n}\operatorname{tr}[\bD_{j}^{-1}(z)\bx_{p}\bx_{p}'\bD_{j}^{-1}(z)]=\frac{1}{n}\bx_{p}'\bD_{j}^{-2}(z)\bx_{p}=\frac{1}{np}\boldsymbol{1}_{p}'\bD_{j}^{-2}(z)\boldsymbol{1}_{p}=O(n^{-1})$.

Write
\begin{align}
	&\;\bx_{p}'\bD^{-1}(z)\bx_{p}-\bx_{p}'\Expe\bD^{-1}(z)\bx_{p}\nonumber\\
 =&\;\sum_{j=1}^{n}\Bigl(\bx_{p}'\Expe_{j}\bD^{-1}(z)\bx_{p}-\bx_{p}'\Expe_{j-1}\bD^{-1}(z)\bx_{p}\Bigr)\nonumber\\
	=&\;\sum_{j=1}^{n}\Bigl(\bx_{p}'\Expe_{j}\bigl(\bD^{-1}(z)-\bD_{j}^{-1}(z)\bigr)\bx_{p}-\bx_{p}'\Expe_{j-1}\bigl(\bD^{-1}(z)-\bD_{j}^{-1}(z)\bigr)\bx_{p}\Bigr)\nonumber\\
	=&\;-\sum_{j=1}^{n}(\Expe_{j} -\Expe_{j-1})\beta_{j}(z)\br_{j}'\bD_{j}^{-1}(z)\bx_{p}\bx_{p}'\bD_{j}^{-1}(z)\br_{j}.\nonumber
\end{align}
By Burkholder's inequality, \eqref{eq:rdxxdr}, and $|\beta_{j}(z)|\leq \frac{|z|}{v_0}$, we have
\begin{align}\label{xdx}
	&\ \ \ \ \Expe|\bx_{p}'\bD^{-1}(z)\bx_{p}-\bx_{p}'\Expe\bD^{-1}(z)x_{p}|^{2}\nonumber\\
	&\leq K\sum_{j=1}^{n}\Expe|(\Expe_{j} -\Expe_{j-1})\beta_{j}(z)\br_{j}'\bD_{j}^{-1}(z)\bx_{p}\bx_{p}'\bD_{j}^{-1}(z)\br_{j}|^{2}\nonumber\\
	&\leq K\sum_{j=1}^{n}\Expe|\beta_{j}(z)\br_{j}'\bD_{j}^{-1}(z)\bx_{p}\bx_{p}'\bD_{j}^{-1}(z)\br_{j}|^{2}\nonumber\\
	&\leq Kn^{-1}.
\end{align}
Thus, we have 
\begin{equation}\label{eq:xdx_xEDx}
\Expe|\bx_{p}'\bD^{-1}(z)\bx_{p}-\bx_{p}'\Expe\bD^{-1}(z)\bx_{p}|^2\to 0.
\end{equation}
% If $\Im z \ge v_0>0$, then $|\beta_j(z)|\le \frac{|z|}{v_0}$, so \eqref{xdx} can get a sharper bound
% \begin{eqnarray}\label{xdx2}
% 	\Expe|\bx_{p}'\bD^{-1}(z)\bx_{p}-\bx_{p}^{*}\Expe\bD^{-1}(z)\bx_{p}|^{2}=O\left(n^{-1}\right).
% \end{eqnarray}
% From the proof above, one should note that \eqref{xdx} and \eqref{xdx2} hold for $\bD_j(z)$ and any unit vector.

Note that $\bD(z)+z\widetilde{\bQ}_p(z)=\sum_{j=1}^{n}\br_{j}\br_{j}'+z\frac{\sigma^2}{\mu^2}\Expe\underline{m}_{p}(z)\bI_p$, where $\widetilde{\bQ}_p(z)=\frac{\sigma^2}{\mu^2}\Expe\underline{m}_{p}(z)\bI_p+\bI_p$.
Recalling $\underline{m}_p(z)=-\frac{1}{nz}\sum_{j=1}^n\beta_j(z)$ and using the identity
$\br_{j}'\bD^{-1}(z)=\beta_{j}(z)\br_{j}' \bD_{j}^{-1}(z)$,  we obtain
% \begin{eqnarray*}
% 	&&(-z\widetilde{\bQ}_p(z))^{-1}-\Expe\bD^{-1}(z)\nonumber\\
% 	&=& -\left[(z\widetilde{\bQ}_p(z))^{-1}(\Expe\bD(z)+z\widetilde{\bQ}_p(z))\Expe\bD^{-1}(z)\right]\nonumber\\
% 	&=&-z^{-1}\widetilde{\bQ}^{-1}_p(z)(z)\Expe\big[\big(\sum^{n}_{j=1}\br_j\br_j'+z\frac{\sigma^2}{\mu^2}\Expe\underline{m}_p(z)\bI_p\big)\bD^{-1}(z)\big]\nonumber\\
% 	&=&-z^{-1}\sum^{n}_{j=1}\Expe\beta_j(z)\big[\widetilde{\bQ}^{-1}_p(z)\br_j\br_j'\bD_j^{-1}(z)\big]\nonumber\\
% 	&&-z^{-1}\Expe\big[\widetilde{\bQ}^{-1}_p(z)\big(z\frac{\sigma^2}{\mu^2}\Expe\underline{m}_p(z)\big)\bI_p\bD^{-1}(z)\big]\nonumber\\
% 	&=&-z^{-1}\sum^{n}_{j=1}\Expe\beta_j(z)\big[\widetilde{\bQ}^{-1}_p(z)\br_j\br_j'\bD_j^{-1}(z)\big]\nonumber\\
% 	&&-z^{-1}\Expe\big[\widetilde{\bQ}^{-1}_p(z)\big(- \frac{1}{n}\frac{\sigma^2}{\mu^2}\sum_{j=1}^n\Expe\beta_j(z)\big)\bI_p\bD^{-1}(z)\big]\nonumber\\
% 	&=&-z^{-1}\sum^{n}_{j=1}\Expe\beta_j(z)\big[\widetilde{\bQ}^{-1}_p(z)\br_j\br_j'\big(\underline{\bB}_{(j)}^{p}-z\bI_p\big)^{-1}-\frac{1}{n}\frac{\sigma^2}{\mu^2}\widetilde{\bQ}^{-1}_p(z)\Expe\bD^{-1}(z)\big]\nonumber\\
% 	&=&-z^{-1}n\Expe\beta_1(z)\big[\widetilde{\bQ}^{-1}_p(z)\br_1\br_1'\bD^{-1}_1(z)
% 	-\frac{1}{n}\frac{\sigma^2}{\mu^2}\widetilde{\bQ}^{-1}_p(z)\Expe\bD^{-1}(z)\big].
% \end{eqnarray*}
\begin{align*}
	&\;\bigl(-z\widetilde{\bQ}_p(z)\bigr)^{-1}-\bD^{-1}(z)\nonumber\\
	=&\; -\bigl(z\widetilde{\bQ}_p(z)\bigr)^{-1}\bigl(\bD(z)+z\widetilde{\bQ}_p(z)\bigr)\bD^{-1}(z)\nonumber\\
	=&\;-z^{-1}\widetilde{\bQ}^{-1}_p(z)\biggl(\sum^{n}_{j=1}\br_j\br_j'+z\frac{\sigma^2}{\mu^2}\Expe\underline{m}_p(z)\bI_p\biggr)\bD^{-1}(z)\nonumber\\
    =&\;-z^{-1}\sum^{n}_{j=1}\beta_j(z)\widetilde{\bQ}^{-1}_p(z)\br_j\br_j'\bD_j^{-1}(z)
    -z^{-1}\widetilde{\bQ}^{-1}_p(z)\Bigl(- \frac{\sigma^2}{\mu^2}\Expe\beta_1(z)\Bigr)\bD^{-1}(z).
\end{align*}
Taking expectation of the above identity yields that
\begin{align*}
    &\; \bigl(-z\widetilde{\bQ}_p(z)\bigr)^{-1}-\Expe\bD^{-1}(z)\\
    =&\;-z^{-1}n\Expe\beta_1(z)\Bigl[\widetilde{\bQ}^{-1}_p(z)\br_1\br_1'\bD^{-1}_1(z)
	-\frac{1}{n}\frac{\sigma^2}{\mu^2}\widetilde{\bQ}^{-1}_p(z)\Expe\bD^{-1}(z)\Bigr].
\end{align*}
Multiplying by $-\bx_{p}'$ on the left and $\bx_{p}$ on the right, we have
\begin{align*}
	&\;\bx_{p}'\Expe\bD^{-1}(z)\bx_{p}-\bx_{p}' \bigl(-z\widetilde{\bQ}_p(z)\bigr)^{-1}\bx_{p}\nonumber\\
	=&\;z^{-1}n\Expe\Bigl\{\beta_1(z)\Bigl[\bx_{p}' \widetilde{\bQ}^{-1}_p(z)\br_1\br_1'\bD^{-1}_1(z)\bx_{p}
	-\frac{1}{n}\frac{\sigma^2}{\mu^2}\bx_{p}' \widetilde{\bQ}^{-1}_p(z)\Expe\bD^{-1}(z)\bx_{p} \Bigr]\Bigr\}\nonumber\\
	=:&\;\delta_1 +\delta_2 +\delta_3 ,
\end{align*}
where 
\begin{align*}
\delta_1 &=\frac{n}{z}\Expe(\beta_1 (z) \alpha_1 (z)),\\
\alpha_1 (z)&=\bx_{p}' \widetilde{\bQ}^{-1}_p(z)\br_1\br_1'\bD^{-1}_1(z)\bx_{p}
-\frac{1}{n}\frac{\sigma^2}{\mu^2}\bx_{p}^* \widetilde{\bQ}^{-1}_p(z)\bD^{-1}_1(z)\bx_{p},\\
\delta_2 &=\frac{1}{z}\frac{\sigma^2}{\mu^2}\Expe\beta_1 (z)\bx_{p}'\widetilde{\bQ}^{-1}_p(z)\bigl(\bD_{1}^{-1}(z)-\bD^{-1}(z)\bigr)\bx_{p},\\
\delta_3 &=\frac{1}{z}\frac{\sigma^2}{\mu^2}\Expe\beta_1 (z)\bx_{p}' \widetilde{\bQ}^{-1}_p(z)\bigl(\bD^{-1}(z)-\Expe\bD^{-1}(z)\bigr)\bx_{p}.
\end{align*}
Recalling the notations defined above and the following equalities:
\begin{align*}
    \delta_1 &=\frac{n}{z}\Expe\bar \beta_1 (z) \alpha_1 (z)-\frac{n}{z}\Expe\left[\beta_1 (z) \bar \beta_1 (z) \varepsilon_1 (z)\alpha_1 (z)\right],\\
    \bar{\beta}_1 (z)&=b_p(z)-\frac{\nu_2}{n}b_p(z)\bar \beta_{1}(z)tr(\bD_1^{-1}(z)-\Expe\bD_{1}^{-1}(z)),\\
    \Expe\alpha_1(z) &=-(\nu_2\Expe\underline{m}_p(z)+1\big)^{-1}\frac{\nu_2}{n}\bigl[\Expe \bx'_{p}\bD_{1}^{-1}(z)\bx_{p}+o(1)\bigr]. 
\end{align*}
From Lemma \ref{estsqu} and \eqref{covmatrix},
it is easy to see that
\[
n\Expe\beta_1 (z) \alpha_1 (z)=\biggl[\frac{1}{1+c\nu_2\Expe m_{p}(z)}+o(1)\biggr]n\Expe \alpha_1 (z).
\]
Therefore, $\delta_1 =\frac{\nu_2 z\Expe\underline{m}_p(z)}{\nu_2 z\Expe\underline{m}_p(z)+z}\bx_{p}'\Expe(\bD_{1}^{-1}(z))\bx_{p} +o(1)$.
Similarly to \cite{Bai2007On}, one may have $\delta_2 =o(1)$ and $\delta_3 =o(1)$.
Hence, we obtain 
$$\Big(1- \frac{\nu_2 z\Expe\underline{m}_p(z)}{\nu_2 z\Expe\underline{m}_p(z)+z}\Big)\bx_{p}'\Expe(\bD^{-1}(z))\bx_{p}+\frac{1}{\sigma^2/\mu^2 z\Expe\underline{m}_p(z)+z} \to 0,$$
which implies that
$$\bx_{p}'\Expe\bD^{-1}(z)\bx_{p} \to -\frac{1}{z}.$$
This, together with \eqref{eq:xdx_xEDx}, completes the proof of Lemma \ref{xdxz}.

% \subsection{Proof of Theorem \ref{clt}}\label{prfclt}

\subsection{Proof of Corollary \ref{corpoly}}\label{prfcorpoly}

For ease of presentation, we denote $\lambda=\sigma^2/\mu^2$ in this section.
Note that 
\begin{align}
    z &= -\frac{1}{\underline{m}}+\frac{c\lambda}{1+\lambda\underline{m}},\label{zeq}\\%\label{zexpr}
    \dif z &= \Bigl[\underline{m}^{-2}-c\lambda^{2}(1+\lambda \underline{m})^{-2}]\dif \underline{m}\nonumber\\
    &=\left[1-c\lambda^{2}\underline{m}^{2}(1+\lambda \underline{m}^{-2})\right]\underline{m}^{-2}\dif\underline{m},\label{dzexpr}\\
    m&=\frac{1}{c}\Bigl(\underline{m}+\frac{1-c}{z}\Bigr),\label{mexpr}\\
    m'&=\frac{1}{c}\Bigl(\underline{m}'-\frac{1-c}{z^2}\Bigr).\label{dmexpr}
\end{align}

For $1\leq r_1, r_2 \in\bbN^{+}$, 
\begin{align}%\label{polycovprfeq}
  &\;\Cov(X_{x^{r_1}}, X_{x^{r_2}}) \nonumber\\
  =&\;2 \times (\lambda c)^{r_{1}+r_{2}} \sum_{k_{1}=0}^{r_{1}-1} \sum_{k_{2}=0}^{r_{2}}\binom{r_1}{k_1}\binom{r_2}{k_2}\left(\frac{1-c}{c}\right)^{k_{1}+k_{2}} \nonumber\\
&\qquad\times \sum_{\ell=1}^{r_{1}-k_{1}} \ell\binom{2 r_{1}-1-\left(k_{1}+\ell\right)}{r_{1}-1}
\binom{2 r_{2}-1-k_{2}+\ell}{r_{2}-1}\label{polycovprfeq1}\\
&\ \ \ \ +\frac{1}{\lambda^2} \left(\beta+h_2-2\lambda h_1 \right)(\lambda c)^{r_{1}+r_{2}}\nonumber\\
&\ \ \ \ \ \ \ \ \times\sum_{k_{1}=0}^{r_{1}} \sum_{k_{2}=0}^{r_{2}}\binom{r_1}{k_1}\binom{r_2}{k_2}\left(\frac{1-c}{c}\right)^{k_{1}+k_{2}}
\binom{2r_1-k_1}{r_1-1}\binom{2r_2-k_2}{r_2-1}.\label{polycovprfeq2}
\end{align}
The proof of \eqref{polycovprfeq1} is exactly analogous with \cite{bai2004clt}, it is then omitted. Next, we prove \eqref{polycovprfeq2}. The contours $\mathcal{C}, \mathcal{C}_1, \mathcal{C}_2$ are closed and taken in the positive direction in the complex plane, each enclosing the support of $F^{c,H}$.
Note that
\begin{align}  \oint_{\mathcal{C}_1}\oint_{\mathcal{C}_2}\frac{z_1^{r_1}z_2^{r_2}}{(1+\lambda \underline{m}_1)^2(1+\lambda \underline{m}_2)^2}\dif \underline{m}_1\dif \underline{m}_2&=\oint_{\mathcal{C}_1}\frac{z_1^{r_1}}{(1+\lambda \underline{m}_1)^2}\dif \underline{m}_1\times\oint_{\mathcal{C}_2}\frac{z_2^{r_2}}{(1+\lambda \underline{m}_2)^2}\dif \underline{m}_2.
\end{align}
By \eqref{zeq},
\begin{align}
    \oint_{\mathcal{C}_1}\frac{z_1^{r_1}}{(1+\lambda \underline{m}_1)^2}\dif \underline{m}_1&=\oint_{\mathcal{C}_1}\frac{\Bigl(-\frac{1}{\underline{m}_1}+\frac{c\lambda}{1+\lambda\underline{m}_1}\Bigr)^{r_1}}{(1+\lambda \underline{m}_1)^2}\dif \underline{m}_1\nonumber\\
    &=(\lambda c)^{r_1}\oint_{\mathcal{C}_1}\Bigl(\frac{1}{1+\lambda \underline{m}_1}+\frac{1-c}{c}\Bigr)^{r_1}(1-(1+\lambda \underline{m}_1))^{-r_1}(1+\lambda \underline{m}_1)^{-2}\dif \underline{m}_1\nonumber\\
    &=(\lambda c)^{r_1}\oint_{\mathcal{C}_1}\sum_{k_1=0}^{r_1}\binom{r_{1}}{k_{1}}
    \Bigl(\frac{1-c}{c}\Bigr)^{k_1}(1+\lambda \underline{m}_1)^{k_1-r_1}\nonumber\\
&\ \ \ \ \times\sum_{j=0}^{\infty}\binom{
r_{1}+j-1}{j}(1+\lambda \underline{m}_1)^{j}(1+\lambda \underline{m}_1)^{-2}\dif \underline{m}_1\nonumber\\
&=(\lambda c)^{r_1}\sum_{k_1=0}^{r_1}\binom{
r_{1}}{k_{1}}\Bigl(\frac{1-c}{c}\Bigr)^{k_1}\nonumber\\
&\ \ \ \ \times\oint_{\mathcal{C}_1}\sum_{j=0}^{\infty}\binom{r_{1}+j-1}{j}(1+\lambda \underline{m}_1)^{k_1-r_1+j-2}\dif \underline{m}_1,
\end{align}
by substitution $\widetilde{\underline{m}}_1=\lambda \underline{m}_1$, we get
\begin{align}
    \oint_{\mathcal{C}_1}\frac{z_1^{r_1}}{\Bigl(1+\lambda \underline{m}_1\Bigr)^2}\dif \underline{m}_1&=\frac{1}{\lambda}\Bigl(\lambda c\Bigr)^{r_1}\sum_{k_1=0}^{r_1}\left(\begin{array}{l}
r_{1} \\
k_{1}
\end{array}\right)\Bigl(\frac{1-c}{c}\Bigr)^{k_1}\nonumber\\
&\ \ \ \ \times\oint_{\mathcal{C}_1}\sum_{j=0}^{\infty}\left(\begin{array}{l}
r_{1}+j-1 \\
\ \ \ \ j
\end{array}\right)\Big(1+\underline{m}_1\Big)^{k_1-r_1+j-2}\dif \underline{m}_1.\label{covinttran}
\end{align}
For integral \eqref{covinttran}, the pole is $-1$, we have by residual theorem
\begin{align}
    \oint_{\mathcal{C}_1}\frac{z_1^{r_1}}{(1+\lambda \underline{m}_1)^2}\dif \underline{m}_1&=\frac{2\pi i}{\lambda}(\lambda c)^{r_1}\sum_{k_1=0}^{r_1}\binom{r_1}{k_1}\Bigl(\frac{1-c}{c}\Bigr)^{k_1}\binom{2r_{1}-k_1}{r_{1}-1}.\label{z1term}
\end{align}
Similarly, we get
\begin{align}
    \oint_{\mathcal{C}_2}\frac{z_2^{r_2}}{(1+\lambda \underline{m}_2)^2}\dif \underline{m}_2&=\frac{2\pi i}{\lambda}(\lambda c)^{r_2}\sum_{k_2=0}^{r_2}\binom{r_2}{k_2}\Bigl(\frac{1-c}{c}\Bigr)^{k_2}\binom{2r_{2}-k_2}{r_{2}-1}.\label{z2term}
\end{align}
By \eqref{z1term} and \eqref{z2term}, \eqref{polycovprfeq2} is derived.

For $f_1=x$,  
\begin{align}\label{f1expe0}
    \Expe X_x
    &= \frac{1}{2\pi i}\oint_{\mathcal{C}} c\lambda^{2}z\underline{m}^{3}(z)\left[1-c\lambda^{2}\underline{m}^{2}(z)(1+\lambda \underline{m}(z))^{-2}\right]^{-1}\nonumber\\
    &\qquad \times\left(1+\lambda \underline{m}(z)\right)^{-3}\left[1-c\lambda^{2}\underline{m}^{2}(z)(1+\lambda \underline{m}(z))^{-2}\right]^{-1}\dif z\nonumber\\
    &\quad - \frac{1}{2\pi i}\oint_{\mathcal{C}} z\underline{m}(z)\left[1-c\lambda^{2}\underline{m}^{2}(z)(1+\lambda \underline{m}(z))^{-2}\right]^{-1}\nonumber\\
    &\qquad\times z\underline{m}(z)\left(1+\lambda \underline{m}(z)\right)^{-1} \times \left(h_{1}m(z)+\lambda m(z)+\lambda/z\right)\dif z\nonumber\\
    &\quad-\frac{1}{2\pi i}\oint_{\mathcal{C}} cz^{3}\underline{m}^{3}(z)\left[1-c\lambda^{2}\underline{m}^{2}(z)(1+\lambda \underline{m}(z))^{-2}\right]^{-1}\nonumber\\
    &\qquad\times\left(1+\lambda \underline{m}(z)\right)^{-1}\left[\beta m^{2}(z)+(h_2-2\lambda h_{1})m^{2}(z)+2\lambda^{2}m'(z)\right]\dif z\nonumber\\
    &=:I_{1}(f_1) + I_{2}(f_1) + I_{3}(f_1).
\end{align}
For $I_1(f_1)$, we get
\begin{align}%\label{f1integr01}
    I_1(f_1)   &= \frac{1}{2\pi i}\oint_{\mathcal{C}} c\lambda^{2}z\underline{m}^{3}(z)\left[1-c\lambda^{2}\underline{m}^{2}(z)(1+\lambda \underline{m}(z))^{-2}\right]^{-1}\nonumber\\
    &\qquad \times\left(1+\lambda \underline{m}(z)\right)^{-3}\left[1-c\lambda^{2}\underline{m}^{2}(z)(1+\lambda \underline{m}(z))^{-2}\right]^{-1}\dif z\nonumber\\
      &= \frac{1}{2\pi i}\oint_{\mathcal{C}} c\lambda^{2}z\underline{m}^{3}(z)\left[1-c\lambda^{2}\underline{m}^{2}(z)(1+\lambda \underline{m}(z))^{-2}\right]^{-1}\nonumber\\
    &\qquad \times\left(1+\lambda \underline{m}(z)\right)^{-3}\left[1-c\lambda^{2}\underline{m}^{2}(z)(1+\lambda \underline{m}(z))^{-2}\right]^{-1}\nonumber\\
    &\qquad \times[\underline{m}^{-2}-c\lambda^2(1+\lambda \underline{m})^{-2}]\dif\underline{m}\nonumber\\
    &= \frac{1}{2\pi i}\oint_{\mathcal{C}}\frac{c\lambda^2(-c\lambda\underline{m}+\lambda\underline{m}+1)}{(1+\lambda\underline{m})^2\big[c(\lambda\underline{m})^2-(1+\lambda\underline{m})^2\big]}\dif\underline{m}.\nonumber
    % &=\lambda (1+c)-\frac{\lambda}{2}(1+\sqrt{c})^2-\frac{\lambda}{2}(1-\sqrt{c})^2\nonumber\\
    % &=0
\end{align}   
For the first integral $I_1(f_1)$, the poles are $-\lambda^{-1}$, $-\frac{1}{(1\pm\sqrt{c})\lambda}$, we have by the residue theorem
\begin{equation}\label{f1integr1}
    I_1(f_1)  =\lambda (1+c)-\frac{\lambda}{2}(1+\sqrt{c})^2-\frac{\lambda}{2}(1-\sqrt{c})^2=0.
\end{equation} 
For the second integral $I_2(f_1)$, the poles are $-\lambda^{-1}$, we have, by the residue theorem,
\begin{align}\label{f1integr2}
    I_2(f_1) &= - \frac{1}{2\pi i}\oint_{\mathcal{C}} z\underline{m}(z)\left[1-c\lambda^{2}\underline{m}^{2}(z)(1+\lambda \underline{m}(z))^{-2}\right]^{-1}\nonumber\\
    &\qquad\times z\underline{m}(z)\left(1+\lambda \underline{m}(z)\right)^{-1} \times \left(h_{1}m(z)+\lambda m(z)+\lambda/z\right)\dif z\nonumber\\
&=- \frac{1}{2\pi i}\oint_{\mathcal{C}} z\underline{m}(z)\left[1-c\lambda^{2}\underline{m}^{2}(z)(1+\lambda \underline{m}(z))^{-2}\right]^{-1}\nonumber\\
    &\qquad\times z\underline{m}(z)\left(1+\lambda \underline{m}(z)\right)^{-1} \times \left(h_{1}\underline{m}(z)+\lambda \underline{m}(z)+ \frac{\lambda}{z}+(h_1+\lambda )\frac{1-c}{cz}\right)\nonumber\\
     &\qquad \times[\underline{m}^{-2}-c\lambda^2(1+\lambda \underline{m})^{-2}]\dif\underline{m}\nonumber\\
    &= - \frac{1}{2\pi i}\oint_{\mathcal{C}} \frac{(\lambda^2\underline{m}-h_1)(-c\lambda \underline{m}+\lambda \underline{m}+1)}{\underline{m}(1+\lambda \underline{m})^3}\dif\underline{m}\nonumber\\
    &=h_1.
\end{align}
For $I_3(f_1)$, we get
\begin{align}%\label{f1integr3}
    I_3(f_1) &= -\frac{1}{2\pi i}\oint_{\mathcal{C}} czz^{2}\underline{m}^{3}(z)\left[1-c\lambda^{2}\underline{m}^{2}(z)(1+\lambda \underline{m}(z))^{-2}\right]^{-1}\nonumber\\
  &\qquad\times\left[1+\lambda \underline{m}(z)\right]^{-1}\left[\beta m^{2}(z)+(h_2-2\lambda h_{1})m^{2}(z)+2\lambda^{2}m'(z)\right]\dif z\nonumber\\
  &= -\frac{1}{2\pi i}\oint_{\mathcal{C}} cz^{3}\underline{m}^{3}(z)\left[1-c\lambda^{2}\underline{m}^{2}(z)(1+\lambda \underline{m}(z))^{-2}\right]^{-1}\nonumber\\
  &\qquad\times\left[1+\lambda \underline{m}(z)\right]^{-1}\Big[(\beta+h_2-2\lambda h_{1})\Bigl(\frac{1}{c^2}\underline{m}^2z^2+\frac{2(1-c)}{c^2}\underline{m}z+\frac{(1-c)^2}{c^2}\Bigr)\nonumber\\
  &\qquad +2\lambda^{2}\Bigl(\frac{1}{c}z^2\underline{m}'-\frac{1-c}{c}\Bigr)\Big]\times\Big[\underline{m}^{-2}-c\lambda^2(1+\lambda \underline{m})^{-2}\Big]\dif\underline{m}\nonumber\\
    &=-\frac{1}{2\pi i}\oint_{\mathcal{C}} \Big[c(\lambda \underline{m} + 1)^4\Bigl(c\lambda^2\underline{m}^2 - (\lambda \underline{m} + 1)^2\Bigr)\Big]^{-1}\nonumber\\
    &\ \ \ \ \times\Big[2c\lambda^2\big(\lambda \underline{m} + 1\big)^2\big((c - 1)[c\lambda^2\underline{m}^2 - (\lambda \underline{m} + 1)^2] - (-c\lambda \underline{m} + \lambda \underline{m} + 1)^2\big)\nonumber\\
    &\ \ \ \ \ \ \ \ \ + \big(c\lambda^2\underline{m}^2 - (\lambda \underline{m} + 1)^2\big)\big(-2\lambda h_1 + \beta + h_2\big)\nonumber\\
    &\ \ \ \ \ \ \ \ \ \ \ \ \ \ \ \ \ \ \ \ \times \big(-c\lambda \underline{m} + \lambda \underline{m} + (c - 1)(\lambda \underline{m} + 1) + 1\big)^2\Big]\nonumber\\
    &\ \ \ \ \ \ \ \ \times(-c\lambda \underline{m} + \lambda \underline{m} + 1)\dif\underline{m}\nonumber
    %  &=-2\lambda (1+c)+\lambda (1+\sqrt{c})^2+\lambda (1-\sqrt{c})^2
    % \nonumber\\
    % &=0
\end{align}
For the third integral $I_3(f_1)$, the poles are $-\lambda^{-1}$, $-\frac{1}{(1\pm\sqrt{c})\lambda}$, we have by the residue theorem
\begin{equation}\label{f1integr3}
I_3(f_1) =(1+c)+\lambda (1+\sqrt{c})^2+\lambda (1-\sqrt{c})^2=0.
\end{equation}
Thus, by \eqref{f1expe0}, \eqref{f1integr1}--\eqref{f1integr3}, we get 
\begin{align}
    \Expe X_x = h_1.
\end{align}
For $f_2=x^2$, we have
    \begin{align}\label{f2expe0}
    \Expe X_{x^2}
    &= \frac{1}{2\pi i}\oint_{\mathcal{C}} c\lambda^{2}z^2\underline{m}^{3}(z)\left[1-c\lambda^{2}\underline{m}^{2}(z)(1+\lambda \underline{m}(z))^{-2}\right]^{-1}\nonumber\\
    &\qquad \times\left(1+\lambda \underline{m}(z)\right)^{-3}\left[1-c\lambda^{2}\underline{m}^{2}(z)(1+\lambda \underline{m}(z))^{-2}\right]^{-1}\dif z\nonumber\\
    &\quad - \frac{1}{2\pi i}\oint_{\mathcal{C}} z^2\underline{m}(z)\left[1-c\lambda^{2}\underline{m}^{2}(z)(1+\lambda \underline{m}(z))^{-2}\right]^{-1}\nonumber\\
    &\qquad\times z\underline{m}(z)\left(1+\lambda \underline{m}(z)\right)^{-1} \times \left(h_{1}m(z)+\lambda m(z)+\lambda/z\right)\dif z\nonumber\\
    &\quad-\frac{1}{2\pi i}\oint_{\mathcal{C}} cz^4\underline{m}^{3}(z)\left[1-c\lambda^{2}\underline{m}^{2}(z)(1+\lambda \underline{m}(z))^{-2}\right]^{-1}\nonumber\\
    &\qquad\times\left(1+\lambda \underline{m}(z)\right)^{-1}\left[\beta m^{2}(z)+(h_2-2\lambda h_{1})m^{2}(z)+2\lambda^{2}m'(z)\right]\dif z\nonumber\\  
    &=:I_1(f_2)+I_2(f_2)+I_3(f_2).
\end{align}
For the first integral $I_1(f_2)$, the poles are $-\lambda^{-1}$, $-\frac{1}{(1\pm\sqrt{c})\lambda}$, we have by the residue theorem
\begin{align}\label{f2integr1}
    I_1(f_2)   &= \frac{1}{2\pi i}\oint_{\mathcal{C}} c\lambda^{2}z^2\underline{m}^{3}(z)\left[1-c\lambda^{2}\underline{m}^{2}(z)(1+\lambda \underline{m}(z))^{-2}\right]^{-1}\nonumber\\
    &\qquad \times\left(1+\lambda \underline{m}(z)\right)^{-3}\left[1-c\lambda^{2}\underline{m}^{2}(z)(1+\lambda \underline{m}(z))^{-2}\right]^{-1}\dif z\nonumber\\
      &= \frac{1}{2\pi i}\oint_{\mathcal{C}} c\lambda^{2}z^2\underline{m}^{3}(z)\left[1-c\lambda^{2}\underline{m}^{2}(z)(1+\lambda \underline{m}(z))^{-2}\right]^{-1}\nonumber\\
    &\qquad \times\left(1+\lambda \underline{m}(z)\right)^{-3}\left[1-c\lambda^{2}\underline{m}^{2}(z)(1+\lambda \underline{m}(z))^{-2}\right]^{-1}\nonumber\\
    &\qquad \times[\underline{m}^{-2}-c\lambda^2(1+\lambda \underline{m})^{-2}]\dif\underline{m}\nonumber\\
    &= \frac{1}{2\pi i}\oint_{\mathcal{C}}\frac{-c\lambda^2(-c\lambda \underline{m}+\lambda \underline{m}+1)^2}{\underline{m}(1+\lambda \underline{m})^3\big[c(\lambda \underline{m})^2-(1+\lambda \underline{m})^2\big]}\dif\underline{m} \nonumber\\
    &=\lambda^2(1+5c+c^2)-\frac{\lambda^2}{2}(1+\sqrt{c})^4-\frac{\lambda^2}{2}(1-\sqrt{c})^4\nonumber\\
    &=-c\lambda^2.
\end{align}
For the second integral $I_2(f_2)$, the poles are $-\lambda^{-1}$, we have by the residue theorem,
\begin{align}\label{f2integr2}
    I_2(f_2) 
    &= - \frac{1}{2\pi i}\oint_{\mathcal{C}} z^2\underline{m}(z)\left[1-c\lambda^{2}\underline{m}^{2}(z)(1+\lambda \underline{m}(z))^{-2}\right]^{-1}\nonumber\\
    &\qquad\times z\underline{m}(z)\left(1+\lambda \underline{m}(z)\right)^{-1} \left[h_{1}m(z)+\lambda m(z)+\lambda/z\right]\dif z\nonumber\\
&=- \frac{1}{2\pi i}\oint_{\mathcal{C}} z^2\underline{m}(z)\left[1-c\lambda^{2}\underline{m}^{2}(z)(1+\lambda \underline{m}(z))^{-2}\right]^{-1}\nonumber\\
    &\qquad\times z\underline{m}(z)\left(1+\lambda \underline{m}(z)\right)^{-1} \times \left(h_{1}\underline{m}(z)+\lambda \underline{m}(z)+\lambda/z+(h_1+\lambda )\frac{1-c}{cz}\right)\nonumber\\
     &\qquad \times[\underline{m}^{-2}-c\lambda^2(1+\lambda \underline{m})^{-2}]\dif\underline{m}\nonumber\\
    &= - \frac{1}{2\pi i}\oint_{\mathcal{C}} \frac{-(\lambda^2\underline{m}-h_1)[-c\lambda \underline{m}+\lambda \underline{m}+1]^2}{\underline{m}^2(1+\lambda \underline{m})^4}\dif\underline{m}\nonumber\\
    &=\lambda (\lambda +2ch_1+2h_1).
\end{align}
For the third integral $I_3(f_2)$, the poles are $-\lambda^{-1}$, $-\frac{1}{(1\pm\sqrt{c})\lambda}$, we have by the residue theorem
\begin{align}
    I_3(f_2) &= -\frac{1}{2\pi i}\oint_{\mathcal{C}} cz^2z^2\underline{m}^3(z)\left[1-c\lambda^{2}\underline{m}^2(z)(1+\lambda \underline{m}(z))^{-2}\right]^{-1}\nonumber\\
  &\qquad\times\left[1+\lambda \underline{m}(z)\right]^{-1}\left[\beta m^2(z)+(h_2-2\lambda h_1)m^2(z)+2\lambda^{2}m'(z)\right]\dif z\nonumber\\
  &= -\frac{1}{2\pi i}\oint_{\mathcal{C}} cz^2z^2\underline{m}^{3}(z)\left[1-c\lambda^{2}\underline{m}^{2}(z)(1+\lambda \underline{m}(z))^{-2}\right]^{-1}\nonumber\\
  &\qquad\times\left[1+\lambda \underline{m}(z)\right]^{-1}\Big[(\beta+h_2-2\lambda h_{1})\left(\frac{1}{c^2}\underline{m}^2z^2+\frac{2(1-c)}{c^2}\underline{m}z+\frac{(1-c)^2}{c^2}\right)\nonumber\\
  &\qquad +2\lambda^{2}\Bigl(\frac{1}{c}z^2\underline{m}'-\frac{1-c}{c}\Bigr)\Big]\times\Big[\underline{m}^{-2}-c\lambda^2(1+\lambda \underline{m})^{-2}\Big]\dif\underline{m}\nonumber\\
 &=-2\lambda^2(1+5c+c^2)+c(\beta+h_2-2\lambda h_1)+\lambda^2(1+\sqrt{c})^4+\lambda^2(1-\sqrt{c})^4\nonumber\\
 &=2c\lambda^2+c(\beta+h_2-2\lambda h_1).\label{f2integr3}
\end{align}
Thus, by \eqref{f2expe0}, \eqref{f2integr1}--\eqref{f2integr3}, we get
\begin{align}
    \Expe X_{x^2}&=-c\lambda^2+\lambda (\lambda +2ch_1+2h_1)+2c\lambda^2+c(\beta+h_2-2\lambda h_1)\nonumber\\
 &=  (c+1)\lambda^2+2(c+1)\lambda h_1+c(\beta+h_2-2\lambda h_1).\nonumber
\end{align}

For $f_3=x^3$, we have
\begin{align}\label{f3expe0}
    \Expe X_{x^3}
    &= \frac{1}{2\pi i}\oint_{\mathcal{C}} c\lambda^{2}z^3\underline{m}^{3}(z)\left[1-c\lambda^{2}\underline{m}^{2}(z)(1+\lambda \underline{m}(z))^{-2}\right]^{-1}\nonumber\\
    &\qquad \times\left(1+\lambda \underline{m}(z)\right)^{-3}\left[1-c\lambda^{2}\underline{m}^{2}(z)(1+\lambda \underline{m}(z))^{-2}\right]^{-1}\dif z\nonumber\\
    &\quad - \frac{1}{2\pi i}\oint_{\mathcal{C}} z^3\underline{m}(z)\left[1-c\lambda^{2}\underline{m}^{2}(z)(1+\lambda \underline{m}(z))^{-2}\right]^{-1}\nonumber\\
    &\qquad\times z\underline{m}(z)\left(1+\lambda \underline{m}(z)\right)^{-1} \times \left[h_{1}m(z)+\lambda m(z)+\lambda/z\right]\dif z\nonumber\\
    &\quad-\frac{1}{2\pi i}\oint_{\mathcal{C}} cz^5\underline{m}^{3}(z)\left[1-c\lambda^{2}\underline{m}^{2}(z)(1+\lambda \underline{m}(z))^{-2}\right]^{-1}\nonumber\\
    &\qquad\times\left(1+\lambda \underline{m}(z)\right)^{-1}\left[\beta m^{2}(z)+(h_2-2\lambda h_{1})m^{2}(z)+2\lambda^{2}m'(z)\right]\dif z\nonumber\\  
    &=:I_1(f_3)+I_2(f_3)+I_3(f_3).
\end{align}
For the first integral $I_1(f_3)$, the poles are $-\lambda^{-1}$, $-\frac{1}{(1\pm\sqrt{c})\lambda}$, we have by the residue theorem
\begin{align}\label{f3integr1}
    I_1(f_3)   &= \frac{1}{2\pi i}\oint_{\mathcal{C}} c\lambda^{2}z^3\underline{m}^{3}(z)\left[1-c\lambda^{2}\underline{m}^{2}(z)(1+\lambda \underline{m}(z))^{-2}\right]^{-1}\nonumber\\
    &\qquad \times\left(1+\lambda \underline{m}(z)\right)^{-3}\left[1-c\lambda^{2}\underline{m}^{2}(z)(1+\lambda \underline{m}(z))^{-2}\right]^{-1}\dif z\nonumber\\
      &= \frac{1}{2\pi i}\oint_{\mathcal{C}} c\lambda^{2}z^3\underline{m}^{3}(z)\left[1-c\lambda^{2}\underline{m}^{2}(z)(1+\lambda \underline{m}(z))^{-2}\right]^{-1}\nonumber\\
    &\qquad \times\left(1+\lambda \underline{m}(z)\right)^{-3}\left[1-c\lambda^{2}\underline{m}^{2}(z)(1+\lambda \underline{m}(z))^{-2}\right]^{-1}\nonumber\\
    &\qquad \times[\underline{m}^{-2}-c\lambda^2(1+\lambda \underline{m})^{-2}]\dif\underline{m}\nonumber\\
    &= \frac{1}{2\pi i}\oint_{\mathcal{C}}\frac{c\lambda^2(-c\lambda \underline{m}+\lambda \underline{m}+1)^3}{\underline{m}^2(1+\lambda \underline{m})^4\big[c(\lambda \underline{m})^2-(1+\lambda \underline{m})^2\big]}\dif\underline{m} \nonumber\\
    &=\lambda^3(1+12c+12c^2+c^3)-\frac{\lambda^3}{2}(1+\sqrt{c})^6-\frac{\lambda^3}{2}(1-\sqrt{c})^6\nonumber\\
    &=-3c(c+1)\lambda^3.
\end{align}
For the second integral $I_2(f_3)$, the poles are $-\lambda^{-1}$, we have by the residue theorem
\begin{align}\label{f3integr2}
    I_2(f_3) 
    &= - \frac{1}{2\pi i}\oint_{\mathcal{C}} z^3\underline{m}(z)\left[1-c\lambda^{2}\underline{m}^{2}(z)(1+\lambda \underline{m}(z))^{-2}\right]^{-1}\nonumber\\
    &\qquad\times z\underline{m}(z)\left(1+\lambda \underline{m}(z)\right)^{-1} \times \left[h_{1}m(z)+\lambda m(z)+\lambda/z \right]\dif z\nonumber\\
&=- \frac{1}{2\pi i}\oint_{\mathcal{C}} z^3\underline{m}(z)\left[1-c\lambda^{2}\underline{m}^{2}(z)(1+\lambda \underline{m}(z))^{-2}\right]^{-1}\nonumber\\
    &\qquad\times z\underline{m}(z)\left(1+\lambda \underline{m}(z)\right)^{-1} \times \left(h_{1}\underline{m}(z)+\lambda \underline{m}(z)+\lambda/z+(h_1+\lambda )\frac{1-c}{cz}\right)\nonumber\\
     &\qquad \times[\underline{m}^{-2}-c\lambda^2(1+\lambda \underline{m})^{-2}]\dif\underline{m}\nonumber\\
    &= - \frac{1}{2\pi i}\oint_{\mathcal{C}} \frac{[\lambda^2\underline{m}-h_1][-c\lambda \underline{m}+\lambda \underline{m}+1]^3}{\underline{m}^3(1+\frac{\sigma^2}{\mu^3}\underline{m})^5}\dif\underline{m}\nonumber\\
    &=\lambda^2(3c\lambda +2\lambda +3c^2h_1+9ch_1+3h_1).
\end{align}
For the third integral $I_3(f_3)$, the poles are $-\lambda^{-1}$, $-\frac{1}{(1\pm\sqrt{c})\lambda}$, we have by the residue theorem
\begin{align}\label{f3integr3}
    I_3(f_3) &= -\frac{1}{2\pi i}\oint_{\mathcal{C}} cz^3z^{2}\underline{m}^{3}(z)\left[1-c\lambda^{2}\underline{m}^{2}(z)(1+\lambda \underline{m}(z))^{-2}\right]^{-1}\nonumber\\
 &\qquad\times\left[1+\lambda \underline{m}(z)\right]^{-1}\left[\beta m^{2}(z)+(h_2-2\lambda h_{1})m^{2}(z)+2\lambda^{2}m'(z)\right]\dif z\nonumber\\
  &= -\frac{1}{2\pi i}\oint_{\mathcal{C}} cz^3z^{2}\underline{m}^{3}(z)\left[1-c\lambda^{2}\underline{m}^{2}(z)(1+\lambda \underline{m}(z))^{-2}\right]^{-1}\nonumber\\
  &\qquad\times\left[1+\lambda \underline{m}(z)\right]^{-1}\Big[(\beta+h_2-2\lambda h_{1})\Bigl(\frac{1}{c^2}\underline{m}^2z^2+\frac{2(1-c)}{c^2}\underline{m}z+\frac{(1-c)^2}{c^2}\Bigr)\nonumber\\
  &\qquad +2\lambda^{2}\Bigl(\frac{1}{c}z^2\underline{m}'-\frac{1-c}{c}\Bigr)\Big]\times\Big[\underline{m}^{-2}-c\lambda^2(1+\lambda \underline{m})^{-2}\Big]\dif\underline{m}\nonumber\\
 &=-2\lambda^3(1+12c+12c^2+c^3)+3c(c+1)\lambda (\beta+h_2-2\lambda h_1)\nonumber\\
 &\ \ \ \ \ \ +\lambda^3(1+\sqrt{c})^6+\lambda^3(1-\sqrt{c})^6\nonumber\\
 &=3c(c+1)\lambda (\beta+h_2-2\lambda h_1)+6c(c+1)\lambda^3.
\end{align}
Thus, by \eqref{f3expe0}, \eqref{f3integr1}--\eqref{f3integr3}, we get
\begin{align}
   \Expe X_{x^3} 
    &= -3c(c+1)\lambda^3+\lambda^2(3c\lambda +2\lambda +3c^2h_1+9ch_1+3h_1)\nonumber\\
    &\ \ \ \ \ \ +3c(c+1)\lambda (\beta+h_2-2\lambda h_1)+6c(c+1)\lambda^3\nonumber\\
    &= (3c^2+6c+2)\lambda^3+3(c^2+3c+1)\lambda^2h_1 +3c(c+1)\lambda (\beta+h_2-2\lambda h_1).\nonumber
\end{align}
Thus,
\begin{align}
    \mu_1&:=\Expe X_{f_1} = h_1,\nonumber\\
    \mu_2&:=\Expe X_{f_2} %&=
   % -c\lambda^2+\lambda (\lambda +2ch_1+2h_1)+c\Big(\beta+h_2-2\lambda h_1\Big)+2c\lambda^2\nonumber\\
    =  (1+c)\lambda^2+2(1+c)\lambda h_1+c(\beta+h_2-2\lambda h_1),\nonumber\\
    \mu_3&:=\Expe X_{f_3} %&= -3c(c+1)\lambda^3+\lambda^2\Big(3c\lambda +2\lambda +3c^2h_1+9ch_1+3h_1\Big)\nonumber\\
    %&\ \ \ \ \ \ +3c(c+1)\lambda \Big(\beta+h_2-2\lambda h_1\Big)+6c(c+1)\lambda^3\nonumber\\
    = (2+6c+3c^2)\lambda^3+3(1+3c+c^2)\lambda^2h_1+3c(1+c)\lambda (\beta+h_2-2\lambda h_1),\nonumber
\end{align}
and
\begin{align}
    V_{(1,1)} &= 2c\lambda^2+c(\beta+h_2-2\lambda h_1),\nonumber\\%(\beta + h_2-2\lambda h_1)
    V_{(2,2)} &=4c(2+c)(1+2c)\lambda^4+4c(1+c)^2\lambda^2(\beta+h_2-2\lambda h_1),\nonumber\\%&=4c(2+5c+2c^2)\lambda^4+4c(1+c)^2\lambda^2\xi\nonumber\\%(\beta + h_2-2\lambda h_1)
    V_{(3,3)} &=6c(1+6c+3c^2)(3+6c+c^2)\lambda^6 +9c(1+3c+c^2)^2\lambda^4(\beta+h_2-2\lambda h_1).\nonumber%(\beta + h_2-2\lambda h_1).
\end{align}

\subsection{Tightness of $M_{p}^{(1)}(z)$}\label{prftight}

The tightness of $M_p^{(1)}(z)$ is similar to that provided in \cite{bai2004clt}. It is sufficient to prove the moment condition (12.51) of \cite{billingsley1968convergence}, i.e.
\begin{eqnarray}\label{momentbound-supp}
\sup_{n; z_1, z_2\in\mathcal{C}_n}\frac{\Expe|M_p^{(1)}(z_1)-M_p^{(1)}(z_2)|^2}{|z_1-z_2|^2}
\end{eqnarray}
is finite.

Before proceeding, we provide some results needed in the proof later. First, moments of $\|\bD^{-1}(z)\|$, $\|\bD^{-1}_j(z)\|$ and $\|\bD_{ij}^{-1}(z)\|$ are bounded in $p$ and $z\in\mathcal{C}_p$. It is easy to see that it is true for $z\in\mathcal{C}_u$ and for $z\in\mathcal{C}_{\ell}$ if $x_{\ell}<0$. For $z\in\mathcal{C}_r$ or, if $x_{\ell}>0$, $z\in\mathcal{C}_{\ell}$, we have from Proposition \ref{extreigen} that
\begin{eqnarray*}
\Expe\|\bD^{-1}_j(z)\|^{m}&\leq&K_1+v^{-m}P(\|\bB_{(j)}\|\geq\eta_r \ or \ \lambda_{\min}(\bB_{(j)})\leq\eta_{\ell})\nonumber\\
&\leq&K_1+K_2n^{m}\varepsilon^{-m}n^{-\ell}\leq K
\end{eqnarray*}
for large $\ell$, where $\bB_{(j)}=\bB_p-\br_j\br_j'$. Here $\eta_r$ is any number between $\frac{\sigma^2}{\mu^2}(1+\sqrt{c})^2$ and $x_r$; if $x_{\ell}>0$, $\eta_{\ell}$ is any number between $x_{\ell}$ and $\frac{\sigma^2}{\mu^2}(1-\sqrt{c})^2$ and if $x_{\ell}<0$, $\eta_{\ell}$ can be any negative number. So for any positive integer $m$,
\begin{eqnarray}
\max\Big(\Expe\|\bD^{-1}(z)\|^{m}, \Expe\|\bD_j^{-1}(z)\|^{m}, \Expe\|\bD_{ij}^{-1}(z)\|^{m}\Big)\leq K.
\end{eqnarray}
By the argument above, we can extend Lemma \ref{bound} and get
\begin{eqnarray}\label{boundextend}
\Big|\Expe\Big(a(v)\prod^{q}_{\ell=1}\big(\br_1'\bB_{(\ell)}(v)\br_1-n^{-1}\tr\bB_{(\ell)}(v)\big)\Big)\Big|\leq Kn^{-1}\delta_n^{2q-4},
\end{eqnarray}
where $\bB_{\ell}(v)$ is independent of $\br_1$ and
\begin{eqnarray*}
\max\bigl(|a(v)|, \|\bB_{(\ell)}(v)\|\bigr)\leq K\bigl(1+n^{s}I(\|\bB_p\|\geq\eta_r \ \text{ or } \ \lambda_{\min}(\widetilde{\bB})\leq\eta_{\ell})\bigr),
\end{eqnarray*}
with $\widetilde{\bB}$ being $\bB_{(j)}$ or $\bB_p$.
By \eqref{boundextend}, we have
\begin{eqnarray}\label{epsibound}
\Expe|\varepsilon_j(z)|^{m}\leq K_mn^{-1}\delta_n^{2m-4}.
\end{eqnarray}

Let $\gamma_j(z)=\br_j'\bD_j^{-1}(z)\br_j-n^{-1}\nu_2\Expe \tr\bD^{-1}_j(z)$.
By Lemma \ref{burk}, \eqref{boundextend} and H\"{o}lder's inequality, with similar derivation on page 580 of \cite{bai2004clt}, we have
\begin{eqnarray}\label{gamepsbound}
\Expe|\gamma_j(z)-\varepsilon_j(z)|^{m}\leq\frac{K_m}{n^{m/2}}.
\end{eqnarray}
It follows from \eqref{epsibound} and \eqref{gamepsbound} that
\begin{eqnarray}\label{gambound}
\Expe|\gamma_j|^{m}\leq K_mn^{-1}\delta_n^{2m-4}, \ \ m\geq 2.
\end{eqnarray}
Next, we prove that $b_p(z)$ is bounded. With \eqref{boundextend}, we have for any $m\geq 1$,
\begin{eqnarray}\label{beta1}
\Expe|\beta_1(z)|^{m}\leq K_m.
\end{eqnarray}
Since $b_p(z)=\beta_1(z)+\beta_1(z)b_p(z)\gamma_1(z)$, it is derived from \eqref{gambound} and \eqref{beta1} that $|b_p(z)|\leq K_1+K_2|b_p(z)|n^{-1/2}$.

Then
\begin{eqnarray}\label{bpbound}
|b_p(z)|\leq\frac{K_1}{1-K_2n^{-1/2}}\leq K.
\end{eqnarray}

With \eqref{boundextend}--\eqref{bpbound} and the same approach on page 581-583 of \cite{bai2004clt}, we can obtain that \eqref{momentbound-supp} is finite.

 \section{Simulation of CLT for $M_p(z)$}\label{simuforcltmpz}

In this section, we compare the empirical mean and covariance of $M_p(z)=\tr(\bB_p^0-z\bI_p)^{-1}-pm_{F^{c_n}}(z)$ with their theoretical limits as stated in Proposition \ref{mpzclt}. This proposition is a key step for the proof of our main result, Theorem \ref{clt}. Readers are referred to Section \ref{prfcltall} for more details of $M_p(z)$. 
We consider two types of data distribution of $w_{ij}$ as follows:
\begin{enumerate}
    \item $w_{ij}$ follows the exponential distribution with rate parameter $5$;
    \item $w_{ij}$ follows the Chi-square distribution with degree of freedom $1$. 
\end{enumerate}

Empirical values of $\Expe M_p(z)$ and $\Cov(M_p(z_1),M_p(z_2))$ are calculated for various combinations of $(p,n)$ with $p/n=3/4$ or $p/n=1$.
For each pair of $(p,n)$, $2000$ independent replications are used to obtain the empirical values. Table \ref{tab:Mp-mean} reports the empirical mean of $M_p(z)$ with $z=\pm 3 + 2i$ for both $\mathrm{Exp}(5)$ population and $\chi^2(1)$ population. The empirical results of $\Cov(M_p(z_1),M_p(z_2))$ are reported in Tanle \ref{tab:Mp-Cov}.  As shown in Tables \ref{tab:Mp-mean} -- \ref{tab:Mp-Cov}, the empirical values of $\Expe M_p(z)$ and $\Cov(M_p(z_1),M_p(z_2))$ closely match their respective theoretical limits under all scenarios.

\begin{table}[htbp]
	\centering
	\caption{Empirical mean of $M_p(z)$ with $z=\mp 3 + 2i$.}
        \begin{tabular}{cccccccc}
		\toprule
		&       &       & \multicolumn{2}{c}{Exp(5)} &       & \multicolumn{2}{c}{$\chi^2(1)$} \\
		\cmidrule{4-5}\cmidrule{7-8}          & $p/n$ & $n$   & -3+2$i$ & 3+2i  &       & -3+2$i$ & 3+2$i$ \\
		\midrule
		\multirow{4}[2]{*}{Emp} & \multirow{4}[2]{*}{$3/4$} & 100   & 0.0586+0.0857$i$ & -0.0373-0.249$i$ &       & 0.1405+0.1628$i$ & -0.55-0.2732$i$ \\
		&       & 200   & 0.0582+0.0858$i$ & -0.0311-0.2526$i$ &       & 0.1459+0.1697$i$ & -0.5761-0.3089$i$ \\
		&       & 300   & 0.0567+0.0844$i$ & -0.0336-0.2566$i$ &       & 0.1465+0.1712$i$ & -0.5705-0.3212$i$ \\
		&       & 400   & 0.0596+0.0878$i$ & -0.0352-0.2528$i$ &       & 0.1463+0.172$i$ & -0.5631-0.3465$i$ \\
		\cmidrule{2-8}    Theo  &       &       & \textbf{0.0587+0.0872i} & \textbf{-0.029-0.2529i} &       & \textbf{0.15+0.1768i} & \textbf{-0.5792-0.3764i} \\
		\midrule
		\multirow{4}[2]{*}{Emp} & \multirow{4}[2]{*}{$5/4$} & 100   & 0.0547+0.0766$i$ & -0.1069-0.2671$i$ &       & 0.1366+0.1473$i$ & -0.5458-0.1545$i$ \\
		&       & 200   & 0.0572+0.0793$i$ & -0.1109-0.2757$i$ &       & 0.1395+0.1518$i$ & -0.5847-0.1787$i$ \\
		&       & 300   & 0.0587+0.0808$i$ & -0.1074-0.2752$i$ &       & 0.1382+0.1511$i$ & -0.5747-0.1934$i$ \\
		&       & 400   & 0.0559+0.0778$i$ & -0.0949-0.2733$i$ &       & 0.1434+0.1553$i$ & -0.5751-0.1933$i$ \\
		\cmidrule{2-8}    Theo  &       &       & \textbf{0.0578+0.0804$i$} & \textbf{-0.0919-0.2764$i$} &       & \textbf{0.1432+0.1569$i$} & \textbf{-0.6025-0.2149$i$} \\
		\bottomrule
	\end{tabular}%
	\label{tab:Mp-mean}%
\end{table}%

\begin{table}[htbp]
	\centering
	\caption{Empirical covariance between $M_p(z_1)$ and $M_p(z_2)$.}
	\begin{tabular}{cccccccc}
		\toprule
		&       &       & \multicolumn{2}{c}{Exp(5)} &       & \multicolumn{2}{c}{$\chi^2(1)$} \\
		\cmidrule{4-5}\cmidrule{7-8}          & $p/n$ & $n$   & (-3+2$i$,-1+1$i$)\tabnoteref{tab1} & (3+2$i$,5+1$i$) &       & (-3+2$i$,-1+1$i$) & (3+2$i$,5+1$i$) \\
		\midrule
		\multirow{4}[2]{*}{Emp} & \multirow{4}[2]{*}{$3/4$} & 100   & -0.0038+0.0147$i$ & -0.04+0.0035$i$ &       & 0+0.0304$i$ & 0.089+0.014$i$ \\
		&       & 200   & -0.0041+0.0163$i$ & -0.0418+0.0022$i$ &       & 0.0004+0.0326$i$ & 0.117+0.0284$i$ \\
		&       & 300   & -0.0043+0.0171$i$ & -0.0446+0.0011$i$ &       & 0+0.0335$i$ & 0.1372+0.0294$i$ \\
		&       & 400   & -0.0043+0.0168$i$ & -0.0465-0.0003$i$ &       & 0.0002+0.0356$i$ & 0.1273+0.036$i$ \\
		\cmidrule{2-8}    Theo  &       &       & \textbf{-0.0044+0.0172$i$} & \textbf{-0.0448-0.0002$i$} &       & \textbf{0.0006+0.0363$i$} & \textbf{0.1491+0.0524$i$} \\
		\midrule
		\multirow{4}[2]{*}{Emp} & \multirow{4}[2]{*}{$5/4$} & 100   & -0.0032+0.0197$i$ & -0.0483+0.0765$i$ &       & 0.0025+0.0349$i$ & 0.0931-0.0373$i$ \\
		&       & 200   & -0.0032+0.0196$i$ & -0.0545+0.0763$i$ &       & 0.0032+0.035$i$ & 0.0991-0.0406$i$ \\
		&       & 300   & -0.0036+0.0212$i$ & -0.0566+0.0708$i$ &       & 0.0026+0.0336$i$ & 0.0955-0.0209$i$ \\
		&       & 400   & -0.0032+0.02$i$ & -0.0594+0.0742$i$ &       & 0.0038+0.0374$i$ & 0.1138-0.0297$i$ \\
		\cmidrule{2-8}    Theo  &       &       & \textbf{-0.0034+0.0206$i$} & \textbf{-0.0624+0.0743$i$} &       & \textbf{0.0035+0.0388$i$} & \textbf{0.1099-0.0323$i$} \\
		\bottomrule
	\end{tabular}%
	\label{tab:Mp-Cov}%
        \begin{tabnotes}
			\tabnotetext[id=tab1,mark=a]{This row denotes different combinations of $(z_1,z_2)$.}
		\end{tabnotes}
\end{table}%

\bibliographystyle{imsart-nameyear} % Style BST file (imsart-number.bst or imsart-nameyear.bst)
\bibliography{reference}       % Bibliography file (usually '*.bib')

%% or include bibliography directly:
% \begin{thebibliography}{}
% \bibitem{b1}
% \end{thebibliography}